\documentclass[a4paper, 11pt]{article}
\title{Quotients of real algebraic sets and $\mathcal{AS}$-sets equipped with a linearizable compact group action}
\author{Fabien Priziac}

\date{}

\usepackage[utf8]{inputenc}

\usepackage{hyperref}

\usepackage{amsfonts}
\usepackage{amsmath}
\usepackage{amsthm}
\usepackage{amssymb}
\usepackage[all]{xy}
\usepackage{graphicx}

\newtheorem{de}{Definition}[subsection]
\newtheorem{theo}[de]{Theorem}
\newtheorem{prop}[de]{Proposition}
\newtheorem{cor}[de]{Corollary}

\newtheorem{lem}[de]{Lemma}

\newcommand{\N}{\mathbb{N}}
\newcommand{\Nstar}{\mathbb{N} \setminus \{0\}}
\newcommand{\R}{\mathbb{R}}
\newcommand{\C}{\mathbb{C}}
\newcommand{\OO}{{\rm O}}
\newcommand{\GL}{{\rm GL}}
\newcommand{\SO}{{\rm SO}}
\newcommand{\M}{{\rm M}}
\newcommand{\SL}{{\rm SL}}
\newcommand{\hookuparrow}{\mathrel{\rotatebox[origin=t]{-90}{\reflectbox{$\hookrightarrow$}}}}
\newcommand{\gq}{/\!\!/}

\theoremstyle{remark}
\newtheorem{rem}[de]{Remark}
\newtheorem{rems}[de]{Remarks}
\newtheorem{ex}[de]{Example}
\newtheorem{exs}[de]{Examples}

\begin{document}

\maketitle

\begin{abstract} We present the full-detailed construction of new geometric quotient structures for affine real algebraic varieties equipped with a linearizable action of a compact Lie group $G$. These quotients are functors from the category of arc-symmetric/$\mathcal{AS}$-sets equipped with a linearizable action of $G$ and equivariant continuous maps with $\mathcal{AS}$-graph to the category of~$\mathcal{AS}$-sets. We furthermore provide a complete review of the results employed in the constructions, including properties of equivariant real algebraic geometry with respect to polynomial group actions, as well as an introduction to semialgebraic arc-symmetric sets and $\mathcal{AS}$-sets of compact affine Nash manifolds.
\end{abstract}

\footnote{ \hspace{-0.6cm}{\bf Fabien Priziac} (fabien.priziac@univ-ubs.fr) Université Bretagne Sud, UMR CNRS 6205, LMBA, F-56000 Vannes, France



Research partially supported by Projet ANR NewMIRAGE (ANR-23-CE40-0002).
}

\tableofcontents

\section{Introduction}

In this text, we deal with real algebraic sets equipped with some action of a compact Lie group $G$: a real algebraic set of a real affine space $\R^n$ is the zero set of a polynomial with real coefficients. We consider more generally the semialgebraic arc-symmetric sets of $\R^n$ or of a compact Nash submanifold of $\R^n$, as well as their finite Boolean combinations called $\mathcal{AS}$-sets. A semialgebraic subset of $\R^n$ is a finite union of sets of $\R^n$ defined by polynomial equalities and inequalities, a Nash submanifold of $\R^n$ is a semialgebraic analytic submanifold of $\R^n$ having analytic charts with semialgebraic data, while an arc-symmetric set of $\R^n$ (this notion has been defined by K. Kurdyka in \cite{KurAR}) is a subset $S$ of $\R^n$ such that for any analytic arc with values in $\R^n$, either the arc meets $S$ at isolated points, either the arc is entirely contained in $S$.

Our goal is to construct and study some functorial quotients, with respect to the considered actions of~$G$, with a semialgebraic arc-symmetric structure, generalizing one part of our work in \cite{FP-QIAS}. In order to give the quotient by $G$ this geometric structure, our main tool will be the real algebro-geometric quotient by $G$. That is why we consider compact real algebraic groups with polynomial group structure (have in mind that a compact Lie group has, up to an isomorphism of Lie groups, such a representative: see remark \ref{remcompliegpreprescompralggp} of the present document and section 2.5 of Chapter 5 of \cite{OV}) and their polynomial actions on real algebraic sets (these actions are all linearizable up to equivariant polynomial isomorphism, hence the title of this document). Indeed, if $G$ is such a compact real algebraic group, Hilbert's finiteness theorem is valid and allows, together with Real Nullstellensatz, to associate to any $G$-real algebraic set a~$G$-invariant polynomial map which restricts into the topological quotient map (contrary to the complex case, this real algebro-geometric quotient map is not surjective in general).

We chose to explain the construction of our $\mathcal{AS}$-quotients in full details, as a self-contained pedagogical introduction to the subject for any reader being a little familiar with real algebraic sets and semialgebraic sets. In particular, we will show most part of the employed results.

A key new result is theorem \ref{theoquotientfreearcsymrnisarcsym} which asserts notably that, if $G$ is a compact real algebraic group, the semialgebraic quotient set of a $G$-real algebraic set on which $G$ acts \emph{freely} is an arc-symmetric set of a real affine space. From theorem \ref{theoquotientfreearcsymrnisarcsym} and its generalization to $\mathcal{AS}$-sets of $\R^n$ (theorem \ref{theoquotientfreearrnisar}), all statements of section \ref{sectionasquotientsofgassets} are original, while the first part of this document is dedicated to the proofs of preliminaries that lead to theorem \ref{theoquotientfreearcsymrnisarcsym} as well as surrounding properties, that are interesting on their own. This includes for instance a complete demonstration of Hilbert's finiteness theorem in the real case and a review of specific features of compact real algebraic groups. Even if statements of sections \ref{sectionralggp} to \ref{sectionfreepolyactionscompralggp} may reveal classical, well-known or else true in the complex case, some of them were hard or impossible to find in the literature: this was another motivation to write down those arguments, for future reference.

In order to preserve the flow of the arguments, some statements and proofs have been gathered in an appendix. In the last part of this appendix, we give in particular an introduction to arc-symmetric sets and $\mathcal{AS}$-sets. These classes of objects generalize in some sense real algebraic sets and closed analytic sets, keeping enough rigidity while allowing more flexibility. Following the seminal papers \cite{KurAR} and \cite{KurPar}, we first state definitions and properties about arc-symmetric sets of $\R^n$ together with full-detailed proofs. As for arc-symmetric sets in more general spaces, we decided to place ourselves in the framework of Nash manifolds. In a second time, we restrict to compact affine Nash manifolds and justify properties and characterizations of section 3 of \cite{KurPar} about arc-symmetric sets and $\mathcal{AS}$-sets of $\mathbb{P}^n(\R)$. Again, our intention is to provide a self-contained exhaustive review of the results that we need in sections \ref{sectionfreepolyactionscompralggp} and \ref{sectionasquotientsofgassets}, as well as to emphasize the importance of arc-symmetric sets and $\mathcal{AS}$-sets in real algebraic geometry.
\\

Let us now describe, in a more developed way, the different parts of this document, with this aim of constructing functorial $\mathcal{AS}$-quotients for any real algebraic set equipped with a polynomial action of a compact group.

As we said above, the fundamental tool that we use is the real geometric quotient which results from a combination of Hilbert's finiteness theorem and Real Nullstellensatz. These theorems impose to consider groups that are real algebraic sets and have polynomial group structure (we call such a group a \emph{real algebraic group}: definition \ref{defrealalggroup}) and, if $G$ is such a real algebraic group, real algebraic sets $X$ equipped with a polynomial action $\alpha : G \times X \rightarrow X$ of $G$ (we call such a set a \emph{$G$-real algebraic set}: definition \ref{degrealalgset}). Section \ref{sectionralggp} is dedicated to the former, while section \ref{sectionpolyactralggpralgsets} is dedicated to the latter. Notice that the polynomial action of $G$ on a $G$-real algebraic set $X$ can always be linearized (proposition \ref{proplinearizgrealalgset}): there exists an equivariant polynomial isomomorphism from $X$ to $G$-real algebraic subset of a polynomial representation $\R^N$ of $G$ (definition \ref{defrealpolyrepres}: the linear action of $G$ on $\R^N$ is given by a matrix with polynomial coordinates in $g \in G$). 
\\

Let $X$ be a $G$-real algebraic set. The $\R$-algebra $\mathcal{P}(X)$ of polynomial functions on $X$ can be equipped with an induced linear action of $G$, by right composition with the action, and we can consider the sub-algebra $\mathcal{P}(X)^G$ of $G$-invariant polynomial functions on $X$. If we assume that the real algebraic group $G$ is linearly reductive (definition \ref{defreallinearreductivity}, see also the characterizations of proposition \ref{propequivlinred}), according to Hilbert's finiteness theorem (theorem \ref{theohilbertfiniteness}), the $\R$-algebra $\mathcal{P}(X)^G$ is finitely generated: we dedicate section \ref{linredralggphilbfinth} to these notion and result. The equivalence of categories provided by Real Nullstellensatz (theorem \ref{theorealnullstellensatzequivcategories}) then asserts that the inclusion morphism $\mathcal{P}(X)^G \hookrightarrow \mathcal{P}(X)$ corresponds to a polynomial map $\pi_X : X \rightarrow X\gq G$ that we call the \emph{geometric quotient} map of $X$ (definition \ref{derealgeoquotient}). This geometric quotient is functorial with respect to equivariant polynomial maps between $G$-real algebraic sets and is compatible with equivariant inclusion: see lemmas \ref{lemfunctorrealalgquotient} and \ref{lemrestrictionrealalgquotient} of section \ref{sectionrealgeoquotient}. 

The above construction is analog to the construction of the geometric quotient map in the complex case and, due to the polynomial hypothesis, the real geometric quotient map $\pi_X$ of $X$ can be complexified into the complex geometric quotient $\pi_{X_{\C}}$ of its complexification $X_{\C}$ by the complexified action of $G_{\C}$ (see lemma \ref{lemcomplexificationpolyaction} and subsection \ref{subsectioncomplexrealgeomquotient}). From another point of view, $\pi_X$ is induced by $\pi_{X_{\C}}$:  we will also take advantage of this facts.

On the other hand, differences have to be noticed between the real and complex versions of the above notions and constructions. For instance, contrary to the complex case, the image of a morphism of real algebraic groups is not necessarily a real algebraic group (remark \ref{remsimagemorphrealalggr} 1) and the image of the real geometric quotient map is not necessarily surjective (see examples \ref{exrealalgquotientplaneorthoaction}). Nevertheless, real algebraic geometry comes with an extra feature, namely \emph{semialgebraic} sets. This allows to consider the image of a morphism of real algebraic groups as a semialgebraic group (definition \ref{defsemialggroup}), while the image of a real geometric quotient map is always a semialgebraic set.
\\

From section \ref{sectpolyactionscomprealalggroups}, we focus on polynomial (or, equivariantly, linearizable) actions of \emph{compact} real algebraic groups. If $G$ is a compact real algebraic group, the associated Haar measure allows to define a Reynolds operator for $G$ (definition \ref{dereynoldsop} and example \ref{exlinearlyredrealalggr} 2) and~$G$ is therefore a linearly reductive real algebraic group (corollary \ref{corequivlinearlyredreynoldsopgroup}). As a consequence, we can consider the functorial geometric quotient of $G$-real algebraic sets. 

A fundamental property of separation of orbits (lemma \ref{lemseporbitscompactsubg}) also allows to show that the image of a compact real algebraic group by a morphism of real algebraic groups is a (compact) real algebraic group. Moreover, polar decomposition can be used to prove that a bijective morphism of real algebraic groups between compact real algebraic groups (in fact any bijective continuous group homomorphism between compact real algebraic groups) is an isomorphism (corollary \ref{corbijmorphcomprealalggroupisiso} and example \ref{exsapplicationsofcontgpmorphbtwcompralggpsispoly} 3, see also theorem \ref{thpolardecompositionorthogonalcompactrealalggroup} as well as proposition \ref{propcompactrealalggroupisorthog} that asserts that any compact real algebraic group is isomorphic to a real algebraic subgroup of an orthogonal group $\OO_N(\R)$). 

This draws a parallel between the behaviors of compact real algebraic groups and of complex algebraic groups. But we have actually more: if $X$ is a $G$-real algebraic set, the orbit of any point of $X$ by the action of $G$ is a real algebraic set (proposition \ref{proporbitcompactsubgrouprealalggroup}) and the restriction~$\varpi_X : X \rightarrow \pi_X(X)$ of the geometric quotient of $X$ onto its semialgebraic image (called the \emph{semialgebraic quotient} of $X$: definition \ref{defsemialgquot}) is the topological quotient of $X$ by $G$ (proposition~\ref{propquotientmapproperclosedopen}). We thereafter show the functoriality of semialgebraic quotient with respect to semialgebraic (resp. continuous, resp. proper) maps (proposition \ref{propfunctsemialgquotientsemialgcontprop}).
\\

As we can see through examples \ref{exrealalgquotientplaneorthoaction}, the semialgebraic quotient set $X/G := \pi_X(X) \subset \R^r$ of a $G$-real algebraic set $X$ is not even an $\mathcal{AS}$-set of $\R^r$ in general. Nevertheless, if the action of~$G$ on $X$ is free, then the set $X/G$ is arc-symmetric in $\R^r$ (theorem \ref{theoquotientfreearcsymrnisarcsym}). This is a consequence of theorem \ref{theoimagenonsingularpointtrivstab} which asserts that a nonsingular point $x$ of $X$ \emph{with trivial stabilizer} is sent onto a nonsingular point in dimension $\dim X - \dim G$ of $X\gq G$. This fact implies the existence of a local Nash section for the geometric quotient map $\pi_X$ (theorem \ref{theolocalnashsectionifnonsingandfree}) which allows to lift some analytic arcs on $X \gq G$ passing through $\pi_X(x)$. The same technique allows to show further functorialities of the semialgebraic quotient, for instance its functoriality with respect to equivariant arc-analytic maps defined on $G$-real algebraic sets on which $G$ acts freely: see proposition \ref{propfunctarcanalyticmapfromfreegsaset}.

Let us also highlight a key property that is used to prove theorem \ref{theoimagenonsingularpointtrivstab}. If $x$ is a point of~$X$ such that $G_x = \{e\}$, then the bijective polynomial map which associates to any $g \in G$ the element $g \cdot x$ of the real algebraic orbit of $x$ is actually a polynomial isomorphism. The latter then induces an isomorphism of $\R$-algebras $\mathcal{P}(G \cdot x) \rightarrow \mathcal{P}(G)$ which reveals crucial in the proof of theorem \ref{theoimagenonsingularpointtrivstab}.
\\

In section \ref{sectionasquotientsofgassets}, we consider more generally $\mathcal{AS}$-sets of compact affine Nash manifolds. In view of section \ref{sectionfreepolyactionscompralggp}, we consider, in a first time (subsections \ref{subsectcompnashsubraff} and \ref{subsectglingassets}), $\mathcal{AS}$-sets that can be seen as a $G$-stable $\mathcal{AS}$-set of a compact Nash submanifold $N$ of a polynomial representation $\R^d$ of $G$ with the extra hypothesis that the action of $G$ on the $\mathcal{AS}$-closure in $N$ is free: see definition \ref{degassetfreegasset} for the precise definition of a \emph{free $G$-$\mathcal{AS}$-set}. 

Any such free $G$-$\mathcal{AS}$-set gives rise to a well-defined $\mathcal{AS}$-quotient: with the above notations, if $S$ is a $G$-stable $\mathcal{AS}$-set of $N$ such that the action of $G$ on $\overline{S}^{\mathcal{AS}}$ is free, the set $S/G := \pi_{\R^d}(S) \subset \R^r$ is an $\mathcal{AS}$-set of $\mathbb{P}^r(\R)$ and the restriction $\varrho_S : S \rightarrow S/G$ of $\pi_{\R^d}$ is a (continuous) $\mathcal{AS}$-map, i.e. a map with $\mathcal{AS}$-graph. This $\mathcal{AS}$-quotient (definition \ref{deasquotientoffreegassets}) is functorial with respect to equivariant continuous maps with $\mathcal{AS}$-graph (theorem \ref{theofunctorialityasquotientfreegas}) and compatible with equivariant closed inclusions (proposition \ref{propasquotfreeassetspreservesclosedincl}).

Let $X$ be a $G$-real algebraic set of a real affine space $\R^n$. Though the action of $G$ on~$X$ is linearizable by proposition \ref{proplinearizgrealalgset} (and even orthogonalizable since $G$ is compact: see proposition~\ref{proporthogonalizationactioncompactrealalggroup}), $X$ is not a free $G$-$\mathcal{AS}$-set of $\mathbb{P}^n(\R)$ in general. However, if $S$ is a free $G$-$\mathcal{AS}$-set, the Cartesian product $X \times S$ (equipped with the diagonal action of $G$) is a free $G$-$\mathcal{AS}$-set as well (proposition \ref{propproductwithfreegassetisgasset}, see also proposition \ref{propgralgsetisgasset}) and we can consider the $\mathcal{AS}$-quotient $\varrho_{X \times S} : X \times S \rightarrow (X \times S)/G$, that we call the \emph{relative $\mathcal{AS}$-quotient of $X$ with respect to $S$} (definition \ref{derelasquotientwithrespecttoafreegasset}). This relative quotient by $G$ is a bifunctor with respect, in both entries, to equivariant continuous $\mathcal{AS}$-maps.

The latter construction will be the starting point of a subsequent paper where we will generalize the second part of \cite{FP-QIAS}. More precisely, if $G$ is a compact real algebraic group, we will apply an, in some sense, additive homological functor (namely the Nash-constructible filtered chain complex defined and studied by C. McCrory and A. Parusiński in section 3 of \cite{MCP}) to a sequence of relative $\mathcal{AS}$-quotients by $G$ associated to any $G$-$\mathcal{AS}$-set (definition \ref{degassetfreegasset}). We will then deduce additive numerical invariants of $G$-$\mathcal{AS}$-sets, in particular of $G$-real algebraic sets and $G$-stable Zariski-open subsets of $G$-real algebraic sets. Because the relative $\mathcal{AS}$-quotient of $G$ is (bi)functorial with respect to equivariant continuous $\mathcal{AS}$-maps, the invariants that we will define will be invariant with respect to equivariant homeomorphisms with $\mathcal{AS}$-graph.
\\

Let us finally introduce notations that we will use all along this text. 

We will first denote by $\R[x_1,\ldots,x_n]$, resp. $\C[x_1,\ldots,x_n]$, the $\R$-algebra, resp. $\C$-algebra, of polynomials in $n$ indeterminates $x_1,\ldots,x_n$ with real, resp. complex, coefficients. If $E$ is any subset of the ring of functions from $\R^n$ to $\R$, resp. $\C^n$ to $\C$, we denote by $V(E)$, resp.~$\mathsf{V}(E)$, the subset of points $x$ of $\R^n$, resp. $\C^n$, satisfying $f(x) = 0$ for all $f \in E$. If $S$ is any subset of $\R^n$, resp. $\C^n$, we denote by $I(S)$, resp. $\mathsf{I}(S)$, the ideal of polynomials $P$ of the ring $\R[x_1,\ldots,x_n]$, resp. $\C[x_1,\ldots,x_n]$, such that $P(x) = 0$ for any $x \in S$. The Zariski closure of $S$, i.e. the smallest real, resp. complex, algebraic set of $\R^n$, resp. $\C^n$, containing $S$, will be denoted by~$\overline{S}^{\mathcal{Z}}$, resp.~$\overline{S}^{\mathfrak{Z}}$. 

In addition, if $\mathbb{K}$ is a commutative field, the zero vector of a $\mathbb{K}$-vector space $V$ will be in general denoted by $\bf{0}$, while $\mathcal{L}(V)$, resp. $\mathcal{GL}(V)$, will denote the $\mathbb{K}$-vector space of linear endomorphisms, resp. automorphisms, of $V$. The zero, resp. identity, matrix of $\M_n(\mathbb{K})$ will be denoted by $0_n$, resp. $I_n$, and, when this is relevant, we will identify any square matrix $\left(a_{ij}\right)_{1\leq i,j\leq n}$ of $\M_n(\mathbb{K})$ as a vector of $\mathbb{K}^{n^2}$. Moreover, if $i,j \in \{1,\ldots,n\}$, we will use the Kronecker delta symbol $\delta_{ij}$ which is equal to $1$ if $i=j$ and $0$ otherwise, and if $\lambda_1,\ldots,\lambda_n$ are scalars, we will write ${\rm Diag}(\lambda_1,\ldots,\lambda_n)$ to denote the diagonal matrix of $\M_n(\mathbb{K})$ with diagonal coefficients $\lambda_1,\ldots,\lambda_n$ in this order.

To conclude, let us also mention that the $\mathcal{C}^{\infty}$, analytic or Nash manifolds we consider are all Hausdorff, that a smooth map means a $\mathcal{C}^{\infty}$ map and that a (sub)manifold means a~$\mathcal{C}^{\infty}$~(sub)manifold. In addition, a regular map between Zariski-open subsets of algebraic sets (i.e. set-theoritic differences of algebraic sets of affine spaces) means a map whose coordinate functions are well-defined quotients of polynomial functions, and a biregular map is a bijective regular map with regular inverse (contrary to the real case, a regular map defined on a complex algebraic set is always the restriction of a polynomial map: see theorem \ref{theoregularfunctiononcompalgsetispoly} of the appendix).
\\

{\bf Acknowledgements.} The author is grateful to many persons for fruitful discussions and comments, including Olivier Benoist, Jean-Baptiste Campesato, Michel Coste, Goulwen Fichou, Julien Grivaux, Francesco Guaraldo, Jimmy Guillou, Evelyne Hubert, Johannes Huisman, Tomohiro Kawakami, Boris Kolev, Jean-Philippe Monnier and Marc Olive.

\section{Real algebraic groups with polynomial structure maps} \label{sectionralggp}

\subsection{First definitions and properties of real algebraic groups} \label{subsecfirstdefproprealalggr}

In the first part of this text, we will be considering real algebraic groups with \emph{polynomial} group structure i.e. whose structure maps are given by restriction of polynomial maps (and not maps that are just regular). Our motivation to make this hypothesis is our subsequent use of Hilbert's finiteness theorem (theorem \ref{theohilbertfiniteness} below) in order to consider algebro-geometric quotient of real algebraic sets acted by compact real algebraic group action. Indeed, Hilbert's finiteness theorem requires that both the involved group and its action are given by polynomial maps. 

Let $m$ be a positive integer.

\begin{de} \label{defrealalggroup} Let $G$ be a real algebraic subset of $\mathbb{R}^m$ and let $\mu : G \times G \rightarrow G$ be a binary operation such that the pair $(G,\mu)$ is a group. Let us furthermore denote by $\omega$ the inverse map $G \rightarrow G~;~g \mapsto g^{-1}$.

We say that $G$ is a \emph{real algebraic group} if both maps $\mu$ and $\omega$ are (restrictions of) polynomial maps.
\end{de}

If $G$ is a real algebraic group and $H$ is a real algebraic subset of $\mathbb{R}^m$ as well as a subgroup of $G$, then we say that $H$ is a \emph{real algebraic subgroup} of $G$ (in this case, $H$ is a real algebraic group too). As first examples, let us mention the groups that will follow us along this paper: let $n$ be a positive integer.

\begin{exs} \label{exspolyrealalggroups} 
~
\begin{enumerate}
		\item Consider the real orthogonal group $\OO_n(\mathbb{R})$ of orthogonal matrices of~$\M_n(\R)$ as being the real algebraic subset of points $\left(a_{ij}\right)_{1\leq i,j \leq n}$ of $\mathbb{R}^{n^2}$ such that $$\sum_{k =1}^n a_{ki} a_{kj} = \delta_{ij}$$
for any $i,j \in \{1,\ldots,n\}$. The group $\OO_n(\mathbb{R})$ is a real algebraic subgroup of the real algebraic group $\Omega_n(\R) := \left\{A \in \mathbb{R}^{n^2}~|~\det(A) = \pm 1\right\}$ equipped with matrix multiplication (if $A \in \Omega_n(\R)$, the inverse of $A$ is given by $\det(A)$ times the transpose of the cofactor matrix of $A$, which is a polynomial map in the coefficients of $A$). In addition, the special orthogonal group $\SO_n(\R)$ of real orthogonal matrices of size $n$ and determinant $1$ is a real algebraic subgroup of $\OO_n(\mathbb{R})$. Recall also that any finite group can be realized as a group of permutation matrices and can then be considered as a real algebraic subgroup of an orthogonal group. 
		\item In view of our definition, the general linear group $\GL_n(\R)$ is not a real algebraic set in itself: it is a Zariski open subset of $\M_n(\R)$, considered as $\R^{n^2}$, biregularly isomorphic to the real algebraic group $\left\{(A,y) \in \R^{n^2+1}~|~y \, \det(A) = 1\right\}$. Nevertheless, any \emph{real algebraic group} which is a subgroup of a general linear group $\GL_n(\R)$, considered as a subset of $\R^{n^2}$, will be called a \emph{linear real algebraic group}.
		\item The real algebraic set $K := \left\{(x,y) \in \mathbb{R}^2~|~xy=1\right\}$ ($K$ is biregularly isomorphic to the Zariski open subset $\mathbb{R}^*$ of $\mathbb{R}$), together with the polynomial binary operation which associates to any points $(x,y)$ and $(x',y')$ of $K$ the couple $(xx',yy')$, is a real algebraic group (the inverse map of the group is $(x,y) \in K \mapsto (y,x) \in K$).
		\item If $G \subset \mathbb{R}^{m}$ and $H \subset \mathbb{R}^{M}$ are real algebraic groups, the direct product group $G \times H \subset \mathbb{R}^{m+M}$ is a real algebraic group.
	\end{enumerate}
\end{exs}

Before stating some basic geometrical results on real algebraic groups, let us make precise the notions of morphism and isomorphism that we will be considering on these objects:

\begin{de} Let $G \subset \mathbb{R}^m$ and $H \subset \mathbb{R}^{M}$ be real algebraic groups and let $\varphi : G \rightarrow H$ be a map. We say that $\varphi$ is a \emph{morphism of real algebraic groups} if $\varphi$ is the restriction of a polynomial map $\mathbb{R}^m \rightarrow \mathbb{R}^M$ as well as a group homomorphism. 

If $\varphi$ is a bijective morphism of real algebraic groups and its inverse $\varphi^{-1}$ is also a morphism of real algebraic groups (i.e. the restriction of a polynomial map $\mathbb{R}^M \rightarrow \mathbb{R}^m$, since the inverse map of a bijective group homomorphism is a group homomorphism as well), we will say that $\varphi$ is an \emph{isomorphism of real algebraic groups} and that $G$ and $H$ are \emph{isomorphic as real algebraic groups}. 
\end{de}

In the following, a bijective restriction of a polynomial map between real algebraic sets whose inverse is also the restriction of a polynomial map will also be called a \emph{polynomial isomorphism} (and we will say that the source and target of a polynomial isomorphism are \emph{polynomially isomorphic} real algebraic sets): an isomorphism of real algebraic groups is a group homomorphism between real algebraic groups which is a polynomial isomorphism.  

\begin{exs} \label{exmorphrealalggroups} \begin{enumerate}
	\item The unit circle $S^1 := \left\{(x,y) \in \mathbb{R}^2~|~x^2+y^2 = 1\right\}$ of $\mathbb{R}^2$, equipped with the binary operation which associates to any points $(x,y)$ and $(x',y')$ of $S^1$ the couple $(xx'-yy',xy'+x'y)$, is a real algebraic group (the inverse map of the group is $(x,y) \in S^1 \mapsto (x,-y) \in S^1$), and the map 
	$$\begin{array}{ccc}S^1 & \rightarrow & \SO_2(\mathbb{R})\\(x,y) & \mapsto & \begin{pmatrix}x&-y\\y&x\end{pmatrix}\end{array}$$
is an isomorphism of real algebraic groups.  
	\item Let $G \subset \mathbb{R}^m$ and $H \subset \mathbb{R}^M$ be {\it finite} real algebraic groups of same cardinality $N \in \Nstar$ and suppose that $G$ and $H$ are isomorphic as groups. Then $G$ and $H$ are isomorphic as real algebraic groups. Indeed, if $\varphi : G \rightarrow H$ is an isomorphism of groups and if we write $G = \{g_1,\ldots,g_N\}$ and $H = \{h_1,\ldots,h_N\}$ so that, for all $k \in \{1,\ldots,N\}$, $h_k = \varphi(g_k)$, then~$\varphi$ is the restriction of the polynomial map 
$$\begin{array}{ccc}\mathbb{R}^m & \rightarrow & \mathbb{R}^M\\g & \mapsto & \displaystyle{\sum_{k=1}^N\left(\prod_{1\leq l\leq N, \, l \neq k} \frac{\|g-g_l\|^2}{\|g_k - g_l\|^2}\right) h_k.}\end{array}$$
and is therefore an isomorphism of real algebraic groups (the inverse map $\varphi^{-1}$ is the restriction of an analogous polynomial map).
\end{enumerate}
\end{exs}

\begin{rems} \label{remsimagemorphrealalggr}
~	
	\begin{enumerate}
		\item If $\varphi : G \rightarrow H$ is a morphism of real algebraic groups, then the kernel of~$\varphi$, or more generally the inverse image by $\varphi$ of any real algebraic subgroup of $H$, is always a real algebraic subgroup of $G$, since the inverse image of a real algebraic set by a polynomial map is a real algebraic set. However, the image of a morphism of real algebraic groups is not a real algebraic group in general, contrary to the complex case (see for instance Proposition 21.2.4 of \cite{TY}). Indeed, consider for example the real algebraic group $K := \left\{(x,y) \in \mathbb{R}^2~|~xy=1\right\}$ (see example \ref{exspolyrealalggroups} 3) and the endomorphism $\psi : (x,y) \in K \mapsto (x^2,y^2) \in K$ of $K$: the image of $\psi$ is the subgroup $\left\{(x,y) \in \mathbb{R}^2~|~xy=1, \, x > 0, \, y > 0\right\}$ of $K$, which is not a real algebraic set.   

		\item Contrary to the complex case (see Proposition 21.2.6 of \cite{TY}), a bijective morphism of real algebraic groups is not an isomorphism of real algebraic groups in general: if we consider again the above real algebraic group $K$, the inverse map of the bijective endomorphism of real algebraic group $(x,y) \in K \mapsto (x^3,y^3) \in K$ is not the restriction of a polynomial map.
	\end{enumerate}
\end{rems}

In the following, let $G \subset \mathbb{R}^m$ be a real algebraic group and denote by $\mu$ its polynomial structure map, by $\omega$ its associated polynomial inverse map and by $e$ its identity element. We are going to state two general geometrical properties of real algebraic groups, which are analogous to the complex case. We will give their proofs, for the sake of self-containedness.

\begin{prop} \label{propcosetsirredcomprealalggroup} There is a unique irreducible component $G_0$ of the real algebraic set $G$ containing $e$. Furthermore, $G_0$ is a real algebraic subgroup of $G$, which is normal and of finite index as a subgroup of $G$, and the cosets of $G$ modulo $G_0$ are the irreducible components of the real algebraic set $G$. 
\end{prop}
  
\begin{proof} We follow the proofs of Lemma 21.1.4 and Theorem 21.1.6 of \cite{TY}. Let first $G_1$ and~$G_2$ be two irreducible components of the real algebraic set $G$ containing $e$ and let us show that $G_1 = G_2$, following the proof of Lemma 21.1.4 of \cite{TY}. The product $G_1 \times G_2$ is an irreducible real algebraic subset of $\mathbb{R}^m \times \mathbb{R}^m$, so that $\mu(G_1 \times G_2)$ is an irreducible subset of~$G$ equipped with its Zariski topology since the polynomial map $\mu$ is in particular continuous with respect to Zariski topology (see Proposition 1.1.7 (i) of \cite{TY}). But $\mu(G_1 \times G_2)$ contains the irreducible components $G_1 = \mu(G_1 \times \{e\})$ and $G_2 = \mu(\{e\} \times G_2)$ of $G$, and therefore, by maximality, $G_1 = \mu(G_1 \times G_2) = G_2$.

So now denote by $G_0$ the unique irreducible component of $G$ containing $e$. By the above arguments, $\mu(G_0 \times G_0) = G_0$ and, furthermore, $\omega(G_0) = G_0$ because, since $\omega : G \rightarrow G$ is a polynomial isomorphism and $\omega(e) = e$, $\omega(G_0)$ is an irreducible component of $G$ containing $e$. As a consequence, $G_0$ is a real algebraic subgroup of $G$.

In order to show that $G_0$ is a normal subgroup of $G$, if $g \in G$, denote by $l_g$ the left translation map $h \in G \mapsto \mu(g,h) \in G$ and by $r_g$ the right translation map $h \in G \mapsto \mu(h,g) \in G$: $l_g$ and~$r_g$ are both polynomial isomorphisms (of respective inverse maps $l_{\omega(g)}$ and $r_{\omega(g)}$). For all~$g \in G$, the image $l_g \circ r_{\omega(g)}(G_0)$ is therefore an irreducible component of $G$ containing $e$ so that $l_g \circ r_{\omega(g)}(G_0) = G_0$, and $G_0$ is then a normal subgroup of $G$. 

Finally, the group $G$ is the disjoint union of its cosets $l_g(G_0)$, $g \in G$, but, for all $g \in G$, $l_g(G_0)$ is an irreducible component of the real algebraic set $G$ containing $g$ and any real algebraic set is the finite union of its irreducible components, hence the last statements of the proposition. 
\end{proof}  

\begin{rems} 
~	
	\begin{enumerate}
		\item Any irreducible real algebraic set is Zariski-connected but the converse is not true. However, by previous proposition \ref{propcosetsirredcomprealalggroup}, any real algebraic group is irreducible if and only if it is Zariski-connected (because the irreducible components of a real algebraic group do not intersect). 
		\item A real algebraic group which is connected with respect to Euclidean topology is necessarily irreducible as a real algebraic set (since the irreducible components of a real algebraic group do not intersect and any real algebraic set is a closed set with respect to Euclidean topology), while the converse is not true, as illustrated by the real algebraic group $\left\{(x,y) \in \mathbb{R}^2~|~xy=1\right\}$ (example \ref{exspolyrealalggroups} 3).
	\end{enumerate}
\end{rems}

We remarked that the irreducible components of the real algebraic group $G$ do not intersect (notice also that they are all polynomially isomorphic to one another via the translation maps of $G$). Actually, $G$ does not have any singular point :

\begin{lem} \label{lemrealalggroupnonsing}
The real algebraic group $G$ is a nonsingular real algebraic set. 
\end{lem}

\begin{proof} Let $g$ be a nonsingular point of $G$ (such a point always exists by, for instance, Proposition 3.3.14 of \cite{BCR}) and let $h \in G$. We have $h = l_{\mu(h,\omega(g))}(g)$ and $l_{\mu(h,\omega(g))}$ is a polynomial isomorphism (in particular, a biregular isomorphism), so that $h$ is a nonsingular point of $G$ as well.
\end{proof}

In particular, $G$ is a Nash submanifold of $\mathbb{R}^m$ (see Proposition 3.3.11 and Definition 2.9.9 of \cite{BCR}) and then, since the polynomial structure maps $\mu$ and $\omega$ of the group $G$ are in particular Nash maps (i.e. $\mathcal{C}^{\infty}$ maps with semialgebraic graphs), $G$ is an affine Nash group (see for instance \cite{Gua20}). This implies that $G$ is a Lie group whose dimension (the dimension of $G$ as a $\mathcal{C}^{\infty}$ submanifold of~$\R^m$ coincides with its dimension as a real algebraic variety: see Proposition 2.8.14 of \cite{BCR}) can be computed as the dimension, as an $\mathbb{R}$-vector space, of the associated Lie algebra. For instance, the Lie groups $\OO_n(\R)$ and $\SO_n(\R)$ have the same Lie algebra $\{M \in \M_n(\R)~|~{}^t\!M = - M\}$, so that the real algebraic groups $\OO_n(\R)$ and $\SO_n(\R)$ have dimension $\frac{n(n-1)}{2}$.

\subsection{Semialgebraic groups with polynomial structure maps}

In this part, motivated by the example of the images of morphisms of real algebraic groups (cf. remark \ref{remsimagemorphrealalggr} 1), we want to extend to the larger context of semialgebraic subgroups of real algebraic groups:

\begin{de} \label{defsemialggroup} Let $G$ be a semialgebraic subset of $\mathbb{R}^m$ and let $\mu : G \times G \rightarrow G$ be a binary operation such that the pair $(G,\mu)$ is a group. If $\omega$ denotes furthermore the inverse map $G \rightarrow G~;~g \mapsto g^{-1}$, we say that $G$ is a \emph{semialgebraic group} if both maps $\mu$ and $\omega$ are restrictions of polynomial maps.
\end{de}

\begin{exs} \label{exssemialggps}
~	
	\begin{enumerate}
		\item As first examples, real algebraic groups are of course semialgebraic groups, as well as the image of any morphism of real algebraic groups (since the image of a semialgebraic set by a polynomial map is a semialgebraic set: see Proposition 2.2.7 of \cite{BCR}): for instance, the image $\left\{(x,y) \in \mathbb{R}^2~|~xy=1, \, x > 0, \, y > 0\right\}$ of the endomorphism of real algebraic group $(x,y) \in K \mapsto (x^2,y^2) \in K$, where $K = \left\{(x,y) \in \mathbb{R}^2~|~xy=1\right\}$, is a semialgebraic group.
		\item Remark that, more generally, the kernel and the image of a map between semialgebraic groups which is both a group homomorphism and a semialgebraic map (i.e. which has a semialgebraic graph) are semialgebraic groups, since inverse image or image of any semialgebraic subset by a semialgebraic map is always a semialgebraic set (see again Proposition 2.2.7 of \cite{BCR}). 
	\end{enumerate}
\end{exs}

Let $G \subset \mathbb{R}^m$ be a semialgebraic group. As we suggested in the introduction of this subsection, the semialgebraic group structure (as defined in the latter definition) of $G$ is actually induced from a real algebraic group structure on the Zariski closure $\overline{G}^{\mathcal{Z}}$ of $G$.

Recall that, if $S$ is a subset of $\mathbb{R}^m$, then $\overline{S}^{\mathcal{Z}} = V(I(S))$ (i.e. $\overline{S}^{\mathcal{Z}}$ is the set of points of $\mathbb{R}^m$ at which vanishes any polynomial of $\mathbb{R}[x_1,\ldots,x_m]$ that vanishes at all points of $S$) and that, if~$T$ is a subset of $\mathbb{R}^M$ and $f : S \rightarrow T$ is the restriction of a polynomial map $\widetilde{f} : \mathbb{R}^m \rightarrow \mathbb{R}^M$, then $\widetilde{f}\left(\overline{S}^{\mathcal{Z}}\right) \subset \overline{T}^{\mathcal{Z}}$ so that we can extend~$f$ into the restriction $\overline{f}^{\mathcal{Z}} : \overline{S}^{\mathcal{Z}} \rightarrow \overline{T}^{\mathcal{Z}}$ of $\widetilde{f}$ (this extension is unique in the sense that two polynomial maps coinciding on $S$ coincide on $\overline{S}^{\mathcal{Z}}$).

Now consider the semialgebraic group $G \subset \mathbb{R}^m$, with structure map $\mu : G \times G \subset \mathbb{R}^m \times \mathbb{R}^m \rightarrow G \subset \mathbb{R}^m$, inverse map $\omega : G \subset \mathbb{R}^m \rightarrow G \subset \mathbb{R}^m$ and neutral element $e$. Because they are polynomial, we can extend the maps $\mu$ and $\omega$ to the Zariski closure $\overline{G}^{\mathcal{Z}}$ of $G$ in $\mathbb{R}^m$ in order to make $\overline{G}^{\mathcal{Z}}$ into a real algebraic group:  

\begin{lem} \label{lemzariskiclrealalggroup}The couple $\left(\overline{G}^{\mathcal{Z}}, \overline{\mu}^{\mathcal{Z}}\right)$ is a real algebraic group.
\end{lem}

\begin{proof} We need to show that $\left(\overline{G}^{\mathcal{Z}}, \overline{\mu}^{\mathcal{Z}}\right)$ is a group with inverse map $\overline{\omega}^{\mathcal{Z}}$ and neutral element~$e$. Suppose that $\mu$, resp. $\omega$, is the restriction of a polynomial map $\widetilde{\mu} : \mathbb{R}^m \times \mathbb{R}^m \rightarrow \mathbb{R}^m$, resp. $\widetilde{\omega} : \mathbb{R}^m \rightarrow \mathbb{R}^m$. First, since we have the polynomial equality $\widetilde{\mu}( \cdot , \widetilde{\mu}( \cdot , \cdot )) = \widetilde{\mu}(\widetilde{\mu}(\cdot , \cdot ),\cdot)$ on $G \times G \times G$, we have the same equality on $\overline{G \times G \times G}^{\mathcal{Z}} = \overline{G}^{\mathcal{Z}} \times \overline{G}^{\mathcal{Z}} \times \overline{G}^{\mathcal{Z}}$, i.e. the map $\overline{\mu}^{\mathcal{Z}}$ is associative. Similarly, the polynomial equalities $\widetilde{\mu}( \cdot, e) = \widetilde{\mu}(e, \cdot) = \cdot$ and  $\widetilde{\mu}(\cdot , \widetilde{\omega}(\cdot)) = \widetilde{\mu}(\widetilde{\omega}(\cdot), \cdot) = e$ on $G$ extend to $\overline{G}^{\mathcal{Z}}$. As a consequence, the couple $\left(\overline{G}^{\mathcal{Z}}, \overline{\mu}^{\mathcal{Z}}\right)$ is a group, and therefore a real algebraic group. 
\end{proof}

\begin{rem} \label{remfunctzariskicl} We wanted to place emphasis on how the polynomial group structure on $G$ extends to $\overline{G}^{\mathcal{Z}}$, but this is actually just an illustration of the functoriality of Zariski closure: with the above notations, $\overline{{\rm Id}_{S}}^{\mathcal{Z}} = {\rm Id}_{\overline{S}^{\mathcal{Z}}}$ and, if $S_1 \subset \mathbb{R}^{m_1}$, $S_2 \subset \mathbb{R}^{m_2}$, $S_3 \subset \mathbb{R}^{m_3}$ and $f : S_1 \rightarrow S_2$, $h : S_2 \rightarrow S_3$ are restrictions of polynomial maps, then $\overline{h \circ f}^{\mathcal{Z}} = \overline{h}^{\mathcal{Z}} \circ \overline{f}^{\mathcal{Z}}$
\end{rem}

As a conclusion, a subset of $\R^m$ is a semialgebraic group if and only if it is a subgroup of a real algebraic group which is semialgebraic.

Let us carry on this discussion. Actually, the semialgebraic group $G$ is not any semialgebraic subset of its Zariski closure $\overline{G}^{\mathcal{Z}}$: we will show that $G$ is necessarily a union of connected (for the Euclidean topology on $\mathbb{R}^m$) components of $\overline{G}^{\mathcal{Z}}$. Let us first state the following ``semialgebraic version'' of proposition \ref{propcosetsirredcomprealalggroup}: recall that any semialgebraic set of $\mathbb{R}^m$ has a finite number of connected components, which are all semialgebraic (cf. Theorem 2.4.5 of \cite{BCR}). 

\begin{lem} \label{lempropcosetsconncompsemialggroup} There is a unique connected component $G_c$ of the semialgebraic set $G$ containing~$e$. Furthermore, $G_c$ is a semialgebraic subgroup of $G$, which is normal and of finite index as a subgroup of $G$, and the cosets of $G$ modulo $G_c$ are the (semialgebraic) connected components of the semialgebraic set $G$.
\end{lem}

\begin{proof} Just adapt the proof of proposition \ref{propcosetsirredcomprealalggroup} by replacing the word ``irreducible'' with the word ``connected'' (with respect to the Euclidean topology of $\mathbb{R}^m$): a polynomial map, being continuous, carries connected sets to connected sets and a polynomial isomorphism, being a homeomorphism, carries connected components to connected components.
\end{proof}

\begin{rem} Apart from the semialgebraicity and the finiteness of the index of $G_c$ as a subgroup of $G$, the other properties of lemma \ref{lempropcosetsconncompsemialggroup} just come from the fact that a semialgebraic group is a topological group.
\end{rem}

We are going to show that the connected component $G_c$ of $G$ is also a connected component of $\overline{G}^{\mathcal{Z}}$, which will lead us to the following announced result:

\begin{prop} \label{propsemialggroupunionofccofzarclos} The semialgebraic group $G$ is a union of connected components (with respect to Euclidean topology) of $\overline{G}^{\mathcal{Z}}$.
\end{prop}

\begin{proof} Notice that $d:= \dim G_c = \dim G = \dim \overline{G}^{\mathcal{Z}}$ and consider a stratification of $\overline{G}^{\mathcal{Z}}$ adapted to $G_c$ (we refer to Corollary 3.8 of \cite{CosteSA}): $\overline{G}^{\mathcal{Z}}$ is the finite disjoint union of semialgebraic subsets $C_1,\ldots, C_N$ such that $G_c$ is a union of some $C_k$'s, $k \in \{1,\ldots,N\}$, and, if $k \in \{1,\ldots,N\}$, the Euclidean closure $\overline{C_k}$ of $C_k$ in $\overline{G}^{\mathcal{Z}}$ is the (disjoint) union of $C_k$ and of some $C_l$'s, $l \in \{1,\ldots, N\} \setminus \{k\}$, of smaller dimension.

Let $k \in \{1,\ldots,N\}$ such that $C_k \subset G_c$ and $\dim C_k = d$, and denote $C:=C_k$. We have $C = \overline{G}^{\mathcal{Z}} \setminus \left(\bigcup_{1\leq l \leq N,\, l \neq k} \overline{C_l}\right)$ (since for any $l \in \{1,\ldots, N\} \setminus \{k\}$, $\overline{C_l} \setminus C_l$ is a union of strata of dimension smaller than $d = \dim C$) and then $C$ is an (Euclidean) open subset of $\overline{G}^{\mathcal{Z}}$ contained in $G_c$.

Now, when $g \in G$, denote, as in the previous subsection, the left translation map $h \in G \mapsto \mu(g,h) \in G$ by $l_g$ and remark that $G_c = \bigcup_{g \in G_c} l_g(C)$ ($G_c$ is a subgroup of $G$), so that $G_c$ is itself an open subset of $\overline{G}^{\mathcal{Z}}$ (for all $g \in G$, $l_g$ is a polynomial isomorphism in particular a homeomorphism) as well as the cosets $l_g(G_c)$, $g \in \overline{G}^{\mathcal{Z}}$, of the group $\overline{G}^{\mathcal{Z}}$ modulo its subgroup~$G_c$, which form a partition of $\overline{G}^{\mathcal{Z}}$: writing $G_c = \overline{G}^{\mathcal{Z}} \setminus \left(\bigcup_{g \in \overline{G}^{\mathcal{Z}} \setminus \{e\}} l_g(G_c)\right)$, we deduce that $G_c$ is also closed in $\overline{G}^{\mathcal{Z}}$. As a consequence, since $G_c$ is furthermore connected, $G_c$ is a connected component of $\overline{G}^{\mathcal{Z}}$.

Finally, since the other connected components of $G$ are cosets $l_g(G_c)$, with $g \in G$, of $G$ modulo $G_c$, they are also connected components of $\overline{G}^{\mathcal{Z}}$.
\end{proof}

\begin{rem} The second part of the above argument is purely topological and works as soon as we are considering a subgroup of a topological group such that the former contains an open subset of the latter.
\end{rem}

We end this section with the following property :

\begin{prop} \label{propconnectedcompirredcompneutralelt} With the notations of proposition \ref{propcosetsirredcomprealalggroup} and lemma \ref{lempropcosetsconncompsemialggroup}, we have the equality $\overline{G_c}^{\mathcal{Z}} = \left(\overline{G}^{\mathcal{Z}}\right)_0$.
\end{prop}

\begin{proof} The identity component $G_c$ of $G$, being connected, is contained in only one irreducible component of $\overline{G}^{\mathcal{Z}}$ (because the irreducible components of the real algebraic group $\overline{G}^{\mathcal{Z}}$, which are closed for the Euclidean topology, do not intersect), which has to be $\left(\overline{G}^{\mathcal{Z}}\right)_0$. As a consequence, $\overline{G_c}^{\mathcal{Z}} \subset \left(\overline{G}^{\mathcal{Z}}\right)_0$. On the other hand, we have
$$\dim \overline{G_c}^{\mathcal{Z}} =  \dim G_c = \dim \overline{G}^{\mathcal{Z}} = \dim \left(\overline{G}^{\mathcal{Z}}\right)_0$$
so that $\overline{G_c}^{\mathcal{Z}} = \left(\overline{G}^{\mathcal{Z}}\right)_0$ (the real algebraic set $\left(\overline{G}^{\mathcal{Z}}\right)_0$ is irreducible).
\end{proof} 

\subsection{Complexification of real algebraic groups} \label{subseccomplexificationralgps}

Considering real algebraic groups with polynomial group structure allows to consider the complexification of the latter, as we are going to explain below. A real algebraic group structure can then also be seen as being induced from a complex algebraic group structure given by polynomial maps with real coefficients.

Recall that, if $X$ is a real algebraic subset of $\mathbb{R}^m \subset \mathbb{C}^m$, the {\it complexification of $X$}, denoted by $X_{\mathbb{C}}$, is the complex Zariski closure $\overline{X}^{\mathfrak{Z}}$ of $X$ in $\mathbb{C}^m$ i.e. the set $\mathsf{V}(\mathsf{I}(X))$ of points of $\mathbb{C}^m$ at which vanishes any polynomial of $\mathbb{C}[x_1,\ldots,x_m]$ that vanishes at all points of $X$. We have $X = X_{\mathbb{C}} \cap \mathbb{R}^d$ (since $X = V(I(X))$ and, if $x$ is any point of $X_{\mathbb{C}}$ with real coordinates, any polynomial of $\mathbb{R}[x_1,\ldots,x_m]$ vanishing at all points of $X$ vanishes at $x$) and, if $f : X \rightarrow Y$ is the restriction of a polynomial map $\widetilde{f} : \mathbb{R}^m \rightarrow \mathbb{R}^M$ and if $\widetilde{f}_{\mathbb{C}} : \mathbb{C}^m \rightarrow \mathbb{C}^M$ denotes map whose coordinate functions are given by the same polynomial expressions as the coordinate functions of $\widetilde{f}$, then $\widetilde{f}_{\mathbb{C}}\left(X_{\mathbb{C}}\right) \subset Y_{\mathbb{C}}$, so that we can extend~$f$ into the restriction $f_{\mathbb{C}} : X_{\mathbb{C}} \rightarrow Y_{\mathbb{C}}$ of $\widetilde{f}_{\mathbb{C}}$ (such an extension is unique: two polynomial maps coinciding on $X$ coincide on $X_{\mathbb{C}}$). Let us also recall that the complexification of real algebraic sets and (restrictions of) polynomial maps between real algebraic sets is functorial.

Now, let $G \subset \mathbb{R}^m$ be a real algebraic group with structure map $\mu : G \times G \rightarrow G$, inverse map $\omega : G \subset \mathbb{R}^m \rightarrow G \subset \mathbb{R}^m$ and neutral element $e$, and consider the respective complexifications $G_{\mathbb{C}} \subset \mathbb{C}^m$ and $\mu_{\mathbb{C}} : G_{\mathbb{C}} \times G_{\mathbb{C}} \rightarrow G_{\mathbb{C}}$ of $G$ and $\mu$.

\begin{lem} The couple $\left(G_{\mathbb{C}},\mu_{\mathbb{C}}\right)$ is a complex algebraic group, i.e. a complex algebraic set equipped with a polynomial (or, equivalently, regular) group structure, with inverse map $\omega_{\mathbb{C}}$ and neutral element $e$. 
\end{lem}

\begin{proof} The proof is analogous to the proof of lemma \ref{lemzariskiclrealalggroup} (see also remark \ref{remfunctzariskicl}).
\end{proof}

\begin{rems} \label{remgeneralcomplexalgsetscomplexification}
~	
	\begin{enumerate}
		\item By the Hilbert's Nullstellensatz, a regular function (i.e. well-defined quotient of polynomial maps) on a complex algebraic set is necessarily the restriction of a polynomial function (see for instance theorem \ref{theoregularfunctiononcompalgsetispoly} in the appendix). This is not the case for regular functions on real algebraic sets.
		\item With the above notations, if the ideal $I(X)$ of $\R[x_1,\ldots,x_m]$ is generated by polynomials $P_1,\ldots,P_k$ (and $X$ is then the set $V(P_1,\ldots,P_k)$ of points of $\R^m$ at which vanish the polynomials $P_1,\ldots,P_k$), then the ideal $\mathsf{I}(X)$ of $\C[x_1,\ldots,x_m]$ is also generated by the polynomials $P_1,\ldots,P_k$ (because a polynomial $P + i Q \in \C[x_1,\ldots,x_m]$, with $P,Q \in \R[x_1,\ldots,x_m],$ vanishes on $X$ if and only if $P$ and $Q$ vanish on $X$) and then the complexification $X_{\C}$ is the set $\mathsf{V}(P_1,\ldots,P_k)$ of points of $\C^m$ at which vanish the polynomials $P_1,\ldots,P_k$. In particular, $X_{\C}$ is defined by polynomial equations with real coefficients and can then be considered as an affine algebraic variety over $\R$. In our context, $G_{\C}$ can be considered as an affine algebraic group over $\R$. Notice also that if $x$ is a point of $X_{\C}$, then the point~$\overline{x}$ of $\C^m$ whose coordinates are the respective conjugates of the coordinates of $x$ belongs to~$X_{\C}$ as well.
		\item For general properties on complex algebraic sets and complex algebraic groups, we refer for instance to \cite{TY}.
	\end{enumerate}
\end{rems}

Before giving some examples of complexifications of real algebraic groups, let us state a property which is an application of the following fact: if $X$ a real algebraic subset of $\mathbb{R}^m$ and $X = X_1 \cup \cdots \cup X_N$ is its decomposition into irreducible components, then $X_{\mathbb{C}} = (X_1)_{\mathbb{C}}  \cup \cdots \cup (X_N)_{\mathbb{C}}$ is the decomposition of the complex algebraic set $X_{\mathbb{C}}$ into irreducible components (if $X$ is an irreducible real algebraic set and if $X_{\C}$ is a union of complex algebraic sets $Z_1$, $Z_2$ of $\C^m$, then $X = X_{\C} \cap \R^m = \left(Z_1 \cap \R^m\right) \cup \left(Z_2 \cap \R^m\right)$ so that, without loss of generality, $X = Z_1 \cap \R^m$ and then $X_{\C} \subset Z_1$).

\begin{lem} Let $\left(G_{\mathbb{C}}\right)_0$ denote the unique irreducible component of the complex algebraic group $G_{\mathbb{C}}$ containing $e$ (cf. Lemma 21.1.4 of \cite{TY}). Then $\left(G_0\right)_{\mathbb{C}} = \left(G_{\mathbb{C}}\right)_0$.
\end{lem}

\begin{proof} Since $G_0$ is an irreducible component of $G$, $\left(G_0\right)_{\mathbb{C}}$ is an irreducible component of $G_{\mathbb{C}}$: since $\left(G_0\right)_{\mathbb{C}}$ contains $e$, by uniqueness, we have $\left(G_0\right)_{\mathbb{C}} = \left(G_{\mathbb{C}}\right)_0$:
\end{proof}

\begin{exs} \label{excomplrealalggroup}
~
\begin{enumerate}
	\item Let us prove that the complexification of the real algebraic group $\SO_n(\R)$ is the complex algebraic group $\SO_n(\C)$. First remark that, since the real algebraic group $\SO_n(\R)$, resp. complex algebraic group $\SO_n(\C)$, is connected with respect to Euclidean topology in $\R^{n^2}$, resp. to Euclidean topology in $\R^{2 n^2}$ (see for instance \cite{WAY}), it is irreducible (because irreducible components of real, resp. complex, algebraic groups do not intersect). Since $\SO_n(\R) \subset \SO_n(\C)$, we therefore have an inclusion of irreducible complex algebraic sets $\left(\SO_n(\R)\right)_{\C} \subset \SO_n(\C)$. The equality is then given by the equality of dimensions $\dim \SO_n(\C) = \frac{n(n-1)}{2} = \dim \SO_n(\R) = \dim \left(\SO_n(\R)\right)_{\C}$ (Propositions 2.2.1 and 2.2.5 of~\cite{AK}).
	
	Now, considering a real orthogonal matrix $O \in \OO_n(\R)$ of determinant $-1$, the map which associates to any special orthogonal matrix $A \in \SO_n(\C)$ the matrix product $OA$ is a polynomial isomorphism (over $\R)$ from $\SO_n(\C)$ to the complex algebraic set $\OO_n^{-}(\C)$ of orthogonal matrices of $\OO_n(\C)$ of determinant $-1$, which restricts into a (real) polynomial isomorphism from $\SO_n(\R)$ to the real algebraic set $\OO_n^{-}(\R)$ of orthogonal matrices of~$\OO_n(\R)$ of determinant $-1$ (the former map is the complexification of the latter map). As a consequence, the complexification of $\OO_n^{-}(\R)$ is $\OO_n^{-}(\C)$ and the complexification of $\OO_n(\R) = \SO_n(\R) \cup \OO_n^{-}(\R)$ is $\SO_n(\C) \cup \OO_n^{-}(\C) = \OO_n(\C)$.
	\item Using the connectedness, with respect to Euclidean topology, of the special linear groups $\SL_n(\R)$ and $\SL_n(\C)$ (see \cite{WAY}), we can show in a same manner that the complexification of $\SL_n(\R)$ is $\SL_n(\C)$ (the dimensions of $\SL_n(\R)$ and $\SL_n(\C)$ are both equal to $n^2-1$) and that the complexification of the real algebraic group $\Omega_n(\R) = \left\{A \in \mathbb{R}^{n^2}~|~\det(A) = \pm 1\right\}$ is the complex algebraic group $\Omega_n(\C) := \left\{A \in \mathbb{C}^{n^2}~|~\det(A) = \pm 1\right\}$.
\end{enumerate}
\end{exs}

Notice also that the complexification of a linear real algebraic group of ${\rm M}_n(\R)$ is a linear complex algebraic group of ${\rm M}_n(\C)$, i.e. a complex algebraic group which a subgroup of the Zariski open subset $\GL_n(\C)$ of ${\rm M}_n(\C)$ (considered as $\C^{n^2}$), since the complexification of the matrix multiplication of ${\rm M}_n(\R)$ is the matrix multiplication of ${\rm M}_n(\C)$.
\\

Let us then consider the complexification of morphisms of real algebraic groups:

\begin{lem} \label{lemcomplexmorphrealalggrismorphcomplexalggru} Let $\varphi : G \subset \mathbb{R}^m \rightarrow H \subset \mathbb{R}^M$ be a morphism of real algebraic groups. Then the complexification $\varphi_{\C} : G_{\C} \subset \mathbb{C}^m \rightarrow H_{\C} \subset \mathbb{C}^M$ is a morphism of complex algebraic groups (i.e. the restriction of a polynomial map as well as a group homomorphism).
\end{lem}

\begin{proof} Let $\nu : H \times H \rightarrow H$ be the structure map of the real algebraic group $H$. Since $\varphi$ is a morphism of real algebraic groups, we have the polynomial equality $\varphi\left(\mu(\cdot, \cdot)\right) = \nu\left(\varphi(\cdot), \varphi(\cdot)\right)$ on $G \times G$ and then, by functoriality of complexification, the polynomial equality $\varphi_{\C}\left(\mu_{\C} (\cdot, \cdot)\right) =\nu_{\C} \left(\varphi_{\C} (\cdot), \varphi_{\C} (\cdot)\right)$ on $G_{\C} \times G_{\C}$.
\end{proof}

\begin{ex} Consider the isomorphism of real algebraic groups 
$$\phi : \begin{array}{ccc}S^1 \subset \R^2& \rightarrow & \SO_2(\mathbb{R}) \subset \R^4\\(x,y) & \mapsto & \begin{pmatrix}x&-y\\y&x\end{pmatrix}\end{array}$$ of example \ref{exmorphrealalggroups} 1. By previous lemma \ref{lemcomplexmorphrealalggrismorphcomplexalggru} and by functoriality of the complexification, the map 
$$\phi_{\C} : \begin{array}{ccc}\left(S^1\right)_{\C} \subset \C^2& \rightarrow & \SO_2(\mathbb{C}) \subset \C^4\\(x,y) & \mapsto & \begin{pmatrix}x&-y\\y&x\end{pmatrix}\end{array}$$
is an isomorphism of complex algebraic groups i.e. a bijective morphism of complex algebraic groups whose inverse is also a morphism of complex algebraic groups. In particular, $\left(S^1\right)_{\C} = S^1_{\C} := \left\{(x,y) \in \mathbb{C}^2~|~x^2+y^2 = 1\right\}$ (we could have proven this equality directly, using Theorem~4.5.1 of \cite{BCR} and remark \ref{remgeneralcomplexalgsetscomplexification} 2 above, since $S^1 = V(x^2+y^2 -1)$ and the sign of the polynomial $x^2+y^2-1 \in \R[x,y]$ changes on $\R^2$, so that $I(S^1)$ is generated by $x^2+y^2-1$ and then $\left(S^1\right)_{\C} = \mathsf{V}(x^2+y^2 -1) = S^1_{\C}$).
\end{ex}

As we already mentioned in remark \ref{remsimagemorphrealalggr}, a bijective morphism of complex algebraic groups is an isomorphism of complex algebraic groups (see Proposition 21.2.6 of \cite{TY}), while this is not the case for morphisms of real algebraic groups. On the other hand, a morphism of real algebraic groups whose complexification is bijective is an isomorphism of real algebraic groups, thanks to the following fact:

\begin{lem} \label{lempolyisomiffcomplexpolyisom} Let $f : X \subset \mathbb{R}^m \rightarrow Y \subset \mathbb{R}^M$ be a (restriction of a) polynomial map between real algebraic sets. Then $f$ is a polynomial isomorphism if and only if $f_{\mathbb{C}} : X_{\mathbb{C}} \subset \mathbb{C}^m \rightarrow Y_{\mathbb{C}} \subset \mathbb{C}^M$ is a polynomial isomorphism.
\end{lem}

\begin{proof} We give the proof of this equivalence for the sake of completeness. The direct implication is given by the functoriality of complexification. As for the converse implication, suppose that~$f_{\mathbb{C}}$ has a polynomial inverse map $g : Y_{\mathbb{C}} \subset \mathbb{C}^M \rightarrow X_{\mathbb{C}} \subset \mathbb{C}^m$. Take a point $y \in Y \subset  Y_{\mathbb{C}}$ and consider its image $x := g(y) \in X_{\mathbb{C}}$ by $g$, as well as the conjugated point $\overline{x} \in X_{\mathbb{C}}$ of $x$. Since $f_{\mathbb{C}}(x) = y$, because the coordinate functions of $f_{\mathbb{C}}$ are given by polynomials with real coefficients and $y$ has real coordinates, we have $f_{\mathbb{C}}(\overline{x}) = f_{\mathbb{C}}(x)$ and then $\overline{x} = x$ since the map~$f_{\mathbb{C}}$ is injective. As a consequence, $x = g(y) \in \mathbb{R}^m \cap X_{\mathbb{C}} = X$. In particular, the imaginary parts of the polynomial coordinate functions of $g$ vanish on $Y$ and then on $Y_{\mathbb{C}}$: the map $g$ is therefore the restriction of a polynomial map $\mathbb{C}^M \rightarrow \mathbb{C}^m$ given by polynomials with real coefficients. The restriction $Y \subset \mathbb{R}^M \rightarrow X \subset \mathbb{R}^m$ of $g$ is consequently a real polynomial inverse map for $f$.
\end{proof}

\begin{cor} \label{corisorealalggroupsiffcomplexificationbij} Let $\varphi : G \rightarrow H$ be a morphism of real algebraic groups. Then $\varphi$ is an isomorphism of real algebraic groups if and only if $\varphi_{\mathbb{C}} : G_{\mathbb{C}} \rightarrow H_{\mathbb{C}}$ is bijective.
\end{cor}

\begin{proof} By lemma \ref{lemcomplexmorphrealalggrismorphcomplexalggru}, $\varphi_{\mathbb{C}}$ is a morphism of complex algebraic groups. Therefore, by Proposition 21.2.6 of \cite{TY}, $\varphi_{\mathbb{C}}$ is a polynomial isomorphism if and only if $\varphi_{\mathbb{C}}$ is bijective, hence the result by previous lemma \ref{lempolyisomiffcomplexpolyisom}. 
\end{proof}

\section{Polynomial actions of real algebraic groups on real algebraic sets} \label{sectionpolyactralggpralgsets}

We are now going to introduce the frame we intend to work with, namely the real algebraic sets on which a real algebraic group (in the sense of definition \ref{defrealalggroup}) acts polynomially. Again, our motivation to restrict ourselves, at first, to polynomial actions is to be able to apply Hilbert's finiteness theorem (theorem \ref{theohilbertfiniteness}).

Let $m$ and $n$ denote positive integers and let $G \subset \mathbb{R}^m$ be a real algebraic group with identity element $e$.

\begin{de} \label{degrealalgset} Let $X$ be a real algebraic subset of $\mathbb{R}^n$. We say that $X$ is a \emph{$G$-real algebraic set} if $X$ is equipped with a group action $\alpha : G \times X \rightarrow X~;~(g,x) \mapsto g \cdot x$ of $G$ such that the map~$\alpha$ is (the restriction of) a polynomial map.
\end{de}

Remark that if $X \subset \mathbb{R}^n$ is a $G$-real algebraic set and if $Y$ is a real algebraic subset of $\mathbb{R}^n$ included in $X$ and globally preserved by the considered action $\alpha$ of $G$ on $X$ (i.e. for all $g \in G$, for all $y \in Y$, $g \cdot y = \alpha(g,y) \in Y$), then $Y$ is a $G$-real algebraic set as well: in this case, we will say that $Y$ is a \emph{$G$-real algebraic subset of $X$}. 

Notice also that if $\varphi : H \rightarrow G$ is morphism of real algebraic groups, then we can consider the polynomial action $H \times X \rightarrow X~;~(h,x) \mapsto \alpha(\varphi(h),x)$ of $H$ on $X$, making $X$ into a $H$-real algebraic set as well.
\\
 
Let us give a first list of examples of $G$-real algebraic sets:

\begin{exs} \label{exgrealalgset}
~
	\begin{enumerate}
		\item Matrix multiplication makes the entire affine space $\mathbb{R}^n$ into a $\Omega_n(\R)$-real algebraic set, as well as a $\OO_n(\R)$-real algebraic set and a $\SO_n(\R)$-real algebraic set.
		\item The cone of equation $x^2 + y^2 = z^2$ in $\R^3$ as well as the set of equation $x^2+y^2 + z^4 = z^2$ are real algebraic subsets of $\R^3$ globally preserved by the polynomial action 
		$$\begin{array}{ccc}S^1 \times \R^3 & \rightarrow & \R^3\\ \big( (a,b), (x,y,z)\big) & \mapsto & (ax-by,bx + ay, z)\end{array}$$
of $S^1$ on $\R^3$.
		\item Let $k$ be a positive integer and denote by $V_k(\R^n)$ the set of $k$-tuples of orthonormal vectors of $\R^n$ (with respect to the canonical scalar product) i.e. the real algebraic subset of points $A = \left(a_{ij}\right)_{1\leq i \leq n, 1 \leq j \leq k}$ of $\M_{n,k}(\R)$, identified with $\mathbb{R}^{nk}$, such that 
$${}^t\!A A = I_k ~\Longleftrightarrow ~ \forall i,j \in \{1,\ldots,k\}, \, \sum_{l =1}^n a_{li} a_{lj} = \delta_{ij}.$$		
The real algebraic set $V_k(\R^n)$ is called the {\it Stiefel manifold} of index $k$ (notice that $V_n(\R^n) = \OO_n(\R)$ and that $V_k(\R^n)$ is empty as soon as $k$ is larger than $n$). We can equip $V_k(\R^n)$ with the polynomial group action $\OO_n(\R) \times V_k(\R^n) \rightarrow V_k(\R^n)~;~(O,A) \mapsto OA$ of~$\OO_n(\R)$, as well as with the polynomial group action $\OO_k(\R) \times V_k(\R^n) \rightarrow V_k(\R^n)~;~(O,A) \mapsto A O^{-1} = A {}^t O$ of $\OO_k(\R)$ (notice that the former action is transitive while the latter one is free).
		\item The real algebraic group $G$ acts on itself by left multiplication and the corresponding action map $G \times G \rightarrow G$ is the binary operation of $G$ hence is polynomial: $G$ is itself a $G$-real algebraic set.
		\item If $X$ is a $G$-real algebraic set, $H$ is a real algebraic group and $Y$ is a $H$-real algebraic set, the polynomial action $(G \times H) \times (X \times Y) \rightarrow X \times Y~;~(g,h),(x,y) \mapsto (g \cdot x, h \cdot y)$ makes~$X \times Y$ into a $(G \times H)$-real algebraic set.
	\end{enumerate}
\end{exs}
 
\begin{rem} \label{remfinitegroupgrealalgsets} Suppose that $G$ is a finite group $\left\{g_1,\ldots,g_N\right\}$ and let $X$ be a real algebraic subset of $\R^n$ equipped with a group action $\alpha : G \times X \rightarrow X~;~(g,x) \mapsto g \cdot x$ of $G$ such that, for all $g \in G$, the map $\alpha_g : X \rightarrow X~;~x \mapsto g \cdot x$ is polynomial. Then $\alpha$ is the restriction of the polynomial map 
$$\begin{array}{ccc}\R^m \times \R^n & \rightarrow & \R^n\\(g,x) & \mapsto & \displaystyle{\sum_{k=1}^N\left(\prod_{1\leq l\leq N, \, l \neq k} \frac{\|g-g_l\|^2}{\|g_k - g_l\|^2}\right) \alpha_{g_k}(x)}\end{array}$$
and $X$ is therefore a $G$-real algebraic set.
\end{rem}

Let us then emphasize an important class of $G$-real algebraic sets:

\begin{de} \label{defrealpolyrepres} Let $V$ be an $n$-dimensional real representation of $G$ i.e. an $n$-dimensional real vector space with a group homomorphism $\rho : G \rightarrow \mathcal{GL}(V)$. Consider a basis of $V$ and denote by~$(v_1,\ldots,v_n)$ the corresponding coordinate system. If the action $G \times V \rightarrow V~;~(g,v) \mapsto \rho(g)(v)$ is polynomial with respect to the coordinates $(v_1,\ldots,v_n)$ and the coordinates of the affine space~$\R^m$ (i.e. if, for all $g \in G$, the entries of the matrix of $\rho(g)$ in the considered basis of $V$ are polynomial in the coordinates of $g$; notice that this property is independent of the chosen basis of $V$), we say that $V$ is a \emph{polynomial representation of $G$}. 
\end{de}

\begin{exs} \label{exspolynomialrepres}
~
\begin{enumerate}
	\item Considering the inclusion $\rho : \Omega_n(\R) \hookrightarrow \mathcal{GL}(\R^n)$, the affine space $\R^n$ is a polynomial representation of $\Omega_n(\R)$, as well as of $\OO_n(\R)$ and $\SO_n(\R)$. 
	\item Let $V$ be a polynomial representation of $G$ associated to a group homomorphism $\rho : G \rightarrow \mathcal{GL}(V)$. The dual representation $V^* := \mathcal{L}(V,\R)$ of $V$, associated to the group homomorphism $\rho^* : G \rightarrow \mathcal{GL}(V^*)$ which associates to any $g \in G$ the map $f \in V^* \mapsto f \circ \rho(g^{-1}) \in V^*$, is polynomial as well. Indeed, considering a basis of $V$, for all $g \in G$, the matrix of $\rho^*(g)$ in the corresponding dual basis of $V^*$ is the transpose of the matrix of $\rho(g^{-1})$ in the former basis of $V$ (the map $G \rightarrow G~;~g \mapsto g^{-1}$ is polynomial since $G$ is a real algebraic group). Notice that the canonical linear isomorphism 
$$\begin{array}{ccc}V & \rightarrow & \left(V^*\right)^*\\
			v & \mapsto & \left(f \in V^* \mapsto f(v) \in \R\right)
    \end{array}$$	
is furthermore $G$-equivariant for the induced action of $G$ on $\left(V^*\right)^*$.
\end{enumerate}
\end{exs}

Below, we will state that, up to isomorphism, any $G$-real algebraic set can actually be considered as a $G$-real algebraic subset of a polynomial representation of $G$. But we first have to make precise what notion of isomorphism of $G$-real algebraic sets we will be considering:

\begin{de} Let $X \subset \R^n$ and $Y \subset \R^N$ be $G$-real algebraic sets and let $f : X \rightarrow Y$ be a map. We say that $f$ is a \emph{morphism of $G$-real algebraic sets} if $f$ is both equivariant (i.e. for all $g \in G$ and $x \in X$, $f(g \cdot x) = g \cdot f(x)$) and the restriction of a polynomial map $\R^n \rightarrow \R^N$. If~$f$ is a bijective morphism of $G$-real algebraic sets such that $f^{-1}$ is also a morphism of $G$-real algebraic sets (i.e. the restriction of a polynomial map $\mathbb{R}^N \rightarrow \mathbb{R}^n$), we will say that $f$ is an \emph{isomorphism of $G$-real algebraic sets}, or just an \emph{isomorphism} when the context is clear, and that the $G$-real algebraic sets $X$ and $Y$ are \emph{isomorphic}.
\end{de} 

Remark that, with the vocabulary we introduced in subsection \ref{subsecfirstdefproprealalggr}, a map between $G$-real algebraic sets is an isomorphism if and only if it is an equivariant polynomial isomorphism.

\begin{exs} \label{exsmorphismgrealalgsets}
~	
	\begin{enumerate}
		\item Consider the polynomial action of $S^1$ on $\R^3$ given in example \ref{exgrealalgset} 2, as well as the polynomial action
$$\begin{array}{ccc}S^1 \times \R^2 & \rightarrow & \R^3\\ \big( (a,b), (x,y)\big) & \mapsto & (ax-by,bx + ay)\end{array}$$
of $S^1$ on $\R^2$. If $X$ denotes the $S^1$-real algebraic subset of $\R^3$ of equation $x^2+y^2 = z$, the map $X \rightarrow \R^2 ; (x,y,z) \mapsto (x,y)$ is an isomorphism of $S^1$-real algebraic sets.   
		\item Keep the above $S^1$-action on $\R^3$ and consider the $S^1$-real algebraic subset $Y$ of equation $x^2+y^2 + z^4 = z^2$, as well as the sphere $S^2 := \{(x,y,z) \in \R^3~|~x^2+y^2+z^2 = 1 \}$. Then the map $S^2 \rightarrow Y~;~(x,y,z) \mapsto (xz,yz,z)$ (which is the restriction of the blowing-up of $\R^3$ at the origin) is a morphism of $S^1$-real algebraic sets.
	\end{enumerate}
\end{exs}

Let $X \subset \R^n$ be a $G$-real algebraic set associated to a polynomial group action $\alpha : G \times X \rightarrow X~;~(g,x) \mapsto g \cdot x$. As announced above, the action of $G$ can be polynomially linearized: up to an isomorphism of $G$-real algebraic sets, we can always suppose that the polynomial action of~$G$ on $X$ is given by matrices with polynomial entries in the coordinates of $\R^m$.

\begin{prop} \label{proplinearizgrealalgset} The $G$-real algebraic set $X$ is isomorphic to a $G$-real algebraic subset of a polynomial representation $\R^N$ of $G$. 
\end{prop}

\begin{proof} As in \cite{BCR} section 3.2, denote by $\mathcal{P}(X) = \R[x_1,\ldots,x_n]/I(X)$ the $\R$-algebra of polynomial functions from $X$ to $\R$, and consider the linear action of $G$ on $\mathcal{P}(X)$ which associates, to any $g \in G$ and $f \in \mathcal{P}(X)$, the polynomial function 
$$g \cdot f : \begin{array}{ccc}X & \rightarrow & \R\\x & \mapsto & f\big(\alpha(g^{-1},x)\big)\end{array}$$ 
on $X$. 

Now, let $f \in \mathcal{P}(X)$ and remark that, since $\alpha$ is a polynomial map, the composition $f \circ \alpha$ belongs to $\mathcal{P}(G \times X) = \mathcal{P}(G) \otimes_{\R} \mathcal{P}(X)$ (see for instance \cite{BCR} Theorem 2.8.3 (iii)). There then exist polynomial functions $\phi_1,\ldots,\phi_l \in \mathcal{P}(G)$ and $h_1,\ldots,h_l \in \mathcal{P}(X)$ such that $f \circ \alpha = \sum_{k=1}^l \phi_k \otimes h_k$. In particular, for all $g \in G$, $g \cdot f = \sum_{k=1}^l \phi_k(g^{-1}) h_k \subset {\rm Span}\{h_1,\ldots,h_l\}$: this shows that the~$\R$-vector subspace $V_f$ of $\mathcal{P}(X)$ generated by the orbit $\{g \cdot f~|~g \in G\}$ of $f$ is finite-dimensional and is then a polynomial representation of $G$ (see remark \ref{remfdstablesubspacepolyrepres} 1 below).

Consider the generators $\overline{x_1},\ldots,\overline{x_n}$ of the $\R$-algebra $\mathcal{P}(X)$ and denote $V := \sum_{k=1}^n V_{\overline{x_k}}$: $V$ is also a polynomial representation of $G$ which furthermore generates the $\R$-algebra $\mathcal{P}(X)$. As we stated in example \ref{exspolynomialrepres} 2, the dual representation $V^*$ is polynomial as well.

Fix a basis $(f_1,\ldots,f_N)$ for the $\R$-vector space $V$ (the polynomial functions $f_1,\ldots,f_N$ then generate the $\R$-algebra $\mathcal{P}(X)$) as well as the corresponding dual basis for $V^*$ and let $i$ denote the map $X \rightarrow V^*$ which associates to any~$x \in X$ the evaluation linear form $f \in V \mapsto f(x) \in \R$. The map $i$ is polynomial (for all $x \in X$, the coordinates of $i(x)$ are $(f_1(x),\ldots,f_N(x)$), equivariant (if $x \in X$ and $g \in G$, we have, for all~$f \in V$, $i(g \cdot x)(f) = f(g \cdot x) = (g^{-1} \cdot f)(x) = i(x)(g^{-1} \cdot f) = g \cdot i(x) (f)$) and injective since the polynomial maps $f_1,\ldots,f_N$ generates the $\R$-algebra $\mathcal{P}(X)$ (if $x,x' \in X$ have the same value by any polynomial map on $X$, then $x = x'$).   

Finally, write, if $k \in \{1,\ldots,n\}$, $\overline{x_k} = P_k(f_1,\ldots,f_N)$ with $P_k \in \R[y_1,\ldots,y_N]$. The inverse map of the bijective map $\iota : X \rightarrow i(X)~;~x \mapsto i(x)$ is, in the above considered systems of coordinates, the restriction of the polynomial map $\R^N \rightarrow \R^n~;~(y_1,\ldots,y_N) \mapsto \big(P_1(y_1,\ldots,y_N), \ldots, P_n(y_1,\ldots,y_N)\big)$, so that $i(X)$ is a $G$-real algebraic subset (see remark \ref{remfdstablesubspacepolyrepres}~2 below) of the polynomial representation $V^*$ of $G$ and $\iota$ is an isomorphism of $G$-real algebraic sets. 
\end{proof}

\begin{rems} \label{remfdstablesubspacepolyrepres}
~
\begin{enumerate}
	\item With the notations of the beginning of the above proof, if $V$ is a finite-dimensional $\R$-vector subspace of $\mathcal{P}(X)$ which is furthermore stable under the action of~$G$ (i.e. for all~$g \in G$ and $f \in V$, $g \cdot f \in V$), then $V$ is a polynomial representation of $G$. Indeed, consider a basis $(f_1,\ldots,f_N)$ of $V$ and, for any $r \in \{1,\ldots,N\}$, write 
	$$f_r \circ \alpha = \sum_{k=1}^{l_r} \phi_{r,k} \otimes h_{r,k} \in \mathcal{P}(G) \otimes_{\R} \mathcal{P}(X).$$
In particular, we have $f_r = \sum_{k=1}^{l_r} \phi_{r,k}(e) h_{r,k}$ for any $r \in \{1,\ldots,N\}$, so that $V$ is a subspace of $W := {\rm Span}\left\{h_{r,k}~|~r \in \{1,\ldots,N\}, k \in \{1,\ldots,l_r\}\right\}$. Considering a basis $(h_1,\ldots,h_l)$ of $W$, we can write, if $r \in \{1,\ldots,N\}$, $f_r \circ \alpha = \sum_{k=1}^{l} \widetilde{\phi_{r,k}} \otimes h_{k}$ (with $\widetilde{\phi_{r,k}} \in \mathcal{P}(G)$ for all $k \in \{1,\ldots,l\}$) and then, for all $g \in G$, the vector of coordinates of $g \cdot f_r$ with respect to the basis $(f_1,\ldots,f_N)$ is the image of the vector $\left(\widetilde{\phi_{r,1}}(g^{-1}), \ldots, \widetilde{\phi_{r,l}}(g^{-1})\right)$ (which has polynomial coordinates in $g$) by the matrix of any linear projection $W \rightarrow V$ in the respective bases $(h_1,\ldots,h_l)$ and $(f_1,\ldots,f_N)$ of $W$ and $V$.
	\item Let $Y$ be a real algebraic subset of $\R^n$, $Z$ be any subset of $\R^N$ and $j : Y \rightarrow Z$ be a bijective restriction of a polynomial map $\widetilde{j} : \R^n \rightarrow \R^N$. If the inverse map $p : Z \rightarrow Y$ is the restriction of a polynomial map $\widetilde{p} : \R^N \rightarrow \R^n$, then $Z$ is a real algebraic subset of $\R^N$ (and then the map $j$ is a polynomial isomorphism), since we can write $Z = \{z \in \R^N~|~\widetilde{p}(z) \in Y \mbox{ and } \widetilde{j} \circ \widetilde{p}(z) = z\}$ and $Y$ is algebraic.
\end{enumerate}
\end{rems}

The proof of proposition \ref{proplinearizgrealalgset} (which is identical to the complex case) can be adapted to show that any real algebraic group is isomorphic to a real algebraic group which is a subgroup of a general linear group:

\begin{prop} \label{proprealalggroupislinear} The real algebraic group $G$ is isomorphic to a linear algebraic group (example \ref{exspolyrealalggroups} 2).
\end{prop}

\begin{proof} Consider the polynomial action $\mu : G \times G \rightarrow G$ of $G$ on itself by left multiplication. Keeping the notations of the proof of proposition \ref{proplinearizgrealalgset}, there exists an $\R$-vector subspace $V$ of $\mathcal{P}(G) = \R[x_1,\ldots,x_m]/I(G)$ such that $V$ is a polynomial representation of $G$ and $V$ generates the $\R$-algebra $\mathcal{P}(G)$.

Choose a basis $\{f_1,\ldots,f_M\}$ of $V$ and consider the polynomial map $j : G \rightarrow \mathcal{GL}(V) (\cong \GL_M(\R))$ which associates to any $g \in G$ the (matrix in the basis $\{f_1,\ldots,f_M\}$ of the) linear isomorphism $f \in V \mapsto g \cdot f \in V$. The map $j$ is a group homomorphism (if $g_1,g_2 \in G$ then, for all $f \in V$ and $x \in G$, $j(g_1 g_2)(f)(x) = f \big((g_1 g_2)^{-1} x \big) = f \big(g_2^{-1} g_1^{-1} x \big) = j(g_2) \circ j(g_1)(f)(x)$) which is injective: if $g_1,g_2 \in G$ satisfy $j(g_1) = j(g_2)$ then, for all $f \in V$, $j(g_1)(f)(e) = j(g_2)(f)(e)$ i.e. $f(g_1^{-1}) = f(g_2^{-1})$, and then $g_1^{-1} = g_2^{-1}$ since $V$ generates the $\R$-algebra $\mathcal{P}(G)$.

Furthermore, writing, for $k \in \{1,\ldots,m\}$, $\overline{x_k} = P_k(f_1,\ldots,f_M)$ with $P_k \in \R[y_1,\ldots,y_M]$, the inverse map of the bijective map $\kappa : G \rightarrow j(G)~;~g \mapsto j(g)$ is the restriction of the polynomial map 
$$\R^{M^2} \rightarrow \R^m~;~A \mapsto \omega\Big(P_1\big(A(f_1)(e), \ldots, A(f_M)(e)\big), \ldots, P_m\big(A(f_1)(e), \ldots, A(f_M)(e)\big)\Big),$$
if $\omega$ denotes the polynomial inverse map $G \rightarrow G~;~g \mapsto g^{-1}$. Therefore, $j(G)$ is a real algebraic subset of $\R^{M^2}$ (remark \ref{remfdstablesubspacepolyrepres} 2). Since furthermore the map $j(G) \rightarrow j(G)~;~\gamma \mapsto \gamma^{-1}$ is a polynomial map, as being the composition $\kappa \circ \omega \circ \kappa^{-1}$, the subgroup $j(G)$ of $\GL_M(\R)$ is a real algebraic group and $\kappa : G \rightarrow j(G)$ is an isomorphism of real algebraic groups.
\end{proof}

In this text, we will not only be interested in the $G$-real algebraic sets but also in the semialgebraic sets equipped with a polynomial action of $G$:

\begin{de} \label{defGsemialgset} Let $S$ be a semialgebraic subset of $\mathbb{R}^n$. We say that $S$ is a \emph{$G$-semialgebraic set} if $S$ is equipped with a group action $\alpha : G \times S \rightarrow S~;~(g,x) \mapsto g \cdot x$ of $G$ such that the map~$\alpha$ is the restriction of a polynomial map.
\end{de}

Let $S \subset \mathbb{R}^n$ be a $G$-semialgebraic set with associated polynomial group action $\alpha : G \times S \rightarrow S$: $\alpha$ is the restriction of a polynomial map $\widetilde{\alpha} : \R^m \times \R^n \rightarrow \R^n$ which restricts to a map $\overline{\alpha}^{\mathcal{Z}} : G \times \overline{S}^{\mathcal{Z}} \rightarrow \overline{S}^{\mathcal{Z}}$.

\begin{lem} \label{lemzariskiclosgsaisgras} The polynomial map $\overline{\alpha}^{\mathcal{Z}} : G \times \overline{S}^{\mathcal{Z}} \rightarrow \overline{S}^{\mathcal{Z}}$ is an action of $G$ on $\overline{S}^{\mathcal{Z}}$, making~$\overline{S}^{\mathcal{Z}}$ into a $G$-real algebraic set.
\end{lem}

\begin{proof} If $g_1,g_2$ are any two elements of $G$, we have the polynomial equalities $\widetilde{\alpha}(e, \cdot) = \cdot$ and $\widetilde{\alpha}(g_1, \widetilde{\alpha}(g_2, \cdot)) = \widetilde{\alpha}(g_1 g_2, \cdot)$ on $S$, hence on $\overline{S}^{\mathcal{Z}}$.
\end{proof}

This means that the polynomial action of $G$ on $S$ comes from a (unique) polynomial action of~$G$ on the Zariski closure of $S$ in $\mathbb{R}^n$. Therefore, a semialgebraic subset of $\R^n$ is a~$G$-semialgebraic set if and only if it is a $G$-stable subset of a $G$-real algebraic set of $\R^n$.
\\

We dedicate the last part of this section to the complexification of $G$-real algebraic sets: as for real algebraic groups, considering polynomial group actions on real algebraic sets allows to consider their complexification. In other words, polynomial group actions on real algebraic sets are induced by complex polynomial group actions with real coefficients.

\begin{lem} \label{lemcomplexificationpolyaction} Let $X \subset \R^n$ be a $G$-real algebraic set with polynomial group action $\alpha : G \times X \rightarrow X~;~(g,x) \mapsto g \cdot x$. The complexification $\alpha_{\C} : G_{\C} \times X_{\C} \rightarrow X_{\C}$ of $\alpha$ is a (polynomial) group action of the complex algebraic group $G_{\C}$ on the complex algebraic set $X_{\C}$.
\end{lem} 

\begin{proof} Suppose that $\alpha$, resp. $\mu$, is the restriction of a polynomial map $\widetilde{\alpha} : \R^m \times \R^n \rightarrow \R^n$, resp. $\widetilde{\mu} : \R^m \times \R^m \rightarrow \R^m$. Since we have the polynomial equality $\widetilde{\alpha}_{\C}( \cdot, \widetilde{\alpha}_{\C}(\cdot, \cdot)) = \widetilde{\alpha}_{\C}(\widetilde{\mu}_{\C}(\cdot, \cdot), \cdot)$ on $G \times G \times X$, we have the same equality on the complex Zariski closure $G_{\C} \times G_{\C} \times X_{\C}$. Similarly, the polynomial equality $\widetilde{\alpha}_{\C}(e, \cdot) = \cdot$ on $X$ extends on $X_{\C}$: these equalities are illustrations of the functoriality of the complexification. 
\end{proof}

\begin{exs} \label{exscomplexificationpolyaction}
~
\begin{enumerate}
	\item The usual linear action of $\Omega_n(\R)$, resp. $\OO_n(\R)$, resp. $\SO_n(\R)$, on $\R^n$ (example~\ref{exgrealalgset}~1) extends to the usual linear action of $\Omega_n(\C)$, resp. $\OO_n(\C)$, resp. $\SO_n(\C)$, on~$\C^n$ (see example \ref{excomplrealalggroup} 1).
	\item Let $k$ be a positive integer such that $k<n$. We want to complexify the polynomial actions of $\OO_n(\R)$ and $\OO_k(\R)$ on $V_k(\R^n)$ (example \ref{exgrealalgset} 3). First, denote by $V_k^n(\C)$ the complex algebraic subset of points $A$ of $\M_{n,k}(\C) \cong \mathbb{C}^{nk}$ such that ${}^t\!A A = I_k$ (the set~$V_k^n(\C)$ must not be confused with the usual complex Stiefel manifold $V_k(\C^n) = \{A \in \M_{n,k}(\C)~|~{}^t \overline{A} A = I_k\}$): if $\vartheta$ denotes the non-degenerate bilinear form 
	$$\vartheta : \begin{array}{ccc}\C^n \times \C^n & \rightarrow & \C\\ \big((z_1,\ldots,z_n),(z_1',\ldots,z_n')\big) & \mapsto & \sum_{l=1}^n z_l z_l'\end{array}$$
on $\C^n$, $V_k^n(\C)$ is the set of $k$-tuples $(v_1,\ldots,v_k)$ of pairwise $\vartheta$-orthogonal vectors of $\C^n$ such that, for all $l \in \{1,\ldots,k\}$, $\vartheta(v_l,v_l) = 1$. The set $V_k^n(\C)$ is the complexification of $V_k(\R^n)$. Indeed, the complex algebraic set $V_k^n(\C)$ contains $V_k(\R^n)$, is irreducible ($k < n$) and $\dim V_k^n(\C) = kn -\frac{k(k+1)}{2} = \dim V_k(\R^n) = \left(V_k(\R^n)\right)_{\C}$ (see section \ref{subsectappendixstiefel} of the appendix). As a consequence, the real polynomial actions of $\OO_n(\R)$ and $\OO_k(\R)$ on $V_k(\R^n)$ considered in example \ref{exgrealalgset} 3 extend into the complexified polynomial actions $\OO_n(\C) \times V_k(\C^n) \rightarrow V_k(\C^n)~;~(O,A) \mapsto OA$ and $\OO_k(\C) \times V_k(\C^n) \rightarrow V_k(\C^n)~;~(O,A) \mapsto A O^{-1}$.
	\item Let $V$ be a polynomial $n$-dimensional real representation of $G$ associated to a group homomorphism $\rho : G \rightarrow \mathcal{GL}(V)$ (definition \ref{defrealpolyrepres}). Fix a basis for $V$ and denote $V_{\C} := V \otimes_{\R} \C$. Considering the chosen $\R$-basis of $V$ as a $\C$-basis of $V_{\C}$, the polynomial action $G \times V \rightarrow V~;~(g,v) \mapsto g \cdot v$ extends into the complex polynomial action $G_{\C} \times V_{\C} \rightarrow V_{\C}$, making $V_{\C}$ into the complex polynomial representation of $G_{\C}$ with group homomorphism $\rho_{\C} :  G_{\C} \rightarrow \mathcal{GL}(V_{\C})$ which associates to any $g \in G_{\C}$ the linear automorphism $v \in V_{\C} \mapsto g \cdot v \in V_{\C}$ (the entries of the matrix of $\rho_{\C}(g)$ are given by the same polynomials, with real coefficients, in the coordinates of $g$). 
\end{enumerate}
\end{exs}

Finally, the complexification of a morphism of $G$-real algebraic sets is a morphism of $G_{\C}$-complex algebraic sets :

\begin{lem} \label{lemcomplexificationmorphismGras} Let $f : X \rightarrow Y$ be a morphism, resp. an isomorphism, of $G$-real algebraic sets. Then the complexification $f_{\C} : X_{\C} \rightarrow Y_{\C}$ is a morphism, resp. an isomorphism, of $G_{\C}$-complex algebraic sets (i.e. an equivariant restriction of a polynomial map, resp. an equivariant polynomial isomorphism).
\end{lem}

\begin{proof} The equivariance of $f$ implies the equivariance of $f_{\C}$ by functoriality of the complexication (cf. section \ref{subseccomplexificationralgps}).
\end{proof}

\section{Linearly reductive real algebraic groups, invariant polynomial functions and Hilbert's finiteness theorem} \label{linredralggphilbfinth}

In this section, we will give the proof of Hilbert's finiteness theorem for real algebraic sets (theorem \ref{theohilbertfiniteness}), a result that we will use in the next section in order to construct a real algebro-geometric quotient for any polynomial action of a so-called linearly reductive real algebraic group on a real algebraic set. 

Hilbert's finiteness theorem states that, if $G$ is a real algebraic group (in the sense of definition \ref{defrealalggroup}) the $\R$-algebra of invariant polynomial functions of any $G$-real algebraic set (in the sense of definition \ref{degrealalgset}) is finitely generated as soon as $G$ is linearly reductive. In particular, we will see that the hypotheses of polynomial group structure and polynomial group action that we have been considering are needed here.

\subsection{Linearly reductive real algebraic groups} \label{subsectionlinredragr}

We give below the definition of a linearly reductive real algebraic group that justifies the employed terminology. We will thereafter give several characterizations of the notion.

Let $G \subset \R^m$ be a real algebraic group and first remark that if $V$ is a polynomial representation of $G$ (definition \ref{defrealpolyrepres}), associated to a group homomorphism $\rho : G \rightarrow \mathcal{GL}(V)$, and if $W$ is a subrepresentation of $V$ (i.e. a vector subspace of $V$ such that for all $g \in G$, $\rho(g)(W) \subset W$, so that we can consider the induced group homomorphism $\rho_W : G \rightarrow \mathcal{GL}(W)$ which associates to any $g \in G$ the restriction $W \rightarrow W$ of $\rho(g)$), then $W$ is also a polynomial representation of~$G$.

\begin{de} \label{defreallinearreductivity} The real algebraic group $G$ is called \emph{linearly reductive} if for any polynomial representation $V$ of $G$ and every subrepresentation $W$ of $V$, there exists a subrepresentation $W'$ of $V$ such that $V = W \oplus W'$.
\end{de}

Let $X \subset \R^n$ be a $G$-real algebraic set with polynomial action $\alpha : G \times X \rightarrow X$. As in the proof of proposition \ref{proplinearizgrealalgset}, $\alpha$ induces on the $\R$-algebra $\mathcal{P}(X)$ of polynomial functions on $X$ the linear action of $G$ which associates, to any $g \in G$ and $f \in \mathcal{P}(X)$, the polynomial function 
$$g \cdot f : \begin{array}{ccc}X & \rightarrow & \R\\x & \mapsto & f\big(\alpha(g^{-1},x)\big).\end{array}$$ 
We will then say that a polynomial function $f$ on $X$ is \emph{invariant} if for all $g \in G$, $g \cdot f = f$, and we will denote by $\mathcal{P}(X)^G$ the set of invariant polynomial functions on $X$: $\mathcal{P}(X)^G$ is an~$\R$-subalgebra of $\mathcal{P}(X)$.

We will see that $G$ is linearly reductive if and only if any $G$-real algebraic set has a Reynolds operator:

\begin{de} \label{dereynoldsop} A \emph{Reynolds operator} of the $G$-real algebraic set $X$ is an $\R$-linear map $\mathcal{R} : \mathcal{P}(X) \rightarrow \mathcal{P}(X)^G$ such that, for all $f \in \mathcal{P}(X)^G$, $\mathcal{R}(f) = f$ and, for all $g \in G$ and $f \in \mathcal{P}(X)$, $\mathcal{R}(g \cdot f) = \mathcal{R}(f)$.
\end{de}

In below proposition \ref{propequivlinred}, we are actually going to establish the equivalence of linearly reductivity with other important properties of linearly reductive real algebraic groups. The proof of the statement will be a mix of the proofs of Theorem 2.2.5 of \cite{DK} and Proposition~4.14 of \cite{Hos}; in particular, it overcomes the fact that one conclusion of Schur's Lemma over $\C$ (namely the second point of Schur's Lemma 1.7 in \cite{FH}) is not true in the real framework.

When $V$ is any representation of $G$, denote by $V^G$ the subspace $\left\{v \in V~|~\forall g \in G, \, g \cdot v = v\right\}$ of \emph{$G$-invariant} vectors of $V$.

\begin{prop} \label{propequivlinred} The following properties are equivalent :
\begin{enumerate}
	\item the real algebraic group $G$ is linearly reductive,
	\item for any polynomial representations $V$ and $W$ of $G$, for any surjective $G$-equivariant linear map $\omega : V \rightarrow W$, we have $\omega\left(V^G\right) = W^G$,
	\item for any polynomial representation $V$ of $G$, for any $\phi \in \left(V^*\right)^G \setminus \left\{\bf{0}\right\}$, there exists $v \in V^G \setminus \left\{\bf{0}\right\}$ such that $\phi(v) \neq 0$,
	\item for any polynomial representation $V$ of $G$, there exists a subrepresentation $W$ of $V$ such that $V = V^G \oplus W$ and $\left(W^*\right)^G = \left\{\bf{0}\right\}$,
	\item any $G$-real algebraic set has a Reynolds operator.
\end{enumerate}
\end{prop}

\begin{proof} 
1) $\Rightarrow$ 2) : Suppose that $G$ is linearly reductive and let $\omega : V \rightarrow W$ be a surjective $G$-equivariant linear map between polynomial representations of $G$. Since $\omega$ is $G$-equivariant, the kernel of $\omega$ is a subrepresentation of $V$ and, by hypothesis, there exists a subrepresentation $V'$ of~$V$ such that $V = {\rm Ker} \, \omega \oplus V'$: the map $\omega$ then restricts to a $G$-equivariant linear isomorphism $\widetilde{\omega} : V' \rightarrow W$ and, in particular, $W^G = \widetilde{\omega}({V'}^G) \subset \omega(V^G)$.

2) $\Rightarrow$ 1) : Suppose that $G$ satisfies 2) and let $V$ and $W$ be polynomial representations of~$G$, with respective associated group homomorphisms $\rho : G \rightarrow \mathcal{GL}(V)$ and $\varrho : G \rightarrow \mathcal{GL}(W)$. We consider the $\R$-vector space $\mathcal{L}(V,W)$ of $\R$-linear maps from $V$ to $W$ as the finite-dimensional real representation of $G$ associated to the group homomorphism which associates to $g \in G$ the linear isomorphism 
$$\begin{array}{ccc}\mathcal{L}(V,W) & \rightarrow & \mathcal{L}(V,W)\\ \varpi &  \mapsto & \varrho(g) \circ \varpi \circ \rho(g^{-1}).\end{array}$$
This is a polynomial representation of $G$ (consider respective bases for $V$ and $W$ and the usual associated basis for $\mathcal{L}(V,W)$). Notice also that the linear maps of $\big(\mathcal{L}(V,W)\big)^G$ are exactly the $G$-equivariant linear maps from $V$ to $W$. Now suppose that $W$ is a subrepresentation of $V$ and consider the surjective $G$-equivariant linear map $\Pi : \mathcal{L}(V,W) \rightarrow \mathcal{L}(W,W)$ which associates to~$\varphi \in \mathcal{L}(V,W)$ the restriction $w \in W \mapsto \varphi(w) \in W$: because $G$ satisfies 2) and the identity~${\rm Id}_W$ of $W$ is a $G$-invariant linear automorphism of $W$, there exists $\varpi \in \big(\mathcal{L}(V,W)\big)^G$ such that $\Pi(\varpi) = {\rm Id}_W$ i.e. for all $w \in W$, $\varpi(w) = w$ ($\varpi$ is in particular surjective). If we then set $W' := {\rm Ker} \, \varpi$, the subspace~$W'$ is a subrepresentation of $V$ (since $\varpi$ is equivariant) and we have $V = W \oplus W'$.  

2) $\Rightarrow$ 3) : Suppose that $G$ satisfies 2), let $V$ be a polynomial representation of $G$ and let~$\phi$ be a nonzero $G$-invariant linear form of $V^{*}$ (see example \ref{exspolynomialrepres} 2). Equip $\R$ with the trivial action $G \times \R \rightarrow \R~;~(g,x) \mapsto x$ of $G$, so that $\phi$ is a surjective $G$-equivariant linear map $V \rightarrow \R$: since $G$ satisfies 2), there exists a (nonzero) $G$-invariant vector $v$ of $V$ such that $\phi(v) = 1$.

3) $\Rightarrow$ 2) : Suppose that $G$ satisfies 3). First, let $V$ be a polynomial representation of $G$ and $v \in V^G \setminus \left\{\bf{0}\right\}$. Since the corresponding linear form $\Phi_v : \phi \in V^* \mapsto \phi(v) \in \R$ of $\left(V^*\right)^*$ (see example \ref{exspolynomialrepres} 2) is nonzero and $G$-invariant as well, by hypothesis, there exists a nonzero~$G$-invariant linear form $\phi$ of $V^*$ such that $\Phi_v(\phi) \neq 0$ i.e. $\phi(v) \neq 0$. Now, let $\omega : V \rightarrow W$ be a surjective $G$-equivariant linear map between polynomial representations of $G$ and let us show that $\omega\left(V^G\right) = W^G$ by induction on the dimension of $W$ (the case $\dim W = 0$ is clear). So let~$w$ be a nonzero $G$-invariant vector of $W$: by the above discussion, there exists $\phi \in \left(W^*\right)^G \setminus \left\{\bf{0}\right\}$ such that $\phi(w) \neq 0$. But $\omega$ is surjective so that the $G$-invariant composition~$\phi \circ \omega$ of $V^*$ is also nonzero: since $G$ satisfies 3), there exists $v \in V^G \setminus \left\{\bf{0}\right\}$ such that $\phi \circ \omega(v) \neq 0$. As a consequence, there exists $\lambda \in \R^*$ such that $\phi(w) = \lambda \, \phi \circ \omega(v) = \phi \circ \omega(\lambda v)$ and then $w = \omega(\lambda v) + w_0$ with $w_0 \in {\rm Ker}\, \phi$. Since $\phi$ is $G$-invariant, the restriction $\omega^{-1}\left({\rm Ker} \, \phi\right) \rightarrow {\rm Ker} \, \phi$ of $\omega$ is a surjective $G$-equivariant linear map between polynomial representations of $G$: since $\dim {\rm Ker} \, \phi = \dim W - 1$ (because $\phi$ is nonzero), we can apply the induction hypothesis to assert that there exists~$v_0 \in V^G$ such that $w_0 = \omega(v_0)$. Therefore, $w = \omega(\lambda v + v_0) \in \omega\left(V^G\right)$. 

3) $\Rightarrow$ 4) : Suppose that $G$ satisfies 3) and let $V$ be a polynomial representation of $G$. Denote by $W$ the annihilator 
$$\left\{ v \in V~|~\forall \phi \in (V^*)^G, \, \phi(v) = 0\right\}$$
of $(V^*)^G$ in $V$ ($W$ is a subrepresentation of $V$). If $v \in V^G$ is such that for all $\phi \in (V^*)^G$, $\phi(v) = 0$, then, because $G$ satisfies 3), $v$ is necessarily the zero vector of $V$ (see the argument at the beginning of 3) $\Rightarrow$ 2)) and therefore $V^G \cap W = \left\{\bf{0}\right\}$. Furthermore, the linear map $\phi \in (V^*)^G \mapsto \phi_{|V^G} \in \left(V^G\right)^*$ is injective because $G$ satisfies 3), so that $\dim (V^*)^G \leq \dim \left(V^G\right)^* = \dim V^G$ and then 
$$\dim W = \dim V - \dim (V^*)^G \geq \dim V - \dim V^G:$$
as a consequence, $V = V^G \oplus W$. Finally, $\left(W^*\right)^G = \left\{\bf{0}\right\}$ because $W^G = V^G \cap W = \left\{\bf{0}\right\}$ and $G$ satisfies 3).

4) $\Rightarrow$ 5) : Suppose that $G$ satisfies 4), let $X$ be a $G$-real algebraic set and consider a finite-dimensional subspace $V$ of $\mathcal{P}(X)$ stable under the action of $G$: by remark \ref{remfdstablesubspacepolyrepres} 1, $V$ is a polynomial representation of $G$. There then exists, by hypothesis, a subrepresentation $W$ of $V$ such that $V = V^G \oplus W$ and $\left(W^*\right)^G = \left\{\bf{0}\right\}$, and we set $\mathcal{R}_{V}$ to be the linear projection $V \rightarrow V^G$ along $W$ (notice in particular that if $f \in V$ and $g \in G$, $\mathcal{R}_V(g \cdot f) = \mathcal{R}_V(f)$). Let $V'$ be a subrepresentation of $V$ and, using again the hypothesis, $W'$ be a subrepresentation of $V'$ such that $V' = {V'}^G \oplus W'$ and $\left({W'}^*\right)^G = \left\{\bf{0}\right\}$: if $f \in V'$, we have $\mathcal{R}_V \left(\mathcal{R}_{V'}(f)\right) = \mathcal{R}_{V'}(f)$ and, if $f_0 \in W'$, we have $\mathcal{R}_{V}(f_0) = \overline{0}$ (if there was an element $x$ in $X$ such that $\mathcal{R}_V(f_0)(x) \neq 0$, the~$G$-invariant linear form $h \in W' \rightarrow \mathcal{R}_V(h)(x) \in \R$ would be nonzero whereas $\left({W'}^*\right)^G = \left\{\bf{0}\right\}$) so that    
$$\overline{0} = \mathcal{R}_V\left(f -  \mathcal{R}_{V'}(f)\right) = \mathcal{R}_V(f) -  \mathcal{R}_V\left(\mathcal{R}_{V'}(f)\right) = \mathcal{R}_V(f) - \mathcal{R}_{V'}(f).$$
As a consequence, the map $\mathcal{R} : \mathcal{P}(X) \rightarrow \mathcal{P}(X)^G$ which associates to any $f \in \mathcal{P}(X)$ the $G$-invariant polynomial function $\mathcal{R}_V(f) \in V^G$ where $V$ is a $G$-stable finite-dimensional linear subspace of $\mathcal{P}(X)$ containing $f$ (such a subspace always exists as we showed in the proof of proposition \ref{proplinearizgrealalgset}) is well-defined: if $V_1$ and $V_2$ are two such subspaces of $\mathcal{P}(X)$ then so is their sum $V_1 + V_2$ and we have, by the above equality, $\mathcal{R}_{V_1}(f) = \mathcal{R}_{V_1+V_2}(f) = \mathcal{R}_{V_2}(f)$. Moreover, $\mathcal{R}$ is a Reynolds operator of the $G$-real algebraic set $X$.

5) $\Rightarrow$ 3) : Suppose that $G$ satisfies 5), let $V$ be a polynomial representation of $G$ and $\phi$ be a nonzero $G$-invariant linear form of $V^*$. Let $v_0 \in V$ such that $\phi(v_0) \neq 0$ and denote by $\epsilon$ the evaluation linear map $\mathcal{P}(V) \rightarrow \R~;~f \mapsto f(v_0)$ at $v_0$ (considering a basis of $V$ and the corresponding coordinate system). Since $G$ satisfies 5), the $G$-real algebraic set $V$ has a Reynolds operator $\mathcal{R}$ and we can consider the linear composition $\epsilon \circ \mathcal{R} : \mathcal{P}(V) \rightarrow \R$ as well as the restriction $\Phi : V^* \rightarrow \R$ of the latter to $V^* \subset \mathcal{P}(V)$. The linear form $\Phi$ of $\left(V^{*}\right)^*$ is $G$-invariant, since so is $\mathcal{R}$, and $\Phi$ is therefore the evaluation map $\Phi_v$ at a $G$-invariant vector $v$ of~$V$ (see example \ref{exspolynomialrepres} 2 and 3) $\Rightarrow$ 2)). We finally have $\phi(v) = \Phi(\phi) = \epsilon\left(\mathcal{R}(\phi)\right) = \epsilon(\phi) = \phi(v_0) \neq 0$ (the linear form $\phi$ is $G$-invariant).   
\end{proof}

We complete the above proof with additional properties on the objects we constructed:

\begin{rems} \label{remcharaclinearredprop} Suppose that the real algebraic group $G$ is linearly reductive.
\begin{enumerate}
	\item If $V$ is a polynomial representation of $G$, the annihilator $W$ of $\left(V^*\right)^G$ in $V$ (proof of 3) $\Rightarrow$ 4) above) is the unique subrepresentation of $V$ such that $V = V^G \oplus W$ and $\left(W^*\right)^G = \left\{\bf{0}\right\}$. Indeed, if $W'$ is a subrepresentation of $V$ such that $V = V^G \oplus W'$ and $\left({W'}^*\right)^G = \left\{\bf{0}\right\}$, then $\left(V^*\right)^G$ is a subspace of the annihilator ${W'}^0$ of $W'$ in $V^*$ (because the restriction of any $G$-invariant linear form of $V^*$ to $W'$ vanishes on $W'$ by hypothesis on $W'$) and 
	$$\dim \left(V^*\right)^G = \dim V - \dim W = \dim V - \dim W' = \dim {W'}^0,$$
so that ${W'}^0 = \left(V^*\right)^G$ i.e. $W' = W$. Notice also that, if $W'$ is a subrepresentation of $V$ such that $V = V^G \oplus W'$, then necessarily $\left({W'}^*\right)^G = \left\{\bf{0}\right\}$ because $G$ satisfies property 3 of proposition \ref{propequivlinred} and $\left(W'\right)^G = W' \cap V^G =  \left\{\bf{0}\right\}$. Consequently, $W$ is actually the unique subrepresentation of $V$ such that $V = V^G \oplus W$.
	\item If $X$ is a $G$-real algebraic set, the Reynolds operator $\mathcal{R}$ of $X$ defined in the proof of the implication 4) $\Rightarrow$ 5) of proposition \ref{propequivlinred} is actually the only Reynolds operator of $X$. In order to show this, let $\mathcal{R}'$ be a Reynolds operator of $X$, let $f \in \mathcal{P}(X)$, let $V$ be a $G$-stable linear subspace of $\mathcal{P}(X)$ containing $f$ and denote by $W$ the unique subrepresentation of~$V$ such that $V = V^G \oplus W$. If we write $f = h + f_0$ with $h \in V^G$ and $f_0 \in W$, we have $\mathcal{R'}(f_0) = \overline{0}$ (by the same argument as in the proof of 4) $\Rightarrow$ 5) above) so that
$$\mathcal{R'}(f) = \mathcal{R'}(h) = h = \mathcal{R}(h) = \mathcal{R}(f).$$
We will therefore denote by $\mathcal{R}_X$ the unique Reynolds operator of the $G$-real algebraic set~$X$.
	\item The Reynolds operator $\mathcal{R}_X$ of any $G$-real algebraic set $X$ is a morphism of $\mathcal{P}(X)^G$-modules. Indeed, keep the notations of previous point and let $\xi$ be an invariant polynomial function on $X$: $\xi V$ is then a polynomial representation of $G$ and, since we have $\xi V = \xi V^G + \xi W$ as well as $\xi V^G = (\xi V)^G$ (apply property 2 of proposition \ref{propequivlinred} to the surjective equivariant linear map $V \rightarrow \xi V$ which associates to any polynomial function of $V$ its multiplication by $\xi$) and $(\xi V)^G \cap (\xi W) = (\xi W)^G = \xi W^G = \left\{\overline{0}\right\}$, $\xi W$ is the unique subrepresentation of $\xi V$ such that $\xi V = (\xi V)^G \oplus \xi W$. As a consequence,
$$\mathcal{R}_X(\xi f) = \mathcal{R}_X(\xi h + \xi f_0) = \xi h = \xi \mathcal{R}_X(f).$$
	\item If $X$ is a $G$-real algebraic set and $V$ is any $G$-stable linear subspace of $\mathcal{P}(X)$, we have $\mathcal{R}_X(V) = V^G$ because, if $f \in V$, we have, with the notations of the proofs of proposition~\ref{propequivlinred} and proposition \ref{proplinearizgrealalgset}, $\mathcal{R}_X(f) = \mathcal{R}_{V_f}(f) \in V_f^G \subset V^G$.
\end{enumerate}
\end{rems}

The existence of a Reynolds operator for any $G$-real algebraic set is actually equivalent to the existence of a Reynolds operator for the $G$-real algebraic set $G$ equipped with its action on itself by left multiplication (example \ref{exgrealalgset} 4):

\begin{cor} \label{corequivlinearlyredreynoldsopgroup} The real algebraic group $G$ is linearly reductive if and only if there exists a Reynolds operator for $G$.
\end{cor} 

\begin{proof} First remark that $\mathcal{P}(G)^G = \R$ since any polynomial function $f \in \mathcal{P}(G)^G$ is constant: if~$g_1, g_2 \in G$, we have $f(g_1) = f\left(\left(g_1 g_2^{-1}\right) g_2\right) = f(g_2)$. As a consequence, a Reynolds operator of~$G$ is a $\R$-linear map $\mathcal{R}_G : \mathcal{P}(G) \rightarrow \R$ such that $\mathcal{R}_G\left(\overline{1}\right) = 1$ and, for all $g \in G$ and $f \in \mathcal{P}(G)$, $\mathcal{R}_G(g \cdot f) = \mathcal{R}_G(f)$.

Now, suppose that there is such a Reynolds operator $\mathcal{R}_G : \mathcal{P}(G) \rightarrow \R$ of $G$ and let $X$ be a $G$-real algebraic set associated to a polynomial action $\alpha : G \times X \rightarrow X$. Denote by~$\mathcal{R}$ the composition of the linear map $\mathcal{P}(X) \rightarrow \mathcal{P}(G) \otimes_{\R} \mathcal{P}(X)~;~f \mapsto f \circ \alpha$ with the linear map $\mathcal{P}(G) \otimes_{\R} \mathcal{P}(X) \rightarrow \mathcal{P}(X)$ which associates to any tensor $\phi \otimes h$ with $\phi \in \mathcal{P}(G)$ and $h \in \mathcal{P}(X)$ the quantity $\mathcal{R}_G(\phi) h$: if $f \in \mathcal{P}(X)$ and if we write $f \circ \alpha = \sum_{k=1}^l \phi_k \otimes h_k$ with $\phi_1,\ldots,\phi_l \in \mathcal{P}(G)$ and $h_1,\ldots,h_l \in \mathcal{P}(X)$, we have $\mathcal{R}(f) = \sum_{k=1}^l \mathcal{R}_G\left(\phi_k\right) h_k \in \mathcal{P}(X)$. With these notations, when $g_0 \in G$, we have, for all $g \in G$ and $x \in X$, $(g_0 \cdot f) \circ \alpha(g,x) = f \circ \alpha(g_0^{-1} g, x) =  \sum_{k=1}^l \phi_k(g_0^{-1} g) h_k(x)$, so that 
$$\mathcal{R}(g_0 \cdot f) = \sum_{k=1}^l \mathcal{R}_G\left(g_0 \cdot \phi_k\right) h_k = \sum_{k=1}^l \mathcal{R}_G\left(\phi_k\right) h_k = \mathcal{R}(f).$$ Furthermore, if $f \in \mathcal{P}(X)^G$, we have $f \circ \alpha = \overline{1} \otimes f$ and then $\mathcal{R}(f) = \mathcal{R}_G\left(\overline{1}\right)f = f$. As a consequence, $\mathcal{R}$ is a Reynolds operator of $X$ and $G$ is therefore linearly reductive.

Conversely, if $G$ is linearly reductive, there exists, in particular, a Reynolds operator for the $G$-real algebraic set $G$ (which is unique by remark \ref{remcharaclinearredprop} 3).
\end{proof}

\begin{rem} \label{remexpresreynoldsopg} Suppose that $G$ is linearly reductive and let $X$ be a $G$-real algebraic set. Keeping the notations of the previous proof, we have, for all $x \in X$,
$$\mathcal{R}_X(f)(x) =  \sum_{k=1}^l \mathcal{R}_G\left(\phi_k\right) h_k(x) = \mathcal{R}_G\left(\sum_{k=1}^l h_k(x) \phi_k\right) = \mathcal{R}_G\left(f \circ \alpha(\cdot,x)\right).$$
\end{rem}

The characterization of corollary \ref{corequivlinearlyredreynoldsopgroup} allows us to exhibit the following first classes of examples of linearly reductive real algebraic groups:

\begin{exs} \label{exlinearlyredrealalggr}
~
\begin{enumerate}
	\item If $G$ is finite of cardinality $N$, we can consider the Reynolds operator
$$\begin{array}{ccc}\mathcal{P}(G) & \rightarrow & \R\\ \phi & \mapsto & \frac{1}{N}\sum_{g \in G} \phi(g) \end{array}$$
of $G$, which is therefore linearly reductive.
	\item More generally, suppose that $G \subset \R^m$ is compact i.e. closed and bounded. Since $G$ is then a \emph{compact} Lie group, there is a unique normalized Haar measure $\mu$ on $G$ and we can consider the associated integral on $G$ (see for instance section \textbf{8}.1.2 of \cite{Pro}). Any polynomial function on $G$ is continuous, thus we can define the $\R$-linear map 
	$$\mathcal{R}_G : \begin{array}{ccc}\mathcal{P}(G) & \rightarrow & \R \\ f & \mapsto & \int_G f {\rm d}\mu. \end{array}$$ 
Since $\int_G {\rm d}\mu = \mu(G) = 1$ (because $\mu$ is normalized) and, for all $f \in \mathcal{P}(G)$ and $g_0 \in G$, $\int_G f(g_0 g) {\rm d}\mu(g) = \int_G f(g) {\rm d}\mu(g)$ (by left-invariance of the Haar measure $\mu$ and the associated integral), $\mathcal{R}_G$ is a Reynolds operator for $G$ (see the beginning of the proof of corollary \ref{corequivlinearlyredreynoldsopgroup}). As a consequence, any compact real algebraic group is linearly reductive (for instance, if $n \in \Nstar$, the orthogonal group $\OO_n(\R)$ as well as the special orthogonal group $\SO_n(\R)$ are linearly reductive real algebraic groups). From section \ref{sectpolyactionscomprealalggroups}, we will focus on compact real algebraic groups and their polynomial actions on real algebraic sets and semialgebraic sets.
	\item Consider the real algebraic group $K := \{(x,y) \in \R^2~|~xy=1\}$ (example \ref{exspolyrealalggroups} 3) and denote by $I$ the ideal of $\R[x,y]$ generated by the irreducible polynomial $xy-1$. Since the sign of $xy-1$ changes on $\R^2$, by Theorem 4.5.1 of \cite{BCR}, the ideal $I$ is real (Definition~4.1.3 of \cite{BCR}) and then, by the Real Nullstellensatz (Theorem 4.1.4 of \cite{BCR}), $I(K) = I$, so that $\mathcal{P}(K) = \R[x,y]/I$. Now, consider the $\R$-linear map on $\R[x,y]$ which to any monomial~$x^{\alpha} y^{\beta}$ with $\alpha,\beta \in \N$ associates $\delta_{\alpha \beta}$: it induces an $\R$-linear map $\mathcal{R} : \mathcal{P}(K) \rightarrow \R$ ($I$ is included in the kernel of the former map). We have $\mathcal{R}\left(\overline{1}\right) = 1$ and, for any $f \in  \mathcal{P}(K)$ and $(x_0,y_0) \in K$, $\mathcal{R}((x_0,y_0) \cdot f) = \mathcal{R}(f)$ (because $x_0 y_0 = 1$): $\mathcal{R}$ is therefore a Reynolds operator for the real algebraic group $K$, which is then linearly reductive by corollary~\ref{corequivlinearlyredreynoldsopgroup}.
	\item Suppose that $G$ is linearly reductive and let $H$ be a linearly reductive real algebraic group. Then the linear map $\mathcal{R} : \mathcal{P}(G \times H) = \mathcal{P}(G) \otimes_{\R} \mathcal{P}(H) \rightarrow \R$ defined by, if $f \in \mathcal{P}(G)$ and~$h \in \mathcal{P}(H)$, $\mathcal{R}(f \otimes h) := \mathcal{R}_G(f) \mathcal{R}_H(h)$ is a Reynolds operator for $G \times H$, and the real algebraic group $G \times H$ is therefore linearly reductive as well.
\end{enumerate}
\end{exs}

\subsection{Hilbert's finiteness theorem}

We have been considering linearly reductive real algebraic groups in the previous subsection because we have in view Hilbert's finiteness theorem. This result asserts the finite generation of the $\R$-algebra of invariant polynomial functions of any real algebraic set polynomially acted by a linearly reductive real algebraic group. This algebra will then correspond to a real algebraic set that we will call the real \emph{geometric quotient} in the next section. 

For the proof of Hilbert's finiteness theorem, we follow the argumentation of \cite{DK} (Corollary~2.2.11, Theorem 2.2.10 and Corollary 2.2.9).

\begin{theo}{(Hilbert's finiteness theorem)} \label{theohilbertfiniteness}
Let $G$ be a linearly reductive real algebraic group and let $X$ be a $G$-real algebraic set. The $\R$-algebra $\mathcal{P}(X)^G$ is finitely generated.
\end{theo} 

\begin{proof} Using proposition \ref{proplinearizgrealalgset}, we can consider an isomorphism of $G$-real algebraic sets $i : X \rightarrow Y$ onto a $G$-algebraic subset $Y$ of a polynomial representation $\R^N$ of $G$. Denote by~$\iota$ the induced morphism of $\R$-algebras $\mathcal{P}(\R^N) = \R[x_1,\ldots,x_N] \rightarrow \mathcal{P}(X)~;~f \mapsto f_{|Y} \circ i$: $\iota$ is~$G$-equivariant with respect to the induced actions of $G$. For the sake of readability, denote $\mathbb{R}[\underline{x}] := \R[x_1,\ldots,x_N]$ in the following.

Now, let $h \in \mathcal{P}(X)^G$ and let $f \in \mathbb{R}[\underline{x}]$ such that $h \circ i^{-1} = f_{|Y} = \overline{f}$ in $\mathcal{P}(Y)$. Let then~$W$ be the $\R$-vector subspace of $\mathcal{P}(X)$ generated by $h$, and $V$ be the subspace of $\mathbb{R}[\underline{x}]$ generated by the orbit $\{ g \cdot f~|~f \in G\}$ of $f$: the restriction $V \rightarrow W$ of $\iota$ is a surjective $G$-equivariant linear map between polynomial representations of $G$. Therefore, since $G$ is linearly reductive, by characterization~2 of proposition \ref{propequivlinred}, there exists $\widetilde{f} \in V^G \subset \mathbb{R}[\underline{x}]^G$ such that $\iota\left(\widetilde{f}\right) = h$. As a consequence, $\mathcal{P}(X)^G = \iota\left(\mathbb{R}[\underline{x}]^G\right)$. It is thus sufficient to prove that $\mathbb{R}[\underline{x}]^G$ is a finitely generated $\R$-algebra to assert that so is $\mathcal{P}(X)^G$.   

Let $I$ be the ideal of the ring $\mathbb{R}[\underline{x}]$ generated by the non-constant homogeneous $G$-invariant polynomials of $\mathbb{R}[\underline{x}]$. Since the latter is a Noetherian ring, the ideal $I$ is generated by a finite number of polynomials and then, by definition of $I$, by a finite number of non-constant homogeneous $G$-invariant polynomials $f_1,\ldots,f_r$. The rest of the proof will consist in showing that the $\R$-algebra $\mathbb{R}[\underline{x}]^G$ is generated by the polynomials $f_1,\ldots,f_r$ i.e. any $G$-invariant polynomial of $\mathbb{R}[\underline{x}]$ is a polynomial in $f_1,\ldots,f_r$. Since any polynomial is a linear combination of homogeneous polynomials and since the action of $G$ on $\R^N$ is linear and then preserves, if $d \in \mathbb{N}$, the subspace $\mathbb{R}[\underline{x}]_d$ of $\mathbb{R}[\underline{x}]$ of homogeneous polynomials of degree $d$, it suffices to prove the desired property for $\mathbb{R}[\underline{x}]^G_d$.

Reason by induction on $d$. The property being clear for any constant polynomial, let $f$ be a homogeneous $G$-invariant polynomial of $\mathbb{R}[\underline{x}]$ of degree $d \in \mathbb{N} \setminus \{0\}$: in particular $f \in I$, so that there exist $h_1,\ldots,h_r \in \mathbb{R}[\underline{x}]$ such that $f = \sum_{k=1}^r h_k f_k$. Since $f$ is homogeneous of degree $d$, we can furthermore write $f = \sum_{k=1}^r \widetilde{h}_k f_k$ where, for $k \in \{1,\ldots,r\}$, $\widetilde{h}_k$ denotes the homogeneous part of degree $d - {\rm deg} f_k$ of $h_k$ when the latter quantity in non-negative, and $0$ otherwise. Then apply the Reynolds operator $\mathcal{R} := \mathcal{R}_{\R^N} : \mathcal{P}(\R^N) \rightarrow \mathcal{P}(\R^N)^G$ to the latter equality to obtain 
$$f = \mathcal{R}(f) =  \sum_{k=1}^r \mathcal{R}\left(\widetilde{h}_k f_k\right) = \sum_{k=1}^r \mathcal{R}\left(\widetilde{h}_k\right) f_k$$
(use in particular remark \ref{remcharaclinearredprop} 3). Finally, because, for all $k \in \{1,\ldots,r\}$, $\widetilde{h}_k \in \mathbb{R}[\underline{x}]_{d-{\rm deg} f_k}$ and then $ \mathcal{R}\left(\widetilde{h}_k\right) \in \mathbb{R}[\underline{x}]_{d-{\rm deg} f_k}^G$ (see remark \ref{remcharaclinearredprop} 4), we can apply the induction hypothesis to assert that all polynomials $\mathcal{R}\left(\widetilde{h}_k\right)$, $k \in \{1,\ldots,r\}$, are polynomials in $f_1,\ldots,f_r$, hence so is~$f$. 
\end{proof}

\begin{exs} \label{exgeninvariantalgorthogroup}
~
	\begin{enumerate}
		\item For all $n \in \mathbb{N} \setminus \{0\}$, if we consider the usual action of $\OO_n(\R)$ on $\R^n$ (see example \ref{exgrealalgset}~1), the algebras $\mathcal{P}(\R^n)^{\OO_n(\R)} = \R[x_1,\ldots,x_n]^{\OO_n(\R)}$ and $\mathcal{P}(\R^n)^{\SO_n(\R)}= \R[x_1,\ldots,x_n]^{\SO_n(\R)}$ are both generated by the polynomial $\sum_{k=1}^n x_k^2$ (this is a particular case of the First Fondamental Theorem for orthogonal and special orthogonal groups: see for instance~\cite{KP}~10.2).
		\item Let $n \in \mathbb{N} \setminus \{0\}$ and consider the polynomial action $K \times \R^n \rightarrow \R^n~;~(a,b),(x_1,\ldots,x_n) \mapsto (ax_1,\ldots,ax_n)$ of the linearly reductive real algebraic group $K = \{(a,b) \in \R^2~|~ab=1\}$ (example \ref{exlinearlyredrealalggr} 3) on $\R^n$. Let $f \in \mathcal{P}(\R^n)^{K}$ then, for all $(x_1,\ldots,x_n) \in \R^n$ and $m \in \Nstar$, we have $f(x_1,\ldots,x_n) = f\left(\frac{x_1}{m},\ldots, \frac{x_n}{m}\right)$ so that, by continuity of $f$, $f(x_1,\ldots,x_n) = f(0,\ldots,0)$ for all $(x_1,\ldots,x_n) \in \R^n$, i.e. $f = \overline{f(0,\ldots,0)}$. As a consequence, $\mathcal{P}(\R^n)^{K} = \R$ is generated by the constant function $\overline{1}$.
	\end{enumerate}
\end{exs}

\begin{rem} \label{reminvariantalgebraproduct} Let $G$ and $H$ be linearly reductive real algebraic groups. According to example~\ref{exlinearlyredrealalggr} 4, the direct product $G \times H$ is linearly reductive as well. In particular, by Hilbert's finiteness theorem \ref{theohilbertfiniteness}, if $X$ is a $G$-real algebraic set and $Y$ is a $H$-real algebraic set, then the $\R$-algebra $\mathcal{P}(X \times Y)^{G \times H}$ is finitely generated. But remark furthermore that we have the equality 
$$\mathcal{P}(X \times Y)^{G \times H} = \mathcal{P}(X)^G \otimes_{\R} \mathcal{P}(Y)^H,$$
which allows to form generators of the left-hand side algebra with the chosen respective generators of the right-hand side algebras. Let us show the above equality: let $f \in \mathcal{P}(X \times Y)^{G \times H}$ and write $f \circ \alpha = \sum_{k=1}^l f_k \otimes h_k$ with $f_1,\ldots,f_l \in \mathcal{P}(X)$ and $h_1,\ldots,h_l \in \mathcal{P}(Y)$. If $\alpha : G \times X \rightarrow X$ and $\alpha' : H \times Y \rightarrow Y$ denote the respective polynomial actions of $G$ and $H$ on $X$ and $Y$, we have, for any $(x,y) \in X \times Y$,
\begin{eqnarray*}
f(x,y) & = & \mathcal{R}_{X \times Y}(f)(x,y) \\
	& = & \mathcal{R}_{G\times H}\left(\sum_{k=1}^l f_k \circ \alpha( \cdot, x) \otimes h_k \circ \alpha'(\cdot,y)\right) \mbox{ (remark \ref{remexpresreynoldsopg})}\\
	& = & \sum_{k=1}^l \mathcal{R}_G(f_k \circ \alpha(\cdot,x)) \mathcal{R}_H(h_k \circ \alpha'(\cdot,y)) \mbox{ (example \ref{exlinearlyredrealalggr} 4)}\\
	& = & \sum_{k=1}^l \mathcal{R}_X(f_k)(x) \mathcal{R}_Y(h_k)(y) \mbox{ (remark \ref{remexpresreynoldsopg}),}
\end{eqnarray*}
so that $f \in \mathcal{P}(X)^G \otimes_{\R} \mathcal{P}(Y)^H$.
\end{rem}

\subsection{Complexification} \label{subsectcomplexificationlinearlyreductive}

We end this part by studying the complexification of the notions and results that we addressed in the above subsections: this will reveal crucial in the following to understand the complexification of the real geometric quotient of $G$-real algebraic sets when $G$ is a linearly reductive real algebraic group. 

First, the notion of linear reductivity (definition \ref{defreallinearreductivity}) has its complex counterpart for complex algebraic groups and the complex analogs of equivalences of proposition \ref{propequivlinred} and corollary \ref{corequivlinearlyredreynoldsopgroup} remain valid, as well as remarks \ref{remcharaclinearredprop}. 

Now, let $G$ be any real algebraic group. Have in mind the following fact:

\begin{lem} \label{lemidpolyfunctoncomplexification} Let $X \subset \R^n$ be a real algebraic set and denote by $\mathcal{P}(X_{\C})$ the $\C$-algebra $\C[x_1,\ldots,x_n]/\mathsf{I}(X_{\C})$ of complex polynomial functions on $X_{\C}$. The morphism $\Psi$ of $\C$-algebras which associates to any polynomial $p+iq \in \C[x_1,\ldots,x_n]$, with $p,q \in \R[x_1,\ldots,x_n]$, the element $\overline{p}\otimes 1 + \overline{q} \otimes i$ of $\mathcal{P}(X) \otimes_{\R} \C$ induces an isomorphism $\mathcal{P}(X_{\C}) \rightarrow \mathcal{P}(X) \otimes_{\R} \C$. In particular, any complex polynomial function $F$ on $X_{\C}$ writes uniquely as $f_{\C} + i h_{\C}$ where $f$ and $h$ are real polynomial functions on $X$ (below, we will respectively call $f$ and $h$ the \emph{real} and \emph{imaginary parts} of $F$).
\end{lem}

\begin{proof} The morphism $\Psi$ is surjective and, if $p+iq \in \C[x_1,\ldots,x_n]$ with $p,q \in \R[x_1,\ldots,x_n]$, we have $\Psi(p+iq) = 0$ if and only if $\overline{p} = \overline{q} = \overline{0}$ in $\mathcal{P}(X)$ if and only if $p,q \in I(X)$ if and only if~$p+iq \in \mathsf{I}(X) = \mathsf{I}(X_{\C})$.
\end{proof} 

\begin{prop} \label{propequivlinearredragandcomplexification} The real algebraic group $G$ is linearly reductive if and only if the complex algebraic group $G_{\C}$ is linearly reductive.
\end{prop}

\begin{proof} For the direct implication, we use the characterization 2 of (the complex counterpart of) proposition \ref{propequivlinred}: suppose that $G$ is linearly reductive and consider two complex polynomial representations $V$ and $W$ of $G_{\C}$ as well as a surjective $G_{\C}$-equivariant linear map $\omega : V \rightarrow W$. The sets $V$ and $W$ are in particular $\R$-vector spaces and $\omega : V \rightarrow W$ is a $\R$-linear map as well; if we consider respective $\C$-bases $(v_1,\ldots,v_n)$ and $(w_1,\ldots,w_N)$ for $V$ and $W$, then the tuples $(v_1,i v_1,\ldots,v_n,iv_n)$ and $(w_1,i w_1,\ldots,w_N,i w_N)$ are respective $\R$-bases for $V$ and $W$ and the matrix of the $\R$-linear map $\omega$ in the latter bases is given by the real and imaginary parts of the coefficients of the matrix of the $\C$-map $\omega$ in the former bases. Furthermore, restricting the considered linear actions on $V$ and $W$ to $G$, the finite-dimensional $\R$-vector spaces $V$ and $W$ are real polynomial representations of $G$ (the entries of the matrix of the considered action of $G$ are given by the real and imaginary parts of the polynomial entries of the action of $G_{\C}$, using the identification $\mathcal{P}({G}_{\C}) \cong \mathcal{P}(G) \otimes_{\R} \C$ of lemma \ref{lemidpolyfunctoncomplexification}) and $\omega$ is a surjective $G$-equivariant linear map between polynomial representations of $G$. Now, let $w \in W^{G_{\C}}$: $w$ is in particular a invariant vector of $W^G$ and then, since the real algebraic group $G$ is linearly reductive, there exists $v \in V^G$ such that $w = \omega(v)$. But the equalities $g \cdot v = v$, $g \in G$, are given by polynomial equalities (with coefficients in $\C)$ over the real algebraic set $G$ and are then true over its complexification $G_{\C}$. As a consequence, $v \in V^{G_{\C}}$ so that $w \in \omega(V^{G_{\C}})$. The complex algebraic group $G_{\C}$ is therefore linearly reductive.

As for the converse implication, suppose that $G_{\C}$ is linearly reductive and consider the associated Reynolds operator $\mathcal{R}_{G_{\C}} : \mathcal{P}(G_{\C}) \rightarrow \C$ (cf. the complex counterpart of corollary~\ref{corequivlinearlyredreynoldsopgroup}). The linear map $f \in \mathcal{P}(G) \mapsto {\rm Re}\left(\mathcal{R}_{G_{\C}}(f_{\C})\right) \in \R$ (where, if $z \in \C$, ${\rm Re}(z)$ denotes the real part of $z$) is then a Reynolds operator for the real algebraic group $G$, which is therefore linearly reductive.
\end{proof}

\begin{exs} 
~
	\begin{enumerate}
		\item The complexification of any compact real algebraic group is linearly reductive (see example \ref{exlinearlyredrealalggr} 2).
		\item Let $n \in \mathbb{N} \setminus \{0\}$. The real algebraic groups ${\rm SL}_n(\R)$ and $\Omega_n(\R)$ are linearly reductive since their respective complexifications ${\rm SL}_n(\C)$ and $\Omega_n(\C)$ (cf. example \ref{excomplrealalggroup} 2) are linearly reductive (use for instance \cite{Pro} Chapter 7 paragraph 3.2 Theorem).
	\end{enumerate} 
\end{exs}

Let $X \subset \R^n$ be a $G$-real algebraic set and consider the complex algebraic set $X_{\C}$: $X_{\C}$ is a~$G_{\C}$-complex algebraic set (lemma \ref{lemcomplexificationpolyaction}). Suppose that $G$ is linearly reductive: the complex algebraic group $G_{\C}$ is therefore linearly reductive. By the complex counterpart of Hilbert's finiteness theorem \ref{theohilbertfiniteness} (the proof is identical), the $\C$-algebra $\mathcal{P}(X_{\C})^{G_{\C}}$ is finitely generated. Furthermore, we can choose as generators of the latter the complexifications of the generators of the finitely generated $\R$-algebra $\mathcal{P}(X)^G$, by virtue of the next lemma. This will imply, in the next section, that the complex geometric quotient of $X_{\C}$ by $G_{\C}$ is the complexification of the real geometric quotient of $X$ by $G$. 

\begin{lem} \label{lemidpolyfunctoncomplexificationrestrictiontoinv} The isomorphism of $\C$-algebras $\mathcal{P}(X_{\C}) \rightarrow \mathcal{P}(X) \otimes_{\R} \C$ of lemma \ref{lemidpolyfunctoncomplexification} restricts to an isomorphism $\mathcal{P}(X_{\C})^{G_{\C}} \rightarrow \mathcal{P}(X)^G \otimes_{\R} \C$.
\end{lem}

\begin{proof} Let $\overline{p + iq} \in \mathcal{P}(X_{\C})^{G_{\C}}$ with $p,q \in \mathbb{R}[x_1,\ldots,x_n]$, and denote by $\alpha : G \times X \rightarrow X$ the polynomial action of $G$ on $X$. For all $g \in G$ and all $x \in X$, we have $(p+iq)(\alpha(g,x)) = (p+iq)(x)$ i.e. $p(\alpha(g,x)) = p(x)$ and $q(\alpha(g,x)) = q(x)$. As a consequence, the real polynomial functions~$\overline{p}, \overline{q} \in \mathcal{P}(X)$ on $X$ are $G$-invariant.

On the other hand, let $\overline{p}, \overline{q} \in \mathcal{P}(X)^G$, with $p,q \in \mathbb{R}[x_1,\ldots,x_n]$, and consider the element~$\overline{p+iq}$ of $\mathcal{P}(X_{\C})$. We have the polynomial equality $(p+iq)(\widetilde{\alpha}_{\C}(\cdot,*)) = (p+iq)(*)$ over~$G \times X$ (we use the notations of the beginning of subsection \ref{subseccomplexificationralgps}) hence over $G_{\C} \times X_{\C}$, so that $\overline{p + iq} \in \mathcal{P}(X_{\C})^{G_{\C}}$.
\end{proof}

As a consequence, if $f_1,\ldots,f_r$ are invariant polynomial functions on $X$ generating the $\R$-algebra $\mathcal{P}(X)^G$, the extensions $(f_1)_{\C},\ldots, (f_r)_{\C}$ generate the $\C$-algebra $\mathcal{P}(X_{\C})^{G_{\C}}$: if $F$ is a $G_{\C}$-invariant complex polynomial function on $X_{\C}$, the real and imaginary parts of $F$ are in $\mathcal{P}(X)^G$ and, if $f \in \mathcal{P}(X)$ satisfies $f = p(f_1,\ldots,f_r)$ with $p \in \R[z_1,\ldots,z_r]$, then $f_{\C} = p\big((f_1)_{\C},\ldots, (f_r)_{\C}\big)$ in $\mathcal{P}(X_{\C})$.

\begin{rem} The property stated in lemma \ref{lemidpolyfunctoncomplexificationrestrictiontoinv} remains true even if we do not suppose that~$G$ is linearly reductive.
\end{rem}

\section{Real geometric quotient} \label{sectionrealgeoquotient}

Let $G$ be a \emph{linearly reductive} real algebraic group and let $X \subset \R^n$ be a $G$-real algebraic set. By Hilbert's finiteness theorem \ref{theohilbertfiniteness}, the $\R$-algebra $\mathcal{P}(X)^G$ is finitely generated. Consequently, as we will see below, the inclusion $\mathcal{P}(X)^G \hookrightarrow \mathcal{P}(X)$ corresponds, by Real Nullstellensatz (Theorem~4.1.4 of \cite{BCR}), to a $G$-invariant polynomial map from $X$ to a real algebraic set $X \gq G$: we will call this map the \emph{geometric quotient} of $X$ by $G$.

In this section, we will detail this construction of the real geometric quotient and we will study its functoriality. Also, we will emphasize that in general, contrary to the complex case, the real geometric quotient map is not surjective and its image is not a real algebraic set.

\subsection{Construction of the real geometric quotient via the Real Nullstellensatz} \label{subsectioncontrrealgeomquotient}

As we said above, we will combine Hilbert's finiteness theorem with Real Nullstellensatz in order to construct the real geometric quotient of real algebraic sets polynomially acted by a linearly reductive real algebraic group. The consequence of Real Nullstellensatz we will state below asserts the contrafunctorial correspondence between polynomial maps between real algebraic sets and morphisms of \emph{real} finitely generated $\R$-algebras. 

Let $A$ be a commutative ring with zero element $0_A$. We say that $A$ is {\it real} if for any $a_1,\ldots,a_r \in A$, if $a_1^2 + \cdots + a_r^2 = 0_A$ then $a_1 = \cdots = a_r = 0_A$ (remark in particular that any subring of a real ring is itself real). Recall that an ideal $I$ of $A$ is called {\it real} if for any $a_1,\ldots,a_r \in A$, if $a_1^2 + \cdots + a_r^2 \in I$ then $a_1, \ldots, a_r \in I$: the ideal $I$ is real if and only if the quotient ring $A/I$ is real. 

\begin{theo} \label{theorealnullstellensatzequivcategories}
The operation which associates to any real algebraic set $Y$ its $\R$-algebra $\mathcal{P}(Y)$ of polynomial functions on $Y$ gives an equivalence of categories between the category of (polynomial isomorphism classes of) real algebraic sets and polynomial maps and the category of (isomorphism classes of) real (as commutative rings) finitely generated $\R$-algebras and morphisms of $\R$-algebras.
\end{theo}

\begin{proof} If $Y \subset \R^n$ is a real algebraic set, the $\R$-algebra $\mathcal{P}(Y) = \R[x_1,\ldots,x_n]/I(Y)$, generated by the classes $\overline{x_1},\ldots,\overline{x_d}$, is real as a commutative ring: if $f_1,\ldots,f_r \in \mathcal{P}(Y)$ satisfy $f_1(y)^2 + \cdots + f_r(y)^2 = 0$ for all $y \in Y$, then all the functions $f_1,\ldots,f_r$ are the zero function on $Y$. Furthermore, any polynomial map $F : Y \subset \R^n \rightarrow Z \subset \R^N$ induces the morphism of $\R$-algebras $F^* : \mathcal{P}(Z) \rightarrow \mathcal{P}(Y)~;~f \mapsto f \circ F$ and this operation is contrafunctorial.

On the other hand, if $A$ is a real $\R$-algebra generated by elements $\alpha_1,\ldots,\alpha_r$, consider the ideal $I := \left\{p \in \R[x_1,\ldots,x_r]~|~p(\alpha_1,\ldots,\alpha_r) = 0_A\right\}$ of $\R[x_1,\ldots,x_r]$: $I$ is the kernel of the morphism $\Phi : \R[x_1,\ldots,x_r] \rightarrow A$ of $\R$-algebras which associates to any indeterminate $x_k$, $k \in \{1,\ldots,r\}$, the element $\alpha_k$ of $A$. Since the elements $\alpha_1,\ldots,\alpha_r$ generate the $\R$-algebra $A$, the morphism $\Phi$ is surjective and induces an isomorphism of $\R$-algebras $\R[x_1,\ldots,x_r]/I \rightarrow A$. In particular, the ideal $I$ is real so that, by the Real Nullstellensatz (Theorem 4.1.4 of \cite{BCR}), $I(V(I)) = I$. Therefore, if we denote by $X_A$ the real algebraic subset $V(I)$ of $\R^r$, the $\R$-algebra $\mathcal{P}(X_A)$ is isomorphic to $A$. Notice that the real algebraic set $X_A$ is well-defined up to polynomial isomorphism: if $\beta_1,\ldots,\beta_s$ is another system of generators for the $\R$-algebra~$A$ and $Y$ is the associated real algebraic subset of $\R^s$, the polynomial map $X_A \rightarrow Y$ given by polynomials $p_1,\ldots,p_s \in \R[x_1,\ldots,x_r]$ such that, for $l \in \{1,\ldots,s\}$, $\beta_l = p_l(\alpha_1,\ldots,\alpha_r)$, is a polynomial isomorphism whose inverse map is given by polynomials $q_1,\ldots, q_r \in \R[y_1,\ldots,y_s]$ such that, for $k \in \{1,\ldots,r\}$, $\alpha_k = q_k(\beta_1,\ldots,\beta_s)$.

Finally, a morphism of $\R$-algebras $\Psi : A \rightarrow B$ between real $\R$-algebras $A$ and $B$ respectively generated by elements $\alpha_1,\ldots,\alpha_r$ and $\beta_1,\ldots,\beta_s$ induces, in a contrafunctorial way, a polynomial map $\Psi^* : X_B \subset \R^s\rightarrow X_A \subset \R^r$ given by polynomials $q_1,\ldots, q_r \in \R[y_1,\ldots,y_s]$ such that, for $k \in \{1,\ldots,r\}$, $\Psi(\alpha_k) = q_k(\beta_1,\ldots,\beta_s)$.
\end{proof}

\begin{rem} Via the above equivalence of categories, the tensor product over $\R$ of real finitely generated $\R$-algebras corresponds to the Cartesian product of real algebraic sets: in particular, the tensor product over $\R$ of two real finitely generated $\R$-algebras $A$ and $B$ is a real finitely generated $\R$-algebra. Indeed, there exist real algebraic sets $V$ and $W$ such that $A$, resp. $B$, is (isomorphic to) $\mathcal{P}(V)$, resp. $\mathcal{P}(W)$, so that the tensor product $A \otimes_{\R} B$ is (isomorphic to) $\mathcal{P}(V) \otimes_{\R} \mathcal{P}(W) = \mathcal{P}(V \times W)$ which is a real finitely generated $\R$-algebra.
\end{rem}

Keep the notations of the beginning of the section. We apply theorem \ref{theorealnullstellensatzequivcategories} and its proof to the inclusion morphism of real finitely generated $\R$-algebras $\mathcal{P}(X)^G \hookrightarrow \mathcal{P}(X)$ (the real algebraic group $G$ is linearly reductive so that, by Hilbert's finiteness theorem \ref{theohilbertfiniteness}, the former~$\R$-algebra is finitely generated as well): it corresponds to the restriction $\pi_X : X \rightarrow X_{\mathcal{P}(X)^G}$ of the polynomial map $x \in \R^n \mapsto (f_1(x),\ldots,f_r(x)) \in \R^r$, where $f_1,\ldots,f_r$ are invariant polynomial functions on $X$ generating the $\R$-algebra $\mathcal{P}(X)^G$ and
$$X_{\mathcal{P}(X)^G} = V \big(\{ p \in \R[z_1,\ldots,z_r]~|~p(f_1,\ldots,f_r) = 0\}\big).$$

Notice that the map $\pi_X$ is $G$-invariant: if $x \in X$ and $g \in G$, we have 
$$\pi_X(g \cdot x) = (f_1(g \cdot x),\ldots,f_r(g \cdot x)) = (f_1(x),\ldots,f_r(x)) = \pi_X(x).$$  

\begin{de} \label{derealgeoquotient} We denote $X \gq G := X_{\mathcal{P}(X)^G}$ and we call the polynomial map $\pi_X : X \rightarrow X \gq G$ the \emph{geometric quotient of $X$ by $G$}. By abuse of terminology, the real algebraic set $X \gq G$ will also be called the \emph{geometric quotient of $X$ by $G$}.
\end{de}

\begin{rem} \label{remgequowduptopoliso} The geometric quotient is well-defined up to polynomial isomorphism.
\end{rem}

As we already introduced, contrary to the complex case (see for instance \cite{DK} Lemma 2.3.2), the geometric quotient map $\pi_X$ is not surjective, as illustrated by the following examples:

\begin{exs} \label{exrealalgquotientplaneorthoaction} 
~
	\begin{enumerate}
		\item Let $n \in \Nstar$ and consider the usual action of $\OO_n(\R)$ on $\R^n$. Since the $\R$-algebra $\mathcal{P}(\R^n)^{\OO_n(\R)}$ is generated by the polynomial $x_1^2+\cdots+x_n^2$ of $\R[x_1,\ldots,x_n]$ (see example \ref{exgeninvariantalgorthogroup}), we have
		$$\R^n\gq \OO_n(\R) = V\big(\{p \in \R[z]~|~\forall a_1,\ldots,a_n \in \R,~p(a_1^2+\cdots + a_n^2) = 0\}\big) = V(\{0\}) = \R$$
and the geometric quotient of $\R$ by $\OO_n(\R)$ is the map 
$$\R^n  \rightarrow \R~;~(x_1,\ldots,x_n) \mapsto x_1^2+\cdots+x_n^2,$$
whose image is the {\it semialgebraic} half-line $\{z \in \R~|~z \geq 0\}$.
		\item Consider the polynomial action $S^1 \times \R^3~;~\big( (a,b), (x,y,z)\big) \mapsto (ax-by,bx + ay, z)$ of $S^1$ on $\R^3$. Notice that a polynomial function on $\R^3$ is invariant with respect to the latter action if and only it is invariant with respect to the action of $\SO_2(\R) \times \{0\}$ induced by the usual action of $\SO_2(\R)$ on $\R^2$ and the (trivial) action of the trivial group $\{0\}$ on $\R$. Therefore,
	$$\mathcal{P}(\R^3)^{S^1} = \mathcal{P}(\R^2 \times \R)^{\SO_2(\R) \times \{0\}} = \mathcal{P}(\R^2)^{\SO_2(\R)} \otimes_{\R} \mathcal{P}(\R)^{\{0\}}$$
by remark \ref{reminvariantalgebraproduct}, so that the $\R$-algebra $\mathcal{P}(\R^3)^{S^1}$ is generated by the polynomials $x^2 + y^2$ and $z$ of $\R[x,y,z]$ (see example \ref{exgeninvariantalgorthogroup} 1). As a consequence, we have
$$\R^3\gq S^1 = V\left(\{p \in \R[u,v]~|~\forall a,b,c \in \R,~p(a^2+b^2,c) = 0\}\right) = V(\{0\}) = \R^2$$
and the geometric quotient of $\R^3$ by $S^1$ is the map
$$\R^3 \rightarrow \R^2~;~(x,y,z) \mapsto (x^2+y^2,z),$$
whose image is the {\it semialgebraic} half-plane $\{(u,v) \in \R^2~|~u \geq 0\}$.
	\end{enumerate}
\end{exs}

Actually, the image of any real geometric quotient by a linearly reductive real algebraic group is always a semialgebraic set, since the image of a polynomial map between real algebraic sets is semialgebraic by Tarski-Seidenberg theorem (more precisely by its corollary Proposition~2.2.7 of \cite{BCR}).

\begin{de} Keeping the above notations, we denote by $X/G$ the semialgebraic subset~$\pi_X(X)$ of $\R^r$ and by $\varpi_X$ the restriction $X \rightarrow X/G$ of $\pi_X$. The map $\varpi_X : X \rightarrow X/G$ will be called the \emph{semialgebraic quotient of $X$ by $G$}. By abuse of terminology, the semialgebraic set~$X/G$ will also be called the semialgebraic quotient of $X$ by $G$.  
\end{de}

\begin{lem} \label{lemzariskiclossemialgquotisgeomquot} The Zariski closure of $X/G$ in $\R^r$ is $X \gq G$. In particular, $\pi_X$ is surjective if and only if $X/G$ is a real algebraic set.
\end{lem}

\begin{proof} The Zariski closure of $X/G$ in $\R^r$ is 
$$V\big(I(\pi_X(X))\big) = V\big(\{p \in \R[z_1,\ldots,z_r]~|~\forall x \in X,~p(f_1(x),\ldots,f_r(x)) = 0\}\big) = X /\!\!/ G.$$
\end{proof}

We can more generally define the semialgebraic quotient of any $G$-semialgebraic set: let $S \subset \R^n$ be a $G$-semialgebraic set and denote by $Z$ the Zariski closure of $S$ in $\R^n$. Since, by lemma \ref{lemzariskiclosgsaisgras}, $Z$ is a $G$-real algebraic set, we can consider its geometric quotient $\pi_Z : Z \rightarrow Z \gq G \subset \R^s$. The image $\pi_Z(S)$, as being the image of the semialgebraic set $S$ by a polynomial map, is then a semialgebraic set, again by Tarski-Seidenberg theorem (Proposition 2.2.7 of \cite{BCR}):

\begin{de} \label{defsemialgquot} We denote $S/G := \pi_Z(S)$ and call the restriction $\varpi_S : S \rightarrow S/G$ of $\pi_Z$ the \emph{semialgebraic quotient} of $S$ by $G$. 
\end{de}

\begin{lem} The Zariski closure of $S/G$ in $\R^r$ is $Z /\!\!/ G$.
\end{lem}

\begin{proof} If $W$ is any real algebraic subset of $\R^s$ such that $S/G \subset W \subset Z /\!\!/ G$, then $\pi_Z^{-1}(W)$ is a $G$-real algebraic subset of $Z$ containing $\pi_Z^{-1}(S/G) \supset S$, thus $\pi_Z^{-1}(W) = Z$ (because $Z$ is the Zariski closure of $S$ in $\R^n$). Therefore $\pi_Z(Z) = \pi_Z\left(\pi_Z^{-1}(W)\right) \subset W$, so that $Z\gq G = \overline{\pi_Z(Z)}^{\mathcal{Z}} \subset W$.
\end{proof}

We will see in next subsection that the semialgebraic quotient $\varpi_S$ can actually be considered as the restriction of the geometric quotient of any $G$-real algebraic set of $\R^n$ whose $S$ is a $G$-stable subset. 

\begin{rems} 
~
	\begin{enumerate} 
		\item The semialgebraic quotient is well-defined up to polynomial isomorphism (i.e. bijective restriction of a polynomial map with inverse being the restriction of a polynomial map): see above remark \ref{remgequowduptopoliso}.
		\item If $G$ is a \emph{compact} real algebraic group (example \ref{exlinearlyredrealalggr} 2) and if $\R^n$ is a polynomial representation of $G$, \cite{PS} gives explicit inequalities describing the semialgebraic set $\R^n/G$.
	\end{enumerate}
\end{rems}
 
Let us end this part by considering the irreducible components of $X \gq G$:

\begin{rem} \label{remirredcompgeomquotient} Let $Y$ be an irreducible component of $X$. Then $\pi_X(Y)$ is an irreducible subset of $X /\!\!/ G$ (see, for instance, Proposition 1.1.7 (i) of \cite{TY}: since $\pi_X$ is a polynomial map, $\pi_X$ is in particular continuous with respect to Zariski topology) as well as its Zariski closure $\overline{\pi_X(Y)}^{\mathcal{Z}}$ (Proposition 1.1.5 (i) of \cite{TY}). Now, remark that for all $g \in G$, $g \cdot Y := \left\{ g \cdot y~|~y \in Y\right\}$ is an irreducible component of $X$ and, since $X$ has only a finite number of irreducible components, the union $G \cdot Y := \bigcup_{g \in G} g \cdot Y$ is a finite union of irreducible components of $X$ as well as a $G$-real algebraic subset of $X$. Finally, notice that $\pi_X(G \cdot Y) = \pi_X(Y)$ and then $\overline{\pi_X(Y)}^{\mathcal{Z}} = \left(G \cdot Y\right)/\!\!/ G$.

On the other hand, if $Y_1,\ldots,Y_k$ denote the irreducible components of $X$, we have $\pi_X(X) = \pi_X(Y_1) \cup \cdots \cup  \pi_X(Y_r) \subset \overline{\pi_X(Y_1)}^{\mathcal{Z}} \cup \cdots \cup \overline{\pi_X(Y_r)}^{\mathcal{Z}} \subset X /\!\!/ G$. Hence $X /\!\!/ G = \overline{\pi_X(Y_1)}^{\mathcal{Z}} \cup \cdots \cup \overline{\pi_X(Y_r)}^{\mathcal{Z}}$ (recall that $\overline{\pi_X(X)}^{\mathcal{Z}} = X /\!\!/ G$: cf. lemma \ref{lemzariskiclossemialgquotisgeomquot}): there then exist $i_1,\ldots,i_l \in \{1,\ldots,k\}$ such that the irreducible algebraic sets $\overline{\pi_X(Y_{i_1})}^{\mathcal{Z}} = \left(G \cdot Y_{i_1}\right)/\!\!/ G, \ldots, \overline{\pi_X(Y_{i_l})}^{\mathcal{Z}} = \left(G \cdot Y_{i_l}\right)/\!\!/ G$ are the irreducible components of~$X /\!\!/ G$.
\end{rem} 
 
\subsection{Functoriality of the real geometric quotient} 

We are going to establish the functoriality of the geometric quotient (and thus of the semialgebraic quotient) as well as its compability with inclusion of $G$-stable subsets.

\begin{lem} \label{lemfunctorrealalgquotient} Let $Y \subset \R^N$ be a $G$-real algebraic set and let $f : X \rightarrow Y$ be a morphism of $G$-real algebraic sets. The morphism $f$ induces, in a functorial way, (a restriction of) a polynomial map $f_{/\!\!/ G} : X /\!\!/ G \rightarrow Y /\!\!/ G$ such that the diagram 
$$\begin{array}{ccc}X & \overset{f}{\longrightarrow} & Y\\~~~\downarrow{\scriptstyle \pi_X}& &~~~\downarrow{\scriptstyle \pi_Y}\\X /\!\!/ G & \xrightarrow{f_{/\!\!/ G}} & Y /\!\!/ G\end{array}$$
commutes and $f_{/\!\!/ G}$ restricts to a map $f_{/G} : X/G \rightarrow Y/G$.
\end{lem} 

\begin{proof} Since $f$ is equivariant, the corresponding morphism of $\R$-algebras $f^* : \mathcal{P}(Y) \rightarrow \mathcal{P}(X)$ (see the beginning of the proof of theorem \ref{theorealnullstellensatzequivcategories}) restricts to a morphism $f^*_G : \mathcal{P}(Y)^G \rightarrow \mathcal{P}(X)^G$ which corresponds to a polynomial map $f_{/\!\!/ G} : X /\!\!/ G \rightarrow Y /\!\!/ G$ (by theorem \ref{theorealnullstellensatzequivcategories}). Moreover, the commutative diagram
$$\begin{array}{ccc}\mathcal{P}(X) & \overset{f^*}{\longleftarrow} & \mathcal{P}(Y)\\\hookuparrow& &\hookuparrow \\\mathcal{P}(X)^G & \overset{f^*_G}{\longleftarrow} & \mathcal{P}(Y)^G\end{array}$$
corresponds to the commutative diagram of real algebraic sets
$$\begin{array}{ccc}X & \overset{f}{\longrightarrow} & Y\\~~~\downarrow{\scriptstyle \pi_X}& &~~~\downarrow{\scriptstyle \pi_Y}\\X /\!\!/ G & \xrightarrow{f_{/\!\!/ G}} & Y /\!\!/ G\end{array}$$
Finally, let $z \in X/G$ and $x \in X$ such that $z = \pi_X(x)$, we have $f_{/\!\!/ G}(z) = f_{/\!\!/ G}(\pi_X(x)) = \pi_Y(f(x)) \in Y/G$, so that $f_{/\!\!/ G}$ restricts to a map $f_{/G} : X/G \rightarrow Y/G$.
\end{proof}

\begin{rems} \label{remafterfunctorialityrealalgquotient}
~
\begin{enumerate}
	\item Keep the notations of lemma \ref{lemfunctorrealalgquotient} and its proof. The morphism $f^*$ commutes with the Reynolds operator (cf. subsection \ref{subsectionlinredragr}): we have
	$$\mathcal{R}_X \circ f^* = f^*_G \circ \mathcal{R}_Y.$$  
Indeed, if $\alpha : G \times X \rightarrow X$ and $\alpha' : G \times Y \rightarrow Y$ denote the respective polynomial actions of $G$ on $X$ and $Y$, we have, for any $h \in \mathcal{P}(Y)$ and $x \in X$, by remark \ref{remexpresreynoldsopg}, 
$$\mathcal{R}_X \circ f^*(h)(x) = \mathcal{R}_G(h \circ f \circ \alpha(\cdot, x)) = \mathcal{R}_G(h \circ \alpha'(\cdot, f(x))) = \mathcal{R}_Y(h)(f(x)) = f^* \circ \mathcal{R}_Y(h)(x).$$
	\item The geometric quotient is also a functor with respect to the acting group: if $H$ is a linearly reductive real algebraic group and $\varphi : H \rightarrow G$ is a morphism of real algebraic groups, then $X$ is also a $H$-real algebraic set and $\varphi$ induces a polynomial map $\Phi : X /\!\!/ H \rightarrow X /\!\!/ G$ such that $\Phi \circ \pi_X = \pi'_X$, where $\pi'_X$ denotes the geometric quotient map of $X$ by the induced action of $H$ (notice that $\Phi$ then restricts to a map $X/G \rightarrow X/H$). Indeed, the identity morphism $\mathcal{P}(X) \rightarrow \mathcal{P}(X)$ restricts to a morphism $\mathcal{P}(X)^G \hookrightarrow \mathcal{P}(X)^H$, and to the commutative diagram of real finitely generated $\R$-algebras
$$\begin{array}{ccc}\mathcal{P}(X) & \xleftarrow{{\rm Id}_{\mathcal{P}(X)}} & \mathcal{P}(X)\\\hookuparrow& &\hookuparrow \\\mathcal{P}(X)^H & \hookleftarrow & \mathcal{P}(X)^G\end{array}$$
corresponds a commutative diagram of real algebraic sets
$$\begin{array}{ccc}X & \overset{{\rm Id}_X}{\longrightarrow} & X\\~~~\downarrow{\scriptstyle \pi'_X}& &~~~\downarrow{\scriptstyle \pi_X}\\X /\!\!/ H & \xrightarrow{\Phi} & X /\!\!/ G\end{array}$$	
\end{enumerate}
\end{rems}

\begin{lem} \label{lemrestrictionrealalgquotient} Let $Y$ be a $G$-real algebraic subset of $X$. The geometric quotient of $Y$ is (up to polynomial isomorphism) a restriction of the geometric quotient of $X$.
\end{lem}

\begin{proof} The inclusion morphism of $G$-real algebraic sets $i : Y \hookrightarrow X$ corresponds to the surjective restriction morphism $i^* : \mathcal{P}(X) \rightarrow \mathcal{P}(Y)$, which restricts to a morphism of $\R$-algebras $\mathcal{P}(X)^G \rightarrow \mathcal{P}(Y)^G$. The latter morphism is surjective as well: if $f \in \mathcal{P}(Y)^G$, there exists $h \in \mathcal{P}(X)$ such that $f = i^*(h)$ and we have
$$f = \mathcal{R}_X(f) = \mathcal{R}_X(i^*(h)) = i^*(\mathcal{R}_Y(h)) \in i^*\left(\mathcal{P}(X)^G\right)$$
by remark \ref{remafterfunctorialityrealalgquotient} 1. As a consequence, if $h_1,\ldots,h_s$ are invariant polynomial functions on $X$ generating the $\R$-algebra $\mathcal{P}(X)^G$, their images $i^*(h_1),\ldots,i^*(h_s)$ generate $\mathcal{P}(Y)^G$ and then, if $y \in Y$, 
$$\pi_Y(y) = \big(i^*(h_1)(y),\ldots,i^*(h_s)(y)\big) = \left(h_1(y),\ldots,h_s(y)\right) = \pi_X(y).$$ 
\end{proof}

In particular, if $S \subset \R^n$ is a $G$-semialgebraic set which is a $G$-stable subset of $X$, the semialgebraic quotient $\varpi_S : S \rightarrow S/G$ can be considered (up to polynomial isomorphism) as the restriction $S \rightarrow \pi_X(S)$ of the geometric quotient $\pi_X : X \rightarrow X \gq G$ of $X$.

\begin{exs} 
~
	\begin{enumerate}
		\item Consider the action of $S^1$ on $\R^3$ of example \ref{exrealalgquotientplaneorthoaction} 2, as well as the $S^1$-real algebraic subsets $C := \{x^2+y^2 = z^2\}$ and $Y := \{x^2+y^2+z^4 = z^2\}$ of $\R^3$ (examples~\ref{exgrealalgset}~2). The geometric quotient map of $\R^3$ is 
		$$\pi_{\R^3} : \R^3 \rightarrow \R^2~;~(x,y,z) \mapsto (x^2+y^2,z),$$
so that the semialgebraic quotient of $C$ is $C/S^1 = \pi_{\R^3}(C) = \{u = v^2\} = C /\!\!/ S^1$ and the semialgebraic quotient of $Y$ is $Y/S^1 = \{u + v^4 = v^2,~u \geq 0\}$ (the geometric quotient of~$Y$ is $Y /\!\!/ S^1 = \overline{Y/S^1}^{\mathcal{Z}} =  \{u + v^4 = v^2\}$).
		\item Consider the usual action of $\OO_2(\R)$ on $\R^2$ (example \ref{exrealalgquotientplaneorthoaction} 1) with geometric quotient 
		$$\pi_{\R^2} :\R^2 \rightarrow \R~;~(x,y) \mapsto x^2+y^2,$$
as well as the open disc $D := \{ (x,y) \in \R^2~|~x^2+y^2 < 1\}$ which is a $\OO_2(\R)$-semialgebraic subset of $\R^2$. We have $D/\OO_2(\R) = \pi_{\R^2}(D) = [0;1[$.
	\end{enumerate}
\end{exs}

\begin{rem} Previous lemmas also allow to reduce any geometric quotient to the restriction of the geometric quotient of a polynomial representation of $G$ (definition \ref{defrealpolyrepres}), using proposition~\ref{proplinearizgrealalgset}: there exists an isomorphism of $G$-real algebraic sets $f$ from $X$ to a $G$-real algebraic subset $Y$ of a polynomial representation $\R^N$ of $G$, which induces an isomorphism between the geometric quotient $\pi_X$ and $\pi_Y$ (by lemma \ref{lemfunctorrealalgquotient}), and $\pi_Y$ is the restriction of the geometric quotient $\pi_{\R^N}$ (by lemma \ref{lemrestrictionrealalgquotient}). 
\end{rem}

\begin{rem} \label{remgeomquotientproduct} Let $H$ be a linearly reductive real algebraic group and let $Y$ be a $H$-real algebraic set. Then, by remark \ref{reminvariantalgebraproduct}, the geometric quotient of $X \times Y$ by $G \times H$ is, up to polynomial isomorphism, the polynomial map 
$$\pi_X \times \pi_Y : \begin{array}{ccc}X \times Y & \rightarrow & X \gq G \times Y \gq H \\ (x,y) & \mapsto & (\pi_X(x),\pi_Y(y))\end{array},$$ 
where $\pi_X$, resp. $\pi_Y$, denotes the geometric quotient of $X$ by $G$, resp. $Y$ by $H$.
\end{rem}

\subsection{Geometric quotient of a real algebraic group by a linearly reductive normal subgroup} \label{subsectiongeomquotientofrealalggrbylrnormal}

We dedicate one part of this section to the application of the previous construction to the case of the geometric quotient of a real algebraic group by a linearly reductive normal subgroup: this quotient turns out to be a real algebraic group as well, as we will prove below. 

Let $K$ be a real algebraic group with polynomial binary operation $\mu : K \times K \rightarrow K$ and let $H$ be a {\it normal} real algebraic subgroup of $K$ that we furthermore suppose to be linearly reductive. The real algebraic group $K$ polynomially acts on itself by left multiplication (cf. example~\ref{exgrealalgset}~4) and we restrict this action into the polynomial action $H \times K \rightarrow (h,k) \mapsto \mu(h,k)$, making~$K$ into a $H$-real algebraic set. Since $H$ is linearly reductive, we can consider the geometric quotient of~$K$ by $H$.

\begin{prop} \label{propdefrealalgquotientgroup} The geometric quotient $K /\!\!/ H$ is a real algebraic group and the geometric quotient map $\pi_K : K \rightarrow K /\!\!/ H$ is a morphism of real algebraic groups.
\end{prop}  

\begin{proof} Consider the morphism of $\R$-algebras $\mu^* : \mathcal{P}(K) \rightarrow \mathcal{P}(K \times K) = \mathcal{P}(K) \otimes_{\R} \mathcal{P}(K)$ associated to $\mu$. We are going to show that $\mu^*$ restricts to a map $\mathcal{P}(K)^H \rightarrow \mathcal{P}(K)^H \otimes_{\R} \mathcal{P}(K)^H$ and that the corresponding polynomial map $K /\!\!/ H \times K /\!\!/ H \rightarrow K /\!\!/ H$ provides a real algebraic group structure on the geometric quotient $K /\!\!/ H$.

So let $f \in \mathcal{P}(K)^H$ and write $\mu^*(f) = f \circ \mu = \sum_{r=1}^s f_r \otimes f'_r$ with $f_1,\ldots,f_s, f_1',\ldots,f_s' \in \mathcal{P}(K)$. Fix $k' \in K$. Since the polynomial function $f \circ \mu(\cdot, k')$ on $K$ is also $H$-invariant, we have, denoting by $\mathcal{R}'_K$ the Reynolds operator of the $H$-real algebraic set $K$ ($H$ is linearly reductive), 
\begin{eqnarray*}
f \circ \mu(\cdot, k') & =  & \mathcal{R}'_K\left(f \circ \mu(\cdot, k')\right)\\
			      & = & \mathcal{R}'_K\left(\sum_{r=1}^s f'_r(k') f_r \right)\\
			      & = & \sum_{r=1}^s f'_r(k') \mathcal{R}'_K(f_r) \\
			      & = & \left(\sum_{r=1}^s \mathcal{R}'_K(f_r) \otimes f'_r\right)(\cdot,k'),  
\end{eqnarray*}
so that we can suppose that $f_1,\ldots,f_s \in \mathcal{P}(K)^H$. Furthermore, if $k \in K$, because $H$ is a normal subgroup of $K$, the polynomial function $f \circ \mu(k, \cdot)$ is $H$-invariant as well and then
\begin{eqnarray*}
f \circ \mu(k, \cdot) & =  & \mathcal{R}'_K\left(f \circ \mu(k, \cdot)\right)\\
			      & = & \mathcal{R}'_K\left(\sum_{r=1}^s f_r(k) f'_r \right)\\
			      & = & \sum_{r=1}^s f_r(k) \mathcal{R}'_K(f'_r) \\
			      & = & \left(\sum_{r=1}^s f_r \otimes \mathcal{R}'_K(f'_r)\right)(k, \cdot).
\end{eqnarray*}
As a consequence, we can also suppose that $f'_1,\ldots,f'_s \in \mathcal{P}(K)^H$ and then $\mu^*(f) \in \mathcal{P}(K)^H \otimes_{\R} \mathcal{P}(K)^H$. 

Denote by $\overline{\mu}$ the polynomial map $K /\!\!/ H \times K /\!\!/ H \rightarrow K /\!\!/ H$ corresponding to the restriction $\mu^*_H : \mathcal{P}(K)^H \rightarrow \mathcal{P}(K)^H \otimes_{\R} \mathcal{P}(K)^H$ of $\mu^*$. First, remark that $\overline{\mu}$ is associative: the associativity of~$\mu$ translates into a commutative diagram of morphisms of $\R$-algebras that restricts to $\mathcal{P}(K)^H$.

Now, denote by $e$ the neutral element of $K$ and consider the polynomial map $\rho : K \rightarrow K \times K~;~k \mapsto (e,k)$. Notice that, if $f,f' \in \mathcal{P}(K)^H$, $\rho^*(f \otimes f') = f(e) f' \in \mathcal{P}(K)^H$, so that $\rho^*$ restricts to a morphism $\rho^*_H : \mathcal{P}(K)^H \otimes_{\R} \mathcal{P}(K)^H \rightarrow \mathcal{P}(K)^H$ which corresponds to a polynomial map $\overline{\rho} : K /\!\!/ H \rightarrow K /\!\!/ H \times K /\!\!/ H$. On the one hand, we have $\rho^*_H \circ \mu^*_H = {\rm Id}_{\mathcal{P}(K)^H}$ i.e. $\overline{\mu} \circ \overline{\rho} = {\rm Id}_{K /\!\!/ H}$. On the other hand, if $f_1,\ldots,f_s$ are generators of the $\R$-algebra $\mathcal{P}(K)^H$, the tensors $f_1\otimes \overline{1},\ldots,f_s\otimes \overline{1}, \overline{1} \otimes f_1,\ldots,\overline{1} \otimes f_r$ generate $\mathcal{P}(K)^H \otimes_{\R} \mathcal{P}(K)^H$ and we have, for all $r \in \{1,\ldots,s\}$, $\rho^*_H(f_r \otimes \overline{1}) = f_r(e)$ and $\rho^*_H( \overline{1} \otimes f_r) = f_r$. Consequently, for any $g \in K /\!\!/ H$, we have $\overline{\rho}(g) = (\pi_K(e), g)$ (see the end of the proof of theorem \ref{theorealnullstellensatzequivcategories}) and then $\overline{\mu}(\pi_K(e),g) = \overline{\mu} \circ \overline{\rho}(g) = g$. We can prove similarly that, for all $g \in K /\!\!/ H$, $\overline{\mu}(g, \pi_K(e)) = g$. 

Finally, denote by $\omega$ the inverse map $K \rightarrow K~;~k \mapsto k^{-1}$, notice that $\omega^* : \mathcal{P}(K) \rightarrow \mathcal{P}(K)$ restricts to $\omega^*_H : \mathcal{P}(K)^H \rightarrow \mathcal{P}(K)^H$ (because $H$ is a normal subgroup of $K$) and denote by~$\overline{\omega}$ the corresponding polynomial map $K /\!\!/ H \rightarrow K /\!\!/ H$. Denote also by $\epsilon$ the constant map $K \rightarrow K~;~k \mapsto e$: the morphism of $\R$-algebras $\epsilon^*$ restricts to the morphism $\epsilon^*_H : \mathcal{P}(K)^H \rightarrow \mathcal{P}(K)^H~;~f \mapsto f(e)$ which itself corresponds to the map $\overline{\epsilon} : K /\!\!/ H \rightarrow K /\!\!/ H~;~g \mapsto \pi_K(e)$. Since $\mu(\cdot, \omega(\cdot)) = \epsilon$, we have, by the correspondence of theorem \ref{theorealnullstellensatzequivcategories} and restriction to $\mathcal{P}(K)^H$, the equality $\overline{\mu}(\cdot, \overline{\omega}(\cdot)) = \overline{\epsilon}$ i.e. for all $g \in K /\!\!/ H$, $\overline{\mu}(g, \overline{\omega}(g)) = \pi_K(e)$. We can prove similarly that, for all $g \in K /\!\!/ H$, $\overline{\mu}(\overline{\omega}(g), g) = \pi_K(e)$. 

As a conclusion, the couple $(K /\!\!/ H, \overline{\mu})$ is a (real algebraic) group. Moreover, the geometric quotient map $\pi_K$ is a group homomorphism by commutativity of the diagram
$$\begin{array}{ccc}\mathcal{P}(K) \otimes_{\R} \mathcal{P}(K) & \xleftarrow{\mu^*} & \mathcal{P}(K)\\\hookuparrow& &\hookuparrow \\ \mathcal{P}(K)^H \otimes_{\R} \mathcal{P}(K)^H &  \xleftarrow{\mu^*_H} & \mathcal{P}(K)^H \end{array}$$
\end{proof}

\begin{cor} The semialgebraic quotient $K/H$ is a semialgebraic group (definition \ref{defsemialggroup}) and the semialgebraic quotient map $\varpi_K : K \rightarrow K/H$ is a group homomorphism.
\end{cor}

\begin{proof} The semialgebraic quotient $K/H$ is a semialgebraic group as the image of the morphism of real algebraic groups $\pi_K$ (see examples \ref{exssemialggps}) and the semialgebraic quotient map $\varpi_K : K \rightarrow K/H$ is a group homomorphism as a restriction of the latter.
\end{proof}

\begin{ex} Consider the commutative real algebraic group $K := \{(x,y) \in \R^2~|~xy=1\}$ (example \ref{exspolyrealalggroups} 3) and its (normal) real algebraic subgroup $H := \{(1,1),(-1,-1)\}$. Since $H$ is finite, $H$ is a linearly reductive real algebraic group (cf. example \ref{exlinearlyredrealalggr} 1). The polynomial action of $H$ on $K$ by left multiplication is a restriction of the polynomial action $H \times \R^2 \rightarrow \R^2~;~\big((a,b),(x,y)\big) \mapsto (ax,by)$ of $H$ on $\R^2$. The $\R$-algebra $\mathcal{P}(\R^2)^H$ being generated by the polynomials $x^2$, $y^2$ and $xy$, the geometric quotient $\pi_{\R^2}$ is then the polynomial map 
$$(x,y) \in \R^2 \mapsto (x^2,y^2,xy) \in \R^3.$$
Therefore, we have (using lemma \ref{lemrestrictionrealalgquotient}) $K/H = \pi_{\R^2}(K) = \{(u,v,1) \in \R^3~|~uv=1,~u \geq 0,~v \geq 0\}$ and $K /\!\!/ H = \overline{K/H}^{\mathcal{Z}} = \{(u,v,1) \in \R^3~|~uv=1\}$ (the latter real algebraic group is isomorphic to $K$).
\end{ex}

\subsection{Complexification of the real geometric quotient} \label{subsectioncomplexrealgeomquotient}

Complexify the $G$-real algebraic set $X$ into the $G_{\C}$-complex algebraic set $X_{\C}$ (cf. lemma \ref{lemcomplexificationpolyaction}) and consider the inclusion morphism of $\C$-algebras $\mathcal{P}(X_{\C})^{G_{\C}} \hookrightarrow \mathcal{P}(X_{\C})$. Both algebras $\mathcal{P}(X_{\C})$ and $\mathcal{P}(X_{\C})^{G_{\C}}$ are finitely generated (since $G$ is a linearly reductive real algebraic group, $G_{\C}$ is a linearly reductive complex algebraic group and the $\C$-algebra $\mathcal{P}(X_{\C})^{G_{\C}}$ is finitely generated: cf. subsection \ref{subsectcomplexificationlinearlyreductive}) and reduced: recall that a commutative ring $A$ is said {\it reduced} if for any $a \in A$, the existence of a positive integer~$k$ such that $a^k = 0_A$ implies that $a = 0_A$ (in particular, any subring of a reduced ring is reduced as well).

By Hilbert's Nullstellensatz, the operation which associates to any complex algebraic set its algebra of complex polynomial functions gives a contrafunctorial equivalence between the category of complex algebraic sets and polynomial maps and the category of reduced finitely generated $\C$-algebras and morphisms of $\C$-algebras (see theorem \ref{theoequivcatcomphilbertsnullstel} of the appendix: the proof is similar to the proof of theorem \ref{theorealnullstellensatzequivcategories}). 

Therefore, the inclusion morphism $\mathcal{P}(X_{\C})^{G_{\C}} \hookrightarrow \mathcal{P}(X_{\C})$ corresponds to a polynomial map $\pi : X_{\C} \rightarrow Z$ given, in the same way as in the real case, by generators of $\mathcal{P}(X_{\C})^{G_{\C}}$. But, as we saw at the end of subsection \ref{subsectcomplexificationlinearlyreductive}, those generators can actually be chosen to be the complexifications of generators of $\mathcal{P}(X)^G$: if we keep the above considered generators $f_1,\ldots,f_r$ of $\mathcal{P}(X)^G$, then the extensions $(f_1)_{\C},\ldots, (f_r)_{\C}$ generate the $\C$-algebra $\mathcal{P}(X_{\C})^{G_{\C}}$. 

Furthermore:

\begin{lem}
The complex algebraic set $Z$ is (up to polynomial isomorphism) the complexification of $X /\!\!/ G$ in $\C^r$.
\end{lem}

\begin{proof}
Taking the elements $(f_1)_{\C},\ldots, (f_r)_{\C}$ as generators of the $\C$-algebra $\mathcal{P}(X_{\C})^{G_{\C}}$, the set~$Z$ is the complex algebraic set $\mathsf{V}(J)$ of $\C^r$ where $J$ is the ideal
$$\{P \in \C[z_1,\ldots,z_r]~|~\forall x \in X_{\C},~P\big((f_1)_{\C}(x),\ldots, (f_r)_{\C}(x)\big) = 0\}$$
of $\C[z_1,\ldots,z_r]$. But, for all $p,q \in \R[z_1,\ldots,z_r]$, $p+iq \in J$ if and only if $p,q$ belong to the ideal 
$$\{P \in \R[z_1,\ldots,z_r]~|~\forall x \in X,~P(f_1(x),\ldots,f_r(x)) = 0\} = I(X /\!\!/ G)$$
of $\R[z_1,\ldots,z_r]$, which then happens if and only if $p+i q$ vanishes on $X /\!\!/ G$. As a consequence,
$$Z = \mathsf{V}(J) = \mathsf{V}(\mathsf{I}(X \gq G)) = \left(X \gq G\right)_{\C}.$$
\end{proof} 

The polynomial map $\pi : X_{\C} \rightarrow Z~;~x \mapsto \big((f_1)_{\C}(x),\ldots, (f_r)_{\C}(x)\big)$ is therefore the complexification of the real geometric quotient $\pi_X : X \rightarrow X /\!\!/ G$. 

\begin{de}
We denote $X_{\C} /\!\!/ G_{\C} := Z$ and $\pi_{X_{\C}} := \pi$, and call the map $\pi_{X_{\C}} : X_{\C} \rightarrow X_{\C} /\!\!/ G_{\C}$ the (complex) \emph{geometric quotient of $X_{\C}$ by $G_{\C}$}.
\end{de}

\begin{ex} Consider the usual action of $\OO_n(\R)$ on $\R^n$ which complexifies into the usual action of $\OO_n(\C)$ on $\C^n$ (see example \ref{exscomplexificationpolyaction} 1). Since the $\R$-algebra $\mathcal{P}(\R^n)^{\OO_n(\R)}$ is generated by the polynomial $\sum_{k=1}^n x_k^2$ of $\R[x_1,\ldots,x_n]$ (see example \ref{exgeninvariantalgorthogroup}), the $\C$-algebra $\mathcal{P}(\C^n)^{\OO_n(\C)}$ is generated by the polynomial~$\sum_{k=1}^n x_k^2$ of $\C[x_1,\ldots,x_n]$ and the complex geometric quotient of~$\C^n$ by $\OO_n(\C)$ is the map
$$\C^n \rightarrow \C~;~(x_1,\ldots,x_n) \mapsto \sum_{k=1}^n x_k^2,$$
whose image is $\C$ and which restricts to the real geometric quotient map $\R^n \rightarrow \R~;~(x_1,\ldots,x_n) \mapsto \sum_{k=1}^n x_k^2$ of $\R^n$ by $\OO_n(\R)$ (example \ref{exrealalgquotientplaneorthoaction} 1).
\end{ex}

Contrary to the real geometric quotient $\pi_X$ of $X$ by $G$, the complex geometric quotient $\pi_{X_{\C}}$ of $X_{\C}$ by $G_{\C}$ is always surjective: we recall the proof of this result (following the argumentation of \cite{DK} Lemma 2.3.2) for the sake of completeness and in order to highlight the fact that it uses Hilbert's Nullstellensatz (which is not true for real algebraic sets).

\begin{lem} The map $\pi_{X_{\C}} : X_{\C} \rightarrow X_{\C} /\!\!/ G_{\C}$ is surjective.
\end{lem}

\begin{proof} For the sake of readability, let us use the first considered notation $\pi : X_{\C} \rightarrow Z$ for the complex geometric quotient of $X_{\C}$ by $G_{\C}$. Let $Y$ be a non-empty complex algebraic subset of~$Z$: there exist polynomial functions $h_1,\ldots,h_s$ of $\mathcal{P}(Z)$ such that $Y$ is the set of points of $Z$ at which vanish (the polynomial functions of the ideal of $\mathcal{P}(Z)$ generated by) the functions $h_1,\ldots,h_s$. For any $x \in X_{\C}$, $x$ is in $\pi^{-1}(Y)$ if and only if $\widetilde{h}_1(x) = \cdots = \widetilde{h}_s(x) = 0$, where, if $l \in \{1,\ldots,s\}$,~$\widetilde{h}_l$ is the $G_{\C}$-invariant polynomial function $x \in X_{\C} \mapsto h_l \big((f_1)_{\C}(x),\ldots, (f_r)_{\C}(x)\big) \in \C$: the set~$\pi^{-1}(Z)$ is therefore the set of points of $X_{\C}$ at which vanish the elements of the ideal $\mathcal{I}$ of $\mathcal{P}(X_{\C})$ generated by the polynomial functions $\widetilde{h}_1,\ldots,\widetilde{h}_s$.

Suppose that $\pi^{-1}(Y)$ is an empty set. By Hilbert's Nullstellensatz (cf. for instance Theorem~1.1 of \cite{Bump}), the ideal $\mathcal{I}$ of $\mathcal{P}(X_{\C})$ is then the entire ring $\mathcal{P}(X_{\C})$: there exist polynomial functions $\rho_1,\ldots,\rho_s$ on $X_{\C}$ such that $\overline{1} = \sum_{l=1}^s \rho_l \widetilde{h}_l$ in $\mathcal{P}(X_{\C})$. We then apply the Reynolds operator of $X_{\C}$ to obtain the equality $\overline{1} = \sum_{l=1}^s \mathcal{R}_{X_{\C}}(\rho_l) \widetilde{h}_l$ in $\mathcal{P}(X_{\C})^{G_{\C}}$ (the functions $\widetilde{h}_1,\ldots,\widetilde{h}_s$ are $G_{\C}$-invariant polynomial functions on $X_{\C}$). As a consequence, by applying the isomorphism of $\C$-algebras $\mathcal{P}(X_{\C})^{G_{\C}} \rightarrow \mathcal{P}(Z)$ (which associates to any function $(f_k)_{\C}$, $k \in \{1,\ldots, r\}$, the function $\overline{z_k}$ on $Z$: see the proof of theorem \ref{theorealnullstellensatzequivcategories}), the function $\overline{1}$ of $\mathcal{P}(Z)$ belongs to the ideal of $\mathcal{P}(Z)$ generated by the polynomial functions $h_1,\ldots,h_s$. Consequently, $Y$ is empty, which is a contradiction. 

In particular, if $y \in Z$, the fiber $\pi^{-1}(\{y\})$ is non-empty and $\pi$ is then surjective.
\end{proof}

\section{Polynomial actions of compact real algebraic groups on real algebraic and semialgebraic sets} \label{sectpolyactionscomprealalggroups}

In the following of this text, we will focus on compact real algebraic groups and their polynomial actions on real algebraic sets. As we have seen in example \ref{exlinearlyredrealalggr} 2, a compact real algebraic group is linearly reductive, so that we can consider the associated geometric and semialgebraic quotients. But the polynomial action of a compact real algebraic group on a real algebraic set has a further property of separation of orbits. This has many consequences, in particular the fact that two points have the same image by the geometric quotient if and only of they have the same orbit. We will see that the semialgebraic quotient is actually (up to homeomorphism) the topological quotient and we will then establish several topological and geometric properties of the semialgebraic and geometric quotients associated to a compact real algebraic group.

The first part of this section will be dedicated to compact real algebraic groups themselves. In fact, our starting point will be to consider the orbits of the polynomial action of any compact subgroup of a real algebraic group on a real algebraic set: due to their separation property, those orbits are always real algebraic sets. This will imply some rigidity and stability of the class of compact real algebraic groups: for instance, any closed subgroup of a compact real algebraic group is a compact real algebraic group and the image of a compact real algebraic group by a morphism of real algebraic groups is a compact real algebraic group. 

In a subsequent subsection, we will establish further properties in this direction, considering the complexification of compact (Lie) subgroups of orthogonal groups. 

\subsection{Compact subgroups of real algebraic groups}

Let $G \subset \R^m$, $m \in \Nstar$, be a real algebraic group with binary operation $\mu : G \times G \rightarrow G$ and let $H$ be a \emph{compact} subgroup of $G$ i.e. a subgroup of $G$ which is a closed bounded subset of~$\R^m$. Let $X \subset \R^n$, $n \in \Nstar$, be a $G$-real algebraic set with polynomial action $\alpha : G \times X \rightarrow X$ of $G$ on $X$. We are going to state the crucial fact that all the orbits of the induced action of~$H$ on $X$ are real algebraic sets: this will be our starting point to establish some properties of compact real algebraic groups and their polynomial actions on real algebraic sets.

First consider the following important property of separation of orbits by an invariant polynomial function: if $x \in X$, we denote by $H \cdot x$ the orbit of $x$ under the action of $H$ on $X$ induced by the action of $G$.

\begin{lem} \label{lemseporbitscompactsubg} Let $x \in X$ and $y \in X \setminus H \cdot x$. There exists a polynomial function $p$ on $X$ which is $H$-invariant (i.e. for all $h \in H$ and $z \in X$, $p(h \cdot z) = p(z)$) and such that $p(x) \neq p(y)$.
\end{lem}

\begin{proof} Consider the continuous function $F : (H \cdot x) \cup (H \cdot y) \rightarrow \R$ on the compact set $(H \cdot x) \cup (H \cdot y)$ which takes value $1$ on $H \cdot x$ and value $-1$ on $H \cdot y$ (we have $(H \cdot x) \cap (H \cdot y) = \emptyset$ and $H \cdot x$ and $H \cdot y$ are compact). By Stone-Weiertrass theorem, there exists a polynomial of~$\R[x_1,\ldots,x_n]$ restricting into a polynomial function $f$ on $X$ such that for all $z \in (H \cdot x) \cup (H \cdot y)$, $|F(z) - f(z)| < \frac{1}{2}$.

Now, because $H$ is a compact topological group (as a compact subgroup of the topological group $G$), we can consider the unique normalized Haar measure $\mu$ and associated integral on~$H$ (see for instance section \textbf{8}.1.2 of \cite{Pro}) and we then define the function 
$$p : X \rightarrow \R~;~z \mapsto \int_H f \circ \alpha(h,z) {\rm d}\mu(h)$$
on $X$ (if $z \in X$, the function $f \circ \alpha(\cdot,z)$ is continuous on $H$). Using the existence of polynomial functions $\phi_1,\ldots,\phi_l \in \mathcal{P}(G)$ and $h_1,\ldots,h_l \in \mathcal{P}(X)$ such that $f \circ \alpha = \sum_{k=1}^l \phi_k \otimes h_k$ (see for instance the proof of proposition \ref{proplinearizgrealalgset}), we have, if $z \in X$,
$$p(z) = \int_H \sum_{k=1}^l \phi_k(h) h_k(z) {\rm d}\mu(h) = \sum_{k=1}^l  \left(\int_H \phi_k(h) {\rm d}\mu(h)\right) h_k(z),$$
so that $p$ is, in particular, a polynomial function on $X$. Moreover, $p$ is $H$-invariant since, for all $z \in X$ and $k \in H$,
$$p(k \cdot z) = \int_H f \circ \alpha(h,k \cdot z) {\rm d}\mu(h) = \int_H f \circ \alpha(hk ,z) {\rm d}\mu(h) = \int_H f \circ \alpha(h ,z) {\rm d}\mu(h) = p(z)$$
(by invariance of the Haar measure $\mu$). 

Finally, $p(x) = \int_H f \circ \alpha(h,x) {\rm d}\mu(h) > \int_H \left( F(h \cdot x) - \frac{1}{2}\right) {\rm d}\mu(h) = \frac{1}{2}$ and $p(y) = \int_H f \circ \alpha(h,y) {\rm d}\mu(h) < \int_H \left( F(h \cdot y) + \frac{1}{2}\right) {\rm d}\mu(h) = - \frac{1}{2}$, so that $p(x) \neq p(y)$.
\end{proof}

\begin{prop} \label{proporbitcompactsubgrouprealalggroup} For all $x \in X$, the orbit $H \cdot x := \{ h \cdot x~|~h \in H\}$ of the action of $H$ on $X$ is a compact real algebraic set.
\end{prop}

\begin{proof} Let $x \in X$ and first remark that the set $H \cdot x$ is compact, as the image of the compact set $H$ by the continuous map $G \rightarrow X~;~g \mapsto g \cdot x = \alpha(g,x)$. We are now going to show that~$H \cdot x$ is a real algebraic subset of $\R^n$, by showing that $H \cdot x = V(I(H \cdot x))$.

Let $y \in X \setminus H \cdot x$ (we have $H \cdot x \subset X$). By previous lemma \ref{lemseporbitscompactsubg}, there exists an $H$-invariant polynomial function $p$ on $X$ such that $p(x) \neq p(y)$. If we denote $q := p - \overline{p(x)} \in \mathcal{P}(X)$, we have, for all $h \in H$, $q(h \cdot x) = 0$ and $q(y) \neq 0$. Therefore, there exists a polynomial $Q \in \R[x_1,\ldots,x_n]$ such that $Q \in I(H \cdot x)$ and $Q(y) \neq 0$, that is $y \notin V(I(H \cdot x))$. Consequently, $V(I(H \cdot x)) \subset H \cdot x$ and then $H \cdot x = V(I(H \cdot x))$.
\end{proof}

An immediate consequence of the former proposition is that any compact subgroup of a real algebraic group is itself a real algebraic group:

\begin{cor} \label{corcompactsubgroupofrealalggrisrealalg} $H$ is a real algebraic set and is then a real algebraic subgroup of $G$.
\end{cor}

\begin{proof} Consider the action of $G$ on itself by left multiplication: $G$ is a $G$-real algebraic set. Therefore, by previous proposition \ref{proporbitcompactsubgrouprealalggroup}, $H = H \cdot e$ (where $e$ is the neutral element of $G$) is a real algebraic set.
\end{proof}

From this result, we deduce several classes of examples of compact real algebraic groups:

\begin{exs} \label{excompactsubgroupofrealalggrisrealalg}
~
	\begin{enumerate}
		\item Any closed (with respect to the ambient Euclidean topology) subgroup of a compact real algebraic group is a real algebraic group. 
		\item Any compact semialgebraic group is a real algebraic group (since it is a subgroup of a real algebraic group: see lemma \ref{lemzariskiclrealalggroup}).
		\item Any compact subgroup of ${\rm GL}_n(\R) = \left\{A \in \R^{n^2}~|~{\rm det}(A) \neq 0\right\}$ is a real algebraic group. Indeed, let $K$ be a compact subgroup of ${\rm GL}_n(\R)$ then, considering the (polynomial) determinant map ${\rm det} : \R^{n^2} \rightarrow \R$, ${\rm det}(K)$ is a compact subset of $\R$. If $A$ was a matrix of the group $K$ such that $|{\rm det}(A)| > 1$, then ${\rm det}(K)$ would be an unbounded subset of~$\R$, and if $A$ was a matrix of $K$ such that $|{\rm det}(A)| < 1$, then the closed set ${\rm det}(K)$ would contain $0$, which is impossible since $K \subset {\rm GL}_n(\R)$. As a consequence, ${\rm det}(K) \subset \{-1;1\}$, that is $K$ is a compact subgroup of the real algebraic group $\Omega_n(\R)$ (see example \ref{exspolyrealalggroups} 1).
	\end{enumerate}
\end{exs}

This also allows to show the stability of the class of compact real algebraic groups under morphisms of real algebraic groups:

\begin{cor} \label{corimagecompactrealalggroupisrealalg} The image of a compact real algebraic group by a morphism of real algebraic groups is a compact real algebraic group. 

In particular, the image of the geometric quotient of a compact real algebraic group by a normal closed subgroup (see subsection \ref{subsectiongeomquotientofrealalggrbylrnormal} and recall that a compact real algebraic group is linearly reductive by example \ref{exlinearlyredrealalggr} 2) is a compact real algebraic group (the associated semialgebraic quotient coincides with the geometric quotient).
\end{cor}

\begin{proof} The image of a compact real algebraic group by a morphism of real algebraic groups, which is in particular continuous, is a compact subgroup of a real algebraic group, hence a compact real algebraic group by corollary \ref{corcompactsubgroupofrealalggrisrealalg}.
\end{proof}

Let us also prove that the irreducible components of a compact real algebraic groups coincide with its (semialgebraic) connected components :

\begin{prop} Suppose that the real algebraic group $G$ is compact. Then the irreducible components of $G$ are exactly the connected components (with respect to Euclidean topology) of~$G$. In particular, $G$ is irreducible if and only if $G$ is connected.
\end{prop}

\begin{proof} Because $G$ is compact, the connected component $G_c$ (see lemma \ref{lempropcosetsconncompsemialggroup}) is a compact semialgebraic group and is therefore a real algebraic group. As a consequence, by proposition \ref{propconnectedcompirredcompneutralelt}, we have $G_c = G_0$, hence the result since the connected components of $G$ are the cosets of $G$ modulo $G_c$ and the irreducible components of $G$ are the cosets of $G$ modulo $G_0$ (see proposition \ref{propcosetsirredcomprealalggroup}). 
\end{proof}

\subsection{Orthogonalizability of the polynomial action of a compact real algebraic group} \label{subsecrealalgsetscomprealalggroupactionquotients}

Now suppose that $G \subset \mathbb{R}^m$ is a {\it compact} real algebraic group and, as above, let $X \subset \mathbb{R}^n$ be a~$G$-real algebraic set associated to a polynomial action map $\alpha : G \times X \rightarrow X$. Before studying the geometric quotient of $X$ by $G$ under this compactness hypothesis on $G$, let us establish that, up to an isomorphism of $G$-real algebraic sets, the action of $G$ on $X$ can be assumed to be given by orthogonal matrices with polynomial entries in the coordinates of $\R^m$ (recall that we have already asserted in proposition \ref{proplinearizgrealalgset} that the action of $G$ on $X$ can be polynomially linearized).

If $V$ is a real representation of a group $H$ associated to a group homomorphism $\rho : H \rightarrow \mathcal{GL}(V)$ and if $V$ is equipped with a Euclidean inner product, we say that the considered representation of $H$ is \emph{orthogonal} if for all $h \in H$, the linear isomorphism $\rho(h)$ is orthogonal with respect to the considered inner product. If $H$ is a real algebraic group, we say that an orthogonal representation of $H$ is a \emph{polynomial orthogonal representation of $H$} if it is a polynomial representation of $H$ (definition \ref{defrealpolyrepres}).

\begin{prop} \label{proporthogonalizationactioncompactrealalggroup} The $G$-real algebraic set $X$ is isomorphic to a $G$-real algebraic subset of a polynomial orthogonal representation $\mathbb{R}^N$ of $G$, where $\mathbb{R}^N$ is equipped with its canonical Euclidean inner product.
\end{prop}

\begin{proof} Thanks to proposition \ref{proplinearizgrealalgset}, up to an isomorphism of $G$-real algebraic sets, we can suppose that the action of $G$ on $X$ is the restriction of an action $\varrho : G \times \R^n \rightarrow \R^n$ of $G$ on $\mathbb{R}^n$ such that, for any $g \in G$, the map $\rho(g) : v \in \R^n \mapsto \varrho(g,v) \in \R^n$ is a linear automorphism given by matrices whose entries are polynomial in $g \in \R^m$. 

Now, denote by $\langle \cdot, \cdot \rangle$ the canonical Euclidean inner product of $\R^n$ and let $v,w$ be vectors of~$\R^n$: since the function $f_{v,w} : G \rightarrow \R~;~g \mapsto \langle \rho(g)(v), \rho(g)(w)\rangle$ is polynomial hence continuous, we can define the quantity
$$\langle v,w \rangle_G := \int_G f_{v,w}(g) {\rm d}\mu(g) = \int_G \langle \rho(g)(v), \rho(g)(w) \rangle {\rm d}\mu(g),$$
where $\mu$ is the unique normalized Haar measure on the compact group $G$. The map $\langle \cdot, \cdot \rangle_G : \R^n \times \R^n \rightarrow \R$ is then a Euclidean inner product and, if $h \in G$, we have 
$$\hspace{-0.3cm}\langle \rho(h)(v),\rho(h)(w) \rangle_G = \int_G \langle \rho(g h)(v), \rho(g h)(w) \rangle {\rm d}\mu(g) = \int_G f_{v,w}(g h) {\rm d}\mu(g) = \int_G f_{v,w}(g) {\rm d}\mu(g) = \langle v,w \rangle_G$$
by invariance of Haar measures of compact groups. As a consequence, for all $g \in G$, the linear isomorphism $\rho(g)$ is orthogonal with respect to $\langle \cdot, \cdot \rangle_G$.

Finally,  consider an orthonormal basis for $\langle \cdot, \cdot \rangle_G$ as well as the associated linear isomorphism $\Psi : \R^n \rightarrow \R^n$ which sends any vector of $\R^n$ onto the vector of its coordinates in the latter basis. If we equip the latter copy of $\R^n$ with the polynomial action
$$\kappa : \begin{array}{ccc} G \times \R^n & \rightarrow & \R^n\\ (g,x) & \mapsto & \Psi \circ \rho(g) \circ \Psi^{-1}(x)\end{array}$$
of $G$, the linear map $\kappa(g,\cdot) : \R^n \rightarrow \R^n$ is orthogonal with respect to the canonical Euclidean inner product on $\R^n$ for all $g \in G$ (since, if $v,w \in \R^n$, we have $\langle v,w \rangle_G = \langle \Psi(v), \Psi(w)\rangle$) and the linear isomorphism $\Psi$ is equivariant, hence the desired result.
\end{proof}

We can also adapt the previous proof to assert that any compact real algebraic group can be seen, up to an isomorphism of real algebraic groups, as a subgroup of an orthogonal group: 

\begin{prop} \label{propcompactrealalggroupisorthog}
The compact real algebraic group $G$ is isomorphic to a real algebraic subgroup of an orthogonal group.
\end{prop}

\begin{proof}
Thanks to proposition \ref{proprealalggroupislinear}, we can first suppose that $G$ is, up to isomorphism of real algebraic groups, a subgroup of $\left\{A \in \R^{M^2}~|~{\rm det}(A) \neq 0\right\}$, with $M \in \Nstar$. As in the previous proof, define, if $v,w \in \R^M$, $\langle v,w \rangle_G := \int_G \langle A v, A w \rangle {\rm d}\mu(A)$ where $\langle \cdot, \cdot \rangle$ denotes the canonical Euclidean inner product of $\R^M$ and $\mu$ is the normalized Haar measure of $G$. Then the matrices of $G$ are representative matrices in the canonical basis of $\R^M$ of orthogonal endomorphisms with respect to the Euclidean inner product $\langle \cdot, \cdot \rangle_G : \R^n \times \R^n \rightarrow \R$. Therefore, the conjugation by the change-of-basis matrix from the canonical basis of $\R^M$ to a fixed orthonormal basis of $\R^M$ for $\langle \cdot, \cdot \rangle_G$ sends $G$ onto a subgroup of the orthogonal group $\OO_M(\R) = \left\{A \in \R^{M^2}~|~{}^t\!A A = I_M\right\}$, which is a real algebraic subgroup by corollary~\ref{corimagecompactrealalggroupisrealalg}.
\end{proof}

\subsection{Compact subgroups of real orthogonal groups, Lie algebras and complexification} \label{subsectcompsubgrorthogpliealgcomplexif}
 
Before studying the geometric and semialgebraic quotients of real algebraic sets by polynomial actions of compact real algebraic groups, we take advantage of above proposition \ref{propcompactrealalggroupisorthog} to assert further rigidity properties of morphisms of compact real algebraic groups. For instance, as in the complex case but contrary to the general real case (see example \ref{remsimagemorphrealalggr} 2), a bijective morphism of real algebraic groups between compact real algebraic groups is an isomorphism of real algebraic groups.

In order to show such a result, we will consider compact real algebraic subgroups of orthogonal groups, the associated Lie algebras, their respective complexifications as well as the matrix exponential map.
\\ 

Let us first recall some basic facts about the differentiation of polynomial maps between nonsingular real algebraic sets and their complexification.
\\

Let $n \in \Nstar$ and let $X \subset \R^n$ be a nonsingular real algebraic set: $X$ is in particular a smooth (actually Nash) submanifold of $\R^n$. Let $x$ be a point of $X$. If $X'$ denotes the unique irreducible component of $X$ containing $x$ and if the vanishing ideal $I(X')$ of $X'$ is generated by polynomials $P_1,\ldots,P_r \in \R[x_1,\ldots,x_n]$, then the (smooth) tangent space $T_x X$ of $X$ at $x$ (which is by definition the vector subspace of $\R^n$ of vectors $v$ such that there exists a smooth curve $\gamma : I \rightarrow X$ defined on an open interval $I$ of $\R$ containing $0$ such that $\gamma(0) = x$ and $\gamma'(0) = v$) is the kernel of the matrix $\left(\frac{\partial P_i}{\partial x_j}(x)\right)_{1\leq i \leq r,\, 1\leq j \leq n}$. 

Moreover, if $Y \subset \R^N$ is a nonsingular real algebraic set and if $f : X \rightarrow Y$ is a restriction of a polynomial map $\widetilde{f} : \R^n \rightarrow \R^N$, then the differential ${\rm d}_x f : T_x X \rightarrow T_{f(x)} Y$ (which is by definition the unique linear map $F : T_x X \rightarrow T_{f(x)} Y$ such that, for any smooth curve $\gamma : I \rightarrow X$ defined on an open interval $I$ of $\R$ containing $0$ such that $\gamma(0) = x$, $F(\gamma'(0)) = (f \circ \gamma)'(0)$) is the restriction of the differential ${\rm d}_x \widetilde{f} : T_x \R^n = \R^n \rightarrow T_{f(x)} \R^N = \R^N$ of $\widetilde{f}$ and is therefore given by the Jacobian matrix of $\widetilde{f}$ at $x$.
\\

Consider the complexification $X_{\C}$ of $X$ and suppose that $X_{\C}$ is a nonsingular complex algebraic set: $X_{\C}$ is in particular a complex submanifold of $\C^n$. The complex vanishing ideal~$\mathsf{I}(X'_{\C})$ of the complexification $X'_{\C}$ of $X'$ is generated by the same polynomials $P_1,\ldots, P_r$ (see remark~\ref{remgeneralcomplexalgsetscomplexification}~2). As a consequence, the complex tangent space $T^{\C}_x X_{\C}$ of $X_{\C}$ at $x$ (which is by definition the $\C$-vector subspace of $\C^n$ of vectors $v$ such that there exists a holomorphic curve $\gamma : U \rightarrow X$ defined on a domain $U$ of $\C$ containing $0$ such that $\gamma(0) = x$ and $\gamma'(0) = v$) is the kernel (in $\C^n$) of the same matrix $\left(\frac{\partial P_i}{\partial x_j}(x)\right)_{1\leq i \leq r,\, 1\leq j \leq n}$. Therefore, if $v,w \in \R^n$, $v+i \, w \in T^{\C}_x X_{\C}$ if and only if $v,w \in T_x X$. 

Keep the above notations and suppose that the complexification $Y_{\C}$ is nonsingular as well. The complex differential ${\rm d}^{\C}_x f_{\C} : T^{\C}_x X_{\C} \rightarrow T^{\C}_{f(x)} Y_{\C}$ (which is by definition the unique $\C$-linear map $F : T^{\C}_x X_{\C} \rightarrow T^{\C}_{f(x)} Y_{\C}$ such that, for any holomorphic curve $\gamma : U \rightarrow X$ defined on a domain $U$ of $\C$ containing $0$ such that $\gamma(0) = x$, $F(\gamma'(0)) = (f \circ \gamma)'(0)$) is the restriction of the differential ${\rm d}^{\C}_x \widetilde{f}_{\C} : T_x \C^n = \C^n \rightarrow T_{f(x)} \C^N = \C^N$ of $\widetilde{f}_{\C}$ and is then given by the Jacobian matrix of the polynomial map $\widetilde{f}_{\C}$ at $x$, which is equal to the real Jacobian matrix of $\widetilde{f}$ at $x$. Therefore, if $v,w \in T_x X$, ${\rm d}^{\C}_x f_{\C}(v+i \, w) = {\rm d}_x f(v) + i \, {\rm d}_x f(w)$.
\\

Now, let $G$ be a compact subgroup of the real orthogonal group $\OO_n(\R)$: by corollary \ref{corcompactsubgroupofrealalggrisrealalg},~$G$ is a compact real algebraic subgroup of $\OO_n(\R)$. In particular, $G$ a nonsingular algebraic set of $\R^{n^2}$ (cf. lemma \ref{lemrealalggroupnonsing}) and a compact real Lie group. The Lie algebra of $G$ is by definition the tangent space ${\rm Lie}(G) := T_{I_n} G$ to $G$ at its neutral element $I_n$. 

Let ${\rm exp}$ denote the matrix exponential ${\rm M}_n(\R) \rightarrow {\rm M}_n(\R)$ on ${\rm M}_n(\R)$. We have ${\rm exp}\left({\rm Lie}(G)\right) \subset G$ and, if $\varphi : G \rightarrow H$ is a morphism of real algebraic groups to a compact subgroup of $\OO_N(\R)$, ${\rm exp} \circ {\rm d}_{I_n} \varphi = \varphi \circ {\rm exp}$: we give the proof of this fact with some additional remarks in section \ref{subsectsubgprealrotholiealgexpmap} of the appendix.
 
Consider the complexification $G_{\C}$ of $G$: by the same arguments as in the proof of lemma~\ref{lemrealalggroupnonsing}, the complex algebraic group $G_{\C}$ is a nonsingular complex algebraic set and then a complex Lie group. We can then consider its complex Lie algebra ${\rm Lie}^{\C}(G_{\C}) := T^{\C}_{I_n} G_{\C}$ and, considering the matrix exponential map ${\rm exp} : {\rm M}_n(\C) \rightarrow {\rm M}_n(\C)$, we can show by the same arguments as in the proof of lemma \ref{lemexpliealgcommwithdiff} (using holomorphic curves on domains of $\C$ containing $0$) that ${\rm exp}\left({\rm Lie}^{\C}(G_{\C})\right) \subset G_{\C}$ and, with above notations, that ${\rm exp} \circ {\rm d}^{\C}_{I_n} \varphi_{\C} = \varphi_{\C} \circ {\rm exp}$.
\\

We are going to show that the polar decomposition of invertible complex matrices allows to recover $G_{\C}$ from $G$ and ${\rm Lie}(G)$ via the matrix exponential map. We will need the following lemma, where we denote by ${\rm H}_n(\C) := \left\{H \in {\rm M}_n(\C)~|~{}^t\!\overline{H} = H\right\}$ the $\C$-vector space of Hermitian matrices of ${\rm M}_n(\C)$ and by ${\rm H}_n^+(\C)$ the subset of positive definite Hermitian matrices.

\begin{lem} \label{complexalgebraicmatgrouprealpowerofpositiveherm} Let $L$ be a subgroup of ${\rm GL}_n(\C)$ and suppose that $L \subset \C^{n^2}$ is a complex algebraic group. Let $H$ be a matrix of $L \cap {\rm H}_n^+(\C)$. Then, for any $t \in \R$, ${\rm exp}\left(t \, {\rm ln}(H)\right) \in L \cap {\rm H}_n^+(\C)$, where ${\rm ln}$ denotes the inverse map of the bijective restriction ${\rm H}_n(\C) \rightarrow {\rm H}_n^+(\C)$ of ${\rm exp}$.
\end{lem} 

\begin{proof} There exist a unitary matrix $U \in {\rm U}_n(\C)$ and ordered positive real numbers $\lambda_1,\ldots,\lambda_n$ such that ${}^t \overline{U} H U$ is the diagonal matrix ${\rm Diag}(\lambda_1,\ldots,\lambda_n)$. Then 
$${\rm ln}(H) = U \, {\rm Diag}\big({\rm ln}(\lambda_1),\ldots,{\rm ln}(\lambda_n)\big) {}^t \overline{U}$$
and, if $t \in \R$,
$${\rm exp}\left(t \, {\rm ln}(H)\right) = U \, {\rm Diag}\left(\lambda_1^t,\ldots,\lambda_n^t\right) {}^t \overline{U} \in {\rm H}_n^+(\C).$$
We then follow the arguments of (hints to) Problem 1.32 of Chapter 5 of \cite{OV} to prove that ${\rm exp}\left(t \, {\rm ln}(H)\right) \in L$: we show that ${\rm exp}\left(t \, {\rm ln}(H)\right) \in \mathsf{V}(\mathsf{I}(L))$. 

So let $p$ be a polynomial of the ideal $\mathsf{I}(L)$ of $\C[z_{k,l}, 1 \leq k,l \leq n]$. From the above expression of ${\rm exp}\left(t \, {\rm ln}(H)\right)$, we can assert that there are complex numbers $\alpha_1,\ldots,\alpha_M \in \C$ and positive real numbers $\mu_1 > \ldots > \mu_M$ such that $p\left({\rm exp}\left(t \, {\rm ln}(H)\right)\right) = \sum_{m=1}^M {\alpha}_m \mu_m^t$. Denote by $\mu(t)$ the latter expression and notice that, if $r \in \Nstar$, since ${\rm exp}\left(r \, {\rm ln}(H)\right) = H^r \in L$ ($L$ is a subgroup of ${\rm GL}_n(\C)$), we have $\mu(r) = p\left({\rm exp}\left(r \, {\rm ln}(H)\right)\right) = 0$ (because $p \in \mathsf{I}(L)$), so that 
$$\left|\left|\alpha_1 \mu_1^r\right| - \left|\sum_{m=2}^M {\alpha}_m \mu_m^r\right|\right| \leq |\mu(r)| = 0 \mbox{~\, i.e. } \left|\alpha_1 \mu_1^r\right| = \left|\sum_{m=2}^M {\alpha}_m \mu_m^r\right|.$$ 
Suppose by absurd that $\alpha_1 \neq 0$. Because 
$$\frac{\left|\sum_{m=1}^M {\alpha}_m \mu_m^r\right|}{\left|\alpha_1 \mu_1^r\right|} \leq \sum_{m=2}^M \left| \frac{\alpha_m}{\alpha_1}\right| \left(\frac{\mu_m}{\mu_1}\right)^r$$
tends to $0$ when $r$ tends to infinity, we have $\left|\alpha_1 \mu_1^r\right| > \left|\sum_{m=2}^M {\alpha}_m \mu_m^r\right|$ for $r$ big enough, which is in contradiction with the above equality. Therefore $\alpha_1 = 0$ and, by induction, $\alpha_2 = \ldots = \alpha_M = 0$: as a consequence, for any $t \in \R$, $\mu(t) = 0$ i.e. $p\left({\rm exp}\left(t \, {\rm ln}(H)\right)\right) = 0$, and then ${\rm exp}\left(t \, {\rm ln}(H)\right) \in \mathsf{V}(\mathsf{I}(L)) = L$ for any $t \in \R$.
\end{proof}

\begin{rems} \label{remafterlemcomplexalgebraicmatgrouprealpowerofpositiveherm}
~
\begin{enumerate}
	\item The arguments of the second part of the above proof show that ${\rm exp}\left(z \, {\rm ln}(H)\right) \in L$ for all~$z \in \C$, i.e. ${\rm ln}(H) \in {\rm Lie}^{\C}(L)$ (consider the complex analog of remark \ref{remafterlemexpliealgcommwithdiff} 1).
	\item Similarly, we can prove that, denoting by ${\rm S}_n(\R)$ the $\R$-vector space $\left\{S \in {\rm M}_n(\R)~|~{}^t\!S = S\right\}$ of symmetric matrices of ${\rm M}_n(\R)$ and by ${\rm S}_n^+(\R)$ the subset of positive definite symmetric matrices, if $K$ is a subgroup of ${\rm GL}_n(\R)$ which is a real algebraic group and $S$ is a matrix of $K \cap {\rm S}_n^+(\R)$, then ${\rm exp}\left(t \, {\rm ln}(S)\right) \in K \cap {\rm S}_n^+(\R)$ for all $t \in \R$ (in particular ${\rm ln}(S) \in {\rm Lie}(K)$: see remark \ref{remafterlemexpliealgcommwithdiff} 1).
\end{enumerate}
\end{rems} 
 
\begin{theo} \label{thpolardecompositionorthogonalcompactrealalggroup} The map $\Psi : G \times {\rm Lie}(G) \rightarrow G_{\C}~;~(R, T) \mapsto R \, {\rm exp}(iT)$ is a bijection. 
\end{theo} 
 
\begin{proof} First notice that if $T \in {\rm Lie}(G)$, then $iT \in {\rm Lie}^{\C}(G_{\C})$ so that ${\rm exp}(iT) \in G_{\C}$ (consider again the complex analog of remark \ref{remafterlemexpliealgcommwithdiff} 1). Therefore, if we furthermore consider a matrix $R \in G \subset G_{\C}$, we have $R \, {\rm exp}(iT) \in G_{\C}$.

Now, let $A \in G_{\C}$ (notice that the conjugate $\overline{A}$ also belongs to $G_{\C}$, as well as the inverse $A^{-1}$ since $G_{\C}$ is a group) and consider the polar decomposition of $A$: if we denote by $Q$ the unique square root of the positive definite Hermitian matrix ${}^t\!\overline{A}A$ (i.e. the unique matrix $Q \in {\rm H}_n^+(\C)$ such that $Q^2 = {}^t\!\overline{A}A$) and by $R$ the unitary matrix $A Q^{-1}$, the equality $A = R Q$ is the unique polar decomposition of $A$. Remark that $Q$ is actually the matrix ${\rm exp}\left(\frac{1}{2} {\rm ln}\left({}^t\!\overline{A}A\right)\right)$: the latter matrix is a positive definite Hermitian matrix and ${\rm exp}\left(\frac{1}{2} {\rm ln}\left({}^t\!\overline{A}A\right)\right)^2 = {\rm exp}\left({\rm ln}\left({}^t\!\overline{A}A\right)\right) = {}^t\!\overline{A}A$. In particular, according to previous lemma \ref{complexalgebraicmatgrouprealpowerofpositiveherm}, $Q = {\rm exp}\left(\frac{1}{2} {\rm ln}\left({}^t\!\overline{A}A\right)\right) \in G_{\C}$. As a consequence, the unitary matrix $R = A Q^{-1}$ belongs to $G_{\C}$. But $G_{\C} \subset \OO_n(\C)$, since $G \subset \OO_n(\R)$ (see also example \ref{excomplrealalggroup} 1), and ${\rm U}_n(\C) \cap \OO_n(\C) = \OO_n(\R)$, so that $R \in G$.

On the other hand, we have $\frac{1}{2}{\rm ln}\left({}^t\!\overline{A}A\right) \in {\rm Lie}^{\C}(G_{\C}) \cap {\rm H}_n(\C)$ (remark \ref{remafterlemcomplexalgebraicmatgrouprealpowerofpositiveherm}) while ${\rm Lie}^{\C}(G_{\C}) \cap {\rm H}_n(\C) = i {\rm Lie}(G)$: if $B+iC \in {\rm Lie}^{\C}(G_{\C}) \cap {\rm H}_n(\C)$ with $B,C \in {\rm Lie}(G)$, we have $B+iC = {}^t\!B - i \, {}^t C = - B + iC$ (since $G \subset \OO_n(\R)$, the elements of ${\rm Lie}(G)$ are anti-symmetric matrices: see remark \ref{remafterlemexpliealgcommwithdiff} 2) i.e. $B$ is the zero matrix $0_n$. 

Consequently, $A = R \, {\rm exp}(iT)$ with $R \in G$ and $T \in {\rm Lie}(G)$, hence the surjectivity of $\Psi$. As for the injectivity of $\Psi$, it is provided by the uniqueness of the polar decomposition and the injectivity of the restriction ${\rm H}_n(\C) \rightarrow {\rm H}_n^+(\C)$ of ${\rm exp}$ (notice that if $M$ is a real anti-symmetric matrix, then $iM$ is an Hermitian matrix).
\end{proof} 

We will use this theorem to prove the announced result on morphisms of compact real algebraic groups, through the following corollary:

\begin{cor} \label{cormorphismcompactsuborthoginjdiff} Let $H$ be a compact subgroup of a real orthogonal group $\OO_N(\R)$ and let $\varphi : G \rightarrow H$ be a morphism of real algebraic groups. 
\begin{enumerate}
	\item If the linear map ${\rm d}_{I_n} \varphi$ is injective, then ${\rm Ker} \, \varphi_{\C} = {\rm Ker} \, \varphi$.
	\item We have ${\rm Im} \, \varphi_{\C} = \left({\rm Im} \, \varphi\right)_{\C}$ (recall that ${\rm Im} \, \varphi$ is, in this context, a real algebraic group: cf. corollary \ref{corimagecompactrealalggroupisrealalg}).
\end{enumerate}
\end{cor}
 
\begin{proof}
	\begin{enumerate}
		\item Suppose that the map ${\rm d}_{I_n} \varphi$ is injective and let $A \in {\rm Ker} \, \varphi_{\C}$: by previous theorem~\ref{thpolardecompositionorthogonalcompactrealalggroup}, there exist $R \in G$ and $T \in {\rm Lie}(G)$ such that $A = \Psi(R,T) = R \, {\rm exp}(iT)$ and we then have
\begin{eqnarray*}
\Psi(I_n, 0_n) = I_n = \varphi_{\C}(A) = \varphi_{\C}(R) \varphi_{\C}({\rm exp}(iT)) = \varphi(R) {\rm exp} \left({\rm d}^{\C}_{I_n} \varphi_{\C}(iT)\right) & = & \varphi(R) {\rm exp} \left(i \, {\rm d}_{I_n}\varphi(T)\right) \\ & = &\Psi\left(\varphi(R), {\rm d}_{I_n}\varphi(T)\right).
\end{eqnarray*}
By the injectivity of $\Psi$, we have $\varphi(R) = I_n$ and ${\rm d}_{I_n}\varphi(T) = 0_n$, and then $T = 0_n$ by hypothesis. Therefore $A = R \in {\rm Ker} \, \varphi$.

		\item Since ${\rm Im} \, \varphi_{\C}$ is a complex algebraic group (cf. for instance Proposition 21.2.4 of \cite{TY}) containing ${\rm Im} \, \varphi$, we have $\left({\rm Im} \, \varphi\right)_{\C} \subset {\rm Im} \, \varphi_{\C}$. Conversely, if $A \in G_{\C}$ and if $R \in G$ and $T \in {\rm Lie}(G)$ are such that $A = R \, {\rm exp}(iT)$, we have $\varphi_{\C}(A) = \Psi\left(\varphi(R), {\rm d}_{I_n}\varphi(T)\right)$ (by the above equalities) and ${\rm d}_{I_n}\varphi(T) \in {\rm Lie}({\rm Im} \, \varphi)$ since, for all $t \in \R$, 
$${\rm exp}\left(t \, {\rm d}_{I_n}\varphi(T)\right) = {\rm exp}\left({\rm d}_{I_n}\varphi(t \, T)\right) = \varphi \left( {\rm exp}( t \, T)\right) \in \varphi(G)$$
because $T \in {\rm Lie}(G)$ (see lemma \ref{lemexpliealgcommwithdiff} and remark \ref{remafterlemexpliealgcommwithdiff} 1). As a consequence, $\varphi_{\C}(A) \in \left({\rm Im} \, \varphi\right)_{\C}$ (notice that ${\rm Im} \, \varphi$ is a compact subgroup of $\OO_N(\R)$).
	\end{enumerate}
\end{proof}
 
\begin{rem} Keeping the notations of above corollary \ref{cormorphismcompactsuborthoginjdiff}, if ${\rm d}_{I_n} \varphi$ is injective, then ${\rm Ker} \, \varphi_{\C} = {\rm Ker} \, \varphi$ is a finite set of (real) points. Indeed, in this case, since ${\rm Ker} \, \varphi_{\C}$ is a complex algebraic set containing ${\rm Ker} \, \varphi$, we have ${\rm Ker} \, \varphi \subset \left({\rm Ker} \, \varphi\right)_{\C} \subset {\rm Ker} \, \varphi_{\C} = {\rm Ker} \, \varphi$, hence $\left({\rm Ker} \, \varphi\right)_{\C} = {\rm Ker} \, \varphi$. This implies that $\dim {\rm Ker} \, \varphi = 0$, because the dimension of $\left({\rm Ker} \, \varphi\right)_{\C}$ considered as a real algebraic subset of $\R^{2 N^2}$ is $2 \dim {\rm Ker} \, \varphi$ (Propositions 2.2.1 and 2.2.5 of \cite{AK} and Proposition 3.1.1 of \cite{BCR}).
\end{rem} 

\begin{cor} \label{corbijmorphcomprealalggroupisiso} Let $\varphi : K \rightarrow H$ be a morphism of real algebraic groups between any compact real algebraic groups. If $\varphi$ is bijective, then $\varphi^{-1}$ is (the restriction of) a polynomial map i.e. $\varphi$ is an isomorphism of real algebraic groups.
\end{cor}
 
\begin{proof} Using proposition \ref{propcompactrealalggroupisorthog}, we can suppose that $\varphi$ is a morphism of real algebraic groups between compact subgroups of orthogonal groups. Suppose now that $\varphi$ is bijective: $\varphi$ is then a bijective morphism of Lie groups and is therefore an isomorphism of Lie groups. As a consequence, the differential ${\rm d}_{I_n}\varphi$ is bijective. According to above corollary \ref{cormorphismcompactsuborthoginjdiff} 1, the morphism of complex algebraic groups $\varphi_{\C}$ is then injective as well, and $\varphi_{\C}$ is also surjective by corollary \ref{cormorphismcompactsuborthoginjdiff} 2. Consequently, by corollary \ref{corisorealalggroupsiffcomplexificationbij}, $\varphi$ is an isomorphism of real algebraic groups.
\end{proof}

\begin{cor} \label{cormorphismcompactrealalggrinducesisopuq} Let $\varphi : K \rightarrow H$ be a morphism of real algebraic groups between compact real algebraic groups. Then $\varphi$ induces an isomorphism of real algebraic groups $\overline{\varphi} : K/ {\rm Ker} \, \varphi \rightarrow {\rm Im} \, \varphi$.
\end{cor}

\begin{proof} First, by corollary \ref{corimagecompactrealalggroupisrealalg}, the image of $\varphi$ is a compact real algebraic group and $\varphi$ then induces a morphism of real algebraic groups $\widetilde{\varphi} : K \rightarrow {\rm Im} \, \varphi~;~g \mapsto \varphi(g)$. Now, consider the (polynomial) action of the (compact) normal real algebraic subgroup ${\rm Ker} \, \varphi$ of $K$ on $K$ by left multiplication, as well as the trivial action of ${\rm Ker}\, \varphi$ on ${\rm Im} \, \varphi$. By lemma \ref{lemfunctorrealalgquotient}, the morphism $\widetilde{\varphi}$ induces the polynomial map $\overline{\varphi} := \widetilde{\varphi}_{\gq {\rm Ker} \, \varphi} : K\gq {\rm Ker} \, \varphi \rightarrow {\rm Im} \, \varphi$: by corollary~\ref{corimagecompactrealalggroupisrealalg}, $K\gq {\rm Ker} \, \varphi = K/ {\rm ker} \, \varphi$ and $\overline{\varphi}$ is a morphism of real algebraic groups between compact real algebraic groups. Moreover, $\overline{\varphi}$ is a group homomorphism since so are $\varphi$ and the real algebraic quotient $\pi_K$ of $K$ by ${\rm Ker} \, \varphi$ (cf. proposition \ref{propdefrealalgquotientgroup}): for all $g \in K$, we have $\overline{\varphi}\left(\pi_K(g)\right) = \varphi(g)$. Finally, $\overline{\varphi}$ is both surjective and injective, hence the result by corollary \ref{corbijmorphcomprealalggroupisiso}.
\end{proof} 
 
Corollary \ref{corbijmorphcomprealalggroupisiso} can also be used to show that any \emph{continuous} group homomorphism between compact real algebraic groups is necessarily the restriction of a polynomial map:

\begin{cor} \label{corcontinuousgpmorphismbtwcompralggpispoly} Let $\varphi : K \rightarrow H$ be a group homomorphism between compact real algebraic groups which is furthermore continuous. Then $\varphi$ is (the restriction of) a polynomial map.
\end{cor}

\begin{proof} We follow the arguments of the proof of Theorem 2.11 of Chapter 5 of \cite{OV}. Consider the graph $\Gamma_{\varphi} = \left\{(x,y) \in K \times H~|~y = \varphi(x)\right\}$ of $\varphi$: since $\varphi$ is continuous, $\Gamma_{\varphi}$ is a closed subset of $K \times H$. Being also a subgroup of the compact real algebraic group $K \times H$ (because $\varphi$ is a group homomorphism), $\Gamma_{\varphi}$ is itself a compact real algebraic group (corollary \ref{corcompactsubgroupofrealalggrisrealalg}). 

Now, consider the projections $\pi : \Gamma_{\varphi} \rightarrow K~;~(x,y) \mapsto x$ and $\varpi : \Gamma_{\varphi} \rightarrow H~;~(x,y) \mapsto y$, and remark that $\pi$ is bijective and $\varphi = \varpi \circ \pi^{-1}$. But $\pi$ is also a morphism of real algebraic groups between the compact real algebraic groups $\Gamma_{\varphi}$ and $K$, so that $\pi^{-1}$ is a morphism of real algebraic groups as well by corollary \ref{corbijmorphcomprealalggroupisiso}. Since $\varpi$ is also a morphism of real algebraic groups, so is the composition $\varphi = \varpi \circ \pi^{-1}$. In particular, $\varphi$ is the restriction of a polynomial map.
\end{proof}

\begin{exs} \label{exsapplicationsofcontgpmorphbtwcompralggpsispoly}
~
	\begin{enumerate}
		\item Any continuous group homomorphism between compact subgroups of general linear groups is polynomial in matrix entries (example \ref{excompactsubgroupofrealalggrisrealalg} 3).
		\item Any continuous (real finite-dimensional) representation of a compact real algebraic group~$G$ is a polynomial representation of $G$. Indeed, the image of a continuous representation map of $G$ is a compact subgroup of a general linear group and is therefore a compact real algebraic group as well (example \ref{excompactsubgroupofrealalggrisrealalg} 3).
		\item Any bijective continuous group homomorphism between compact real algebraic groups is an isomorphism of real algebraic groups (combine corollary \ref{corcontinuousgpmorphismbtwcompralggpispoly} and corollary \ref{corbijmorphcomprealalggroupisiso}).
	\end{enumerate}
\end{exs}

\begin{rem} \label{remcompliegpreprescompralggp} Any compact Lie group is isomorphic, as a Lie group, to a compact real algebraic group, and the latter is unique up to an isomorphism of real algebraic groups (see section 2.5 of Chapter 5 of \cite{OV}).
\end{rem}

\subsection{Geometric, semialgebraic and topological quotients by the polynomial action of a compact real algebraic group} \label{subsectgeosatopquotient}

Let $G$ be a compact real algebraic group and let $X \subset \R^n$ be a $G$-real algebraic set. Recall that, since $G$ is compact, $G$ is a linearly reductive real algebraic group (example \ref{exlinearlyredrealalggr} 2) so that the~$\R$-algebra $\mathcal{P}(X)^G$ is, by Hilbert's finiteness theorem \ref{theohilbertfiniteness}, finitely generated: if $f_1,\ldots,f_r$ are generators of $\mathcal{P}(X)^G$, the geometric quotient of $X$ by $G$ is then the polynomial map 
$$\pi_X : X \rightarrow X /\!\!/ G~;~x \mapsto (f_1(x),\ldots,f_r(x))$$
where $X \gq G = V \left(\left\{ p \in \R[z_1,\ldots,z_r]~|~p(f_1,\ldots,f_r) = 0\right\}\right)$ (see subsection \ref{subsectioncontrrealgeomquotient}).

The fact that $G$ is compact implies the following important fact that $\pi_X$ separates the orbits of the action of $G$ on $X$. If $x,y \in X$, recall that $G \cdot x = G \cdot y$ if and only if $(G \cdot x) \cap (G \cdot y) \neq \emptyset$ if and only if $y \in G \cdot x$.

\begin{prop} \label{propseparationorbits} For all $x,y \in X$, we have $G \cdot x = G \cdot y$ if and only if $\pi_X(x) = \pi_X(y)$. In other words, for all $x \in X$, $G \cdot x = \pi_X^{-1}(\pi_X(x))$.
\end{prop}

\begin{proof} The direct implication is given by the invariance of the polynomial functions $f_1,\ldots,f_r$. As for the converse one, let $x,y \in X$ and suppose that $(G \cdot x) \cap (G \cdot y) = \emptyset$. In particular, $y \notin G \cdot x$ and then, by lemma \ref{lemseporbitscompactsubg}, there exists an invariant polynomial function $p$ of $\mathcal{P}(X)^G$ such that $p(x) \neq p(y)$. As a consequence, $\pi_X(x) \neq \pi_X(y)$ since $\pi_X(x) = \pi_X(y)$ if and only if for all $k \in \{1,\ldots,r\}$, $f_k(x) = f_k(y)$ if and only if for all $f \in \mathcal{P}(X)^G$, $f(x) = f(y)$ (because the functions $f_1,\ldots,f_r$ generates the $\R$-algebra $\mathcal{P}(X)^G$).
\end{proof}

\begin{rem} 
The previous result is not true for general linearly reductive real algebraic groups. As a real counterpart of Example 2.3.1 of \cite{DK}, suppose that $n \geq 2$ and consider the polynomial action $K \times \R^n \rightarrow \R^n~;~(a,b),(x_1,\ldots,x_n) \mapsto (ax_1, \ldots,ax_n)$ of the linearly reductive real algebraic group $K = \{(a,b) \in \R^2~|~ab=1\}$ (example \ref{exlinearlyredrealalggr} 3) on $\R^n$. According to example~\ref{exgeninvariantalgorthogroup}~2, the real algebraic quotient map of $\R^n$ by $K$ is constant, therefore does not separate the infinitely many orbits of the action of $K$ on $\R^n$.
\end{rem}

\begin{cor} \label{corstabilityofimageintersectiondifference} Let $A$ be any $G$-stable subset of $X$ (i.e. for all $x \in A$, for all $g \in G$, $g \cdot x \in A$) and $B$ be any subset of $X$. Then $\pi_X(A \cap B) = \pi_X(A) \cap \pi_X(B)$ and $\pi_X(X \setminus A) = \pi_X(X) \setminus \pi_X(A)$.
\end{cor}

\begin{proof} Let $x \in A$ and $y \in B$ such that $\pi_X(x) = \pi_X(y)$: by proposition \ref{propseparationorbits}, there exists $g \in G$ such that $y = g \cdot x$, so that $y \in B \cap A$ (because $A$ is $G$-stable) and $\pi_X(y) \in \pi_X(A \cap B)$. As for the second equality, since $A$ is a $G$-stable subset of $X$ (notice that so is $X \setminus A$), we have $\pi_X(A) \cap \pi_X(X \setminus A) = \pi_X\big( A \cap (X \setminus A)\big) = \emptyset$.
\end{proof}

In the following of this subsection, we will be considering a $G$-stable semialgebraic subset~$S$ of $X$ and the semialgebraic quotient of $S$ by~$G$, that is the restriction $\varpi_S : S \rightarrow S/G$ of~$\pi_X$ (definition \ref{defsemialgquot}, see also lemma \ref{lemrestrictionrealalgquotient}). The previous properties will allow us to induce from any equivariant map between two $G$-semialgebraic sets a well-defined map between their semialgebraic quotients, in a functorial way. This operation will preserves semialgebraicity as well as continuity.

But let us first show that the semialgebraic quotient of $S$ by $G$ coincides with the topological quotient of $S$ by $G$. In order to prove this result, we will use some general topological properties that we state and prove in the appendix (see section \ref{subsectappendixtopology}).

\begin{prop} \label{propquotientmapproperclosedopen} Equip $S \subset \R^n$ and $S/G \subset \R^r$ with the respective topologies induced by the usual Euclidean topology. The continuous map $\varpi_S : S \rightarrow S/G$ is proper, closed and open, and~$\varpi_S$ is (up to left composition with an homeomorphism) the topological quotient of $S$ by the considered action of $G$.
\end{prop}

\begin{proof} By proposition \ref{proporthogonalizationactioncompactrealalggroup} and lemma \ref{lemfunctorrealalgquotient}, we can suppose that $\R^n$ is a polynomial orthogonal representation of $G$ ($\R^n$ being equipped with its canonical Euclidean inner product) and that $X$ is a $G$-real algebraic subset of $\R^n$. By lemma \ref{lemrestrictionrealalgquotient}, $\pi_X$ and $\varpi_S$ are restrictions of the geometric quotient map of $\R^n$ by the considered polynomial orthogonal action of $G$ and all these maps are then restrictions of the polynomial map $\pi : \mathbb{R}^n \rightarrow \mathbb{R}^r~;~x \mapsto (f_1(x),\ldots,f_r(x))$, where $f_1,\ldots,f_r$ are considered as generators of the $\R$-algebra $\mathcal{P}(\R^n)^G = \R[x_1,\ldots,x_n]^G$. Let us prove that the latter map $\pi$ is proper and closed.

Consider the polynomial map $\phi : \R^n \rightarrow \R~;~x_1^2+\cdots+x_n^2$: $\phi$ is an invariant polynomial of $\mathcal{P}(\R^n)^G$ ($\phi$ is invariant under right composition with any orthogonal isomorphism of $\R^n$), so that there exists $Q \in \R[z_1,\ldots,z_r]$ such that $\phi = Q(f_1,\ldots,f_r) = Q \circ \pi$. But the map $\phi = Q \circ \pi$ is proper (if $K$ is a compact subset of $\R$, $\phi^{-1}(K)$ is a closed subset of an Euclidean closed ball of $\R^n$, which is compact), so that, by lemma \ref{lemtopgenproperopenclosed} 1, $\pi$ is proper as well. By lemma \ref{lemtopgenproperopenclosed} 3, $\pi$ is also a closed map (since $\R^r$ is locally compact).

According to lemma \ref{lemtopgenproperopenclosed} 2, the restriction $\varpi_S : S \rightarrow S/G$ of $\pi$ is therefore proper as well ($\pi^{-1}(S/G) = \pi_{\R^n}^{-1}(S/G) = S$ by proposition \ref{propseparationorbits}). Furthermore, if $F$ is a closed subset of~$\R^n$, we have $\pi(S \cap F) = \pi(S) \cap \pi(F)$ by corollary \ref{corstabilityofimageintersectiondifference} (notice that if $T$ is any subset of~$\R^n$, $\pi(T) = \pi_{\R^n}(T)$) so that, by lemma \ref{lemtopgenproperopenclosed} 4, the restriction $\varpi_S$ of $\pi$ is closed. Finally, if $U$ is any open subset of $S$, denote by $\widetilde{U}$ the union of all open subsets $g \cdot U := \{g \cdot x~|~x \in S\}$, $g \in G$, of $S$: $\widetilde{U}$ is a $G$-stable open subset of $S$ such that $\varpi_S(\widetilde{U}) = \varpi_S(U)$ and we have $(S/G) \setminus \varpi_S(U) = (S/G) \setminus \varpi_S(\widetilde{U}) = \varpi_S(S \setminus \widetilde{U})$ (by corollary \ref{corstabilityofimageintersectiondifference}), so that $\varpi_S(U)$ is an open subset of $S/G$ since $\varpi_S$ is a closed map.

As for the last claim of proposition \ref{propquotientmapproperclosedopen}, denote $\mathcal{Q} := \left\{ G \cdot x~|~x \in S\right\}$ and denote by $\omega_S$ the topological quotient 
$$\begin{array}{ccc}S & \rightarrow & \mathcal{Q} \\x & \mapsto & G \cdot x\end{array}$$
of $S$ by $G$ where, by definition, a subset $V$ of $\mathcal{Q}$ is open if and only if the set $\{x \in S~|~G \cdot x \in V\}$ is an open subset of $S$. The map $\varsigma : \mathcal{Q} \rightarrow S/G~;~ G \cdot x \mapsto \varpi_S(x)$ is then a well-defined homeomorphism such that $\varpi_S = \varsigma \circ \omega_S$. Indeed, for any $x,y \in S$, $\varpi_S(x) = \varpi_S(y)$ if and only if $G \cdot x = G \cdot y$ (by proposition \ref{propseparationorbits}) and, for any subset $U$ of $S/G$, $U$ is open in $S/G$ if and only $\varpi_S^{-1}(U)$ is open in $S$ (because the map $\varpi_S$ is open and surjective) and $\varsigma^{-1}(U) = \{ G \cdot x \in \mathcal{Q}~|~\varpi_S(x) \in U\}$ is open in $\mathcal{Q}$ if and only 
$$\{x \in S~|~G \cdot x \in \varsigma^{-1}(U)\} = \{x \in S~|~\varpi_S(x) \in U\} = \varpi_S^{-1}(U)$$
is open in $S$.
\end{proof}

\begin{rem} The semialgebraic quotient map of a semialgebraic set by a polynomial action of any linearly reductive real algebraic group is not proper in general. Indeed, let $n \in \mathbb{N} \setminus \{0\}$ and consider for instance the polynomial action $K \times \R^n \rightarrow \R^n~;~(a,b),(x_1,\ldots,x_n) \mapsto (ax_1,\ldots,ax_n)$ of the non-compact linearly reductive real algebraic group $K = \{(a,b) \in \R^2~|~ab=1\}$ on~$\R^n$. Since $\mathcal{P}(\R^n)^K = \R$ is generated by $\overline{1}$ (example \ref{exgeninvariantalgorthogroup} 2), the semialgebraic quotient of $\R^n$ by~$K$ is the constant map $\R^n \rightarrow \{1\}~;~x \mapsto 1$ which is not proper (notice also that the semialgebraic quotient of the $K$-stable semialgebraic set $\R^n \setminus \{\bf{0}\}$ is not the topological quotient $\R^n \setminus \{\bf{0}\} \rightarrow \mathbb{P}^n(\R)$ of $\R^n \setminus \{\bf{0}\}$ by the considered action of $K$).
\end{rem}

As a first consequence of the topological properties of the semialgebraic quotient $\varpi_S : S \rightarrow S/G$, let us state that $\varpi_S$ preserves locally compactness in the following sense:

\begin{cor} \label{corsemialgquotientpreslocallycompactness} If $T$ is a $G$-stable locally compact subset of $S$, then $\varpi_S(T)$ is a locally compact subset of $S/G$. On the other hand, if $W$ is a locally compact subset of~$S/G$, then $\varpi_S^{-1}(W)$ is a ($G$-stable) locally compact subset of $S$.
\end{cor}

\begin{proof} Let $x \in T$. Since $T$ is locally compact, there exist an open subset $U$ of $S$ and a compact subset $K$ of $T$ such that $x \in U\cap T \subset K$. We then have $\varpi_S(x) \in \varpi_S(U \cap T) \subset \varpi_S(K)$ and, by corollary \ref{corstabilityofimageintersectiondifference}, $\varpi_S(U \cap T) = \varpi_S(U) \cap \varpi_S(T)$: since $\varpi_S$ is an open and continuous map (proposition \ref{propquotientmapproperclosedopen}), $\varpi_S(U)$ and $\varpi_S(K)$ are respectively an open subset and a compact subset of $S/G$. Consequently, the set $\varpi_S(T)$ is locally compact in $S/G$. 

As for the second part of the statement, let $x \in \varpi_S^{-1}(W)$ and let $V$ and $L$ be respectively an open subset of $S/G$ and a compact subset of $W$ such that $\varpi_S(x) \in V \cap W \subset L$. Then $x \in \varpi_S^{-1}(V) \cap \varpi_S^{-1}(W) \subset \varpi_S^{-1}(L)$: since $\varpi_S$ is continuous, $\varpi_S^{-1}(V)$ is an open subset of~$S$ and, since $L$ is compact and $\varpi_S$ is proper (proposition \ref{propquotientmapproperclosedopen}), $\varpi_S^{-1}(L)$ is compact. As a consequence,~$\varpi_S^{-1}(W)$ is locally compact in $S$.
\end{proof}

\begin{rem} A subset $D$ of a locally compact Hausdorff space $A$ is locally compact if and only if it is locally closed (see lemma \ref{lemequivloccomplocclosedloccomphaus} 4 in the appendix).
\end{rem}

\begin{cor} There exist $G$-stable locally compact semialgebraic subsets $S_1,\ldots,S_l$ of $S$ such that $S = \bigsqcup_{k=1}^l S_k$.
\end{cor} 
 
\begin{proof} By Theorem 2.3.6 of \cite{BCR}, the semialgebraic subset $S/G$ of $\R^r$ is the disjoint union of semialgebraic sets $T_1,\ldots,T_l$, each of them being semialgebraically homeomorphic to an open hypercube $]0;1[^d \subset \R^d$, $d \in \mathbb{N}$: for all $i \in \{1,\ldots,l\}$, $T_i$ is therefore locally compact. We then have $S = \varpi_S^{-1}(S/G) = \bigsqcup_{i=1}^l \varpi_S^{-1}(T_i)$ and, for all $i \in \{1,\ldots,l\}$, $\varpi_S^{-1}(T_i)$ is a $G$-stable locally compact (by previous corollary \ref{corsemialgquotientpreslocallycompactness}) semialgebraic set (Proposition 2.2.7 of \cite{BCR}).
\end{proof} 

Let us now deal with functoriality properties of the semialgebraic quotient. Let $T$ be a~$G$-semialgebraic set and let $f : S \rightarrow T$ be any equivariant map. Proposition \ref{propseparationorbits} allows to define the map
$$f_{/G} : S/G \rightarrow T/G~;~\varpi_S(x) \mapsto \varpi_T(f(x))$$
(the map $f_{/G}$ is well-defined up to left and right compositions with polynomial isomorphisms). This functorial operation preserves semialgebraicity, continuity and properness:

\begin{prop} \label{propfunctsemialgquotientsemialgcontprop} If the map $f$ is semialgebraic (i.e. the graph of $f$ is a semialgebraic set), resp. continuous, resp. proper, then so is $f_{/G}$.
\end{prop}

\begin{proof} Suppose that $f$ is semialgebraic. The graph 
\begin{eqnarray*}
\Gamma_{f_{/G}} & = & \left\{(w,z) \in S/G \times T/G ~|~z = f_{/G}(w)\right\} \\
		            & = & \left\{(w,z) \in S/G \times T/G ~|~\exists x \in S,~\exists y \in T,~\varpi_S(x) = w,~\varpi_T(y) = z,~y = f(x)\right\}
\end{eqnarray*}
of $f_{/G}$ is the image of the graph $\Gamma_f = \left\{(x,y) \in S \times T ~|~y = f(x)\right\}$ of $f$ by the (restriction of a) polynomial map $S \times T \rightarrow S/G \times T/G~;~(x,y) \mapsto (\varpi_S(x), \varpi_T(y))$. Therefore, since $\Gamma_f$ is semialgebraic, so is $\Gamma_{f_{/G}}$ by Tarski-Seidenberg theorem (Proposition 2.2.7 of \cite{BCR}).   

Suppose now that $f$ is continuous and let $U$ be an open subset of $T/G$, then $f_{/G}^{-1}(U) = \varpi_S \left(f^{-1}\left(\varpi_T^{-1}(U)\right)\right)$ is an open subset of $S/G$ since $\varpi_T$ is continuous, $f$ is continuous and $\varpi_S$ is open (proposition \ref{propquotientmapproperclosedopen}). 

Finally, if $f$ is proper and $K$ is a compact subset of $T/G$ then $f_{/G}^{-1}(K) = \varpi_S \left(f^{-1}\left(\varpi_T^{-1}(K)\right)\right)$ is a compact subset of $S/G$ since $\varpi_T$ is proper (proposition \ref{propquotientmapproperclosedopen}), $f$ is proper and $\varpi_S$ is continuous.
\end{proof}

\begin{rem} 
If $f$ is (the restriction of) a polynomial map, the map $f_{/G}$ defined above coincides with the polynomial map $f_{/G} : S/G \rightarrow T/G$ defined in lemma \ref{lemfunctorrealalgquotient}.
\end{rem}

Let us end this part with the general remark that the geometric (resp. semialgebraic) quotient of a real algebraic (resp. semialgebraic) set by a polynomial action of a compact real algebraic group can always be computed, up to polynomial isomorphism, as a restriction of the geometric quotient of a real affine space $\R^N$ by a closed subgroup of $\OO_N(\R)$:

\begin{rem} By proposition \ref{proporthogonalizationactioncompactrealalggroup} and lemma \ref{lemfunctorrealalgquotient}, the geometric quotient of $X$ by $G$ is polynomially isomorphic to the geometric quotient of a $G$-real algebraic subset of a polynomial orthogonal representation $\R^N$ of $G$ (where $\R^N$ is equipped its canonical Euclidean inner product) which is, by lemma \ref{lemrestrictionrealalgquotient}, a restriction of the geometric quotient of $\R^N$ by $G$. 

Let $Y$ be any $G$-real algebraic subset of $\R^N$. If $\rho : \R^N \rightarrow \mathcal{GL}(\R^N)$ is the group homomorphism associated to the representation $\R^N$ of $G$, consider the morphism of real algebraic groups $\varrho : G \rightarrow \OO_N(\R)$ which associates to any element $g$ of $G$ the matrix of $\rho(g)$ in the canonical basis of~$\R^N$: by corollary \ref{corimagecompactrealalggroupisrealalg}, the image $K$ of $\varrho$ is a compact real algebraic subgroup of $\OO_N(\R)$, with respect to which $\R^N$ is a $K$-real algebraic set and $Y$ is a $K$-real algebraic subset of $\R^N$. Since a polynomial function of $\mathcal{P}(Y)$ is invariant with respect to the action of $G$ on $Y$ if and only if it is invariant with respect to the action of $K$, the geometric quotient $Y \rightarrow Y \gq G$ of $Y$ by $G$ is the same map as the geometric quotient $Y \rightarrow Y \gq K$ of $Y$ by $K$, which is a restriction of the geometric quotient of $\R^N$ by $K$.
\end{rem}

\subsection{Dimensions of the orbits}

Keep the previous notations and let $x \in X$: by proposition \ref{proporbitcompactsubgrouprealalggroup}, the orbit $G \cdot x = \{\alpha(g,x),~g \in G\}$ of $x$ is a compact real algebraic set. Notice furthermore that the stabilizer $G_x := \left\{g \in G~|~\alpha(g,x) = x\right\}$ of $x$ is a (compact) real algebraic subgroup of $G$ (since $G$ is a real algebraic set and $\alpha$ is a polynomial map). The dimensions of the real algebraic sets $G \cdot x$, $G$ and $G_x$ are related by the following equality:

\begin{lem} \label{lemdimorbitstabilizer} We have $\dim G \cdot x = \dim G - \dim G_x$.
\end{lem}

\begin{proof} Consider the orbit map $\alpha_x : G \rightarrow G \cdot x~;~g \mapsto \alpha(g,x)$. Since the map $\alpha_x$ is a (restriction of a) polynomial map, it is in particular continuous and semialgebraic and we can apply Hardt's semialgebraic triviality theorem (Theorem 9.3.2 of \cite{BCR}) to $\alpha_x$: there exist a finite partition of~$G \cdot x$ into semialgebraic sets $T_1,\ldots,T_s$, semialgebraic sets $F_1,\ldots,F_s$ and, for each $i \in \{1,\ldots,s\}$, a semialgebraic homeomorphism $\theta_i : T_i \times F_i \rightarrow \alpha_x^{-1}(T_i)$ such that $\alpha_x \circ \theta_i$ is the projection $T_i \times F_i \rightarrow T_i~;~(y,w) \mapsto y$.

Let $i \in \{1,\ldots,s\}$ and $y \in T_i$. The map $\theta_i$ restricts into a semialgebraic homeomorphism $\{y \} \times F_i \rightarrow \alpha_x^{-1}(\{y\})$. On the other hand, if we write $y = \alpha(g,x)$ with $g \in G$, we have 
$$\alpha_x^{-1}(\{y\}) = \{ h \in G~|~ \alpha(h,x) = \alpha(g,x)\} = g G_x,$$
so that $F_i$ is semialgebraically homeomorphic to the stabilizer $G_x$ (left multiplication by $g$ in $G$ induces a polynomial isomomorphism between $G_x$ and $g G_x$) and then $\dim \alpha_x^{-1}(T_i) = \dim T_i + \dim F_i = \dim T_i + \dim G_x$ (Theorem 2.8.8 of \cite{BCR}).

As a consequence, if $i_0 \in \{1,\ldots,s\}$ is such that $\dim T_{i_0} = \underset{i \in \{1,\ldots,s\}}{\rm max} \dim T_i = \dim G \cdot x$ (Proposition 2.8.5 of \cite{BCR}), we have 
$$\dim G = \underset{i \in \{1,\ldots,s\}}{\rm max} \dim \alpha_x^{-1}(T_i) = \dim \alpha_x^{-1}(T_{i_0}) = \dim T_{i_0} + \dim G_x = \dim G \cdot x + \dim G_x.$$
\end{proof}

We can also use Hardt's semialgebraic triviality theorem to bound the dimension of the semialgebraic quotient $S/G$ of the $G$-stable semialgebraic subset $S$ of $X$:

\begin{prop} \label{lembounddimsemialgquotient} We have 
$$(\dim S - \dim G \leq)~\dim S - \underset{x \in S}{\rm max} \dim G \cdot x \leq \dim S/G \leq \dim S - \underset{x \in S}{\rm min} \dim G \cdot x.$$
\end{prop}

\begin{proof} Apply Hardt's semialgebraic triviality theorem (Theorem 9.3.2 of \cite{BCR}) to the semialgebraic quotient $\varpi_S : S \rightarrow S/G$: there exist a finite partition of $S/G$ into semialgebraic sets $T_1,\ldots,T_s$, semialgebraic sets $F_1,\ldots,F_s$ and, for each $i \in \{1,\ldots,s\}$, a semialgebraic homeomorphism $\theta_i : T_i \times F_i \rightarrow \varpi_S^{-1}(T_i)$ such that $\varpi_S \circ \theta_i$ is the projection $T_i \times F_i \rightarrow T_i~;~(z,w) \mapsto z$.

Let $i \in \{1,\ldots,s\}$ and $z \in T_i$, the map $\theta_i$ restricts to a semialgebraic homeomorphism $\{z \} \times F_i \rightarrow \varpi_S^{-1}(\{z\})$. On the other hand, if $x \in \varpi_S^{-1}(\{z\})$, we have $\varpi_S^{-1}(\{z\}) = G \cdot x$ by proposition \ref{propseparationorbits}, therefore $F_i$ is semialgebraically homeomorphic to $G \cdot x$ and then
$$\dim \varpi_S^{-1}(T_i) = \dim T_i + \dim F_i = \dim T_i + \dim G \cdot x.$$

As a consequence, if $i$ is such that $\dim \varpi_S^{-1}(T_i) = \dim S$, we have
$$\dim S = \dim \varpi_S^{-1}(T_i) = \dim T_i + \dim G \cdot x \leq \dim S/G + \underset{x \in S}{\rm max} \dim G \cdot x \leq \dim S/G + \dim G,$$
while, if $i$ is such that $\dim T_i = \dim S/G$, we have
$$\dim S \geq \dim \varpi_S^{-1}(T_i) = \dim T_i + \dim G \cdot x = \dim S/G + \dim G \cdot x \geq  \dim S/G + \underset{x \in S}{\rm min} \dim G \cdot x,$$
hence the result.
\end{proof}

\begin{cor} Suppose that the action of $G$ on $S$ is free i.e. for all $x \in S$, $G_x = \{e\}$ (where~$e$ denotes the neutral element of $G$). Then $\dim S/G = \dim S - \dim G$.
\end{cor}

\begin{proof} For any $x \in S$, since $G_x = \{e\}$, the polynomial map $\alpha_x : G \rightarrow G \cdot x$ is bijective and therefore $\dim G \cdot x = \dim \alpha_x(G) = \dim G$ by Theorem 2.8.8 of \cite{BCR}. As a consequence, 
$$\underset{x \in S}{\rm max} \dim G \cdot x = \underset{x \in S}{\rm min} \dim G \cdot x = \dim G$$
and then, by previous proposition \ref{lembounddimsemialgquotient}, $\dim S/G = \dim S - \dim G$.
\end{proof}

\begin{ex} If $H$ is a normal closed subgroup of $G$, we have $\dim G/H = \dim G - \dim H$ (see corollary \ref{corimagecompactrealalggroupisrealalg}). As a consequence, if $K$ is a compact real algebraic group and $\varphi : G \rightarrow K$ is a morphism of real algebraic groups, we have $\dim G = \dim {\rm Ker} \, \varphi + \dim {\rm Im} \, \varphi$ (see corollary~\ref{cormorphismcompactrealalggrinducesisopuq}).
\end{ex}

Let us also consider the complexifications of the orbits of the $G$-real algebraic set $X$. Because the real algebraic group $G$ is compact, the latter are the orbits of the complexification:

\begin{prop} \label{lemrealcomplexorbits} Let $x \in X$. The orbit $G_{\C} \cdot x = \{h \cdot x~|~h \in G_{\C}\}$ of $x$ under the action of $G_{\C}$ on $X_{\C}$ is a complex algebraic subset of $\C^n$ which is the complexification of $G \cdot x$. In particular, the dimension of the complex algebraic set $G_{\C} \cdot x$ is equal to $\dim G \cdot x$.  
\end{prop}

\begin{proof} In order to deal with the first point, let $f : X \rightarrow Y$ be an isomorphism of $G$-real algebraic sets onto a $G$-real algebraic subset $Y$ of a polynomial representation $\R^N$ of $G$ (proposition~\ref{proplinearizgrealalgset}): the complexification $f_{\C} : X_{\C} \rightarrow Y_{\C}$ is then an isomorphism of $G_{\C}$-complex algebraic sets (lemma~\ref{lemcomplexificationmorphismGras}). Denote $y := f(x)$. Since $G_{\C}$ is a linearly reductive algebraic group over $\R$ (cf. remark \ref{remgeneralcomplexalgsetscomplexification} 2. and proposition \ref{propequivlinearredragandcomplexification}) and because the orbit~$G \cdot y$ is a compact subset of $\R^n$, we can assert by Corollary 5.3 of \cite{Birkes} that the orbit $G_{\C} \cdot y$ of $y$ under the linear action of $G_{\C}$ on $\C^N$ (given by the complexification of the linear action of $G$ on $\R^N$: see example~\ref{exscomplexificationpolyaction} 3) is Zariski-closed. Since $f_{\C}$ is an equivariant polynomial isomorphism, it follows that $G_{\C} \cdot x = f_{\C}^{-1}\left(G_{\C} \cdot y \right)$ is itself a complex algebraic set of $\C^n$.

As for the second point, consider the complexification $\left(\alpha_x\right)_{\C} : G_{\C} \rightarrow \left(G \cdot x\right)_{\C}~;~g \mapsto \alpha_{\C}(g,x)$ of the orbit map $\alpha_x$ (see the proof of lemma \ref{lemdimorbitstabilizer}): we have in particular $G_{\C} \cdot x \subset \left(G \cdot x\right)_{\C}$ and, since $G_{\C} \cdot x$ is a complex algebraic set containing $G \cdot x$, we actually have the equality $\left(G \cdot x\right)_{\C} = G_{\C} \cdot x$.
\end{proof}

\begin{rem} \label{remcomplexorbitnotclosed} Keeping the notations of previous proposition, notice that the points of $X_{\C} \setminus X$ may have a non-closed orbit. Consider for instance the usual linear action of $\SO_2(\R)$ on $\R^2$ and its complexification (see example~\ref{exscomplexificationpolyaction} 1). Any matrix $A$ of $\SO_2(\C)$ writes $\begin{pmatrix}a&-b\\b&a\end{pmatrix}$ with $a,b \in \C$ satisfying $a^2+b^2=1$, and we then have $A(1,i) = (a-ib) (1,i)$. As a consequence, the orbit $\SO_2(\C) \cdot (1,i)$ is the linear line generated by $(1,i)$ minus the origin (consider the biregular isomorphism  of affine complex algebraic varieties $S^1_{\C} = \left\{(a,b) \in \C^2~|~a^2+b^2=1\right\} \rightarrow \C^{*}~;~(a,b) \mapsto a+ib$, with inverse $\C^{*} \rightarrow S^1_{\C}~;~z \mapsto \left(\frac{z+z^{-1}}{2}, \frac{z-z^{-1}}{2i}\right)$). In particular, $\SO_2(\C) \cdot (1,i)$ is not closed.
\end{rem}

\begin{cor} For all $x,y \in X$, we have $G_{\C} \cdot x = G_{\C} \cdot y$ if and only if $\pi_{X_{\C}}(x) = \pi_{X_{\C}}(y)$.
\end{cor}

\begin{proof} Recall that the complex geometric quotient $\pi_{X_{\C}} : X_{\C} \rightarrow X_{\C}\gq G_{\C}$ is the complexification of $\pi_X$ and is given by the complexifications $(f_1)_{\C},\ldots,(f_r)_{\C}$ which are $G_{\C}$-invariant (see subsection \ref{subsectioncomplexrealgeomquotient}), hence the direct implication. As for the converse implication, let $x,y \in X$ such that $\pi_{X_{\C}}(x) = \pi_{X_{\C}}(y)$ (i.e. $\pi_X(x) = \pi_X(y)$): the respective complex Zariski closures of the orbits $G_{\C} \cdot x$ and $G_{\C} \cdot y$ in $\C^n$ then intersect (see for instance Corollary 2.3.8 of \cite{DK}). But, by previous proposition \ref{lemrealcomplexorbits}, the orbits $G_{\C} \cdot x$ and $G_{\C} \cdot y$ are Zariski-closed, so that $G_{\C} \cdot x$ and $G_{\C} \cdot y$ intersect i.e. $G_{\C} \cdot x = G_{\C} \cdot y$.
\end{proof}

\begin{rems} \label{remsafterrealpointssamecomplexorbitsequivsameimages}
~
	\begin{enumerate}
		\item The converse implication of the previous statement is in general not true for orbits that that does not contain real points. For instance, keeping the context of the counterexample of remark \ref{remcomplexorbitnotclosed}, notice that the respective images of $(1,i)$ and $(1,-i)$ by the complex geometric quotient map $\C^2 \rightarrow \C^2 \gq \SO_2(\C) = \C~;~(u,v) \mapsto u^2+v^2$ (which is the complexification of the real geometric quotient $\R^2 \rightarrow \R^2 \gq \SO_2(\R) = \R~;~(x,y) \mapsto x^2+y^2$: see examples \ref{exrealalgquotientplaneorthoaction} 1 and \ref{exgeninvariantalgorthogroup} 1) are both $0$. However, the orbits $\SO_2(\C) \cdot (1,i) = \C^* (1,i)$ and $\SO_2(\C) \cdot (1,-i) = \C^* (1,-i)$ do not intersect (their Zariski closures intersect at the origin of $\C^2$).
		\item We can also show the following fact: if $Y$ is a $G$-real algebraic subset of $X$ and if $x$ is a point of $X$ such that $\pi_X(x) \in Y /\!\!/ G$, then $x \in Y$ (in other words, $\pi_X^{-1}\left(Y \gq G\right) = Y$). Indeed, since
$$Y /\!\!/ G = \left(Y_{\C} /\!\!/ G_{\C}\right) \cap \R^r= \left(\pi_{X_{\C}}\left(Y_{\C}\right)\right) \cap \R^r$$
(see subsection \ref{subsectioncomplexrealgeomquotient}), there exists $y \in Y_{\C}$ such that $\pi_X(x) = \pi_{X_{\C}}(x) = \pi_{X_{\C}}(y)$: the complex Zariski closures of the orbits $G_{\C} \cdot x$ and $G_{\C} \cdot y$ therefore intersect (Corollary 2.3.8 of \cite{DK}). But the orbit $G_{\C} \cdot x$ is Zariski-closed by proposition \ref{lemrealcomplexorbits}, so that there exists $g \in G_{\C}$ such that $g \cdot x$ belongs to the Zariski closure of $G_{\C} \cdot y$, which is contained in $Y_{\C}$. As a consequence, $x \in Y_{\C} \cap \R^n = Y$.
		\item Any $G$-real algebraic subset $Y$ of $X$ can be described as the zero set of a $G$-invariant polynomial function on $X$ since, by the previous item, $Y = \pi_X^{-1}\left(Y \gq G\right)$ (the coordinate functions of $\pi_X$ are $G$-invariant polynomial functions on $X$ and $Y \gq G$, as any real algebraic set, is the zero set of a single polynomial). In particular, if we suppose that $\R^n$ is a $G$-real algebraic set such that $X$ is a $G$-real algebraic subset of $\R^n$, then $X$ is the zero set of a single $G$-invariant polynomial of $\R[x_1,\ldots,x_n]$.
	\end{enumerate}
\end{rems}

\section{Free polynomial actions of compact real algebraic groups} \label{sectionfreepolyactionscompralggp}

We dedicate this part to free polynomial actions of compact real algebraic groups on real algebraic sets and their geometric and semialgebraic quotients. If $G$ is a compact real algebraic group and $X$ is a $G$-real algebraic set on which $G$ acts freely, i.e. for all $x \in X$, $G_x = \{e\}$ (if $e$ denotes the neutral element of $G$), we will notably establish that the semialgebraic quotient~$X/G$ satisfies an additional geometric property, namely symmetry with respect to analytic arcs. This result will be a consequence of the fact that, because the action of $G$ is free, the geometric quotient map sends nonsingular points of $X$ to nonsingular points of $X\gq G$. We will also give a proof of the latter property, using the orbit maps associated to the action of $G$ on the points of $X$: since $G$ acts freely on $X$, all orbit maps are polynomial isomorphisms.

\subsection{Orbit map of a point with trivial stabilizer}

For all this section, fix a compact real algebraic group $G \subset \R^m$ with neutral element $e$ and let $X \subset \R^n$ be any $G$-real algebraic set with polynomial action $\alpha : G \times X \rightarrow X$ of $G$. We begin with the property that the orbit map associated to a point with trivial stabilizer is a polynomial isomorphism:

\begin{prop} \label{proporbitmapisoiftrivialstabilizer} Let $x \in X$ such that $G_x = \{e\}$. The polynomial map $\alpha_x : G \rightarrow G \cdot x~;~g \mapsto g \cdot x$ is a polynomial isomorphism (recall that the orbit $G \cdot x$ is a real algebraic set by proposition~\ref{proporbitcompactsubgrouprealalggroup}).
\end{prop}

\begin{proof} Thanks to proposition \ref{proporthogonalizationactioncompactrealalggroup}, we can suppose that $\R^n$ is a polynomial orthogonal representation of $G$ (with respect to the canonical Euclidean inner product on $\R^n$) and that $X$ is a~$G$-real algebraic subset of $\R^n$. Let $\rho : \R^n \rightarrow \mathcal{GL}(\R^n)$ be the group homomorphism associated to the representation $\R^n$ of $G$ and consider the morphism of real algebraic groups $\varrho : G \rightarrow \OO_n(\R)$ which associates to any element $g$ of $G$ the matrix of $\rho(g)$ in the canonical basis of~$\R^n$. Denote $K := {\rm Im} \, \varrho$: $K$ is a compact real algebraic group (by corollary \ref{corimagecompactrealalggroupisrealalg}) and~$X$ is a $K$-real algebraic set. Furthermore, $\varrho$ is an injective map: if $g \in G$ satisfies $\varrho(g) = I_n$, then $g \cdot x = x$ so that $g = e$ since $G_x = \{e\}$. As a consequence,~$\varrho$ induces a bijective morphism of real algebraic groups $G \rightarrow K$, which is an isomorphism of real algebraic groups by corollary~\ref{corbijmorphcomprealalggroupisiso} because $G$ and $K$ are both compact. Therefore, the orbit map $\alpha_x : G \rightarrow G \cdot x$ is a polynomial isomorphism if and only if the orbit map $\varphi_x : K \rightarrow K \cdot x~;~ A \mapsto Ax$ is a polynomial isomorphism (we have $K \cdot x = G \cdot x$ and $\alpha_x = \varphi_x \circ \varrho$).

Consider the stabilizer $H := \left(K_{\C}\right)_x$ of $x$ for the complexified action of $K_{\C}$ on $X_{\C}$. First remark that, by Proposition 21.4.3 (iii) of \cite{TY}, proposition \ref{lemrealcomplexorbits} and lemma \ref{lemdimorbitstabilizer}, we have 
$$\dim H = \dim K_{\C} - \dim K_{\C} \cdot x = \dim K - \dim K \cdot x = \dim K_x = 0$$
(the stabilizer $K_x$ is reduced to $I_n$), so that $H$ is a finite group.

Now, let $A \in H$. Since $K$ is a compact real algebraic subgroup of $\OO_n(\R)$, $H$ is a complex algebraic subgroup of $\OO_n(\C)$, in particular ${}^t\!A = A^{-1} \in H$. Moreover, $\overline{A}x = \overline{A x} = \overline{x} = x$ ($x$ has real coordinates), therefore $\overline{A} \in H = \left(K_{\C}\right)_x$. Consequently, the positive definite Hermitian matrix ${}^t\!\overline{A} A$ is in $H$. On the other hand, there exist a unitary matrix $U \in {\rm U}_n(\C)$ and a diagonal matrix $D$ with positive real diagonal coefficients such that $U^{-1}  {}^t\!\overline{A} A U = D$. If $N$ denotes the order of the finite group $H$, we have $D^N = U^{-1}  \left({}^t\!\overline{A} A\right)^N U = I_n$, which implies that $D = I_N$ (because the diagonal matrix $D$ has positive real diagonal coefficients) i.e. ${}^t\!\overline{A} A = I_n$ i.e. $A$ is a unitary matrix.

Since ${\rm U}_n(\C) \cap \OO_n(\C) = \OO_n(\R)$, we deduce that $H = \left(K_{\C}\right)_x = K_x = \{I_n\}$. The complex orbit map $K_{\C} \rightarrow K_{\C} \cdot x~;~A \mapsto Ax$, which is the complexification of $\varphi_x$ by proposition \ref{lemrealcomplexorbits}, is therefore bijective. By \cite{TY} Theorem 25.1.2 (iv), we can then assert that the $K_{\C}$-equivariant polynomial map $\left(\varphi_x\right)_{\C} : K_{\C} \rightarrow K_{\C} \cdot x$, between $K_{\C}$-complex algebraic sets on which $K_{\C}$ acts transitively (considering the transitive action of $K_{\C}$ on itself by left multiplication and recalling that $K_{\C} \cdot x$ is a complex algebraic set by proposition \ref{lemrealcomplexorbits}), is a polynomial isomorphism, so that $\varphi_x$ is a polynomial isomorphism by lemma \ref{lempolyisomiffcomplexpolyisom}.
\end{proof}

\begin{rem} Even if the action of $G$ on $X$ is free, the action of $G_{\C}$ on $X_{\C}$ is not free in general. Consider for instance the real algebraic set $Y := \left\{(x,y,z) \in \R^3~|~x^2+y^2 = z^2 + 1\right\}$ and the polynomial action $S^1 \times Y \rightarrow Y~;~(a,b), (x,y,z) \mapsto (ax -by, bx + ay,z)$ of $S^1$ on $Y$: this action is free whereas the stabilizer of the point $(0,0,i)$ of $Y_{\C} = \left\{(x,y,z) \in \C^3~|~x^2+y^2 = z^2 + 1\right\}$ under the complexified action of $S^1_{\C}$ is the entire group $S^1_{\C}$ (notice that, by Theorem 4.5.1 of~\cite{BCR},~$I(Y)$ is generated by the irreducible polynomial $x^2+y^2-z^2-1$ of $\R[x,y,z]$, so that $Y_{\C}$ is given by the same polynomial equation $x^2+y^2-z^2-1$ in $\C^3$: see remark \ref{remgeneralcomplexalgsetscomplexification} 2.)
\end{rem}

\subsection{Nonsingular point with trivial stabilizer}

We are going to use above proposition \ref{proporbitmapisoiftrivialstabilizer} to show that the image of a nonsingular point of~$X$ with trivial stabilizer by the geometric quotient $\pi_X$ is a nonsingular point of $X \gq G$. We will also use Nakayama's lemma (lemma \ref{lemnakayamaandconsq} of the appendix) as well as the following property:

\begin{lem} \label{lemimageonlyoneirreduciblecomponent} Let $x \in X$. If $x$ belongs to only one irreducible component of $X$, then $\pi_X(x)$ belongs to only one irreducible component of $X /\!\!/ G$.
\end{lem}

\begin{proof} Let $Y_1, \ldots,Y_l$ be irreducible components of $X$ such that the irreducible algebraic sets $\left(G \cdot Y_{1}\right)/\!\!/ G, \ldots, \left(G \cdot Y_{l}\right)/\!\!/ G$ are the irreducible components of $X \gq G$: see remark \ref{remirredcompgeomquotient}. Now, suppose that $\pi_X(x)$ belongs to two such different components $\left(G \cdot Y_{i}\right)/\!\!/ G$ and $\left(G \cdot Y_{j}\right)/\!\!/ G$ of~$X /\!\!/ G$, with $i,j \in \{1,\ldots,l\}$ and $i \neq j$. Then, by remark \ref{remsafterrealpointssamecomplexorbitsequivsameimages} 2, $x$ belongs to $G \cdot Y_{i}$ and~$G \cdot Y_{j}$ so that there exist $g,h \in G$ such that $x$ belongs to the two different irreducible components $g \cdot Y_i$ and $h \cdot Y_j$ of $X$ (we have $g \cdot Y_i \neq h \cdot Y_j$ since otherwise $G \cdot Y_i$ would be equal to $G \cdot Y_j$ and then $\left(G \cdot Y_{i}\right)/\!\!/ G$ would be equal to $\left(G \cdot Y_{j}\right)/\!\!/ G$ which is not the case).
\end{proof}

\begin{theo} \label{theoimagenonsingularpointtrivstab} Let $x \in X$. Suppose that $x$ is a nonsingular point in dimension $d \in \mathbb{N}$ of $X$ and that $G_x = \{e\}$. Then $\pi_X(x)$ is a nonsingular point in dimension $d - \dim G$ of $X /\!\!/ G$.
\end{theo}  

\begin{proof}
Let $Z$ be the unique irreducible component of $X$ containing $x$: $x$ is a nonsingular point of $Z$ and $\dim Z = d$ (cf. Proposition 3.3.10 of \cite{BCR}). According to previous lemma \ref{lemimageonlyoneirreduciblecomponent} (see also remark \ref{remirredcompgeomquotient}), $\overline{\pi_X(Z)}^{\mathcal{Z}} = \left(G \cdot Z\right)/\!\!/ G$ is the only irreducible component of $X /\!\!/ G$ containing~$\pi_X(x)$. Using again Proposition 3.3.10 of \cite{BCR}, we are therefore reduced to show that~$\pi_X(x)$ is a nonsingular point (in dimension $d-\dim(G)$) of $\left(G \cdot Z\right)/\!\!/ G$: consequently, let us assume in the following that $X = G \cdot Z$ (and then $d = \dim X$) and let us show that the local ring $\mathcal{R}_{X /\!\!/ G, \pi_X(x)}$ of germs of regular functions at $\pi_X(x)$ of the real algebraic set $X /\!\!/ G$ (which is then irreducible) is regular (Proposition 3.3.6 of \cite{BCR}).

For the sake of simplicity, let us denote $W := X /\!\!/ G$ and $\pi := \pi_X$. Recall that $\mathcal{R}_{W,\pi(x)}$ is the localization $\mathcal{P}(W)_{\mathfrak{m}_{\pi(x)}}$ of the ring $\mathcal{P}(W)$ at its maximal ideal $\mathfrak{m}_{\pi(x)}$ of polynomial functions vanishing at $\pi(x)$. If $m_{\pi(x)}$ denotes the unique maximal ideal $\mathfrak{m}_{\pi(x)} \mathcal{R}_{W,\pi(x)}$ of $\mathcal{R}_{W,\pi(x)}$, we have to show that the Krull dimension of $\mathcal{R}_{W,\pi(x)}$ is equal to the dimension of the vector space $m_{\pi(x)}/m_{\pi(x)}^2$ over $\mathcal{R}_{W,\pi(x)}/m_{\pi(x)} \cong \R$, the latter isomorphism being induced by the evaluation at $\pi(x)$ (recall that if $R$ is a local ring with maximal ideal $m$, we always have that the Krull dimension of $R$ is not greater than the dimension of the $R/m$-vector space~$m/m^2$: see for instance Corollary 15.25 of \cite{Sharp}). Notice that the isomorphism of $\R$-algebras $\mathcal{P}(W) \rightarrow \mathcal{P}(X)^G~;~f \mapsto f \circ \pi$ (see the proof of theorem \ref{theorealnullstellensatzequivcategories}) sends $\mathfrak{m}_{\pi(x)}$ to the maximal ideal $\mathfrak{m}_{x}^G$ of~$\mathcal{P}(X)^G$ of invariant polynomial functions on $X$ vanishing at $x$ (and then at any point of the orbit $G \cdot x$), so that the local ring $\mathcal{R}_{W,\pi(x)}$ is isomorphic to $\mathcal{P}(X)^G_{\mathfrak{m}_{x}^G}$. 

We will adapt and develop arguments of the proof of Corollory 9.52 of \cite{Mukai}. In particular, we are going to consider the ideal $I$ of polynomial functions on $X$ vanishing on the $G$-real algebraic subset $Y = G \cdot x$ (notice that $\mathfrak{m}_{x}^G = I^G$ and that $I$ is a $G$-stable linear subspace of $\mathcal{P}(X)$). Since $I$ is a $\mathcal{P}(X)$-module, $I/I^2$ is a $\mathcal{P}(X)/I$-module and hence a $\mathcal{P}(Y)$-module (the restriction to $Y$ induces an isomorphism of $\R$-algebras $\mathcal{P}(X)/I \rightarrow \mathcal{P}(Y)$). But, because $G_x = \{e\}$, the morphism of $\R$-algebras $\alpha_x^* : \mathcal{P}(Y) \rightarrow \mathcal{P}(G)$ is an isomorphism by proposition~\ref{proporbitmapisoiftrivialstabilizer}, so that we can consider the quotient $I/I^2$ as a $\mathcal{P}(G)$-module: for any $\varphi \in \mathcal{P}(G)$ and $f \in I$, if $\varphi \circ \alpha_x^{-1} \in \mathcal{P}(Y)$ is the restriction to $Y$ of a polynomial function $h$ on $X$, we set 
$$\varphi \cdot \overline{f} := \left(\varphi \circ \alpha_x^{-1}\right) \cdot \overline{f} = \overline{h \cdot f} \in I/I^2.$$
The $\mathcal{P}(G)$-module $I/I^2$ is a free module of finite rank and we can find polynomial functions $f_1,\ldots,f_r \in I^G$ such that the family $\left(\overline{f_1},\ldots,\overline{f_r}\right)$ is a $\mathcal{P}(G)$-basis of $I/I^2$: we postpone the demonstration of this crucial fact in the appendix (lemma \ref{lemcomplementprooftheofreeactionnonsingpoint}) in order to make the reading of this proof more fluent.

We thereafter apply Nakayama's Lemma, precisely the fourth statement of lemma \ref{lemnakayamaandconsq}, to the generators $\overline{f_1}, \ldots,\overline{f_r}$ of the $\mathcal{P}(X)$-module $I/I^2$ (we have $\mathcal{P}(G) \cong \mathcal{P}(X)/I$ and the structure of $\mathcal{P}(G)$-module of $I/I^2$ is induced by the structure of $\mathcal{P}(X)$-module of $I$): there exists $h_0 \in 1+I$ such that the polynomial functions $f_1,\ldots,f_r$ generate the $\mathcal{P}(X)_{h_0}$-module~$I_{h_0}$. There then exists $k \in \mathbb{N}$ such that $h_0^k I$ is included in the ideal $\langle f_1,\ldots,f_r \rangle$ of $\mathcal{P}(X)$ generated by $f_1,\ldots,f_r$ (consider a finite family of generators of the ideal $I$ of the Noetherian ring $\mathcal{P}(X)$), and the element $\widetilde{h} := h_0^k$ therefore belongs to the ideal $J := \{f \in \mathcal{P}(X)~|~f I \subset \langle f_1,\ldots,f_r \rangle\}$. Remark that $J$ is a $G$-stable linear subspace of $\mathcal{P}(X)$ (since so is $I$ and $f_1,\ldots,f_r \in I^G$) and that $\widetilde{h}$ belongs to $1+I$ as well. As a consequence, applying the Reynolds operator of $X$ ($G$ is linearly reductive), we obtain an invariant polynomial function $h := \mathcal{R}_X\left(\widetilde{h}\right) \in J^G \cap \left(1+I^G\right)$ (by remark \ref{remcharaclinearredprop} 4): $h$ is in particular an element of $1+I^G$ such that $h I \subset \langle f_1,\ldots,f_r \rangle$. Using the fact that the Reynolds operator of $X$ is a morphism of $\mathcal{P}(X)^G$-modules (this is remark \ref{remcharaclinearredprop} 3), we can furthermore assert that $h I^G$ is included in the ideal of $\mathcal{P}(X)^G$ generated by $f_1,\ldots,f_r$.

On the one hand, we deduce that, if $m_G$ denotes the maximal ideal $\mathfrak{m}_{x}^G \mathcal{P}(X)^G_{\mathfrak{m}_{x}^G} = I^G \mathcal{P}(X)^G_{I^G}$ of $\mathcal{P}(X)^G_{\mathfrak{m}_{x}^G} = \mathcal{P}(X)^G_{I^G}$, $m_G$ is the ideal of $\mathcal{P}(X)^G_{I^G}$ generated by $f_1,\ldots,f_r \in I^G$. In particular, the vector space $m_G/m_G^2$ over $\mathcal{P}(X)^G_{\mathfrak{m}_{x}^G} / m_G \cong \R$ is generated by the classes of the latter elements. Consequently, we have the inequalities
$$r \geq \dim_{\R} m_G/m_G^2 \geq \dim \mathcal{R}_{W,\pi(x)}$$
(recall that the rings $\mathcal{P}(X)^G_{\mathfrak{m}_{x}^G}$ and $\mathcal{R}_{W,\pi(x)}$ are isomorphic). Furthermore, since the real algebraic set $W$ is irreducible, we have $\dim W = \dim \mathcal{P}(W) = \dim \mathcal{P}(W)_{\mathfrak{m}_{\pi(x)}} = \dim \mathcal{R}_{W,\pi(x)}$ (if $R$ is an integral domain and a finitely generated algebra over a field and $P$ is a prime ideal of~$R$, we have $\dim R = \dim R_P$: cf. Corollary 11.12 and Lemma 11.6 (c) of \cite{GathCom}) while, by lemma~\ref{lembounddimsemialgquotient}, $\dim W \geq \dim X - \dim G$, so that
$$r \geq \dim_{\R} m_G/m_G^2 \geq \dim \mathcal{R}_{W,\pi(x)} \geq \dim X - \dim G.$$ 
 
On the other hand, let us consider the irreducible component $G_0$ of $G$ containing its neutral element $e$ (cf. proposition \ref{propcosetsirredcomprealalggroup}) and denote $Y_0 := G_0 \cdot x$: the orbit map $G_0 \rightarrow Y_0~;~g \mapsto g \cdot x = \alpha(g,x)$ is a polynomial isomorphism as well. Moreover, since $Y_0$ is an irreducible real algebraic set (because so is $G_0$), the ideal $I_0$ of polynomial functions on $X$ vanishing on $Y_0$ is a prime ideal of $\mathcal{P}(X)$, and the ideal $I_0 \mathcal{R}_{X,x}$ of the localization $\mathcal{R}_{X,x}$ of $\mathcal{P}(X)$ is a prime ideal of~$\mathcal{R}_{X,x}$. But $\mathcal{R}_{X,x}$ is a regular local ring (of dimension $\dim X$) since $x$ is a nonsingular point (in dimension $\dim X$) of $X$, so that its localization $\left(\mathcal{R}_{X,x}\right)_{I_0 \mathcal{R}_{X,x}}$ at the prime ideal $I_0 \mathcal{R}_{X,x}$ is a regular local ring as well (see for instance Theorem 6.30 of \cite{Bump}). Now, denote by $m_0$ the maximal ideal $I_0 \left(\mathcal{R}_{X,x}\right)_{I_0 \mathcal{R}_{X,x}}$ of the local ring $\left(\mathcal{R}_{X,x}\right)_{I_0 \mathcal{R}_{X,x}}$ and remark that the quotient ring $\left(\mathcal{R}_{X,x}\right)_{I_0 \mathcal{R}_{X,x}}/m_0$ is isomorphic to the field of fractions $\mathcal{K}(Y_0)$ of $\mathcal{P}(Y_0)$ (via restriction of functions to $Y_0$). Let us prove that the classes of the elements $f_1,\ldots,f_r$ of $I \subset I_0$ in $m_0/m_0^2$ are linearly independent over $\left(\mathcal{R}_{X,x}\right)_{I_0 \mathcal{R}_{X,x}}/m_0 \cong\mathcal{K}(Y_0)$. 

For any $i \in \{1,\ldots,r\}$, let $(p_i/q_i)/(s_i/t_i) \in \left(\mathcal{R}_{X,x}\right)_{I_0 \mathcal{R}_{X,x}}$, with $p_i,q_i,s_i,t_i \in \mathcal{P}(X)$ and $q_i(x) \neq 0$, $t_i(x) \neq 0$, $s_i \notin I_0$, such that $\sum_{i=1}^r \left[(p_i/q_i)/(s_i/t_i)\right] f_i \in m_0^2$. There then exists, for any $i \in \{1,\ldots,r\}$, an element $\alpha_i \notin I_0$ such that $\sum_{i=1}^r \alpha_i p_i f_i \in I_0^2$. If $G$ is not irreducible, let $p$ denote a polynomial function on $X$ vanishing on $Y \setminus Y_0$ but not on $Y_0$ (recall that the irreducible components of $G_0$, and then the irreducible components of $Y$, are disjoint) and, if~$G$ is irreducible, set $p$ to be the constant function on $X$ equal to $1$. In both cases, we have $\sum_{i=1}^r (p^2 \alpha_i p_i) f_i \in I^2$. Since the (classes of the) elements $f_1,\ldots,f_r$ form a basis of the free $\mathcal{P}(X)/I$-module $I/I^2$, we deduce that, for any $i \in \{1,\ldots,r\}$, $p^2 \alpha_i p_i \in I \subset I_0$ and then $p_i \in I_0$ ($I_0$ is a prime ideal and $p^2 \alpha_i \notin I_0$). Therefore, for any $i \in \{1,\ldots,r\}$, $(p_i/q_i)/(s_i/t_i) \in m_0$ and the (classes of the) elements $f_1,\ldots,f_r$ of $I \subset I_0$ form a $\mathcal{K}(Y_0)$-linearly independent family of~$m_0/m_0^2$.

As a consequence, we have the inequality
$$\dim \mathcal{R}_{X,x} - \dim \left(\mathcal{R}_{X,x}/I_0 \mathcal{R}_{X,x}\right) \geq  \dim \left(\mathcal{R}_{X,x}\right)_{I_0 \mathcal{R}_{X,x}} = \dim_{\mathcal{K}(Y_0)} m_0/m_0^2 \geq r$$
(recall that the local ring $\left(\mathcal{R}_{X,x}\right)_{I_0 \mathcal{R}_{X,x}}$ is regular): the left-hand side inequality is provided by the general fact that, if $R$ is a ring and $I$ is a prime ideal of $R$, then $\dim R \geq \dim R/P + \dim R_P$ (we refer to Remark 11.5 (b) in \cite{GathCom}). But $\dim \mathcal{R}_{X,x} = \dim X$ since $x$ is a nonsingular point in dimension $\dim X$ of $X$ and the ring $\mathcal{R}_{X,x}/I_0 \mathcal{R}_{X,x}$ is isomorphic to $\mathcal{R}_{Y_0,x}$ (via restriction to $Y_0$) since $Y_0$ is irreducible, so that 
$$\dim \left(\mathcal{R}_{X,x}/I_0 \mathcal{R}_{X,x}\right) = \dim \mathcal{R}_{Y_0,x} = \dim Y_0 = \dim G_0 = \dim G.$$
As a conclusion, 
$$\dim X - \dim G \geq r \geq \dim_{\R} m_G/m_G^2 \geq \dim \mathcal{R}_{W,\pi(x)} \geq \dim X - \dim G$$
and then $\dim_{\R} m_G/m_G^2 = \dim \mathcal{R}_{W,\pi(x)} = \dim X - \dim G$: the local ring $\dim \mathcal{R}_{W,\pi(x)}$ is therefore regular of dimension $\dim X - \dim G$ (again, recall that the rings $\mathcal{P}(X)^G_{\mathfrak{m}_{x}^G}$ and $\mathcal{R}_{W,\pi(x)}$ are isomorphic) i.e. $\pi(x)$ is a nonsingular point in dimension $\dim X - \dim G$ of $W$. Since $W$ is irreducible, we can moreover infer that $\dim W = \dim X - \dim G$.
\end{proof}

Using previous theorem \ref{theoimagenonsingularpointtrivstab} and its proof, we will show that if $X$ is irreducible and $x$ is a nonsingular point of $X$ with trivial stabilizer, then the geometric quotient $\pi_X$ has a local Nash section at $x$. This will allow us, in the next subsection, to lift analytic arcs on $X/G$ to analytic arcs on $X$ if the action of $G$ on $X$ is free, in order to prove the arc-symmetry of the semialgebraic quotient $X/G$ in this case.

Let us first prove the following property given by the semialgebraic version of inverse function theorem. We refer to Definition 2.9.9 of \cite{BCR} for the definition of a Nash submanifold of an affine space and the definition of a Nash map between two such objects. If $A$ is any subset of~$\R^n$ and $y \in A$, an \emph{open semialgebraic neighborhood of $y$ in $A$} is the intersection of $A$ with an open (with respect to Euclidean topology) semialgebraic subset of $\R^n$ containing $y$.

\begin{lem} \label{lemlocalsectionifdiffsurj} Let $M$ be a Nash submanifold of $\R^n$ of dimension $d \in \mathbb{N}$, $N$ be a Nash submanifold of $\R^r$ of dimension $D \in \mathbb{N}$, $f : M \rightarrow N$ be a Nash map and $y \in M$. If the differential ${\rm d}_y f : T_y M \rightarrow T_{f(y)} N$ is surjective, then there exist an open semialgebraic neighborhood $U$ of $y$ in~$M$, an open semialgebraic neighborhood $V$ of $f(y)$ in $N$ and a Nash map $h : V \rightarrow U$ such that $f \circ h = {\rm Id}_V$.
\end{lem}

\begin{proof} By definition, up to restrictions and left and right compositions with Nash diffeomorphisms between open semialgebraic sets of affine spaces sending respectively $y$ to ${\bf 0} \in \R^n$ and~$f(y)$ to ${\bf 0} \in \R^r$, we can suppose that $f$ restricts to a Nash map $F = (F_1,\ldots,F_D)$ from an open semialgebraic neighborhood $U_0$ of the origin in $\R^d$ to~$\R^D$ such that $F(\bf{0}) = \bf{0}$, and the differential ${\rm d}_y f : T_y M \rightarrow T_{f(y)} N$ is then the differential ${\rm d}_{\bf{0}} F : \R^d \rightarrow \R^D$ of $F$ at the origin of~$\R^d$.

Suppose that ${\rm d}_y f$ is surjective i.e. ${\rm d}_{\bf{0}} F = \left({\rm d}_{\bf{0}} F_1,\ldots,{\rm d}_{\bf{0}} F_D\right)$ is surjective. The family $\left({\rm d}_{\bf{0}} F_1,\ldots,{\rm d}_{\bf{0}} F_D\right)$ of $\left(\R^d\right)^*$ is therefore linearly independent: let $\phi_{D+1},\ldots, \phi_d$ be linear forms of~$\left(\R^d\right)^*$ such that the family $\left({\rm d}_{\bf{0}} F_1,\ldots,{\rm d}_{\bf{0}} F_D, \phi_{D+1},\ldots,\phi_d\right)$ is a basis of $\left(\R^d\right)^*$ and consider the Nash map $\phi = \left(F_0,\ldots,F_D,\phi_{D+1},\ldots,\phi_d\right) : U_0 \rightarrow \R^{d}$. By construction, the differential ${\rm d}_{\bf{0}} \phi$ of $\phi$ is invertible so that, by the semialgebraic version of inverse function theorem (see Proposition 2.9.7 of \cite{BCR}), there exist open semialgebraic neighborhoods $U$ and $W$ of the origin in $\R^d$ such that $U \subset U_0$ and $\phi$ restricts to a Nash diffeomorphism from $U$ to $W$.

Denote by $V$ the open semialgebraic subset $\left(\R^D \times \{\bf{0}\}\right) \cap W$ of $\R^D$ and consider the restriction $H : V \rightarrow U$ of the Nash inverse of the latter restriction of $\phi$: we have $F \circ H = {\rm Id}_V$, hence the result.
\end{proof} 

\begin{theo} \label{theolocalnashsectionifnonsingandfree} Suppose that the real algebraic set $X$ is irreducible and let $x$ be a nonsingular point of $X$ such that $G_x = \{e\}$. Then there exist an open semialgebraic neighborhood $X'$ of~$x$ in $X$ which is a Nash submanifold of $\R^n$, an open semialgebraic neighborhood $W'$ of $\pi_X(x)$ in~$X /\!\!/ G$ which is a Nash submanifold of $\R^s$ (where $X\gq G \subset \R^s$) and a Nash map $\psi : W' \rightarrow X'$ such that $\pi_X \circ \psi = {\rm Id}_{W'}$.
\end{theo}

\begin{proof} First, since $X$ is irreducible, so is the real algebraic set $X /\!\!/ G$ (cf. remark \ref{remirredcompgeomquotient}). Moreover, by theorem \ref{theoimagenonsingularpointtrivstab}, because $x$ is a nonsingular point of $X$ such that $G_x = \{e\}$, $\pi_X(x)$ is a nonsingular point of $X /\!\!/ G$, which is of dimension $\dim X - \dim G$.

For simplicity, as in the proof of theorem \ref{theoimagenonsingularpointtrivstab}, denote $W := X /\!\!/ G$ and $\pi := \pi_X$. By Proposition 3.3.11 of \cite{BCR}, because $x$, resp. $\pi(x)$, is a nonsingular point of $X$, resp. $W$, there is an open semialgebraic neighborhood $M$, resp. $N$, of $x$, resp. $\pi(x)$, in $X$, resp. $W$, which is a Nash submanifold of dimension $\dim X$, resp. $\dim X - \dim G$, of $\R^n$, resp. $\R^r$, and we can suppose that $\pi(M) \subset N$. Furthermore, the tangent space $T_x M = T_x X$, resp. $T_{\pi(x)} N = T_{\pi(x)} W$, coincides with the Zariski tangent space of $X$, resp. $W$, at $x$, resp. $\pi(x)$, and the tangent map ${\rm d}_x \pi : T_x X \rightarrow T_{\pi(x)} W$ is dual to the map $\mathfrak{m}_{W,\pi(x)}/\mathfrak{m}_{W,\pi(x)}^2 \rightarrow \mathfrak{m}_{X,x}/\mathfrak{m}_{X,x}^2$ induced by the composition with $\pi$ ($\mathfrak{m}_{X,x}$ denotes the maximal ideal of $\mathcal{P}(X)$ of polynomial functions on $X$ vanishing at $x$). We are going to show the injectivity of the latter map, which will then prove the surjectivity of the tangent map ${\rm d}_x \pi : T_x M \rightarrow T_{\pi(x)} N$ (where $\pi$ denotes here, by abuse of notations, the restriction $M \rightarrow N$ of $\pi$) and allow us to apply previous lemma \ref{lemlocalsectionifdiffsurj} in order to conclude.

As in the proof of theorem \ref{theoimagenonsingularpointtrivstab}, denote by $I$ the ideal of $\mathcal{P}(X)$ of polynomial functions vanishing on $Y := G \cdot x$. Using the identification $I^G/\left(I^G\right)^2 \cong \mathfrak{m}_{W,\pi(x)}/\mathfrak{m}_{W,\pi(x)}^2$, we want to prove that the $\R$-linear map $\iota : I^G/\left(I^G\right)^2 \rightarrow \mathfrak{m}_{X,x}/\mathfrak{m}_{X,x}^2$ induced by the inclusion $\mathcal{P}(X)^G \hookrightarrow \mathcal{P}(X)$ (recall that $\mathfrak{m}_{X,x}^G = I^G$) is injective. 

Consider the sequence of $\R$-vector spaces
$$I^G/\left(I^G\right)^2 \xrightarrow{\iota} \mathfrak{m}_{X,x}/\mathfrak{m}_{X,x}^2 \xrightarrow{\kappa} \mathfrak{m}_{Y,x}/\mathfrak{m}_{Y,x}^2,$$
where $\kappa$ denotes the surjective map $\mathfrak{m}_{X,x}/\mathfrak{m}_{X,x}^2 \rightarrow \mathfrak{m}_{Y,x}/\mathfrak{m}_{Y,x}^2$ induced by the restriction $\mathcal{P}(X) \rightarrow \mathcal{P}(Y)$ to $Y$. Let $f_1,\ldots,f_r \in I^G$ such that the classes $\overline{f_1},\ldots,\overline{f_r} \in I/I^2$ form a basis of the free $\mathcal{P}(X)/I$-module $I/I^2$ (see the proof of theorem \ref{theoimagenonsingularpointtrivstab} and lemma \ref{lemcomplementprooftheofreeactionnonsingpoint}) and let $f \in \mathfrak{m}_{X,x}$ such that $f_{|Y} \in \mathfrak{m}_{Y,x}^2$ i.e. there exist $p \in \mathfrak{m}_{X,x}^2$ and $q \in I$ such that $f = p + q$. Let $q_1,\ldots,q_r \in \mathcal{P}(X)$ and $h \in I^2$ such that $q = h + \sum_{i=1}^r q_i f_i$, we then have 
$$f \equiv \sum_{i=1}^r q_i(x) f_i \mbox{~~mod~~$\mathfrak{m}_{X,x}^2$}$$ 
(we have $I \subset \mathfrak{m}_{X,x}$ and, for any $i \in \{1,\ldots,r\}$, $(q_i-q_i(x))f_i \in \mathfrak{m}^2_{X,x}$), so that ${\rm Ker} \, \kappa = {\rm Im} \, \iota$.

We thereafter compute the dimension of the $\R$-vector space ${\rm Ker} \, \iota$, using the rank theorem and the fact that $\pi(x)$ is a nonsingular point of $W$ and $x$ is a nonsingular point of $X$ as well as a nonsingular point of $Y$: we have
\begin{eqnarray*}
\dim_{\R} {\rm Ker} \, \iota & = & \dim_{\R}  I^G/\left(I^G\right)^2 - \dim_{\R} {\rm Im} \, \iota \\
 					    & = & \dim_{\R} \mathfrak{m}_{W,\pi(x)}/\mathfrak{m}_{W,\pi(x)}^2 - \dim_{\R} {\rm Ker} \, \kappa \\
					    & = & \dim W - \left(\dim_{\R} \mathfrak{m}_{X,x}/\mathfrak{m}_{X,x}^2 - \dim_{\R} \mathfrak{m}_{Y,x}/\mathfrak{m}_{Y,x}^2 \right) \\
					    & = & \dim X - \dim G - \left(\dim X - \dim Y\right) = 0.
\end{eqnarray*}
As a consequence, the map $\iota$ is injective and the tangent map ${\rm d}_x \pi : T_x M \rightarrow T_{\pi(x)} N$ is therefore surjective. Lemma \ref{lemlocalsectionifdiffsurj} then provides a local Nash section for $\pi$ at $x$.
\end{proof}

\begin{rems} \label{remsaftertheonashsectionfreenonsingpoint}
~
\begin{enumerate}
 	\item With the notations of the above proof, the injectivity of the linear map $\iota$ means that any polynomial function of $I^G$ belonging to $\mathfrak{m}_{X,x}^2$ can be written as an element of $\left(I^G\right)^2$.
	\item With the notations of theorem \ref{theolocalnashsectionifnonsingandfree}, notice that the set $W' = \pi_X(\psi(W'))$ is included in $X/G$ and is an open semialgebraic neighborhood of $\pi_X(x) = \varpi_X(x)$ in $X/G$.
\end{enumerate}
\end{rems}

Let us state a first consequence of theorem \ref{theolocalnashsectionifnonsingandfree}, namely the further functoriality of semialgebraic quotient (cf. proposition \ref{propfunctsemialgquotientsemialgcontprop}) with respect to equivariant differentiable maps defined on open $G$-semialgebraic sets on which $G$ acts freely: 

\begin{cor} Let $S \subset \R^n$ be an open (with respect to the Euclidean topology of $\R^n$) $G$-semialgebraic set, let $T \subset \R^N$ be any $G$-semialgebraic set and let $f : S \rightarrow T$ be an equivariant $\mathcal{C}^k$ differentiable map with $k \in \Nstar \cup \{ \infty, \omega\}$ (by a $\mathcal{C}^{\omega}$ differentiable map, we mean an analytic map). Suppose that for all $x \in S$, $G_x = \{e\}$. Then the semialgebraic quotient $S/G \subset \R^s$ is a Nash submanifold of dimension $n-\dim G$ of $\R^s$ and the map $f_{/G} : S/G \rightarrow T/G \subset \R^s$ (see section~\ref{subsectgeosatopquotient}) is a $\mathcal{C}^k$ differentiable map (of $\mathcal{C}^k$ manifolds).
\end{cor}

\begin{proof} Since $S$ is a (non-empty) open semialgebraic subset of $\R^n$, the Zariski closure of $S$ in~$\R^n$ is $\R^n$ (see for instance the proof of Proposition 2.8.4 of \cite{BCR}: we have $\overline{S}^{\mathcal{Z}} = V(I(S)) = V(\{0\}) = \R^n$) so that $\R^n$ is a $G$-real algebraic set and $S$ is a $G$-semialgebraic subset of $\R^n$ (see lemma~\ref{lemzariskiclosgsaisgras}). For simplicity, denote by $\varpi$ the semialgebraic quotient map $\varpi_{\R^n} : \R^n \rightarrow \R^n/G \subset \R^r$ and let $x \in S$: since $G_x = \{e\}$, by above theorem \ref{theolocalnashsectionifnonsingandfree}, there exist an open semialgebraic set $U$ of $\R^n$ containing $x$, an open semialgebraic neighborhood $V$ of $\varpi(x) = \varpi_S(x)$ in $\R^n/G$ which is a Nash submanifold of dimension $n - \dim G$ of $\R^r$ and a Nash map $\psi : V \rightarrow U$ such that~$\varpi \circ \psi = {\rm Id}_V$ (see also remark \ref{remsaftertheonashsectionfreenonsingpoint} 2).

On the one hand, since $S$ is an open semialgebraic subset of $\R^n$ and $\varpi$ is a semialgebraic and open map (lemma \ref{propquotientmapproperclosedopen}), $V \cap \varpi(S) = V \cap S/G$ is an open semialgebraic subset of $\R^n/G$ contained in $V$. The set $V \cap S/G$ is therefore an open semialgebraic subset of the Nash submanifold $V$ of $\R^r$ and is then itself a Nash submanifold of dimension $n - \dim G$ of $\R^r$. But $V \cap S/G$ is also an open semialgebraic subset of $S/G$ containing $\varpi(x)$. As a consequence, since for any point $\varpi(y)$, $y \in S$, of $S/G$, there exists an open semialgebraic neighborhood of~$\varpi(y)$ in~$S/G$ which a Nash submanifold of dimension $n - \dim G$, the set $S/G$ is itself a Nash submanifold of dimension~$n - \dim G$ of $\R^r$.

On the other hand, on the open semialgebraic neighborhood $V \cap S/G$ of $\varpi(x)$, $f/G$ writes as $\varpi_T \circ f \circ \psi$ and $f/G$ is therefore a $\mathcal{C}^k$ map on $V \cap S/G$ ($\omega_T$ is the restriction of a polynomial map from $\R^N$ to $\R^s$ and recall also, in the case $k = \omega$, that a Nash map between Nash submanifolds of affine spaces is analytic: cf. Proposition 8.1.8 of \cite{BCR}). Consequently, the map $f_{/G} : S/G \rightarrow T/G$ is a $\mathcal{C}^k$ map.
\end{proof}

\begin{rem} 
Combining the previous result with proposition \ref{propfunctsemialgquotientsemialgcontprop}, the semialgebraic quotient is functorial with respect to equivariant Nash maps defined on open $G$-semialgebraic sets on which $G$ acts freely.
\end{rem}

Adapting the above proof also allows to show the functoriality of semialgebraic quotient with respect to equivariant arc-analytic maps defined on $G$-semialgebraic sets on which $G$ acts freely: a map $f : A \rightarrow B$ between subsets of affines spaces is said \emph{arc-analytic} if for any analytic arc $\gamma : ]-1;1[ \rightarrow A$, the composition $f \circ \gamma : ]-1;1[ \rightarrow B$ is analytic (see definition \ref{defarcanalyticmap}).

\begin{prop} \label{propfunctarcanalyticmapfromfreegsaset} Let $S \subset \R^n$ and $T \subset \R^N$ be $G$-semialgebraic sets and let $f : S \rightarrow T$ be an equivariant arc-analytic map. Suppose that for all $x \in S$, $G_x = \{e\}$. Then $f_{/G} : S/G \rightarrow T/G \subset \R^s$ is an arc-analytic map.
\end{prop}

\begin{proof} Thanks to proposition \ref{proplinearizgrealalgset} (see also proposition \ref{proporthogonalizationactioncompactrealalggroup}), we can suppose that $\R^n$ is a polynomial representation of $G$ and that $S$ is a $G$-stable semialgebraic subset of $\R^n$. Consider the polynomial map $\pi : \R^n \rightarrow \R^r$ induced by $\pi_{\R^n} : \R^n \rightarrow \R^n \gq G \subset \R^r$ and let $\gamma : ]-1;1[ \rightarrow S/G$ be an analytic arc. Let $t_0 \in ]-1;1[$ and let $x \in S$ such that $\pi(x) = \gamma(t_0)$. Since $G_x = \{e\}$, we can apply theorem \ref{theolocalnashsectionifnonsingandfree}: there exist an open semialgebraic neighborhood $U$ of $x$ in $\R^n$, an open semialgebraic neighborhood $V$ of $\pi(x) = \gamma(t_0)$ in $\R^n /\!\!/G$ which is a Nash submanifold of $\R^r$ and a Nash map $\psi : V \rightarrow U$ such that $\pi \circ \psi = {\rm Id}_{V}$. 

Since $\gamma$ has values in $S/G \subset \R^n /\!\!/G$, $\gamma^{-1}(S/G \cap V)$ is then an open subset of $]-1;1[$ containing~$t_0$: there exists $\epsilon >0$ such that $\gamma(]t_0-\epsilon;t_0+\epsilon[) \subset S/G \cap V$ and therefore, using proposition \ref{propseparationorbits}, $\psi \circ \gamma (]t_0-\epsilon;t_0+\epsilon[) \subset S \cap U$ (we have $\pi \circ \psi = {\rm Id}_{V}$ and $S$ is stable under the action of $G$). We then write, for any $t \in ]t_0-\epsilon;t_0+\epsilon[$,
$$f_{/G} \circ \gamma(t) = \varpi_T \circ f \circ \psi \circ \gamma(t),$$
which shows that $f_{/G} \circ \gamma$ is analytic at $t_0$, because $\psi$ is analytic, $f$ is arc-analytic and $\varpi_T$ is the restriction of a polynomial map. As a consequence, the map $f_{/G} : S/G \rightarrow T/G$ is arc-analytic.
\end{proof}

\subsection{Arc-symmetric semialgebraic quotient of a real algebraic set by a free polynomial action of a compact real algebraic group}

In this part, we will achieve one goal of this paper: we prove that the semialgebraic quotient of a $G$-real algebraic set on which the compact real algebraic group $G$ acts freely is an arc-symmetric set of the ambient affine space.

\begin{de} Let $A$ be a subset of $\R^n$. We say that $A$ is an \emph{arc-symmetric set of $\R^n$} if for any analytic arc $\gamma : ]-1;1[ \rightarrow \R^n$, if $\gamma(]-1;0[) \subset A$, then $\gamma(]-1;1[) \subset A$ (see also definition \ref{defarcsymsetrn}).
\end{de}

Notice that the zero set of any analytic function form $\R^n$ to $\R$ is an arc-symmetric set of~$\R^n$ by the identity theorem for analytic functions. In particular, a real algebraic set of $\R^n$ is an arc-symmetric set of $\R^n$ (see also example \ref{exrealalgsetisarcsym}).

This important notion of arc-symmetric set was introduced by K. Kurdyka in \cite{KurAR}. In section \ref{subsectarcsymassets} of the appendix, for the sake of self-completeness, we chose to give an introduction to arc-symmetric sets of affine spaces and real analytic manifolds, deeply inspired by \cite{KurAR} and~\cite{KurPar}, including the detailed proofs of properties of semialgebraic arc-symmetric sets of affine spaces and real projective spaces stated in these two articles, in particular the fact that the semialgebraic arc-symmetric sets of $\R^n$ form the closed sets of a Noetherian topology on~$\R^n$ (theorem \ref{theotopar}), called the \emph{$\mathcal{AR}$ topology} of $\R^n$.

Below, we are actually going to show that the semialgebraic quotient of any arc-symmetric $G$-semialgebraic set on which $G$ acts freely is an arc-symmetric set of the ambient affine space:

\begin{theo} \label{theoquotientfreearcsymrnisarcsym} Let $S$ be a $G$-stable semialgebraic subset of $X$. Suppose that $S$ is arc-symmetric in $\R^n$ and that for any $x \in S$, $G_x = \{e\}$. Then $S/G \subset \R^r$ is a semialgebraic arc-symmetric subset of $\R^r$.
\end{theo}

\begin{proof} Let $f : X \rightarrow Y$ be an isomorphism of $G$-real algebraic sets onto a $G$-real algebraic subset of a polynomial orthogonal representation $\R^N$ of $G$ (proposition \ref{proporthogonalizationactioncompactrealalggroup}). It induces the polynomial isomorphism $f_{\gq G} : X\gq G \subset \R^r \rightarrow Y \gq G \subset \R^s$ which satisfies $f_{\gq G} \circ \pi_X = \pi_Y \circ f$ (lemma~\ref{lemfunctorrealalgquotient}). In particular, we have $f_{\gq G}\left(S/G\right) = f_{\gq G}\left(\pi_X(S)\right) = \pi_Y(f(S)) = f(S)/G$, so that~$S/G$ is arc-symmetric in $\R^r$ if and only if $f(S)/G$ is arc-symmetric in $\R^s$ by proposition~\ref{proparcanalyticmaps}~2: the polynomial maps $f_{\gq G}$ and $f_{\gq G}^{-1}$ are in particular arc-analytic. As a consequence, we can suppose that $\R^n$ is a polynomial representation of $G$ and that $X$ is a $G$-real algebraic set of $\R^n$. In particular, $S/G = \pi_{\R^n}(S)$.

Consider the polynomial map $\pi : \R^n \rightarrow \R^r$ induced by $\pi_{\R^n} : \R^n \rightarrow \R^n \gq G \subset \R^r$. The semialgebraic set $\pi(S)$ is (Euclidean) closed in $\R^r$ since the map $\pi$ is closed (see the proof of lemma \ref{propquotientmapproperclosedopen}) and since any semialgebraic arc-symmetric set of an affine space is closed (proposition \ref{propsemialgarcsymisclosed}). Now, let $\gamma : ]-1;1[ \rightarrow \R^r$ be an analytic arc such that $\gamma(]-1;0[) \subset \pi(S)$: since $\pi(S)$ is closed and $\gamma$ is continuous, we have $\gamma(]-1;0]) \subset \pi(S)$. Denote $w := \gamma(0)$ and let~$x \in S$ such that $\pi(x) = w$. Since $G_x = \{e\}$, we can apply theorem \ref{theolocalnashsectionifnonsingandfree}: there exist an open semialgebraic neighborhood $U$ of $x$ in $\R^n$, an open semialgebraic neighborhood $V$ of $w= \pi(x)$ in $\R^n /\!\!/G$ which is a Nash submanifold of $\R^r$ and a Nash map $\psi : V \rightarrow U$ such that $\pi \circ \psi = {\rm Id}_{V}$. 

Remark also that, because~$\R^n\gq G$ is a real algebraic set and then an arc-symmetric subset of $\R^r$, $\gamma(]-1;1[) \subset \R^n\gq G$: since $\gamma$ has values in $\R^n\gq G$, $\gamma^{-1}(V)$ is an open subset of $]-1;1[$ containing $0$. There then exists $\epsilon >0$ such that $\gamma(]-\epsilon;\epsilon[) \subset V$. We can therefore consider the analytic arc $\widetilde{\gamma} := \psi \circ \gamma : ]-\epsilon;\epsilon[ \rightarrow U \subset \R^n$. Because $\pi\big(\widetilde{\gamma}(]-\epsilon;0])\big) = \gamma(]-\epsilon;0]) \subset \pi(S)$, we have, by proposition \ref{propseparationorbits} and because $S$ is $G$-stable, $\widetilde{\gamma}(]-\epsilon;0]) \subset S$.

But $S$ is arc-symmetric in $\R^n$, so that $\widetilde{\gamma}(]-\epsilon;\epsilon[) \subset S$. As a consequence, $\gamma(]-\epsilon;\epsilon[) = \pi\big(\widetilde{\gamma}(]-\epsilon;\epsilon[)\big) \subset \pi(S)$ and $\pi(S)$ is therefore an arc-symmetric semialgebraic subset of $\R^r$.

\end{proof}

In particular, if the action of $G$ on $X$ is free, then the semialgebraic quotient $X/G$ is arc-symmetric in $\R^r$.

\begin{rem} In general, the arc-symmetric semialgebraic quotient of a $G$-real algebraic set by a free action of $G$ is not a real algebraic set. Consider for instance the polynomial action of~$S^1$ on $\R^3$ given in example \ref{exrealalgquotientplaneorthoaction} 2, as well as the compact $\SO_2(\R)$-stable real algebraic set
$$Z := \left\{(x,y,z) \in \R^3~|~z^2 + (x^2+y^2-2)(x^2+y^2 - 1)(x^2+y^2+1) = 0\right\}.$$
The action of $S^1$ on $Z$ is free and 
$$Z/S^1 = \pi_{\R^3}(Z) = \left\{(u,z) \in \R^2~|~z^2 + (u-2)(u-1)(u+1) = 0, u\geq 0\right\}$$
is an arc-symmetric semialgebraic non-algebraic set (it is the compact connected component of the real algebraic set given by the equation $z^2 + (u-2)(u-1)(u+1) = 0$ in $\R^2$).
\end{rem}

More generally, let us consider the $\mathcal{AS}$-sets of $\R^n$:

\begin{de} Let $A$ be a subset of $\R^n$. We say that $A$ is an \emph{$\mathcal{AS}$-set of $\R^n$} if $A$ is semialgebraic and for any analytic arc $\gamma : ]-1;1[ \rightarrow \R^n$, if $\gamma(]-1;0[) \subset A$, then there exists $\epsilon \in ]0;1[$ such that $\gamma(]-1;1[) \subset A$ (see subsection \ref{subsubsecqarcsymassets}).
\end{de}

The $\mathcal{AS}$-sets of $\R^n$ are exactly the finite Boolean combinations of semialgebraic arc-symmetric subsets of $\R^n$ (proposition \ref{propquasiarcsymrniffasrn}).

The semialgebraic quotient of a $G$-semialgebraic $\mathcal{AS}$-set of $\R^n$ is an $\mathcal{AS}$-set of $\R^r$ as soon as the action of $G$ is free \emph{on the $\mathcal{AR}$-closure} of the latter: 

\begin{theo} \label{theoquotientfreearrnisar} Let $S$ be a $G$-stable semialgebraic subset of $X$ which is an $\mathcal{AS}$-set of $\R^n$. Then $\overline{S}^{\mathcal{AR}}$ is a $G$-stable semialgebraic subset of $X$ and if for any $x \in \overline{S}^{\mathcal{AR}}$, $G_x = \{e\}$, then $S/G \subset \R^r$ is an $\mathcal{AS}$-set of $\R^r$.
\end{theo}

\begin{proof} First notice that since $X$, as a real algebraic set, is an $\mathcal{AR}$-closed set (i.e. a semialgebraic arc-symmetric set) of $\R^n$ containing $S$, we have $\overline{S}^{\mathcal{AR}} \subset X$. Now, let $g \in G$, the set $g^{-1} \cdot \overline{S}^{\mathcal{AR}} := \alpha\left(g^{-1},\overline{S}^{\mathcal{AR}}\right) = \alpha(g,\cdot)^{-1}\left(\overline{S}^{\mathcal{AR}}\right)$ is also an $\mathcal{AR}$-closed set of $\R^n$, as the inverse image by a semialgebraic arc-analytic map of an $\mathcal{AR}$-closed set: cf. proposition \ref{proparcanalyticmapsandarclosedsets} 2. Because $g^{-1} \cdot \overline{S}^{\mathcal{AR}}$ contains $S$ (the set $S$ is stable under the action of $G$), we then have $\overline{S}^{\mathcal{AR}} \subset g^{-1} \cdot \overline{S}^{\mathcal{AR}}$ i.e. $g \cdot \overline{S}^{\mathcal{AR}} \subset \overline{S}^{\mathcal{AR}}$. We also have $\overline{S}^{\mathcal{AR}} \subset \left(g^{-1}\right)^{-1} \cdot \overline{S}^{\mathcal{AR}} = g \cdot \overline{S}^{\mathcal{AR}}$, so that $g \cdot \overline{S}^{\mathcal{AR}} = \overline{S}^{\mathcal{AR}}$ and~$\overline{S}^{\mathcal{AR}}$ is therefore a $G$-stable semialgebraic subset of $X$ as well.

In order to prove the second part of the statement, suppose that the action of $G$ on $\overline{S}^{\mathcal{AR}}$ is free: by theorem \ref{theoquotientfreearcsymrnisarcsym}, $\overline{S}^{\mathcal{AR}}/G = \pi_X\left(\overline{S}^{\mathcal{AR}}\right) \subset \R^r$ is a semialgebraic arc-symmetric set of $\R^r$. We then write $S = \overline{S}^{\mathcal{AR}} \setminus \left(\overline{S}^{\mathcal{AR}} \setminus S\right)$ and, by corollary \ref{corstabilityofimageintersectiondifference},
$$\pi_X(S) = \pi_X\left(\overline{S}^{\mathcal{AR}}\right) \setminus \pi_X\left(\overline{S}^{\mathcal{AR}} \setminus S\right).$$
Because $\overline{S}^{\mathcal{AR}} \setminus S$ is a $G$-stable $\mathcal{AS}$-set of $\R^n$ whose $\mathcal{AR}$-closure is contained in $\overline{S}^{\mathcal{AR}}$ and whose dimension is lower than $\dim S$ (see for instance the proof of proposition \ref{propquasiarcsymrniffasrn}), we can assert that $\pi\left(\overline{S}^{\mathcal{AR}} \setminus S\right)$ is an $\mathcal{AS}$-set of $\R^r$ by induction on dimension (if $\dim S =0$, we have $\dim \pi_X(S) \leq \dim S = 0$ by Theorem 2.8.8 of \cite{BCR} and $\pi_X(S)$ is therefore a finite number of points). As a consequence, $\pi_X(S) = S/G$ is an $\mathcal{AS}$-set of $\R^r$ as well.
\end{proof}

\begin{rems} \label{remsquotientfreeassetsaffinespace}
~
	\begin{enumerate}
		\item Assuming that the action is free on the sole considered $G$-semialgebraic $\mathcal{AS}$-set is not enough to guarantee the semialgebraic quotient to be an $\mathcal{AS}$-set. Consider for instance the usual orthogonal action of ${\rm O}_2(\R)$ on $\R^2$: the Zariski open set $U := \R^2 \setminus \{(0,0)\}$ of~$\R^2$ is a ${\rm O}_2(\R)$-stable $\mathcal{AS}$-set of $\R^2$ on which ${\rm O}_2(\R)$ acts freely. However, the image of $U$ under the geometric quotient $(x,y) \in \R^2 \mapsto x^2+y^2 \in \R$ (example~\ref{exrealalgquotientplaneorthoaction}~1) is the open interval $]0;+\infty[$ of $\R$ which is not an $\mathcal{AS}$-set of $\R$ (the analytic arc $\gamma : t \in ]-1;1[ \mapsto -t \in \R$ satisfies $\gamma(]-1;0[) \subset ]0;+\infty[$ but for all $\epsilon \in ]0;1[$, $\gamma(\epsilon) \notin ]0;+\infty[$).
		\item With the notations of theorem \ref{theoquotientfreearrnisar}, if the action of $G$ on $\overline{S}^{\mathcal{AR}}$ is free, we have 
		$$\overline{\left(S/G\right)}^{\mathcal{AR}} = \overline{\pi_X(S)}^{\mathcal{AR}} = \pi_X\left(\overline{S}^{\mathcal{AR}}\right) = \overline{S}^{\mathcal{AR}}/G.$$
Indeed, $\overline{S}^{\mathcal{AR}}/G = \pi_X\left(\overline{S}^{\mathcal{AR}}\right)$ is an $\mathcal{AR}$-closed set of $\R^r$ containing $S/G$, hence the direct inclusion. Moreover, $T := \overline{\left(S/G\right)}^{\mathcal{AR}}$ is contained in $X \gq G$ (since the latter is a real algebraic set of $\R^r$ containing $S/G$) and then $\pi_X^{-1}(T)$ is a semialgebraic arc-symmetric subset of $\R^n$ included in $X$ by proposition \ref{proparcanalyticmaps} (the polynomial map $\pi_X$ is in particular arc-analytic and semialgebraic). Furthermore, because $T$ contains $S/G$, $\pi_X^{-1}(T)$ contains~$S$ and therefore $\overline{S}^{\mathcal{AR}} \subset \pi_X^{-1}(T)$, so that $\pi_X\left(\overline{S}^{\mathcal{AR}}\right) \subset \pi_X\left(\pi_X^{-1}(T)\right) \subset T$.
	\end{enumerate}
\end{rems} 

\section{$\mathcal{AS}$ quotients of $G$-$\mathcal{AS}$-sets} \label{sectionasquotientsofgassets}

Let $G$ be a compact real algebraic group. In this last section, we consider arc-symmetric sets and $\mathcal{AS}$-sets of compact affine Nash manifolds, as well as actions of~$G$ on these objects which allow to define quotients with arc-symmetric/$\mathcal{AS}$ structures. From these quotients, we will construct, in a subsequent paper, further invariants with respect to equivariant continuous maps with $\mathcal{AS}$-graphs.

Arc-symmetric sets can actually be defined for any real analytic manifold (definition \ref{defarcsymsetsrealanman}):

\begin{de} Let $M$ be a real analytic manifold and let $A$ be a subset of $M$. We say that $A$ is an \emph{arc-symmetric set of $M$} if for any analytic arc $\gamma : ]-1;1[ \rightarrow M$, if $\gamma(]-1;0[) \subset A$, then $\gamma(]-1;1[) \subset A$.
\end{de}

However, we focus on Nash manifolds: on a Nash manifold $M$, we can consider semialgebraic sets and define a Noetherian topology whose closed sets are the semialgebraic arc-symmetric subsets of $M$, called the $\mathcal{AS}$ topology of $M$. We further restrict our attention on compact affine Nash manifolds, that is compact Nash manifolds that can be embedded, via a Nash diffeomorphism, as a Nash submanifold of an affine space (notice that the real projective spaces are compact affine Nash manifolds): on those objects, the $\mathcal{AS}$-sets are exactly the finite Boolean combinations of semialgebraic arc-symmetric subsets (we refer to the subsections \ref{subsubsecsaofnashmfds}, \ref{subsubsecastop} and \ref{subsubsecqarcsymassets} of the appendix for the details and proofs about these notions and facts).

\begin{de} Let $M$ be a compact affine Nash manifold and let $A$ be a subset of $M$. We say that $A$ is an \emph{$\mathcal{AS}$-set of $M$} if $A$ is a  semialgebraic subset of $M$ and for any analytic arc $\gamma : ]-1;1[ \rightarrow M$, if $\gamma(]-1;0[) \subset A$, then there exists $\epsilon \in ]0;1[$ such that $\gamma(]0;\epsilon[) \subset A$.
\end{de} 

\subsection{Compact Nash submanifolds of real affine spaces} \label{subsectcompnashsubraff}

We begin by dealing with the case of compact Nash submanifolds of real affine spaces. Precisely, let $\R^n$ be a polynomial representation of $G$ with geometric quotient $\pi : \R^n \rightarrow \R^n \gq G \subset \R^r$ and let $M$ be a compact Nash submanifold of $\R^n$.

\begin{prop} Let $S$ be a semialgebraic arc-symmetric subset of $M$. Suppose that $S$ is globally stable under the action of $G$ on $\R^n$ and that for any $x \in S$, $G_x = \{e\}$. Then $S/G = \pi(S) \subset \R^r$ is a semialgebraic arc-symmetric subset of $\mathbb{P}^r(\R)$.
\end{prop}

\begin{proof} By proposition \ref{propcompactnashsubastopcoincide}, $S$ is an $\mathcal{AR}$-closed set of $\R^n$ included in $M$. We then apply theorem \ref{theoquotientfreearcsymrnisarcsym} above to assert that $S/G = \pi(S)$ is an $\mathcal{AR}$-closed set of $\R^r$. Because $S$ is closed in the compact space $M$ (proposition \ref{propsemialgsarcsymnashmanclosed}), $S$ is compact in $\R^n$. Therefore, since~$\pi$ is continuous, $S/G = \pi(S)$ is compact as well, and consequently $S/G$ is a semialgebraic arc-symmetric subset of $\mathbb{P}^r(\R)$ by proposition \ref{propcompactsubrnequivsaarcsympnsaarcsymrn}.
\end{proof}

More generally, we have the following result on $G$-stable $\mathcal{AS}$-sets of $M$:

\begin{prop} \label{propquotientfreeassetsnashsub} Let $S$ be an $\mathcal{AS}$-set of $M$ and suppose that $S$ is globally stable under the action of $G$ on $\R^n$. Then the $\mathcal{AS}$-closure of $S$ in $M$ is $G$-stable as well. Furthermore, if for any $x \in \overline{S}^{\mathcal{AS}}$, $G_x = \{e\}$, then $S/G = \pi(S) \subset \R^r$ is an $\mathcal{AS}$-set of $\mathbb{P}^r(\R)$.
\end{prop}

\begin{proof} First, by lemma \ref{remassetcompactsubnashsaaamap}, $S$ is an $\mathcal{AS}$-set of~$\R^n$ included in $M$. We then apply theorem~\ref{theoquotientfreearrnisar}: the $\mathcal{AR}$-closure $\overline{S}^{\mathcal{AR}}$ of $S$ in $\R^n$ is $G$-stable as well. Furthermore, $\overline{S}^{\mathcal{AR}}$ is also the $\mathcal{AS}$-closure $\overline{S}^{\mathcal{AS}}$ of $S$ in $M$ by proposition \ref{propcompactnashsubastopcoincide} ($M$ is $\mathcal{AR}$-closed in $\R^n$). As a consequence, if the action of $G$ on $\overline{S}^{\mathcal{AS}}$ is free, then $S/G = \pi(S)$ is an $\mathcal{AR}$-set of $\R^r$ by theorem \ref{theoquotientfreearrnisar}.

Moreover, by remark \ref{remsquotientfreeassetsaffinespace} 2, the $\mathcal{AR}$-closure of $S/G$ in $\R^r$ is $\pi\left(\overline{S}^{\mathcal{AR}}\right) = \pi\left(\overline{S}^{\mathcal{AS}}\right)$, and the latter set is compact (as the image of the compact subset $\overline{S}^{\mathcal{AS}}$ of $\R^n$ by the continuous map~$\pi$). Consequently, by corollary \ref{corascompactarcloasinproj}, $S/G$ is an $\mathcal{AS}$-set of $\mathbb{P}^r(\R)$.
\end{proof}

In this framework, we can prove the following important property of functoriality of the semialgebraic quotient with respect to $\mathcal{AS}$-maps between $\mathcal{AS}$-sets of Nash submanifolds of affine spaces which satisfy the hypotheses of previous proposition \ref{propquotientfreeassetsnashsub}. A map $f : X \subset M \rightarrow Y \subset N$ between $\mathcal{AS}$-sets of compact affine Nash manifolds is called an \emph{$\mathcal{AS}$-map} if its graph is an $\mathcal{AS}$-set of the compact affine Nash manifold $M \times N$ (definition \ref{defasmap}).

\begin{theo} \label{propsemialgquotientmapbetweenfreeassetscompnashsubisas} Let $N$ be a compact Nash submanifold of a polynomial representation $\R^m$ of~$G$ with geometric quotient $\pi' : \R^m \rightarrow \R^m /\!\!/G \subset \R^s$. Let $X$ be a $G$-stable $\mathcal{AS}$-set of $M$ and~$Y$ be a $G$-stable $\mathcal{AS}$-set of $N$ such that $G$ acts freely on the respective $\mathcal{AS}$-closures of $X$ and $Y$. Let $f : X \subset M \rightarrow Y \subset N$ be an equivariant $\mathcal{AS}$-map. Then the map
$$f_{/G} : X/G \subset \mathbb{P}^r(\R) \rightarrow Y/G \subset \mathbb{P}^s(\R)$$ 
is an $\mathcal{AS}$-map.
\end{theo} 

\begin{proof} According to corollary \ref{corsufficientconditionasmapbetweenrealassets} and example \ref{exaftercorsufficientconditionasmapbetweenrealassets}, the map $f_{/G}$ is an $\mathcal{AS}$-map if and only if its graph 
$$\Gamma_{f_{/G}} = \left\{(w,z) \in X/G \times Y/G~|~z = f_{/G}(w) \right\}$$
is an $\mathcal{AS}$-set of $\R^r \times \R^s$ (the respective $\mathcal{AR}$-closures $\overline{X/G}^{\mathcal{AR}} = \pi\left(\overline{X}^{\mathcal{AS}}\right) \subset \R^r$ and $\overline{Y/G}^{\mathcal{AR}} = \pi'\left(\overline{Y}^{\mathcal{AS}}\right) \subset \R^s$ of $X/G$ and $Y/G$ are compact: see the proof of proposition \ref{propquotientfreeassetsnashsub}).

Denote $\Gamma := \left\{(x,y) \in X \times Y~|~\exists g \in G,~y = g \cdot f(x)\right\}$ and remark that
$$\Gamma_{f_{/G}} = \left\{(w,z) \in X/G \times Y/G ~|~\exists x \in X,~\exists y \in Y,~\pi(x) = w,~\pi'(y) = z,~y = f(x)\right\}$$
is the image of $\Gamma$ by the polynomial map $(x,y) \in \R^n \times \R^m \mapsto (\pi(x),\pi'(y)) \in \R^n /\!\!/G \times \R^m /\!\!/G \subset \R^r \times \R^s$, i.e by the geometric quotient associated to the induced polynomial representation $\R^n \times \R^m$ of $G \times G$ (see remark \ref{remgeomquotientproduct}). We show that the $(G \times G)$-stable set $\Gamma$ (the map $f$ is equivariant) satisfies the hypotheses of proposition \ref{propquotientfreeassetsnashsub}. 

Towards this goal, notice that $\Gamma$ is the image of the graph $\Gamma_{\widetilde{f}}$ of the map $\widetilde{f} : (g,x) \in G \times X \mapsto g \cdot f(x) \in Y$ by the projection map $\Pi : G \times X \times Y \rightarrow X \times Y~;~(g,x,y) \mapsto (x,y)$. On the one hand, the map $\widetilde{f}$ is an $\mathcal{AS}$-map by corollary \ref{corcompasmapsisasmap}, as the composition of the $\mathcal{AS}$-map~$f$ with a polynomial action map. In other words, considering $G$ as a compact nonsingular real algebraic subset of an affine space $\R^d$, the graph $\Gamma_{\widetilde{f}}$ is an $\mathcal{AS}$-set of the compact Nash submanifold $G \times M \times N$ of $\R^d \times \R^r \times \R^s$. 

On the other hand, consider the restriction $\widetilde{\Pi} : \Gamma_{\widetilde{f}} \rightarrow X \times Y$ of $\Pi$: $\widetilde{\Pi}$ is an $\mathcal{AS}$-map as a restriction of the projection $G \times M \times N \rightarrow M \times N~;~(g,x,y) \mapsto (x,y)$ to $\mathcal{AS}$-sets (example \ref{exsasmaps} 2). The map $\widetilde{\Pi}$ is furthermore injective since the action of $G$ on $Y$ is free (if $(g,x),(g',x') \in G \times X$ satisfy $x =x'$ and $g \cdot f(x) = g' \cdot f(x') = g' \cdot f(x)$, then $g = g'$ since the stabilizer of $f(x)$ is trivial). As a consequence, $\Gamma = \widetilde{\Pi}\left(\Gamma_{\widetilde{f}}\right)$ is an $\mathcal{AS}$-set of $M \times N$ by theorem~\ref{thimageofassetbyinjasmapisas}.

Finally, let us check that $G \times G$ acts freely on the $\mathcal{AS}$-closure of $\Gamma$ in $M \times N$: the latter set is included in $\overline{X}^{\mathcal{AS}} \times \overline{Y}^{\mathcal{AS}}$ and, for all $g,g' \in G$ and $(x,y) \in \overline{X}^{\mathcal{AS}} \times \overline{Y}^{\mathcal{AS}}$, if $(g,g') \cdot (x,y) = (x,y)$, we have $g \cdot x = x$ and $g' \cdot y = y$, so that $g = e$ and $g' = e$, because $G$ acts freely on $\overline{X}^{\mathcal{AS}}$ and~$\overline{Y}^{\mathcal{AS}}$.

As a consequence, by (the proof of) previous proposition \ref{propquotientfreeassetsnashsub}, $\Gamma_{f_{/G}} = \Gamma/(G \times G)$ is an $\mathcal{AS}$-set of $\R^r \times \R^s$, hence the result.
\end{proof}

\subsection{$G$-linearizations and $G$-$\mathcal{AS}$-sets} \label{subsectglingassets}

In this part, using the affine case of previous subsection, we precisely emphasize the objects to which we will be able to associate some $\mathcal{AS}$-quotients well-defined up to $\mathcal{AS}$-homeomorphism. 

\begin{de} \label{degassetfreegasset} Let $M$ be a compact affine Nash manifold and $S$ be an $\mathcal{AS}$-set of $M$ equipped with an action of $G$. A \emph{$G$-linearization of $S$} is an equivariant homeomorphism with $\mathcal{AS}$-graph $\Phi : S \rightarrow S'$ onto a $G$-stable $\mathcal{AS}$-set $S'$ of a compact Nash submanifold $N$ of a polynomial representation $\R^n$ of $G$. Such a $G$-linearization of $S$ is said \emph{free} if the action of $G$ on the $\mathcal{AS}$-closure of $S'$ in $N$ (which is also the $\mathcal{AR}$-closure of $S'$ in $\R^n$) is free. If $S$ has a $G$-linearization, we say that $S$ is a \emph{$G$-$\mathcal{AS}$-set}. If $S$ has a free $G$-linearization, we say that $S$ is a \emph{free $G$-$\mathcal{AS}$-set}.
\end{de}

The above definition of a free $G$-linearization is justified by proposition \ref{propquotientfreeassetsnashsub} and we will be able to consider a well-defined $\mathcal{AS}$-quotient of free $G$-$\mathcal{AS}$-sets thanks to the latter result and theorem \ref{propsemialgquotientmapbetweenfreeassetscompnashsubisas}: this will be done in the next subsection. First, let us give some examples of $G$-$\mathcal{AS}$-sets and free $G$-$\mathcal{AS}$-sets. Notice also that, by definition, if $X$ is a $G$-$\mathcal{AS}$-set, the action $G \times X \rightarrow X~;~(g,x) \mapsto g \cdot x$ is a continuous $\mathcal{AS}$-map and is given by \emph{$\mathcal{AS}$-homemorphisms}, i.e. homeomorphisms with $\mathcal{AS}$-graph.

\begin{prop} \label{propgralgsetisgasset} Any $G$-real algebraic set is a $G$-$\mathcal{AS}$-set.
\end{prop}

\begin{proof} Let $X$ be a $G$-real algebraic set of an affine space $\R^n$. Then $X$ is an $\mathcal{AS}$-subset of the compact affine Nash manifold $\mathbb{P}^n(\R)$ ($X$ is the set-theoretic difference of the projective algebraic set obtained by homogenizing generators of the ideal $I(X)$ and the hyperplane at infinity of equation $x_0 = 0$: $X$ is a Zariski open subset of the real projective space $\mathbb{P}^n(\R)$). Moreover, according to proposition \ref{proporthogonalizationactioncompactrealalggroup}, there is an equivariant polynomial isomorphism $f$ from $X$ to a $G$-real algebraic subset $Y$ of a polynomial orthogonal representation $\R^N$ of $G$. The set $Y$ is an $\mathcal{AS}$-subset of $\mathbb{P}^N(\R)$ and the homeomorphism $f$ has $\mathcal{AS}$-graph in $\mathbb{P}^r(\R) \times \mathbb{P}^s(\R)$ (see remark~\ref{rempolynomialmaphasasgraphinproj}).

We can furthermore extend the orthogonal action of $G$ on $\R^N$ to an action on $\mathbb{P}^N(\R)$. Indeed, if $\alpha : G \times \R^N \rightarrow \R^N$ denotes the polynomial action of $G$ on $\R^N$ and if, for $g \in G$, $\alpha_g = (\alpha_g^1,\ldots,\alpha_g^N)$ is the orthogonal map associated to $g$, we set, for all $[y_0:y_1:\ldots:y_N] \in \mathbb{P}^N(\R)$, 
\begin{eqnarray*}
g \cdot [y_0:y_1:\ldots:y_N] & := & [y_0 : \alpha_g^1(y_1,\ldots,y_N):\ldots:\alpha_g^N(y_1,\ldots,y_N)] \\
					   & = & [\widetilde{\alpha}_g^0(y_0, \ldots, y_N):\ldots:\widetilde{\alpha}_g^N(y_0,\ldots,y_N)].
\end{eqnarray*}
The maps $\widetilde{\alpha}_g := \left(\widetilde{\alpha}^0_g, \ldots, \widetilde{\alpha}^N_g\right) : \R^{N+1} \rightarrow \R^{N+1}$, $g \in G$, are orthogonal maps of $\R^{N+1}$ inducing a polynomial action $\widetilde{\alpha}$ of $G$ on $\R^{N+1}$, which itself induces an action of $G$ on $\mathbb{P}^N(\R)$ (by biregular isomorphisms). We then consider the biregular embedding
$$\Lambda : \begin{array}{ccc}\mathbb{P}^N(\R) & \rightarrow & {\rm M}_{N+1}(\R) = \R^{(N+1)^2}\\ (y_0:\ldots:y_N) & \mapsto & \left(\frac{y_iy_j}{\sum_{i=0}^N y_i^2}\right)_{0\leq i,j\leq N}\end{array}$$
of $\mathbb{P}^N(\R)$ as a compact (irreducible) nonsingular real algebraic set of $\R^{(N+1)^2}$ (see Theorem~3.4.4 of \cite{BCR}): the restriction $\mathbb{P}^N(\R) \rightarrow \Lambda\left(\mathbb{P}^N(\R)\right)$ is in particular a Nash isomorphism between compact affine Nash manifolds.

Let $g \in G$ and let $A(g) = \big(a_{i j}(g)\big)_{0\leq i,j\leq N}$ be the matrix of $\widetilde{\alpha}_g$ in the canonical basis of~$\R^{N+1}$. We consider the Kronecker product 
$$A(g) \otimes A(g) = \big(a_{ik}(g)a_{jl}(g)\big)_{(i,j) \in \{0,\ldots,N\}^2, (k,l) \in \{0,\ldots,N\}^2} \in {\rm M}_{(N+1)^2}(\R)$$
of $A(g)$ with itself and the induced polynomial action $G \times {\rm M}_{(N+1)^2} \rightarrow {\rm M}_{(N+1)^2}$ which associates to $g$ the linear isomorphism 
$$\begin{array}{ccc}{\rm M}_{N+1}(\R) & \rightarrow & {\rm M}_{N+1}(\R)\\M = (u_{kl})_{0\leq k,l\leq N} & \mapsto & A(g) \otimes A(g) M =  \left( \sum_{0 \leq k,l \leq N} a_{ik}(g) a_{jl}(g) u_{kl}\right)_{0 \leq i,j\leq N}.\end{array}$$
Notice that the matrix $A(g) \otimes A(g)$ is an orthogonal matrix of ${\rm M}_{(N+1)^2}(\R)$ since, if $(i,j),(k,l) \in \{0,\ldots,N\}^2$,
$$\sum_{0\leq r,s\leq N} a_{ri}(g) a_{sj}(g) a_{rk}(g) a_{sl}(g) = \sum_{r=0}^N a_{ri}(g) a_{rk}(g) \sum_{s=0}^N a_{sj}(g) a_{sl}(g) = \delta_{ik} \delta_{j l},$$
which is equal to $1$ if $(i,j) = (k,l)$ and $0$ otherwise. Consequently, ${\rm M}_{N+1}(\R) = \R^{(N+1)^2}$ is a polynomial orthogonal representation of $G$.

The biregular isomorphism $\Lambda$ is equivariant with respect to the latter action of $G$: for any $[y_0:\ldots:y_N] \in \mathbb{P}^N(\R)$, we have
\begin{eqnarray*}
\Lambda(g \cdot [y_0:\ldots:y_N]) & = & \left(\frac{\widetilde{\alpha}_g^i(y_0, \ldots, y_N) \widetilde{\alpha}_g^j(y_0, \ldots, y_N)}{\sum_{i=0}^N \widetilde{\alpha}_g^i(y_0, \ldots, y_N)^2}\right)_{0\leq i,j\leq N}\\
					    & = & \frac{1}{\sum_{i=0}^N y_i^2} \left( \sum_{0 \leq k,l \leq N} a_{ik}(g) a_{jl}(g) y_k y_l \right)_{0 \leq i,j\leq N} \\
					  & = & g \cdot \left( \frac{y_i y_j}{\sum_{i=0}^N y_i^2} \right)_{0 \leq i,j\leq N} = g \cdot \Lambda([y_0:\ldots:y_N])
\end{eqnarray*}

As a consequence, the map $\Phi : X \rightarrow \Lambda(Y)$ induced by $\Lambda \circ f$ is an equivariant biregular isomorphism between the $G$-real algebraic set $X$ and the Zariski open subset $\Lambda(Y)$ of the compact nonsingular real algebraic set $\Lambda\left(\mathbb{P}^N(\R)\right)$ of the polynomial orthogonal representation~$\R^{(N+1)^2}$ of $G$, $\Phi$ is a $G$-linearization of $X$ and $X$ is a $G$-$\mathcal{AS}$-set.
\end{proof}

\begin{rem} Any $G$-stable Zariski open subset $U$ of a $G$-real algebraic set $X \subset \R^n$ is a $G$-$\mathcal{AS}$-set as well. Indeed, being the set-theoretic difference of~$X$ and a $G$-stable real algebraic subset of~$X$, the set $U$ is an $\mathcal{AS}$-subset of $\mathbb{P}^n(\R)$. Moreover, if $\Phi : X \rightarrow X' \subset M \subset \R^d$ is a $G$-linearization of $X$, $\Phi(U)$ is an $\mathcal{AS}$-set of $M$ (by theorem \ref{thimageofassetbyinjasmapisas}) and the restriction $U \rightarrow \Phi(U)$ of $\Phi$ is an equivariant $\mathcal{AS}$-homeomorphism (see example \ref{exsasmaps} 2).
\end{rem}

\begin{ex} \label{exstiefelfreegasset} If $n,k \in \Nstar$, the Stiefel manifold $V_k\left(\R^n\right)$ is an $\OO_k(\R)$-$\mathcal{AS}$-set on which~$\OO_k(\R)$ acts freely (see example \ref{exgrealalgset} 3) and then, since $V_k\left(\R^n\right)$ is compact (as a closed subset of the product of $k$ copies of the unit sphere of $\R^n$), $V_k\left(\R^n\right)$ a free $\OO_k(\R)$-$\mathcal{AS}$-set. Moreover, if the compact real algebraic group $G$ is isomorphic to a real algebraic subgroup of $\OO_k(\R)$ (proposition~\ref{propcompactrealalggroupisorthog}), then $V_k\left(\R^n\right)$ a free $G$-$\mathcal{AS}$-set. Notice that the considered action of $\OO_k(\R)$ restricts from a polynomial action of $\OO_k(\R)$ on ${\rm M}_{n,k}(\R) \cong \R^{nk}$ which is given by orthogonal matrices (the canonical scalar product on $\R^{nk}$ corresponds to the scalar product on ${\rm M}_{n,k}(\R)$ which associates to any matrices $A, B$ of ${\rm M}_{n,k}(\R)$ the trace of the product ${}^t\!AB$). Notice also that, since $V_k(\R^n)$ is a compact nonsingular real algebraic set (lemma \ref{lemstiefeliscompactnonsing}), $V_k(\R^n)$ is a compact Nash submanifold of the polynomial orthogonal representation $\R^{nk}$ of $\OO_k(\R)$ (and of~$G$ if $G$ is isomorphic to a real algebraic subgroup of $\OO_k(\R)$).
\end{ex}

Let us also mention the finite group case:

\begin{prop} Suppose that $G$ is a finite group. Then any $\mathcal{AS}$-set equipped with an action, resp. a free action, of $G$ given by homeomorphisms with $\mathcal{AS}$-graph is a $G$-$\mathcal{AS}$-set, resp. a free $G$-$\mathcal{AS}$-set.
\end{prop}

\begin{proof} Denote $G := \{g_1,\ldots,g_k\}$, where $g_1$ is the identity element of $G$, and let $S$ be an $\mathcal{AS}$-set of a compact affine Nash manifold $M$ equipped with an action $G \times S \rightarrow S~;~(g,x) \mapsto \alpha_g(x)$ of $G$ where, if $g \in G$, $\alpha_g : S \rightarrow S$ is an homeomorphism with $\mathcal{AS}$-graph. Consider a Nash diffeomorphism $\Lambda$ from $M$ to a compact Nash submanifold $N$ of an affine space $\R^n$ as well as the map
$$\Phi : \begin{array}{ccc}S & \rightarrow & N \times \cdots \times N\\x & \mapsto & \left(\Lambda(x),\Lambda\left(\alpha_{g_2^{-1}}(x)\right),\ldots,\Lambda\left(\alpha_{g_k^{-1}}(x)\right)\right).\end{array}$$

Remark that the Cartesian product $N \times \cdots \times N$ is a compact Nash submanifold of $\left(\R^n\right)^k$ and that~$\Phi$ is an injective continuous map with $\mathcal{AS}$-graph (use example \ref{exsasmaps} 1, corollary~\ref{corinverseimageassetbyasmapisasset} and remark \ref{remafterthcorasmaps} 1). By theorem \ref{thimageofassetbyinjasmapisas}, the image of $\Lambda$ is an $\mathcal{AS}$-set of $N \times \cdots \times N$ and the restriction $S \rightarrow \Lambda(S)$ of $\Lambda$ is an homeomorphism with $\mathcal{AS}$-graph (see example \ref{exsasmaps} 2). 

Furthermore, if we equip the product $\left(\R^n\right)^k$ with the action of $G$ given by the permutations induced by the product in $G$, $\left(\R^n\right)^k$ is a polynomial representation of $G$ (see remark \ref{remfinitegroupgrealalgsets}) and~$\Lambda$ is equivariant.
\end{proof}

The Cartesian product with a fixed free $G$-$\mathcal{AS}$-set allows to make any $G$-$\mathcal{AS}$-set into a ``relative'' free $G$-$\mathcal{AS}$-set:

\begin{prop} \label{propproductwithfreegassetisgasset} Let $S$ be a free $G$-$\mathcal{AS}$-set and let $A$ be any $G$-$\mathcal{AS}$-set. Then the Cartesian product $A \times S$, equipped with the diagonal action $G \times A \times S \rightarrow A \times S~;~(g,x,y) \mapsto (g \cdot x,g \cdot y)$ of $G$, is a free $G$-$\mathcal{AS}$-set. 
\end{prop}

\begin{proof} Let $\Phi : S \rightarrow S'$ be a free $G$-linearization of $S$ into a compact Nash submanifold $N$ of a polynomial representation $\R^n$ of $G$ and let $\Psi : A \rightarrow A'$ be a $G$-linearization of $A$ into a compact Nash submanifold $L$ of a polynomial representation $\R^l$ of $G$. First, the Cartesian product $A \times S$ is an $\mathcal{AS}$-set and the Cartesian product $A' \times S'$ is an $\mathcal{AS}$-set of the compact Nash submanifold $L \times N$ of $\R^{l+n}$ (remark \ref{remcartprodquasiarcsymsisquasiarcsym}). Furthermore, the product map $\Psi \times \Phi : A \times S \rightarrow A' \times S'~;~(x,y) \mapsto (\Phi(x),\Psi(y))$ is an homeomorphism with $\mathcal{AS}$-graph, since so are the maps $\Psi$ and $\Phi$. Moreover, equipping the affine space $\R^{l} \times \R^n$ with the diagonal action of~$G$ makes~$\R^{l+n}$ into a polynomial representation of $G$ with respect to which the set $A' \times S'$ is~$G$-stable and the map $\Psi \times \Phi$ is equivariant. 

It remains to check that the action of $G$ on the $\mathcal{AS}$-closure of $A' \times S'$ in $L \times N$ (i.e. the $\mathcal{AR}$-closure of $A' \times S'$ in $\R^{l+n}$) is free. But, according to corollary \ref{corarclosureprodisproductarclos}, we have $\overline{A' \times S'}^{\mathcal{AR}} = \overline{A'}^{\mathcal{AR}} \times \overline{S'}^{\mathcal{AR}}$ and, if $g \in G$, $x \in \overline{A'}^{\mathcal{AR}}$ and $y \in \overline{S'}^{\mathcal{AR}}$ satisfy $(x,y) = (g \cdot x, g \cdot y)$, then in particular $g \cdot y = y$, so that $g = e$ (the action of $G$ on $\overline{S'}^{\mathcal{AR}}$ is free by hypothesis).

As a consequence, the product map $\Psi \times \Phi$ is a free $G$-linearization of $A \times S$ and $A \times S$ is a free $G$-$\mathcal{AS}$-set. 
\end{proof}

\begin{ex} If $A$ is any $G$-$\mathcal{AS}$-set, $G$ is isomorphic to a real algebraic subgroup of $\OO_k(\R)$ and $n \in \Nstar$, then the Cartesian product $A \times V_k\left(\R^n\right)$ is a free $G$-$\mathcal{AS}$-set (see example \ref{exstiefelfreegasset}).
\end{ex}

Let us also give the following remark:

\begin{rem} \label{remalllinearizationoffreegassetarenotfree} Even if $S$ is a free $G$-$\mathcal{AS}$-set, a $G$-linearization of $S$ is not necessarily free. Consider for instance the cylinder $X := \{(x,y,z) \in \R^3~|~x^2+y^2 = 1\}$ of $\R^3$ equipped with the restriction of the polynomial action of $S^1$ on $\R^3$ that we considered in example \ref{exgrealalgset} 2. The embedding $(x,y,z) \in X \mapsto \big((x,y),[1:z]\big) \in S^1 \times \mathbb{P}^1(\R)$ makes $X$ into a free $S^1$-$\mathcal{AS}$-set, whereas the $S^1$-linearization of $X$ induced by the map $(x,y,z) \in X \mapsto \left( \frac{x}{\sqrt{z^2+1}}, \frac{y}{\sqrt{z^2+1}}, \frac{z}{\sqrt{z^2+1}}\right) \in S^2$ is not free (the complement of the image of $X$ in $S^2$ consists of the two points $(0,0,1)$ and~$(0,0,-1)$ on which $S^1$ acts trivially).
\end{rem}

\subsection{$\mathcal{AS}$-quotient of free $G$-$\mathcal{AS}$ sets}

We are now going to use proposition \ref{propquotientfreeassetsnashsub} and theorem \ref{propsemialgquotientmapbetweenfreeassetscompnashsubisas} to define an $\mathcal{AS}$-quotient for any free $G$-$\mathcal{AS}$-set, which will be well-defined up to homeomorphism with $\mathcal{AS}$-graph. Let~$S$ be a free $G$-$\mathcal{AS}$-set.

\begin{prop} \label{propwhichimpliesasquotientwelldefined} Let $\Phi : S \rightarrow S'$ and $\Psi : S \rightarrow S''$ be free $G$-linearizations of $S$. The map 
$$\left(\Psi \circ \Phi^{-1}\right)_{/G} : S'/G \rightarrow S''/G$$
is an \emph{$\mathcal{AS}$-homeomorphism}, i.e. an homeomorphism with $\mathcal{AS}$-graph.
\end{prop}

\begin{proof} The map $\Psi \circ \Phi^{-1}$ is an equivariant $\mathcal{AS}$-homeomorphism (cf. remark \ref{remasmapembeddings} 1 and corollary \ref{corcompasmapsisasmap}) between~$\mathcal{AS}$-sets satisfying the hypotheses of theorem \ref{propsemialgquotientmapbetweenfreeassetscompnashsubisas}. Therefore, by the latter result and by proposition \ref{propquotientfreeassetsnashsub} as well as proposition \ref{propfunctsemialgquotientsemialgcontprop}, the induced map $\left(\Psi \circ \Phi^{-1}\right)_{/G} : S'/G \rightarrow S''/G$ is an $\mathcal{AS}$-homeomorphism.
\end{proof}

\begin{de} \label{deasquotientoffreegassets} If $\Phi : S \rightarrow S'$ is a free $G$-linearization of the free $G$-$\mathcal{AS}$-set $S$, we denote by $S/G$ the $\mathcal{AS}$-set $S'/G$ and we define the \emph{$\mathcal{AS}$-quotient} $\varrho_S : S \rightarrow S/G$ of $S$ by $G$ to be the composition of $\Phi$ with the semialgebraic quotient $\varpi_{S'} : S' \rightarrow S'/G$ of $S'$ by $G$. By abuse of terminology, the $\mathcal{AS}$-set $S/G$ will also be called the $\mathcal{AS}$-quotient of $S$ by $G$. 
\end{de}

The $\mathcal{AS}$-quotient $\varrho_S : S \rightarrow S/G$ is a continuous $\mathcal{AS}$-map which is, according to above proposition \ref{propwhichimpliesasquotientwelldefined}, well-defined up to $\mathcal{AS}$-homeomorphism. We can furthermore show the functoriality of this $\mathcal{AS}$-quotient with respect to continuous $\mathcal{AS}$-maps:

\begin{theo} \label{theofunctorialityasquotientfreegas} Let $T$ be a free $G$-$\mathcal{AS}$-set and let $f : S \rightarrow T$ be a continuous map with~$\mathcal{AS}$-graph. The map $f$ induces, in a functorial way, a continuous $\mathcal{AS}$-map $f_{/G} : S/G \rightarrow T/G$ such that $f_{/G} \circ \varrho_S = \varrho_T \circ f$.
\end{theo}

\begin{proof} Let $\Phi : S \rightarrow S'$, resp. $\Psi : T \rightarrow T'$, be a free $G$-linearization of $S$, resp. $T$. We consider the map $\left(\Psi \circ f \circ \Phi^{-1}\right)_{/G} : S'/G \rightarrow T'/G$ and denote it by $f_{/G}$: by theorem \ref{propsemialgquotientmapbetweenfreeassetscompnashsubisas} and proposition \ref{propfunctsemialgquotientsemialgcontprop}, the map $f_{/G}$ is a continuous $\mathcal{AS}$-map. Furthermore, we have
$$\left(\Psi \circ f \circ \Phi^{-1}\right)_{/G} \circ \left(\varpi_{S'} \circ \Phi\right) = \varpi_{T'} \circ \left(\Psi \circ f \circ \Phi^{-1}\right) \circ \Phi = (\varpi_{T'} \circ \Psi) \circ f,$$
i.e. $f_{/G} \circ \varrho_S = \varrho_T \circ f$.
\end{proof}

\begin{rem} With the above notations, if the map $f$ is proper, then so is the induced map~$f_{/G}$ (by proposition \ref{propfunctsemialgquotientsemialgcontprop} and using the fact that the maps $\Phi$ and $\Psi$ are homeomorphisms). 
\end{rem}

The $\mathcal{AS}$-quotient of free $G$-$\mathcal{AS}$-sets also preserves closed inclusions in the following sense:

\begin{prop} \label{propasquotfreeassetspreservesclosedincl} Let $M$ be a compact affine Nash manifold such that $S$ is an $\mathcal{AS}$-set of $M$ and let $T$ be an $\mathcal{AS}$-set of $M$ such that $T$ is a $G$-stable subset of $S$. Then $T$ is a free $G$-$\mathcal{AS}$-set and, if $T$ is closed in $S$, $T/G$ is a closed subset of $S/G$.
\end{prop}

\begin{proof} Let $\Phi : S \rightarrow S'$ be a free $G$-linearization of $S$ into a compact Nash submanifold $N$ of a polynomial representation $\R^n$ of $G$. Then the restriction $T \rightarrow \Phi(T)$ is a free $G$-linearization of~$T$. Indeed, $T' := \Phi(T)$ is an $\mathcal{AS}$-set of $N$ by theorem \ref{thimageofassetbyinjasmapisas} and the latter map is an equivariant $\mathcal{AS}$-homeomorphism (cf. example \ref{exsasmaps} 2). Moreover, the $\mathcal{AS}$-closure of $T'$ in $N$ is included in the $\mathcal{AS}$-closure of $S'$ in $N$ and therefore the action of $G$ on $T'$ is free. Finally, notice that, if $T$ is closed in $S$, since $\Phi$ is an homeomorphism, $T'$ is closed in $S'$, so that $T/G = T'/G = \varpi_{S'}(T')$ is a closed subset of $S'/G = S/G$, because the map $\varpi_{S'} : S' \rightarrow S'/G$ is closed (proposition~\ref{propquotientmapproperclosedopen}).
\end{proof}

\begin{rems} 
~
	\begin{enumerate}
		\item With the above notations, if $T$ is open in $S$, then $T/G$ is an open subset of $S/G$.
		\item The $\mathcal{AS}$-quotient of $T$ by $G$ is (up to $\mathcal{AS}$-homeomorphism) the restriction $T \rightarrow T/G$ of the $\mathcal{AS}$-quotient $\varrho_S : S \rightarrow S/G$ of $S$.
	\end{enumerate}
\end{rems}

\subsection{Relative quotients of $G$-$\mathcal{AS}$-sets}

Taking advantage of proposition \ref{propproductwithfreegassetisgasset}, we now plan to associate to any $G$-$\mathcal{AS}$-set an $\mathcal{AS}$-quotient relatively to a fixed free $G$-$\mathcal{AS}$-set, in a way which would be functorial in both entries. Keep the above free $G$-$\mathcal{AS}$-set $S$ and let $A$ be any $G$-$\mathcal{AS}$-set: according to proposition \ref{propproductwithfreegassetisgasset}, the Cartesian product $A \times S$, equipped with the induced diagonal action of $G$, a free $G$-$\mathcal{AS}$-set. 

\begin{de} \label{derelasquotientwithrespecttoafreegasset} We denote by $\varrho_{A}^S$ the $\mathcal{AS}$-quotient $\varrho_{A \times S} : A \times S \rightarrow \left(A \times S\right)/G$ and call it the \emph{relative $\mathcal{AS}$-quotient of $A$ with respect to $S$}. We will also denote $A^S/G := \left(A \times S\right)/G$.
\end{de}

The relative $\mathcal{AS}$-quotient with respect to the free $G$-$\mathcal{AS}$-set $S$ is functorial with respect to continuous $\mathcal{AS}$-maps, but we also have functoriality with respect to the entry $S$ in the following sense:

\begin{theo} Let $B$ be a $G$-$\mathcal{AS}$-set and $f : A \rightarrow B$ be a continuous $\mathcal{AS}$-map. Let $T$ be a free $G$-$\mathcal{AS}$-set and $h : S \rightarrow T$ be a continuous $\mathcal{AS}$-map. The map $f$ induces, in a functorial way, a continuous $\mathcal{AS}$-map $f^S_{/G} : A^S/G \rightarrow B^S/G$ and the map $h$ induces, in a functorial way, a continuous $\mathcal{AS}$-map $h_{/G}^A : A^S/G \rightarrow A^T/G$ such that the following diagram commutes:

\begin{equation*}
	\xymatrix{A \times S \ar[rd]_{\varrho_{A}^S} \ar[ddd]_{{\rm Id}_A \times h} \ar[rrr]^{f \times {\rm Id_S}}& & & B \times S \ar[ld]^{\varrho^S_{B}} \ar[ddd]^{{\rm Id}_B \times h}\\
		      & A^S/G \ar[r]^{f^S_{/G}} \ar[d]_{h^A_{/G}} &  B^S/G \ar[d]^{h^B_{/G}} & \\
		      & A^T/G \ar[r]_{f^{T}_{/G}} &  B^T/G& \\ 
	               A \times T \ar[ru]^{\varrho_A^{T}} \ar[rrr]_{f \times {\rm Id_T}}& & & B \times T \ar[lu]_{\varrho_B^{T}}
	}
\end{equation*} 

\end{theo}

\begin{proof} We consider and denote by $f^S_{/G}$ the continuous $\mathcal{AS}$-map $\left(f \times {\rm Id}_S\right)_{/G} : A^S/G \rightarrow B^S/G$ given by theorem \ref{theofunctorialityasquotientfreegas} (the sets $A \times S$ and $B \times S$ are free $G$-$\mathcal{AS}$-sets by proposition \ref{propproductwithfreegassetisgasset}) and by $h^A_{/G}$ the continuous $\mathcal{AS}$-map $\left({\rm Id}_A \times h\right)_{/G} : A^S/G \rightarrow A^T/G$ given by the same theorem (the sets $A \times T$ is a free $G$-$\mathcal{AS}$-set as well).   

The equality $f^T_{/G} \circ h^A_{/G} = h^B_{/G} \circ f^S_{/G}$ is provided by functoriality of the semialgebraic quotient by $G$, applied to the equality $\left(f \times {\rm Id_T}\right) \circ \left({\rm Id}_A \times h\right) = \left({\rm Id}_B \times h\right) \circ \left(f \times {\rm Id_S}\right)$.
\end{proof}

The relative $\mathcal{AS}$-quotient also preserves closed inclusions:

\begin{prop} Let $M$ be a compact affine Nash manifold such that $A$ is an $\mathcal{AS}$-set of $M$ and let $B$ be an $\mathcal{AS}$-set of $M$ such that $B$ is a $G$-stable subset of $A$. Then $B$ is a $G$-$\mathcal{AS}$-set and, if $B$ is closed in $A$, $B^S/G$ is a closed subset of $A^S/G$.
\end{prop}

\begin{proof} First, let $\Psi : A \rightarrow A'$ be a $G$-linearization of $A$ into a compact Nash submanifold $L$ of a polynomial representation $\R^l$ of $G$. Then the restriction $B \rightarrow \Psi(B)$ is a $G$-linearization of $B$ since $B' := \Psi(B)$ is an $\mathcal{AS}$-set of $L$ by theorem \ref{thimageofassetbyinjasmapisas} and the latter map is then an equivariant~$\mathcal{AS}$-homeomorphism. Moreover, if $B$ is closed in $A$, then $B \times S$ is a $G$-stable closed subset of $A \times S$ and we can use proposition \ref{propasquotfreeassetspreservesclosedincl} to assert that $B^S/G$ is closed in $A^S/G$.
\end{proof}

In the same way (see also proposition \ref{propasquotfreeassetspreservesclosedincl}), we can show that if $T$ is an $\mathcal{AS}$-set which is a~$G$-stable closed subset of $S$, then $A^T/G$ is a closed subset of $A^S/G$.

\begin{ex} Suppose that $G$ is isomorphic to a real algebraic subgroup of $\OO_k(\R)$ and consider the free $G$-$\mathcal{AS}$-sets $V_k\left(\R^n\right)$, $n \in \Nstar$ (example \ref{exstiefelfreegasset}). For any $n \in \Nstar$, we have an equivariant closed embedding $V_k\left(\R^n\right) \hookrightarrow V_k\left(\R^{n+1}\right)$ consisting in adding a line of zeroes to any matrix of $V_k\left(\R^n\right)$, so that, up to $\mathcal{AS}$-homeomorphism, the $\mathcal{AS}$-set $A^{V_k\left(\R^n\right)}/G$ is a closed subset of the $\mathcal{AS}$-set~$A^{V_k\left(\R^{n+1}\right)}/G$.
\end{ex}

\section{Appendix}

In this appendix, we gather some more or less classical properties that are used in the present document, together with their proofs, in order for this text to be almost self-contained.

\subsection{Around Hilbert's Nullstellensatz}

Hilbert's Nullstellensatz (cf. for instance Theorem 1.2 of \cite{Bump}) asserts that, when $I$ is an ideal of the polynomial ring $\C[x_1,\ldots,x_n]$, then $\mathsf{I}\left(\mathsf{V}(I)\right) = I$ if and only if $I$ is radical (if $B$ is any commutative ring and $J$ is an ideal of $B$, $J$ is said \emph{radical} if for any $b \in B$, the existence of a positive integer $k$ such that $b^k \in J$ implies that $b \in J$). We will give two consequences of this result. First, let us state the complex counterpart of theorem \ref{theorealnullstellensatzequivcategories}:

\begin{theo} \label{theoequivcatcomphilbertsnullstel} The operation which associates to any complex algebraic set $Z$ its $\C$-algebra of polynomial functions on $Z$ gives an equivalence of categories between the category of (polynomial isomorphism classes of) complex algebraic sets and polynomial maps and the category of (isomorphism classes of) reduced (as commutative rings) finitely generated $\C$-algebras and morphisms of $\C$-algebras (recall that a commutative ring $A$ is said \emph{reduced} if for any $a \in A$, the existence of a positive integer~$k$ such that $a^k = 0_A$ implies that $a = 0_A$).
\end{theo}

\begin{proof} The proof is analog to the proof of theorem \ref{theorealnullstellensatzequivcategories}, using the fact that, if $f$ is a polynomial function on a complex algebraic set $Z$ and $k$ is a positive integer such that $f^k(z) = 0$ for all~$z \in Z$, then $f$ is the zero function on $Z$. Furthermore, if $B$ is any commutative ring and $J$ is an ideal of $B$, $J$ is radical if and only if the quotient ring $B/J$ is reduced and we can apply Hilbert's Nullstellensatz in the second part of the proof.
\end{proof}

In the previous statement, we can replace the word ``polynomial'' with the word ``regular'' since any regular map on a complex algebraic set (a \emph{regular map} on a complex algebraic set~$Z$ is a map $Z \rightarrow \C^n$ whose any coordinate function is of the form $z \in Z \mapsto \frac{p(z)}{q(z)} \in \C$ where $p$ and $q$ are polynomial functions on $Z$ such that for all $z \in Z$, $q(z) \neq 0$) is the restriction of a polynomial map: this is the second consequence of Hilbert's Nullstellensatz that we state here.

\begin{theo} \label{theoregularfunctiononcompalgsetispoly} Let $f : Z \rightarrow \C$ be a regular function on a complex algebraic set $Z$ of $\C^m$. Then~$f$ is a polynomial function on $Z$. 
\end{theo} 

\begin{proof} Let $F_1,\ldots, F_k, P,Q \in \C[x_1,\ldots,x_m]$ such that $Z = \mathsf{V}(F_1,\ldots,F_k)$ and, for any $z \in Z$, $Q(z) \neq 0$ and $f(z) = \frac{P(z)}{Q(z)}$. Denote by $I$ the ideal of $\C[x_1,\ldots,x_m]$ generated by the polynomials $F_1,\ldots,F_k,Q$: we have $\mathsf{V}(I) = \emptyset$ and then $\mathsf{I}\left(\mathsf{V}(I)\right) = \C[x_1,\ldots,x_m]$. On the other hand, by Hilbert's Nullstellensatz (see for instance Theorem 1.2 of \cite{Bump}), $\mathsf{I}\left(\mathsf{V}(I)\right)$ is the radical of $I$, so that the constant polynomial $1$ belongs to the latter and therefore to $I$. As a consequence, there exist polynomials $R_1,\ldots,R_k, S \in \C[x_1,\ldots,x_m]$ such that $1 = S Q + \sum_{l=1}^k R_l F_l$ and then, for any $z \in Z$, $1 = S(z) Q(z)$ so that we can write
$$f(z) = 1 \times f(z) = S(z) Q(z) \frac{P(z)}{Q(z)} = S(z) Q(z).$$
\end{proof}

\subsection{About Stiefel manifolds} \label{subsectappendixstiefel}

If $k$ and $n$ are positive integers, an important example of $\OO_k(\R)$-real algebraic set that we considered in the text is the Stiefel manifold $V_k\left(\R^n\right)$ of index $k$ on $\R^n$ (example \ref{exgrealalgset} 2). We prove here some basic properties of $V_k\left(\R^n\right)$ as a real algebraic subset of $\R^{nk}$.

\begin{lem} \label{lemstiefeliscompactnonsing} The subset $V_k\left(\R^n\right)$ of $\R^{nk}$ is a compact nonsingular real algebraic set of dimension $nk - \frac{k(k+1)}{2}$.
\end{lem}

\begin{proof} The set $V_k\left(\R^n\right)$ is the set of $k$-tuples of orthonormal vectors of $\R^n$ i.e. the set of points $A = \left(a_{ij}\right)_{1\leq i \leq n, 1 \leq j \leq k}$ of $\mathbb{R}^{nk}$ such that 
$${}^t\!A A = I_k ~\Longleftrightarrow ~ \forall j,j' \in \{1,\ldots,k\}, \, \sum_{i =1}^n a_{ij} a_{ij'} = \delta_{jj'}.$$	
Denote $P_l := \sum_{i=1}^n a_{il}^2 - 1 = 0$, $1 \leq l \leq k$, and $Q_{j,j'} := \sum_{i=1}^n a_{ij} a_{ij'} = 0$, $1 \leq j < j' \leq k$, so that $V_k\left(\R^n\right) = V(P_l, Q_{j,j'},~1 \leq l \leq k,~1 \leq j < j' \leq k)$.

We have, for $i \in \{1,\ldots,n\}$ and $l \in \{1,\ldots,k\}$, $\frac{\partial P_l}{\partial a_{il}} = 2 a_{il}$ and, for $i \in \{1,\ldots,n\}$ and $j,j' \in \{1,\ldots,k\}$ such that $j<j'$, $\frac{\partial Q_{j,j'}}{\partial a_{ij}} = a_{ij'}$ and $\frac{\partial Q_{j,j'}}{\partial a_{ij'}} = a_{ij}$, and the columns of the associated Jacobian matrix at a point $a = (a_{ij})_{1 \leq i \leq n, 1 \leq j \leq k}$ of $V_k(\R^n)$ are linearly independent. Indeed, if $\lambda_l$, $1 \leq l \leq k$, and $\lambda_{j,j'}$, $1 \leq j < j' \leq k$, are scalars such that for all $i \in \{1,\ldots,n\}$ and $\iota \in \{1,\ldots,k\}$,
$$\sum_{l=1}^k \lambda_l \frac{\partial P_l}{\partial a_{i\iota}}(a) + \sum_{1 \leq j < j' \leq k} \lambda_{j,j'} \frac{\partial Q_{j,j'}}{\partial a_{i\iota}}(a) = 0$$
i.e.
$$2 \lambda_{\iota} a_{i\iota} + \sum_{\iota < j' \leq k} \lambda_{\iota, j'} a_{ij'} + \sum_{1\leq j < \iota} \lambda_{j,\iota} a_{ij} = 0,$$
then all the scalars $\lambda_l$, $1 \leq l \leq k$, and $\lambda_{j,j'}$, $1 \leq j < j' \leq k$, are equal to $0$ since the vectors $(a_{i1})_{1\leq i \leq n},\ldots, (a_{ik})_{1\leq i \leq n}$ form an orthonormal family of $\R^n$. As a consequence, any point of~$V_k(\R^n) \subset \mathbb{R}^{nk}$ is nonsingular in dimension~$nk - \frac{k(k+1)}{2}$ (cf. Proposition 3.3.10 of \cite{BCR}) and the real algebraic set $V_n(\R^k)$ is therefore nonsingular.

Finally, the set $V_k\left(\R^n\right)$ is compact as a closed subset of the product of $k$ copies of the unit sphere of $\R^n$.
\end{proof}

Notice that we can adapt the above proof to show that the complex algebraic set $V_k^n(\C) = \left\{ A \in {\rm M}_{n,k}(\C)~|~{}^t\!A A = k\right\}$ of $\M_{n,k}(\C) \cong \mathbb{C}^{nk}$ that we considered in example \ref{exscomplexificationpolyaction} 2 is nonsingular of dimension $nk - \frac{k(k+1)}{2}$ (see also section 2 of Chapter II of \cite{AK}). We furthermore show that the latter complex algebraic set is irreducible when $k<n$.

\begin{lem} Suppose that $k < n$. Then the complex algebraic set $V_k^n(\C)$ of $\C^{nk}$ is irreducible.
\end{lem}

\begin{proof} The complex algebraic group $\OO_n(\C)$, resp. $\SO_n(\C)$, polynomially acts on $V_k^n(\C)$ by left matrix multiplication: this action is given by the action of $\vartheta$-isometries, resp. $\vartheta$-isometries of determinant $1$, on $\C^n$ (see example \ref{exscomplexificationpolyaction} 2 for the definition of the bilinear form $\vartheta$). Since the bilinear symmetric form $\vartheta$ is non-degenerate, Witt's theorem can be used to assert that the polynomial action of $\SO_n(\C)$ on $V_k^n(\C)$ is transitive (if $A$ and $B$ are any two $k$-tuples of $V_k^n(\C)$, the linear isomorphism from ${\rm Span}(A)$ to ${\rm Span}(B)$ that sends $A$ onto $B$ is an isometry which extends into an isometry of $\C^n$ which can be supposed to be of determinant $1$ since $k<n$).  We deduce that the complex algebraic subset $V_k^n(\C)$ is irreducible, as the image of $\SO_n(\C)$ by the surjective polynomial orbit map $O \in \SO_n(\C) \mapsto OA \in V_k^n(\C)$, where $A \in V_k^n(\C)$. 
\end{proof}

If $k < n$, the real algebraic set $V_k\left(\R^n\right)$ of $\R^{nk}$ is irreducible as well: we can adapt the previous proof considering the canonical scalar product on $\R^n$ or argue that $V_k^n(\C)$ is the complexification of $V_k\left(\R^n\right)$ (the latter fact uses what we showed above: see example \ref{exscomplexificationpolyaction} 2).

\subsection{Continuous, proper, closed, open maps and locally compact sets} \label{subsectappendixtopology}

Let us now assert and prove the following topological properties that we used in the proof of proposition \ref{propquotientmapproperclosedopen}:

\begin{lem} \label{lemtopgenproperopenclosed} Let $\varphi : A \rightarrow B$ be a continuous map between Hausdorff topological spaces $A$ and~$B$.
\begin{enumerate}
	\item If there exists a continuous map $\psi : B \rightarrow C$ between Hausdorff topological spaces such that $\psi \circ \varphi$ is proper, then $\varphi$ is proper.
	\item Let $B'$ be a subset of $B$. If $\varphi$ is proper, then the restriction $\varphi^{-1}(B') \rightarrow B'$ of $\varphi$ is proper as well.
	\item If $\varphi$ is proper and $B$ is locally compact, then $\varphi$ is a closed map.
	\item Let $A'$ be a subset of $A$. If the map $\varphi$ is closed, resp. open, and if for any closed, resp. open, subset $D$ of $A$, $\varphi(D \cap A') = \varphi(D) \cap \varphi(A')$, then the restriction $A' \rightarrow \varphi(A')$ of $\varphi$ is closed, resp. open, as well. 
\end{enumerate}
\end{lem} 

\begin{proof}
\begin{enumerate}
	\item Let $\psi : B \rightarrow C$ be a continuous map between Hausdorff topological spaces such that $\psi \circ \varphi$ is proper and let $K$ be a compact subset of $B$. We have $\varphi^{-1}(K) \subset (\psi \circ \varphi)^{-1}(\psi(K))$. But $\varphi^{-1}(K)$ is a closed subset of $A$ ($\varphi$ is continuous and $K$ is closed in $B$ since $B$ is Hausdorff) and $(\psi \circ \varphi)^{-1}(\psi(K))$ is compact, since $\psi$ is continuous and $\psi \circ \varphi$ is proper. As a consequence, $\varphi^{-1}(K)$ is compact.

	\item Suppose that $\varphi$ is proper. If $K$ is a compact subset of $B'$, then $K$ is a compact subset of~$B$ and, since $\varphi$ is proper, $\varphi^{-1}(K)$ is a compact subset of $A$ included in $\varphi^{-1}(B')$ and then a compact subset of $\varphi^{-1}(B')$.

	\item Suppose that $\varphi$ is proper and $B$ is locally compact, and let $D$ be a closed subset of $A$: let us show that $\overline{\varphi(D)} = \varphi(D)$. Let $y \in \overline{\varphi(D)}$ and let $U$, resp. $K$ be an open, resp. a compact, subset of $B$ such that $y \in U \subset K$ ($B$ is locally compact). Remark that $y \in \overline{K \cap \varphi(D)}$: if $V$ is an open neighborhood of $y$, so is the intersection $U \cap V$ and then $(U \cap V) \cap \varphi(D) \neq \emptyset$, in particular $V \cap (K \cap \varphi(D)) \neq \emptyset$. Denote $L := \varphi^{-1}(K)$:~$L$ is compact since $\varphi$ is proper and then $L \cap D$ is compact as well (since $D$ is closed), so that $\varphi(L \cap D)$ is also compact ($\varphi$ is continuous), in particular closed in $B$ ($B$ is Hausdorff). On the other hand, $\varphi(L \cap D) = \varphi\left( \varphi^{-1}(K) \cap D\right) = K \cap \varphi(D)$, and then $y \in \overline{K \cap \varphi(D)} = \overline{\varphi(L \cap D)} = \varphi(L \cap D) \subset \varphi(D)$.

	\item Suppose that $\varphi$ is closed and that for any closed subset $D$ of $A$, $\varphi(D \cap A') = \varphi(D) \cap \varphi(A')$. Consider a closed subset $F \cap A'$ of $A'$, where $F$ is a closed subset of $A$, then $\varphi(F \cap A') = \varphi(F) \cap \varphi(A')$ is a closed subset of $\varphi(A')$, since $\varphi(F)$ is closed in $B$ because the map $\varphi$ is closed. Finally, notice that the same arguments work if we replace the word ``closed'' with ``open''.
\end{enumerate}
\end{proof} 

Let us also prove the following statement about locally compact subsets of a locally compact Hausdorff space:

\begin{lem} \label{lemequivloccomplocclosedloccomphaus} Let $A$ be a locally compact Hausdorff space.
\begin{enumerate}
	\item If $K$ is any compact subset of $A$ and $x$ is an element of $A \setminus K$, then there exist disjoint open subsets $U$ and $V$ of $A$ such that $x \in U$ and $K \subset V$.
	\item If $x$ is any element of $A$ and $U$ is any open neighborhood of $x$ in $A$, then there exists an open neighborhood $V$ of $x$ such that the closure $\overline{V}$ of $V$ in $A$ is compact and contained in~$U$.
	\item Any open or closed subset of $A$ is locally compact.
	\item If $D$ is any subset of $A$, then $D$ is locally compact if and only $D$ is locally closed.
\end{enumerate}
\end{lem}

\begin{proof} 
\begin{enumerate}
	\item For all $y \in K$, consider an open neighborhood $U_y$ of $x$ in $A$ and an open neighborhood $V_y$ of $y$ in $A$ such that $U_y \cap V_y \neq \emptyset$ ($A$ is Hausdorff). The open sets $V_y$, $y \in K$, cover the compact set $K$, so that there exists $y_1,\ldots,y_n \in K$ such that $K \subset \bigcup_{i=1}^n V_{y_i}$. If we then set $U := \bigcap_{i=1}^n U_{y_i}$ and $V := \bigcup_{i=1}^n V_{y_i}$, $U$ and $V$ are both open in $A$ and satisfy~$U \cap V = \emptyset$.
	\item Since $A$ is locally compact, there exist an open subset $W$ and a compact subset $K$ of~$A$ such that $x \in W \subset K$. Denote by $V$ the open set $U \cap W$, which is contained in $K$, and remark that $K \setminus V$ is then a closed subset of the compact set $K$. The set $K \setminus V$ is therefore compact as well and, since $x \notin K \setminus V$, by the first item, we can find open subsets $U_1$, $U_2$ of $A$ such that $x \in U_1$, $K \setminus V \subset U_2$ and $U_1 \cap U_2 = \emptyset$. If we finally set $V_0 := V \cap U_1 \subset K$, we have $\overline{V_0} \subset K$ and
$$\overline{V_0} \cap U_2 \subset \overline{U_1} \cap U_2 \subset \overline{U_1 \cap U_2} = \emptyset$$
(we have $\overline{U_1} \cap U_2 \subset \overline{U_1 \cap U_2}$ since $U_2$ is open) so that $\overline{V_0} \subset V$. Therefore, we have $x \in V_0 \subset \overline{V_0} \subset V \subset U$ and $\overline{V_0}$, being a closed subset of $K$, is compact.
	\item Let $C$ be a subset of $A$. If $C$ is open and $x \in C$, there exist, by the previous statement, an open neighborhood $V$ of $x$ such that $V \subset \overline{V} \subset C$ and $\overline{V}$ is compact, so that $C$ is locally compact. If $C$ is closed and $x \in C$, let $U$ be an open subset and $K$ be a compact subset of the locally compact space $A$ such that $x \in U \subset K$: we have $x \in U \cap C \subset K \cap C \subset C$, where $U \cap C$ is an open subset of $C$ and $K \cap C$ is a compact subset of $C$ ($C$ is closed), so that $C$ is locally compact in this case as well.
	\item Let us first recall that $D$ is a \emph{locally closed subset of $A$} if $D$ is the intersection of an open set $U$ and a closed set $F$ of $A$, and that $D$ is locally closed in $A$ if and only if $D$ is the intersection of an open set $U$ of $A$ and the closure $\overline{D}$ of $D$ in $A$, if and only if $D$ is an open subset of $\overline{D}$.
	
Now, suppose that $D$ is locally compact and let $x \in D$: there exists an open subset~$U$ and a compact subset $K$ of $A$ such that $x \in U \cap D \subset K \subset D$. We then have $x \in U \cap \overline{D} \subset \overline{U \cap D} \subset K \subset D$ (we have $U \cap \overline{D} \subset \overline{U \cap D}$ since $U$ is open and $\overline{U \cap D} \subset K$ since the compact set $K$ is closed in $A$ because $A$ is Hausdorff). As a consequence, $D$ is open in~$\overline{D}$ i.e. $D$ is a locally closed subset of $A$.

As for the converse implication, since any intersection of two locally compact subsets of~$A$ is locally compact, if $D$ is locally closed i.e. if $D$ is the intersection of an open and a closed subset of $A$, then $D$ is locally compact by the third item.  	
\end{enumerate}
\end{proof}

\subsection{Compact subgroups of real orthogonal groups, Lie algebras and the exponential map} \label{subsectsubgprealrotholiealgexpmap}

Let $G$ be a compact subgroup of the real orthogonal group $\OO_n(\R)$. We show the following property that we used in subsection \ref{subsectcompsubgrorthogpliealgcomplexif}.

\begin{lem} \label{lemexpliealgcommwithdiff} Consider the matrix exponential ${\rm exp} : {\rm M}_n(\R) \rightarrow {\rm M}_n(\R)$. We have ${\rm exp}\left({\rm Lie}(G)\right) \subset G$ and, if $\varphi : G \rightarrow H$ is a morphism of Lie groups to a compact subgroup $H$ of $\OO_N(\R)$, we have ${\rm exp} \circ {\rm d}_{I_n} \varphi = \varphi \circ {\rm exp}$.
\end{lem}

\begin{proof} Let $A \in {\rm Lie}(G)$: there exists a smooth curve $\gamma : I \rightarrow G$, with $I$ an open interval of $\R$ containing $0$, such that $\gamma(0) = I_n$ and $\gamma'(0) = A$. For $t \in I$, write $\gamma(t) = \gamma(0) + t \gamma'(0) + t A(t) = I_n + t A + t A(t)$, where $A(t)$ tends to the zero matrix $0_n$ when $t$ tends to $0$. For $B \in {\rm M}_n(\R)$ with $\| I_n - B\| < 1$ (where $\| \cdot \| : {\rm M}_n(\R) \rightarrow [0;+\infty[$ is a matrix norm on ${\rm M}_n(\R)$), write also the matrix logarithm of $B$ as ${\rm log}(B) = - (I_n - B) + (I_n-B) R(B)$, where $R(B)$ tends to $0_n$ when~$B$ tends to $I_n$. For $t \in I$ close enough to $0$ so that $\| I_n - \gamma(t) \| < 1$, we have
$${\rm log}(\gamma(t)) = t A + t A(t) - t(A + A(t)) R(\gamma(t)) = t(A + B(t))$$
where $B(t)$ tends to $0_n$ when $t$ tends to $0$. Therefore, for $n \in \Nstar$ big enough, we have $n \, {\rm log}\left(\gamma\left(\frac{1}{n}\right)\right) = A + B\left(\frac{1}{n}\right)$ and then 
$$\left(\gamma\left(\frac{1}{n}\right)\right)^n = {\rm exp}\left(A + B\left(\frac{1}{n}\right)\right),$$
so that $\left(\gamma\left(\frac{1}{n}\right)\right)^n \xrightarrow[n \rightarrow + \infty]{} {\rm exp}(A)$ and ${\rm exp}(A) = \underset{n \rightarrow + \infty}{\rm lim} \left(\gamma\left(\frac{1}{n}\right)\right)^n$. But $\gamma(I) \subset G$, $G$ is a matrix group and $G$ is closed in ${\rm M}_n(\R)$: as a consequence, ${\rm exp}(A) \in G$. 
 
As for the second statement, consider the smooth curve $\varphi \circ \gamma : I \rightarrow H$, remark that $\varphi \circ \gamma(0) = I_N$ (since $\varphi$ is a group homomorphism) and recall that ${\rm d}_{I_n} \varphi(\gamma'(0)) = (\varphi \circ \gamma)'(0)$ i.e. $(\varphi \circ \gamma)'(0) = {\rm d}_{I_n} \varphi(A)$. By the previous argument, we have $\left(\varphi \circ \gamma\left(\frac{1}{n}\right)\right)^n \xrightarrow[n \rightarrow + \infty]{} {\rm exp}\left({\rm d}_{I_n} \varphi(A)\right)$ and ${\rm exp}\left({\rm d}_{I_n} \varphi(A)\right) = \underset{n \rightarrow + \infty}{\rm lim} \left(\varphi \circ \gamma\left(\frac{1}{n}\right)\right)^n$, and then, thanks to the continuity of the group homomorphism $\varphi$,
$${\rm exp}\left({\rm d}_{I_n} \varphi(A)\right) = \underset{n \rightarrow + \infty}{\rm lim} \left(\varphi \circ \gamma\left(\frac{1}{n}\right)\right)^n = \underset{n \rightarrow + \infty}{\rm lim} \varphi \left(\gamma\left(\frac{1}{n}\right)^n\right) = \varphi\left(\underset{n \rightarrow + \infty}{\rm lim} \gamma\left(\frac{1}{n}\right)^n\right) = \varphi\left({\rm exp}(A)\right).$$
\end{proof}
 
\begin{rems} \label{remafterlemexpliealgcommwithdiff}
~
	\begin{enumerate}
		\item Thanks to the previous lemma, we can also prove that 
		$${\rm Lie}(G) = \left\{A \in {\rm M}_n(\R)~|~\forall t \in \R,~{\rm exp}(tA) \in G\right\}.$$ 
Indeed, the inclusion ${\rm exp}\left({\rm Lie}(G)\right) \subset G$ shows the inclusion of ${\rm Lie}(G)$ in the right-hand side set (${\rm Lie}(G)$ is a vector subspace of ${\rm M}_n(\R)$). In order to prove the converse inclusion, let $A \in {\rm M}_n(\R)$ such that ${\rm exp}(tA) \in G$ for all $t \in \R$ and consider the smooth curve $\gamma : t \in \R \mapsto {\rm exp}(tA) \in G$. The latter satisfies $\gamma(0) = I_n$ and, for all $t \in \R$, $\gamma'(t) = A \, {\rm exp}(tA)$, in particular $\gamma'(0) = A$: as a consequence, $A \in {\rm Lie}(G)$.
		\item Since $G \subset \OO_n(\R)$, we have ${\rm Lie}(G) = \left\{A \in {\rm M}_n(\R)~|~\forall t \in \R,~{\rm exp}(tA) \in G\right\} \subset {\rm Lie}(\OO_n(\R)) = \left\{A \in {\rm M}_n(\R)~|~{}^t\!A = - A\right\}$.
	\end{enumerate}
\end{rems}

\subsection{Complements for the proof of theorem \ref{theoimagenonsingularpointtrivstab}}

In order for the proof of theorem \ref{theoimagenonsingularpointtrivstab} to be fully complete, we will justify the following fact, keeping the notations that we used in the above-mentioned proof. In a second time, we will also state Nakayama's Lemma and some of its consequences (lemma \ref{lemnakayamaandconsq}). 

\begin{lem} \label{lemcomplementprooftheofreeactionnonsingpoint} Keep the notations of the proof of theorem \ref{theoimagenonsingularpointtrivstab}. The $\mathcal{P}(G)$-module $I/I^2$ is a free module of finite rank and there exist polynomial functions $f_1,\ldots,f_r \in I^G$ such that the family $\left(\overline{f_1},\ldots,\overline{f_r}\right)$ is a $\mathcal{P}(G)$-basis of $I/I^2$.
\end{lem}

This result uses Lemma 8.39 of \cite{Mukai} that we state and show here for the sake of completeness:

\begin{lem} \label{lem839mukai} Let $R$ be a commutative ring containing, as a subring, a field $k$ and let $M$ be an $R$-module. Equip the tensor product $R \otimes_k M$ with the $R$-module structure given by, if $\alpha \in R$, $a \in R$ and $m \in M$, $\alpha \cdot (a \otimes m) = (\alpha \otimes 1) \cdot (a \otimes m) := (\alpha a) \otimes m$. Suppose that there exists a morphism of $R$-modules $\Phi : M \rightarrow R \otimes_k M$ such that 
\begin{enumerate}
	\item if $\mathcal{M}$ denotes the $k$-linear map $R \otimes_k M \rightarrow M$ given by $\mathcal{M}(a \otimes m) := a \cdot m$ if $a \in R$ and $m \in M$, we have $\mathcal{M} \circ \Phi = {\rm Id}_M$, 
	\item if $\mathcal{J}$ denotes the $k$-linear map $R \otimes_k M \rightarrow R \otimes_k R \otimes_k M$ given by $\mathcal{J}(a \otimes m) := a \otimes 1 \otimes m$ if $a \in R$ and $m \in M$, we have $\mathcal{J} \circ \Phi = \left({\rm Id}_M \otimes \Phi \right) \circ \Phi$. 
\end{enumerate}
Then any basis of the $k$-vector subspace $M_0 := \left\{ m \in M~|~\Phi(m) = 1 \otimes m\right\}$ of $M$ is a free basis of the $R$-module $M$. In particular, $M$ is a free $R$-module.
\end{lem}

\begin{proof} We follow the proof of Lemma 8.39 of \cite{Mukai}. First, let $\left(n_i\right)_{i \in I}$ be a $k$-basis of $M_0$, let $i_1,\ldots,i_r$ be distinct indices of $I$ and let $a_{i_1},\ldots, a_{i_r}$ be elements of $R$ such that $\sum_{l=1}^r a_{i_l} \cdot n_{i_l} = 0_M$. We have
$$0_{R} \otimes 0_M = \Phi\left(\sum_{l=1}^r a_{i_l} \cdot n_{i_l}\right) = \sum_{j=1}^r \left(a_{i_l} \otimes 1\right) \cdot \Phi\left(n_{i_l}\right) =  \sum_{l=1}^r \left(a_{i_l} \otimes 1\right) \cdot \left( 1 \otimes n_{i_l}\right) =  \sum_{l=1}^r a_{i_l} \otimes n_{i_l}$$
so that, since $\left(n_i\right)_{i \in I}$ is a $k$-basis of $M_0$, $a_{i_1} = \ldots = a_{i_r} = 0_R$.

Now, consider an element $m$ of $M$. If $\left(\alpha_{j}\right)_{j \in J}$ is a $k$-basis of $R$, there exist indices $j_1,\ldots,j_s$ of $J$ and $m_{j_1},\ldots,m_{j_s} \in M$ such that $\Phi(m) = \sum_{l=1}^s \alpha_{j_l} \otimes m_{j_l}$. By the second hypothesis on $\Phi$, we have 
$$\sum_{l=1}^s \alpha_{j_l} \otimes 1 \otimes m_{j_l} = \mathcal{J}\left(\Phi(m)\right) = \left({\rm Id}_M \otimes \Phi \right) \circ \Phi(m) = \sum_{l=1}^s \alpha_{j_l} \otimes \Phi\left(m_{j_l}\right)$$
so that, since $\left(\alpha_{j}\right)_{j \in J}$ is a $k$-basis of $R$, $\Phi\left(m_{j_l}\right) = 1 \otimes m_{j_l}$ for any $l \in \{1,\ldots,s\}$, i.e. $m_{j_l} \in M_0$ for any $l \in \{1,\ldots,s\}$. On the other hand, by the first hypothesis on $\Phi$, we have
$$m = \mathcal{M}\left(\Phi(m)\right) = \sum_{l=1}^k  \alpha_{j_l} \cdot m_{j_l}.$$

As a consequence, the family $\left(n_i\right)_{i \in I}$ is a free basis of the $R$-module $M$.
\end{proof}

\begin{proof}[Proof of lemma \ref{lemcomplementprooftheofreeactionnonsingpoint}] We are going to construct a map $\Phi : I/I^2 \mapsto \mathcal{P}(G) \otimes_{\R} I/I^2$ that satisfies the hypotheses of lemma \ref{lem839mukai}, considering the structure of $\mathcal{P}(G)$-module of $I/I^2$ defined in the proof of theorem \ref{theoimagenonsingularpointtrivstab}. The latter gives rise to the $\R$-linear map $\mathcal{M} : \mathcal{P}(G) \otimes_{\mathbb{R}} I/I^2\rightarrow I/I^2$ given by $\mathcal{M}(\varphi \otimes \overline{f}) := \varphi \cdot \overline{f}$ if $\varphi \in \mathcal{P}(G)$ and $f \in I$ (we keep the notation of lemma \ref{lem839mukai}).

Let $F \in \mathcal{P}(X)$ and write $\alpha^*(F) = F \circ \alpha = \sum_{i=1}^k \varphi_i \otimes H_i$ with $\varphi_1,\ldots,\varphi_k\in \mathcal{P}(G)$ linearly independent over $\R$ and $H_1,\ldots,H_k \in \mathcal{P}(X)$. If $F \in I$ and $y \in Y$, notice that, for all $g \in G$, $0 = F(\alpha(g,y)) = \sum_{i=1}^k \varphi_i(g) H_i(y)$, i.e. $\sum_{i=1}^k H_i(y) \varphi_i$ is the zero function of $\mathcal{P}(G)$. Since the polynomials functions $\varphi_1,\ldots,\varphi_k$ on $G$ are linearly independent over $\R$, we therefore have $H_1(y) = \cdots = H_k(y) = 0$ and then $H_1,\ldots,H_k \in I$. As a consequence, $F \circ \alpha \in \mathcal{P}(G) \otimes I$. Remark also that if $F \in I^2$, we have $F \circ \alpha \in \mathcal{P}(G) \otimes I^2$ since $\alpha^*$ is a morphism of $\R$-algebras.

Now, consider the map $\mathcal{A} : I \mapsto \mathcal{P}(G) \otimes_{\R} \mathcal{P}(G) \otimes_{\R} I$, resp. $\mathfrak{A} : \mathcal{P}(Y) \mapsto \mathcal{P}(G) \otimes_{\R} \mathcal{P}(G) \otimes_{\R} \mathcal{P}(Y)$ which associates to any polynomial map $F \in I$, resp. $h \in \mathcal{P}(Y)$, its composition with the polynomial map $\alpha(\mu( \bullet ,\omega(\star)), *)$ from $G \times G \times X$ to $X$, resp. from $G \times G \times Y$ to $Y$. The map~$\mathcal{A}$ induces a map $\overline{\mathcal{A}} : I/I^2 \mapsto \mathcal{P}(G) \otimes_{\R} \mathcal{P}(G) \otimes_{\R} I/I^2$: if $F \in I$ and $F \circ \alpha = \sum_{i=1}^k \varphi_i \otimes H_i$ with $\varphi_1,\ldots,\varphi_k\in \mathcal{P}(G)$ and $H_1,\ldots,H_k \in I$, we have $\mathcal{A}(F) = \sum_{i=1}^k \varphi_i \circ \mu(\bullet, \omega(\star)) \otimes H_i$ and
$$\overline{\mathcal{A}}\left(\overline{F}\right) = \sum_{i=1}^k \varphi_i \circ \mu(\bullet, \omega(\star)) \otimes \overline{H_i}$$
(tensor product, over a commutative field $\mathbb{K}$, with a $\mathbb{K}$-vector space is an exact functor). If furthermore $\mathfrak{P}$ denotes the map $\mathcal{P}(G) \rightarrow \mathcal{P}(G) \otimes_{\R} \mathcal{P}(G) \otimes_{\R} \mathcal{P}(G)$ which associates to any polynomial map $\varphi \in \mathcal{P}(G)$ its composition with $\mu(\mu( \bullet ,\omega(\star)), *)$, remark that we have the equality, if $h \in \mathcal{P}(Y)$ and $\varphi := h \circ \alpha_x \in \mathcal{P}(G)$,
$$\overline{\mathcal{A}}\left(\varphi \cdot \overline{F}\right) = \overline{\mathcal{A}}\left(h \cdot \overline{F}\right) = \mathfrak{A}(h) \cdot  \overline{\mathcal{A}}\left(\overline{F}\right) = \mathfrak{P}(\varphi) \cdot  \overline{\mathcal{A}}\left(\overline{F}\right)$$
where, if $\phi,\phi',\psi,\psi' \in \mathcal{P}(G)$ and $p \in \mathcal{P}(Y)$, $H \in I$, 
$$(\phi' \otimes \psi' \otimes p) \cdot \left(\phi \otimes \psi \otimes \overline{H}\right) = (\phi' \cdot \phi) \otimes (\psi' \cdot \psi) \otimes \left(p \cdot \overline{H}\right).$$
We then denote by $\Phi$ the composition $\left({\rm Id}_{\mathcal{P}(G)} \otimes \mathcal{M}\right) \circ \overline{\mathcal{A}} : I/I^2 \mapsto \mathcal{P}(G) \otimes_{\R} I/I^2$: we are going to check that $\Phi$ satisfies the hypotheses of lemma \ref{lem839mukai}. 

First, with the above notations, we have
$$\Phi\left(\overline{F}\right) = \sum_{i=1}^k \left(\varphi_i \circ \mu(\bullet, \omega(\star))\right) \cdot \left( 1 \otimes \overline{H_i}\right),$$
where, if $\phi,\phi' \in \mathcal{P}(G)$ and $\psi \in \mathcal{P}(G)$, $H \in I$,
$$(\phi' \otimes \psi) \cdot \left(\phi \otimes \overline{H}\right) = (\phi' \cdot \phi)  \otimes \left(\psi \cdot \overline{H}\right),$$
and remark that
\begin{eqnarray*}
\mathcal{M} \circ \Phi\left(\overline{F}\right) & = &  \sum_{i=1}^k \left(\varphi_i \circ \mu(\star, \omega(\star))\right) \cdot \overline{H_i} = \overline{\sum_{i=1}^k \varphi_i(e) H_i}  = \overline{F}.
\end{eqnarray*}

Notice also that, if $\mathfrak{M} : \mathcal{P}(G) \otimes \mathcal{P}(G) \rightarrow \mathcal{P}(G)$ is the map given by $\mathfrak{M}(\phi \otimes \psi) := \phi \cdot \psi$ if $\phi,\psi \in \mathcal{P}(G)$, then, for any element $\pi \in \mathcal{P}(G) \otimes_{\R} \mathcal{P}(G) \otimes_{\R} I/I^2$ and any element $\mathfrak{p} \in \mathcal{P}(G) \otimes_{\R} \mathcal{P}(G) \otimes_{\R} \mathcal{P}(G)$, 
$$\left({\rm Id}_{\mathcal{P}(G)} \otimes \mathcal{M}\right)(\mathfrak{p} \cdot \pi) = \left({\rm Id}_{\mathcal{P}(G)} \otimes \mathfrak{M}\right)(\mathfrak{p}) \cdot \left({\rm Id}_{\mathcal{P}(G)} \otimes \mathcal{M}\right)(\pi),$$
where, if $\phi,\phi' \in \mathcal{P}(G)$ and $\psi \in \mathcal{P}(G)$, $H \in I$,
$$(\phi' \otimes \psi) \cdot \left(\phi \otimes \overline{H}\right) = (\phi' \cdot \phi)  \otimes \left(\psi \cdot \overline{H}\right).$$
As a consequence, $\varphi \in \mathcal{P}(G)$ and $F \in I$, we have
\begin{eqnarray*}
\Phi\left(\varphi \cdot \overline{F}\right) & = & \left({\rm Id}_{\mathcal{P}(G)} \otimes \mathcal{M}\right)\left(\mathfrak{P}(\varphi) \cdot  \overline{\mathcal{A}}\left(\overline{F}\right)\right)\\
& = & \left({\rm Id}_{\mathcal{P}(G)} \otimes \mathfrak{M}\right)(\mathfrak{P}(\varphi)) \cdot \left({\rm Id}_{\mathcal{P}(G)} \otimes \mathcal{M}\right)\left(\overline{\mathcal{A}}\left(\overline{F}\right)\right).
\end{eqnarray*}
If we furthermore write $\varphi \circ \mu = \sum_{j =1}^l \phi_j \otimes \psi_j$ with $\phi_1,\ldots,\phi_l, \psi_1,\ldots,\psi_l \in \mathcal{P}(G)$, we have $\mathfrak{P}(\varphi) = \sum_{j =1}^l \phi_j \circ \mu(\bullet, \omega(\star)) \otimes \psi_j(*)$ and then
\begin{eqnarray*}
\left({\rm Id}_{\mathcal{P}(G)} \otimes \mathfrak{M}\right)(\mathfrak{P}(\varphi)) & = & \sum_{j =1}^l \phi_j \circ \mu(\bullet, \omega(*)) \cdot (1 \otimes \psi_j(*))\\
& = &  \varphi \circ \mu(\mu(\bullet, \omega(*)),*) = \varphi \otimes 1
\end{eqnarray*}
so that $\Phi\left(\varphi \cdot \overline{F}\right) = \left(\varphi \otimes 1\right) \cdot \Phi\left(\overline{F}\right)$.

Finally, if $F \in I$, if we write $F \circ \alpha = \sum_{i=1}^k \varphi_i \otimes H_i$ as above and if, for $i \in \{1,\ldots,k\}$, $H_i \circ \alpha = \sum_{\iota_i = 1}^{\kappa_i} \varphi_{i,\iota_i} \otimes H_{i,\iota_i}$ with $\varphi_{i,1}, \ldots,\varphi_{\kappa_i} \in \mathcal{P}(G)$ and $\psi_{i,1}, \ldots,\psi_{\kappa_i} \in I$, we have
\begin{eqnarray*}
\left({\rm Id}_{\mathcal{P}(G)} \otimes \Phi\right)\left(\Phi\left(\overline{F}\right)\right) & = & \left({\rm Id}_{\mathcal{P}(G)} \otimes \Phi\right) \left(\sum_{i=1}^k \left(\varphi_i \circ \mu(\bullet, \omega(\star))\right) \cdot \left( 1 \otimes \overline{H_i}\right)\right)\\
& = & \sum_{i=1}^k \left(\varphi_i \circ \mu(\bullet, \omega(\star)) \otimes 1 \right) \cdot \left(1 \otimes \Phi\left(\overline{H_i}\right)\right) \\
& = & \sum_{i=1}^k \left(\varphi_i \circ \mu(\bullet, \omega(\star)) \otimes 1\right) \cdot \sum_{\iota_i = 1}^{\kappa_i} \left(1 \otimes \varphi_{i,\iota_i} \circ \mu(\star, \omega(*))\right) \cdot \left(1 \otimes 1 \otimes \overline{H_{i,\iota_i}}\right)
\end{eqnarray*}
while
\begin{eqnarray*}
\sum_{i=1}^k \varphi_i \circ \mu(\bullet, \omega(\star)) \sum_{\iota_i = 1}^{\kappa_i} \varphi_{i,\iota_i} \circ \mu(\star,\omega(*)) H_{i,\iota_i}(\blacksquare) & = & F \circ \alpha( \mu(\bullet, \omega(\star)), \alpha( \mu(\star,\omega(*)),\blacksquare)) \\
& = & F \circ \alpha(\mu(\bullet, \omega(*)), \blacksquare) \\
& = & \sum_{i=1}^k \varphi_i \circ \mu(\bullet, \omega(*)) H_i(\blacksquare)
\end{eqnarray*}
so that
$$\left({\rm Id}_{\mathcal{P}(G)} \otimes \Phi\right)\left(\Phi\left(\overline{F}\right)\right) = \sum_{i=1}^k \left(\varphi_i \circ \mu(\bullet, \omega(*))\right) \cdot \left( 1 \otimes 1 \otimes \overline{H_i}\right)$$
is the image of $\Phi\left(\overline{F}\right)$ by the map $\mathcal{J} : \mathcal{P}(G) \otimes_{\R} I/I^2 \rightarrow \mathcal{P}(G) \otimes_{\R} \mathcal{P}(G) \otimes_{\R} I/I^2$ which associates to $\psi \otimes \overline{H}$, with $\psi \in \mathcal{P}(G)$ and $H \in I$, the element $\psi \otimes 1 \otimes \overline{H}$.

Consequently, by lemma \ref{lem839mukai}, any basis of the $\R$-vector space 
$$M_0 := \left\{\overline{f} \in I/I^2~|~\Phi\left(\overline{f}\right) = 1 \otimes \overline{f}\right\}$$
is a free basis of the $\mathcal{P}(G)$-module $I/I^2$. In order to conclude, let us prove that the $\R$-vector space $M_0$ is finite-dimensional and that we can find an $\R$-basis of $M_0$ with representatives in~$I^G$.

First, notice that $M_0$ is a $\R$-vector subspace of $\left(I/I^2\right)^G$: with $F$, $\varphi_1,\ldots,\varphi_k$ and $H_1,\ldots,H_k$ as above, if $\overline{F} \in M_0$ then there exist $\rho_1,\ldots,\rho_s \in \mathcal{P}(G)$ and $h_1,\ldots,h_s \in I^2$ such that $\sum_{i=1}^k \left(\varphi_i \circ \mu(\bullet, \omega(\star))\right) \cdot \left( 1 \otimes H_i\right) = 1 \otimes F + \sum_{j=1}^s \rho_j \otimes h_j$ and then, for any $g \in G$,
\begin{eqnarray*}
\Phi\left(g \cdot \overline{F} \right) = \Phi\left(\overline{g \cdot F}\right) & = & \sum_{i=1}^k \left(\varphi_i \circ \mu(\mu(g^{-1},\bullet), \omega(\star))\right) \cdot \left( 1 \otimes \overline{H_i}\right) \\
& = & \overline{\sum_{i=1}^k \left(\varphi_i \circ \mu(\mu(g^{-1},\bullet), \omega(\star))\right) \cdot \left( 1 \otimes H_i\right)} \\
& = & \overline{1 \otimes F + \sum_{j=1}^s \rho_j\left(\mu(g^{-1},\bullet)\right) \otimes h_j} \\
& = & \overline{1 \otimes F} = 1 \otimes \overline{F} = \Phi\left(\overline{F} \right),
\end{eqnarray*}
so that $g \cdot \overline{F} = \mathcal{M} \circ \Phi\left(g \cdot \overline{F} \right) = \mathcal{M} \circ \Phi\left(\overline{F}\right) = \overline{F}$.

Furthermore, the $\R$-vector space $\left(I/I^2\right)^G$ is finite-dimensional. Indeed, consider the $\R$-linear map $\iota_G : I^G \rightarrow \left(I/I^2\right)^G$ induced by the quotient map $\iota : I \rightarrow I/I^2$, then $\iota_G$ is surjective by the linear reductivity of $G$: take $\overline{f} \in \left(I/I^2\right)^G$ with $f \in I$, then the line $W$ generated by~$\overline{f}$ and the subspace $V_f$ of $I$ generated by the orbit of $f$ (see the proof of proposition \ref{proplinearizgrealalgset}) are both polynomial representations of $G$ and the restriction $V_f \rightarrow W$ of $\iota$ is a surjective equivariant linear map, so that we can apply characterization 2 of proposition \ref{propequivlinred}. We can thereafter consider the surjective map $\mathfrak{m}_{\pi(x)}/\mathfrak{m}_{\pi(x)}^2 \cong I^G/\left(I^G\right)^2 \rightarrow \left(I/I^2\right)^G$ induced by $\iota_G$ (recall that $\mathfrak{m}_{\pi(x)} \cong \mathfrak{m}_{x}^G = I^G$: see the proof of theorem \ref{theoimagenonsingularpointtrivstab}) and, since the $\R$-vector space $\mathfrak{m}_{\pi(x)}/\mathfrak{m}_{\pi(x)}^2$ is finite-dimensional (recall that it is isomorphic to the dual of the Zariski tangent space of $W$ at~$\pi(x)$: we refer to Theorem 2.1 of Chapter 2 section 1.3 of \cite{Shaf} whose proof also works over real numbers), so is $\left(I/I^2\right)^G$.

As a consequence, the $\R$-vector space $M_0$ is finite-dimensional as well and there exists an $\R$-basis $\left(\overline{f_1}, \ldots,\overline{f_r}\right)$ of $M_0$ with $f_1,\ldots,f_r \in I^G$ (the map $\iota_G : I^G \rightarrow \left(I/I^2\right)^G$ is surjective): the family $\left(\overline{f_1}, \ldots,\overline{f_r}\right)$ is therefore a free basis of the $\mathcal{P}(G)$-module $I/I^2$ by lemma \ref{lem839mukai}.
\end{proof}

\begin{lem} \label{lemnakayamaandconsq} Let $R$ be a commutative ring, $I$ be a proper ideal of $R$ and $M$ be a finitely generated $R$-module, and consider the submodule $I M$ of $M$. 
\begin{enumerate}
	\item Suppose that $IM = M$. Then there exists $b \in 1+I$ such that $b M = \{0_M\}$.
	\item Let $N$ be a submodule of $M$ such that $M = IM+N$. Then there exists $b \in 1+I$ such that $bM \subset N$.
	\item Let $\phi : N \rightarrow M$ be a morphism of $R$-modules such that the composition $\overline{\phi} : N \rightarrow M/IM$ of $\phi$ with the quotient morphism $M \rightarrow M/IM~;~x \mapsto \overline{x}$ is surjective. Then there exists $b \in 1+I$ such that the morphism of $R_b$-modules $\phi_b : N_b \rightarrow M_b$ which associates to any fraction $\frac{y}{b^k}$, with $y \in N$ and $k \in \mathbb{N}$, the fraction $\frac{\phi(y)}{b^k}$ is surjective.
	\item Let $r \in \Nstar$ and $x_1,\ldots,x_r \in M$ such that the classes $\overline{x_1},\ldots, \overline{x_r} \in M/IM$ generate the $R$-module $M/IM$. Then there exists $b \in 1+I$ such that $x_1,\ldots,x_r$ generate the $R_b$-module~$M_b$.
\end{enumerate}
\end{lem}

\begin{proof} 
\begin{enumerate}
	\item This first statement is standard Nakayama's Lemma: see for instance Proposition 2.6.7 of \cite{TY} for a proof.
	\item The equality $M = IM + N$ implies the equality $I (M/N) = M/N$ and we can therefore apply previous statement to the finitely generated $R$-module $M/N$: there exists $b \in 1+I$ such that $b (M/N) = \left\{\overline{0_M}\right\}$ i.e. $bM \subset N$.
	\item The surjectivity of the morphism $\overline{\phi}$ means that $M = {\rm Im} \, \phi + IM$.  As a consequence, by above property 2, there exists $b \in 1+I$ such that $b M \subset {\rm Im} \, \phi$: if $x \in M$ and $k \in \mathbb{N}$, there exists $y \in N$ such that $bx = \phi(y)$ and then 
	$$\frac{x}{b^k} = \frac{bx}{b^{k+1}} = \frac{\phi(y)}{b^{k+1}} \in {\rm Im} \, \phi_b.$$
	\item The fact that the elements $\overline{x_1},\ldots, \overline{x_r}$ generate the $R$-module $M/IM$ means that the morphism of $R$-modules $\overline{\phi} : R^r \rightarrow M/IM~;~(a_1,\ldots,a_r) \mapsto \sum_{i=1}^r a_i \overline{x_i}$, induced by the morphism $\phi : R^r \rightarrow M~;~(a_1,\ldots,a_r) \mapsto \sum_{i=1}^r a_i x_i$, is surjective. Consequently, by 3, there exists $b \in 1+I$ such that the morphism $\phi_b$ is surjective, meaning that the elements $x_1,\ldots,x_r$ of $M$ generate the $R_b$-module $M_b$.
\end{enumerate}
\end{proof}

\subsection{Arc-symmetric sets and $\mathcal{AS}$-sets} \label{subsectarcsymassets}

In the final part of this appendix, we chose to consider notions, properties and results of the papers \cite{KurAR} and \cite{KurPar} around the notion of arc-symmetric set which was first introduced by K. Kurdyka in \cite{KurAR}. We will principally follow the outline of \cite{KurPar} and give the detailed proofs of most statements.

\subsubsection{Arc-symmetric sets of real affine spaces}

We begin with the notion of arc-symmetric subset of a real affine space. Let $n \in \Nstar$.

\begin{de} \label{defarcsymsetrn} Let $S$ be a subset of $\R^n$. We say that $S$ is an \emph{arc-symmetric set of $\R^n$} if for any analytic arc $\gamma : ]-1;1[ \rightarrow \R^n$, if $\gamma(]-1;0[) \subset S$, then $\gamma(]-1;1[) \subset S$.
\end{de}

\begin{ex} \label{exrealalgsetisarcsym} Let $f : \R^n \rightarrow \R$ be an analytic function and denote $X := f^{-1}\left(\{0\}\right)$. Then $X$ is arc-symmetric in $\R^n$. Indeed, let $\gamma : ]-1;1[ \rightarrow \R^n$ be an analytic arc such that $\gamma(]-1;0[) \subset X$, i.e. such that the analytic function $f \circ \gamma : ]-1;1[ \rightarrow \R$ vanishes on $]-1;0[$: by the identity theorem (or equivalently the principle of isolated zeros) for analytic functions, the function~$f \circ \gamma$ then vanishes on the entire interval $]-1;1[$ meaning that $\gamma(]-1;1[) \subset X$. In particular, any real algebraic subset of $\R^n$ is arc-symmetric in $\R^n$ (recall that any real algebraic set can be written as the vanishing locus of a single polynomial, using sums of squares).
\end{ex}

Arc-symmetric sets of $\R^n$ can also be characterized as follows:

\begin{lem} \label{lemequivalentchararcsymrn} Let $S$ be a subset of $\R^n$. The following properties are equivalent :
\begin{enumerate}
	\item $S$ is arc-symmetric in $\R^n$,
	\item for any analytic arc $\gamma : ]\epsilon ; \epsilon'[ \rightarrow \R^n$, if there exists $\epsilon'' \in ]\epsilon ; \epsilon'[$ such that $\gamma(]\epsilon; \epsilon''[) \subset S$ or $\gamma(]\epsilon''; \epsilon'[) \subset S$, then $\gamma(]\epsilon; \epsilon'[) \subset S$,
	\item for any analytic arc $\gamma : ]-1;1[ \rightarrow \R^n$, if the interior of $\gamma^{-1}(S)$ in $]-1;1[$ is non-empty, then $\gamma(]-1;1[) \subset S$,
	\item for any analytic arc $\gamma : ]-1;1[ \rightarrow \R^n$, if $\gamma(]-1;0[) \subset S$, then there exists $\epsilon \in ]0;1[$ such that $\gamma(]-1;\epsilon[) \subset S$.
\end{enumerate}
\end{lem}

\begin{proof} 
1) $\Rightarrow$ 2) : Suppose that $S$ is arc-symmetric in $\R^n$ and let $\gamma : ]\epsilon ; \epsilon'[ \rightarrow \R^n$ be an analytic arc. Let $\epsilon'' \in ]\epsilon ; \epsilon'[$ and denote by $\varphi$ the restriction $]-1;1[ \rightarrow ]\epsilon ; \epsilon'[$ of the increasing affine map $\R \rightarrow \R$ which sends $]-1;1[$ onto $]\epsilon ; \epsilon'[$, as well as $\xi := \varphi^{-1}(\epsilon'') \in ]-1;1[$. Now, let $\psi$ be the increasing analytic diffeomorphism $]-1;1[ \rightarrow ]-1;1[~;~x \mapsto \frac{x + \xi}{\xi x+1}$: if $\gamma(]\epsilon; \epsilon''[) \subset S$, the composition $\widetilde{\gamma} := \gamma \circ \varphi \circ \psi : ]-1;1[ \rightarrow \R^n$ is then an analytic arc such that $\widetilde{\gamma}(]-1;0[) = \gamma(\varphi(]-1;\xi[)) = \gamma(]\epsilon; \epsilon''[) \subset S$, so that $\gamma(]\epsilon; \epsilon'[) = \widetilde{\gamma}(]-1;1[) \subset S$ since $S$ has been supposed arc-symmetric. If $\gamma(]\epsilon''; \epsilon'[) \subset S$, consider the further map $\phi : ]-1;1[ \rightarrow ]-1;1[~;~x \mapsto -x$ and the composition $\overline{\gamma} := \gamma \circ \varphi \circ \psi \circ \phi$ is then an analytic arc $]-1;1[ \rightarrow \R^n$ such that $\overline{\gamma}(]-1;0[) = \gamma \circ \varphi \circ \psi(]0;1[) = \gamma(]\epsilon''; \epsilon'[) \subset S$ so that $\gamma(]\epsilon''; \epsilon'[) = \overline{\gamma}(]-1;1[) \subset S$.

2) $\Rightarrow$ 3) : Suppose that $S$ satisfies property 2) and let $\gamma : ]-1;1[ \rightarrow \R^n$ be an analytic arc such that the interior of $\gamma^{-1}(S)$ in $]-1;1[$ is not empty. There then exists an open interval $]a;b[$, with $-1 \leq a < b \leq 1$, such that $]a;b[ \subset \gamma^{-1}(S)$. If we consider the restricted analytic arc $\gamma_{|]a;1[} : ]a;1[ \rightarrow \R^n$, we have $]a,b[ \subset \gamma_{|]a;1[}^{-1}(S)$ so that, by the assumed property 2), $\gamma(]a;1[) = \gamma_{|]a;1[}(]a;1[) \subset S$. We can then assert, again by 2), that $\gamma(]-1;1[) \subset S$.

3) $\Rightarrow$ 4) : Suppose that $S$ satisfies 3) and let $\gamma : ]-1;1[ \rightarrow \R^n$ be an analytic arc such that $\gamma(]-1;0[) \subset S$ i.e. $]-1;0[ \subset \gamma^{-1}(S)$. In particular, the interior of $\gamma^{-1}(S)$ in $]-1;1[$ is non-empty and then, by 3), $\gamma(]-1;1[) \subset S$.

4) $\Rightarrow$ 1) : Suppose that $S$ satisfies last property 4) and let $\gamma : ]-1;1[ \rightarrow \R^n$ be an analytic arc such that $\gamma(]-1;0[) \subset S$. By 4), there then exists $\epsilon \in ]0;1[$ such that $\gamma(]-1;\epsilon[) \subset S$. Let~$\xi$ denote the supremum of the set $\left\{\varepsilon \in ]0;1[~|~\gamma(]-1;\varepsilon[) \subset S\right\}$ in $\R$: in particular, $\gamma(]-1;\xi[) \subset S$. Assume by absurd that $\xi < 1$ and consider the analytic arc $\gamma \circ \psi : ]-1;1[ \rightarrow \R^n$ where $\psi$ is the above increasing analytic diffeomorphism $]-1;1[ \rightarrow ]-1;1[~;~x \mapsto \frac{x + \xi}{\xi x+1}$: we have $\psi(0) = \xi$ and $\gamma \circ \psi(]-1;0[) = \gamma(]-1;\xi[) \subset S$, so that, by 4), there exists $\eta \in ]0;1[$ such that $\gamma \circ \psi (]-1;\eta[) \subset S$, i.e. $\gamma(]-1;\psi(\eta)[) \subset S$. But $\psi(\eta) > \xi$ ($\psi$ is an increasing function) hence a contradiction. As a consequence, $\xi =1$ and $\gamma(]-1;1[) \subset S$.
\end{proof}

Let us use the previous characterizations to highlight some classes of examples of arc-symmetric sets of real affine spaces:

\begin{exs} \label{exanalyticsubsetarcsym} 
~
\begin{enumerate}
	\item Let $X \subset \R^n$ be an \emph{analytic subset of $\R^n$} i.e. for any $x \in X$, there exist an open neighborhood $U$ of $x$ in $\R^n$ and an analytic function $f : U \rightarrow \R$ such that $X \cap U = f^{-1}(0)$. If $X$ is closed in $\R^n$ (with respect to Euclidean topology), then $X$ is arc-symmetric in~$\R^n$. Indeed, suppose that $X$ is closed and let $\gamma : ]-1;1[ \rightarrow \R^n$ be an analytic arc such that $\gamma(]-1;0[) \subset X$. Since $X$ is closed and $\gamma$ is continuous, the point $x := \gamma(0)$ also belongs to~$X$ and we can then consider an open neighborhood $U$ of $x$ in $\R^n$ and an analytic function $f : U \rightarrow \R$ such that $X \cap U = f^{-1}(0)$. The set $\gamma^{-1}(U)$ being an open neighborhood of $0$ in $]-1;1[$, there exists $\epsilon > 0$ such that $]-\epsilon;\epsilon[ \subset \gamma^{-1}(U)$ i.e. $\gamma(]-\epsilon;\epsilon[) \subset U$. Because $\gamma(]-1;0]) \subset X$, we have $\gamma(]-\epsilon;0]) \subset X \cap U = f^{-1}(0)$ i.e. the analytic function $f \circ \gamma_{|]-\epsilon;\epsilon[} : ]-\epsilon;\epsilon[ \rightarrow \R$ vanishes on $]-\epsilon;0]$, and therefore on the entire interval $]-\epsilon;\epsilon[$ (by the identity theorem for analytic functions) meaning that $\gamma(]-\epsilon;\epsilon[) \subset X \cap U$. As a consequence, $\gamma(]-1;\epsilon[) \subset X$ and $X$ is arc-symmetric by characterization 4 of previous lemma \ref{lemequivalentchararcsymrn}.
	\item The connected components (with respect to Euclidean topology) of real algebraic sets of~$\R^n$ are closed in $\R^n$ and analytic, hence arc-symmetric in $\R^n$ by the previous example. Indeed, let $Y = h^{-1}(0)$ be a real algebraic subset of $\R^n$, with $h$ a polynomial function on $\R^n$, and let~$C$ be a connected component of $Y$ (notice that $C$ is closed in $Y$ which is closed in $\R^n$, hence $C$ is closed in $\R^n$). Let $U$ be an open subset of $\R^n$ containing $C$ such that $U \cap (Y \setminus C) = \emptyset$ ($C$ and $Y \setminus C$ are disjoint closed subsets of $\R^n$ since $Y$ is closed in $\R^n$): we have $C = h_{|U}^{-1}(0)$ and $h_{|U} : U \rightarrow \R$ is an analytic function, so that $C$ is an analytic subset of $\R^n$.
	\item Closed analytic submanifolds of $\R^n$ are also instances of closed analytic subsets of $\R^n$, hence of arc-symmetric sets of $\R^n$. Indeed, let $M$ be an analytic submanifold of $\R^n$ of dimension $d \in \N$ and let $x \in M$: there exists an analytic diffeomorphism $\phi$ from an open neighborhood $U$ of $x$ in $\R^n$ onto an open neighborhood $V$ of the origin in $\R^n$ such that $\phi(x) = \bf{0}$ and 
$$\phi(M \cap U) = \left\{(x_1,\ldots,x_n) \in V~|~x_{d+1} = \ldots = x_n = 0\right\}.$$
If, for $i \in \{d+1,\ldots,n\}$, we denote by $\pi_i : \R^n \rightarrow \R$ the projection on the $i$\textsuperscript{th} coordinate, we have $M \cap U = f^{-1}(0)$ where $f$ is the analytic function $U \rightarrow \R~;~y \mapsto \left(\pi^2_{d+1} + \cdots + \pi^2_{n}\right)(\phi(y))$: $M$ is therefore an analytic subset of $\R^n$.
\end{enumerate}
\end{exs}

All the above arc-symmetric sets are analytic sets but the class of arc-symmetric sets is actually wider than the class of analytic sets: we consider below the example given in \cite{KurAR} (Exemple~1.2~2)) and \cite{KurPar} (Example 2.10) of a closed semialgebraic arc-symmetric set which is not an analytic set. As in the latter references, consider also the following statement:

\begin{lem} \label{lemhomogeneousanalyticisalgebraic} Let $S$ be a (non-empty) subset of $\R^n$. Suppose that $S$ is \emph{homogeneous}, i.e. for all $x \in S$ and $t \in \R$, $t x \in S$ (in particular, the origin of $\R^n$ belongs to $S$), and \emph{analytic at the origin of $\R^n$}, i.e. there exist an open neighborhood $U$ of the origin in $\R^n$ and an analytic function $f : U \rightarrow \R$ such that $S \cap U = f^{-1}(0)$. Then $S$ is an algebraic subset of $\R^n$. 
\end{lem}

\begin{proof} We can assume that, on $U$, $f$ is a (convergent) infinite sum of homogeneous polynomial functions $h_k$, $k \in \Nstar$, such that, for $k \in \Nstar$, ${\rm deg}(h_k) = k$. Let $x \in S \cap U$, let $I$ be an open interval containing $0$ such that for all $t \in I$, $t x \in U$ and consider the analytic function $F : I \rightarrow \R~;~t \mapsto f(tx) = \sum_{k=1}^{\infty} h_k(tx) =  \sum_{k=1}^{\infty} t^k h_k(x)$. Since $S$ is homogeneous, for all $t \in I$, $tx \in S \cap U = f^{-1}(0)$ so that $F$ is the zero function on $I$ i.e. for $k \in \Nstar$, $h_k(x) = 0$. As a consequence, $S \cap U \subset \left(\bigcap_{k=1}^{\infty} V(h_k)\right) \cap U$. Conversely, if $x \in \big(\bigcap_{k=1}^{\infty} V(h_k)\big) \cap U$, then $f$ vanishes at $x$ and then $x \in f^{-1}(0) = S \cap U$, so that 
$$S \cap U = \left(\bigcap_{k=1}^{\infty} V(h_k)\right) \cap U.$$
We are going to show that $S$ is equal to the real algebraic set $Z := \bigcap_{k=1}^{\infty} V(h_k)$. Notice that~$Z$ is also an homogeneous subset of $\R^n$ (since the polynomial functions $h_k$, $k \in \Nstar$, are all homogeneous) and let $y \in S$: there exists $t \in \R^*$ such that $t y \in U$ and we then have $t y \in S \cap U$ (since $S$ is homogeneous). In particular, $t y$ belongs to $Z$ (we have $S \cap U = Z \cap U$) and then $y = \frac{1}{t}(ty) \in Z$ because $Z$ is homogeneous as well. The converse inclusion $Z \subset S$ can be showed by exchanging the roles of $S$ and $Z$ in the previous argument. 
\end{proof}

\begin{ex} \label{excartanumbrella} Consider the real algebraic Cartan umbrella $X := \left\{(x,y,z) \in \R^3~|~z(x^2+y^2) = x^3\right\}$ as well as its closed (with respect to Euclidean topology) semialgebraic canopy 
$$S := \left\{(x,y,z) \in \R^3~|~z(x^2+y^2) = x^3,~(x,y) \neq (0,0)\right\} \cup \left\{(0,0,0)\right\},$$
and let us show that $S$ is an arc-symmetric subset of $\R^3$ which is not analytic at the origin of~$\R^3$ (in particular, $S$ is not an analytic subset of $\R^3$). 

First, let $\gamma : ]-1;1[ \rightarrow \R^3$ be a non-constant analytic arc such that $\gamma(]-1;0[) \subset S \subset X$: since $X$ is a real algebraic subset of $\R^3$, we have $\gamma(]-1;1[) \subset X$ (example \ref{exanalyticsubsetarcsym} 2). Suppose by absurd that $\gamma(]-1;1[)$ is not entirely contained in~$S$ i.e. $\gamma^{-1}(X \setminus S) \neq \emptyset$. The set $X \setminus S =  \left\{(x,y,z) \in \R^3~|~x = y = 0,~z \neq 0\right\}$ is an open subset of $X$ and then, since $\gamma(]-1;1[) \subset X$, $\gamma^{-1}(X \setminus S)$ is a non-empty open subset of $]-1;1[$ and contains in particular an open interval $I \subset ]-1;1[$. But $X \setminus S$ is the set-theoritic difference of the handle $T :=  \left\{(x,y,z) \in \R^3~|~x = y = 0\right\}$ of the Cartan umbrella $X$ and the origin of $\R^3$, so that $I \subset \gamma^{-1}(X \setminus S) \subset \gamma^{-1}(T)$. The interior of $\gamma^{-1}(T)$ is therefore non-empty: since $T$ is an algebraic set, we have $\gamma(]-1;1[) \subset T$. As a consequence, $\gamma(]-1;0[) \subset S \cap T = \left\{(0,0,0)\right\}$ i.e. $\gamma$ is constant equal to $(0,0,0)$ on $]-1;0[$, hence on $]-1;1[$ (since $\gamma$ is analytic), which is in contradiction with our hypothesis on $\gamma$.

Now, let us prove that the set $S$ is not analytic at $(0,0,0)$, using previous lemma \ref{lemhomogeneousanalyticisalgebraic}. Suppose by absurd that $S$ is analytic at the origin of $\R^3$ and remark that $S$ is an homogeneous set: by lemma \ref{lemhomogeneousanalyticisalgebraic}, $S$ is then an algebraic subset of $\R^3$. On the other hand, $S$ writes as the disjoint union of $\{(0,0,0)\}$ and the graph of the regular map $\varphi : \R^2 \setminus \{(0,0)\} \rightarrow ~;~(x,y) \mapsto \frac{x^3}{x^2+y^2}$. As a consequence, $\dim S = \dim \R^2 \setminus \{(0,0)\} = 2$ (see for instance Proposition 2.8.5 (i) and Theorem 2.8.8 of \cite{BCR}) and $S$ is then a proper algebraic subset of the irreducible algebraic set $X$ (the polynomial $z(x^2+y^2)-x^3$ of $\R[x,y,z]$ is irreducible and its sign changes on $\R^3$: see Theorem 4.5.1 of \cite{BCR}) such that $\dim S = \dim X = 2$, hence a contradiction (see for instance Proposition 2.2.5 e) of \cite{AK}).
\end{ex}

\begin{rem} \label{remconncomparcsymarcsym} To conclude this part on examples of arc-symmetric sets of real affine spaces, let us notice that the connected components of an arc-symmetric set of $\R^n$ are arc-symmetric as well. Indeed, let $S \subset \R^n$ be an arc-symmetric set of $\R^n$ and let $(S_i)_{i \in I}$ be the family of connected components of $S$. Let $j \in I$ and let $\gamma : ]-1;1[ \rightarrow \R^n$ be an analytic arc such that $\gamma(]-1;0[) \subset S_j$. Then, since $S_j \subset S$ and $S$ is arc-symmetric, we have $\gamma(]-1;1[) \subset S$ and $\gamma(]-1;1[) = \bigsqcup_{i \in I} \gamma(]-1;1[) \cap S_i$. But, for all $i \in I$, $\gamma(]-1;1[) \cap S_i$ is closed and open in $\gamma(]-1;1[)$ and $\gamma(]-1;1[)$ is connected so that, necessarily, $\gamma(]-1;1[) \subset S_j$.
\end{rem}

Before stating a first stability property of arc-symmetric sets of real affine spaces, let us mention that Lemma 2.1 of \cite{KurPar} (see also Lemme 0.1 of \cite{KurAR}) allows to consider only injective analytic arcs in definition \ref{defarcsymsetrn}. Firstly, notice that the arguments of the proof of lemma \ref{lemequivalentchararcsymrn} are still valid if we consider only injective analytic arcs, so that the considered equivalences of lemma \ref{lemequivalentchararcsymrn} remain true if we replace the words ``analytic arc'' with ``injective analytic arc''.

\begin{lem} \label{lemindefarcsymonlyinjectiveanarcs} A subset $S$ of $\R^n$ is arc-symmetric if and only if for any \emph{injective} analytic arc $\gamma : ]-1;1[ \rightarrow \R^n$, if $\gamma(]-1;0[) \subset S$, then $\gamma(]-1;1[) \subset S$.
\end{lem}

\begin{proof} 
Suppose that $S$ satisfies the latter and let $\gamma : ]-1;1[ \rightarrow \R^n$ be any analytic arc such that $\gamma(]-1;0[) \subset S$. By Lemma 2.1 of \cite{KurPar}, there exist open sub-intervals $I$ and $J$ of $]-1;1[$ containing $0$, an analytic diffeomorphism $\phi : I \rightarrow J$ sending $0$ to $0$, a positive integer $k$ and an injective analytic arc $\eta : J \rightarrow \R^n$ such that, for any $t \in J$, $t^k \in J$ and $\gamma(\phi^{-1}(t)) = \eta(t^k)$. Therefore, if $\delta$ denotes the infimum of $I$ and $\psi$ denotes the map $t \in J \mapsto t^k \in J$, we have
$$\eta(\psi\circ \phi(]\delta;0[)) = \gamma(]\delta;0[) \subset S.$$
But $\psi\circ \phi(]\delta;0[)$ is a non-empty open sub-interval of $J$ and $\eta$ is an injective analytic arc, so that, by hypothesis on $S$, we have $\eta(J) \subset S$. As a consequence,
$$\gamma(I) = \gamma\left(\phi^{-1}(J)\right) = \eta \circ \psi(J) \subset \eta(J) \subset S$$
and then the interior of $\gamma^{-1}(S)$ in $]-1;1[$ is non-empty: $S$ is thus arc-symmetric in $\R^n$ by the third characterization of lemma \ref{lemequivalentchararcsymrn}.
\end{proof}

It is also proven in \cite{KurAR} that, in definition \ref{defarcsymsetrn}, if $S$ is a semialgebraic subset of $\R^n$, we can consider only Nash (i.e. analytic and semialgebraic) arcs $]-1;1[ \rightarrow \R^n$ (see Remarque 1.1.1 and Lemme 2.9 of \cite{KurAR}). Remark that, if $S$ is a semialgebraic subset of $\R^n$, the equivalences of lemma \ref{lemequivalentchararcsymrn} remain valid if we consider only Nash arcs, since the analytic diffeomorphisms considered in the proof have semialgebraic graph. Moreover, in definition \ref{defarcsymsetrn}, we can consider only injective Nash arcs:

\begin{lem} \label{lemcharsaarcsymrninjnasharcs} Let $S$ be a semialgebraic subset of $\R^n$. The set $S$ is arc-symmetric in $\R^n$ if and only if for any \emph{injective Nash} arc $\gamma : ]-1;1[ \rightarrow \R^n$, if $\gamma(]-1;0[) \subset S$, then $\gamma(]-1;1[) \subset S$.
\end{lem}

\begin{proof} Suppose that $S$ satisfies the latter and let $\gamma : ]-1;1[ \rightarrow \R^n$ be any Nash arc such that $\gamma(]-1;0[) \subset S$. Use the same arguments and notations as in the proof of lemma \ref{lemindefarcsymonlyinjectiveanarcs}: following the construction given in the proof of Lemma 2.1 of \cite{KurPar}, we can suppose that the analytic diffeomorphism $\phi : I \rightarrow J$ has semialgebraic graph. As a consequence, the injective analytic arc $\eta : J \rightarrow \R^n$ is also Nash. Indeed, consider a coordinate function $\eta_i : J \rightarrow \R$ of $\eta$, $i \in \{1,\ldots,n\}$: if $\xi$ denotes the supremum of $J$, the restriction of $\eta_i$ to $]0;\xi[$ is the composition of the $i$\textsuperscript{th} coordinate function of $\gamma \circ \phi^{-1}$ with the function $s \in ]0;\xi[ \mapsto s^{\frac{1}{k}} \in ]0;\xi[$ and is then Nash. Therefore, there exist polynomials $a_0, \ldots, a_d \in \R[x_1,\ldots,x_n]$, with $a_d \neq 0$, such that, for any $s \in ]0;\xi[$,
$$a_d(s) \left(\eta_i(s)\right)^d + \ldots + a_1(s) \eta_i(s) + a_0(s) = 0$$
(see subsection 8.1 and Proposition 8.1.8 of \cite{BCR}). But, since $\eta_i : J \rightarrow \R$ is an analytic function, by the identity theorem for analytic functions, the above equality is satisfied on the entire interval $J$: consequently,~$\eta_i$ is a Nash function by Proposition 8.1.8 of \cite{BCR} and $\eta : J \rightarrow \R^n$ is an injective Nash arc. The rest of the proof now runs as in the proof of lemma \ref{lemindefarcsymonlyinjectiveanarcs}: by hypothesis on $S$, we have $\eta(J) \subset S$ and $\gamma(I) \subset S$, and therefore $S$ is arc-symmetric in $\R^n$.
\end{proof}

Before restricting our attention to semialgebraic arc-symmetric sets, let us state the following property that make the closed (for Euclidean topology) arc-symmetric subsets of $\R^n$ into the closed sets of a topology on $\R^n$:

\begin{prop} \label{propclosedarcsymtop} Any intersection or finite union of closed arc-symmetric subsets of $\R^n$ is a closed arc-symmetric subset of $\R^n$.
\end{prop}

\begin{proof} Let $S,T \subset \R^n$ be closed (with respect to Euclidean topology) arc-symmetric in $\R^n$ and let $\gamma : ]-1;1[ \rightarrow \R^n$ be an analytic arc such that $\gamma(]-1;0[)$ is contained in the closed set~$S \cup T$. If the interior of $\gamma^{-1}(S)$ in $]-1;1[$ is non-empty, then $\gamma(]-1;1[) \subset S \subset S \cup T$ since~$S$ is arc-symmetric. Suppose that the interior of $\gamma^{-1}(S)$ is empty, then the interior of $\gamma^{-1}(T)$ in $]-1;1[$ is necessarily non-empty. Indeed, if the interior of $\gamma^{-1}(T)$ was also empty, then the interior of $\gamma^{-1}(S \cup T) = \gamma^{-1}(S) \cup \gamma^{-1}(T)$ would be empty as well because $S$ and $T$ (and then $\gamma^{-1}(S)$ and $\gamma^{-1}(T)$) are closed, but $]-1;0[ \subset \gamma^{-1}(S \cup T)$. As a consequence, since $T$ is arc-symmetric, $\gamma(]-1;1[) \subset T \subset S \cup T$.

Now, let $\left(S_i\right)_{i \in I}$ be a family of closed arc-symmetric subsets of $\R^n$ and let $\gamma : ]-1;1[ \rightarrow \R^n$ be an analytic arc such that $\gamma(]-1;0[)$ is contained in the closed set $\bigcap_{i \in I} S_i$ i.e. for all $i \in I$, $\gamma(]-1;0[) \subset S_i$. Since all the sets $S_i$, $i \in I$, are arc-symmetric, we have $\gamma(]-1;1[) \subset S_i$ for all $i \in I$ i.e. $\gamma(]-1;1[) \subset \bigcap_{i \in I} S_i$. 
\end{proof}

\begin{rem} If $S$ is an arc-symmetric set of $\R^n$ and $T$ is an arc-symmetric set of an affine space $\R^m$, then the Cartesian product $S \times T$ is an arc-symmetric set of $\R^{n+m}$. Indeed, if $\gamma : ]-1;1[ \rightarrow \R^{n+m}$ is an analytic arc such that $\gamma(]-1;0[) \subset S \times T$ and if we write $\gamma = (\gamma_1,\gamma_2)$ with $\gamma_1 : ]-1;1[ \rightarrow \R^n$ and $\gamma_2 : ]-1;1[ \rightarrow \R^m$, we have $\gamma_1(]-1;0[) \subset S$ and $\gamma_2(]-1;0[) \subset T$ and then, since $S$ and $T$ are both arc-symmetric, $\gamma_1(]-1;1[) \subset S$ and $\gamma_2(]-1;1[) \subset T$ i.e. $\gamma(]-1;1[) \subset S \times T$.
\end{rem}

\subsubsection{Semialgebraic arc-symmetric sets of real affine spaces and the $\mathcal{AR}$ topology}

From now on, we will focus on \emph{semialgebraic} arc-symmetric sets of $\R^n$. First, remark that all semialgebraic arc-symmetric sets are closed (this is Remarque 1.3 of \cite{KurAR} and Proposition 2.8 of \cite{KurPar}):

\begin{prop} \label{propsemialgarcsymisclosed} Let $S$ be a semialgebraic arc-symmetric subset of $\R^n$. Then $S$ is closed with respect to the Euclidean topology of $\R^n$.
\end{prop}

\begin{proof} Let $x$ be a point of the Euclidean closure $\overline{S}$ of $S$. Since $S$ is semialgebraic, by Nash curve selection lemma (Proposition 8.1.13 of \cite{BCR}), there exists a Nash arc $\gamma : ]-1;1[ \rightarrow \R^n$ such that $\gamma(0) = x$ and $\gamma(]0;1[) \subset S$. Because $S$ is arc-symmetric in $\R^n$, we have $\gamma(]-1;1[) \subset S$, in particular $x = \gamma(0) \in S$.
\end{proof}

One part of this section will be dedicated to the proof of the following result: the semialgebraic arc-symmetric subsets of $\R^n$ are the closed sets of a Noetherian topology on $\R^n$ (Théorème~1.4 of \cite{KurAR} and Theorem 2.12 of \cite{KurPar}).

\begin{theo} \label{theotopar} Any intersection or finite union of semialgebraic arc-symmetric subsets of~$\R^n$ is a semialgebraic arc-symmetric subset of $\R^n$. Furthermore, the topology on $\R^n$ whose closed sets are exactly the semialgebraic arc-symmetric subsets of $\R^n$ is a Noetherian topology.
\end{theo}

Since any semialgebraic arc-symmetric subset of $\R^n$ is closed by previous proposition \ref{propsemialgarcsymisclosed}, any intersection or finite union of semialgebraic arc-symmetric subsets of $\R^n$ is arc-symmetric in $\R^n$ by proposition \ref{propclosedarcsymtop}. Because any finite union of semialgebraic sets is (by definition) a semialgebraic set, the proof of the first statement of theorem \ref{theotopar} amounts to show that any intersection of semialgebraic arc-symmetric subsets of $\R^n$ is a semialgebraic subset of $\R^n$. We are actually going to see that an intersection of semialgebraic arc-symmetric subsets of $\R^n$ is necessarily a finite one: 

\begin{lem} \label{lemprooftopar} Let $M$ be a connected Nash submanifold of $\R^n$.
	\begin{enumerate}
		\item If $S$ is a semialgebraic arc-symmetric subset of $\R^n$, then $M \subset S$ or $\dim (M \cap S) < \dim M$.
		\item If $\left(S_i\right)_{i \in I}$ is a family of semialgebraic arc-symmetric subsets of $\R^n$, then there exist $i_1,\ldots,i_k \in I$ such that
		$$M \cap \bigcap_{i \in I} S_i = M \cap S_{i_1} \cap \cdots \cap S_{i_k}.$$	
	\end{enumerate}
\end{lem}

\begin{proof} \begin{enumerate}
	\item We follow the proofs of Lemme 1.5 of \cite{KurAR} and Lemma 2.14 of \cite{KurPar}. Denote by~$d$ the dimension of $M$ (as a semialgebraic set and a submanifold of $\R^n$: see Proposition~2.8.14 of \cite{BCR}) and suppose that $d = \dim (M \cap S)$. Then, by remark \ref{remafterlemprooftopar} below, there exists a non-empty open subset of $M$ which is contained in $M \cap S$ i.e. the interior $U$ of $M \cap S$ in~$M$ is non-empty. Let us show that~$U$ is furthermore closed in $M$: let $x$ be an element of the closure of $U$ in $M$ (i.e. an element of $M$ and the Euclidean closure of $U$ in $\R^n$). Since~$M$ is a Nash submanifold of dimension $d$ of $\R^n$, there exists a Nash diffeomorphism~$\phi$ from an open semialgebraic neighborhood $W$ of $x$ in $\R^n$ onto an open semialgebraic neighborhood~$V$ of the origin in $\R^n$ such that $\phi(x) = \bf{0}$ and 
$$\phi(M \cap W) = \left\{(x_1,\ldots,x_n) \in V~|~x_{d+1} = \ldots = x_n = 0\right\}.$$
Identify the set $\left\{(x_1,\ldots,x_n) \in \R^n~|~x_{d+1} = \ldots = x_n = 0\right\}$ with $\R^d$, let $r \in ]0;+\infty[$ such that the open ball $B_n$ of radius $r$ centered at the origin in $\R^n$ is contained in $V$ and denote by $B_d$ the open ball of radius $r$ centered at the origin in $\R^d$ ($B_d$ is the intersection of $B_n$ with~$\R^d$). We have 
$$\phi^{-1}(B_d) = \phi^{-1}(\R^d \cap B_n) = M \cap \phi^{-1}(B_n),$$
so that $\phi^{-1}(B_d)$ is an open subset of $M$ containing $x$ and then $\phi^{-1}(B_d) \cap U$ is non-empty~($x$ belongs to the closure of $U$ in $M$).

Let $x_0 \in \phi^{-1}(B_d) \cap U$, denote $y_0 := \phi(x_0) \in B_d$ and let $y$ be any point of $B_d$. Consider the analytic arc $\gamma : ]-\epsilon;1+\epsilon[ \rightarrow \R^d~;~(1-t) y_0 + t y$ with $\epsilon \in ]0;+\infty[$ such that for all $t \in ]-\epsilon;1+\epsilon[$, $\gamma(t) \in B_d$ ($B_d$ is convex and open in $\R^d$). Since $\phi^{-1}(B_d) \cap U$ contains $x_0$ and is the intersection of $M$ with an open set of $\R^n$ contained in $W$, $\phi(\phi^{-1}(B_d) \cap U)$ is an open neighborhood of $y_0$ in $\R^d$. In particular, there exist $\epsilon', \epsilon'' \in ]-\epsilon;1+\epsilon[$ with $\epsilon' < \epsilon''$ such that $\gamma(]\epsilon';\epsilon''[) \subset \phi(\phi^{-1}(B_d) \cap U)$ i.e. $\phi^{-1} \circ \gamma (]\epsilon';\epsilon''[) \subset \phi^{-1}(B_d) \cap U$. But $U \subset S$ and $S$ is arc-symmetric: as a consequence, since the arc $\phi^{-1} \circ \gamma : ]-\epsilon;1+\epsilon[ \rightarrow \R^n$ is analytic, we have $\phi^{-1} \circ \gamma( ]-\epsilon;1+\epsilon[) \subset S$.

In particular, $\phi^{-1}(y) \in S$ and then the open subset $\phi^{-1}(B_d)$ of $M$ containing $x$ is contained in $S$. Consequently, $\phi^{-1}(B_d) \subset U$ (recall that $U$ is the interior of $M \cap S$ in $M$) and $x \in U$. As a conclusion, $U \subset M \cap S$ is an open and closed subset of $M$, which is connected, and therefore $U = M \cap S = M$.
\item We follow the proofs of Lemme 1.5 of \cite{KurAR} and Lemma 2.13 of \cite{KurPar}. Reason by induction on the dimension of $M$ and notice that if $\dim M = 0$, then $M$ is a point since $M$ is connected and the statement is therefore true in this case. Now, let $\left(S_i\right)_{i \in I}$ be a family of semialgebraic arc-symmetric subsets of $\R^n$: since the statement is also trivially true if for all $i \in I$, $M \subset S_i$, suppose that there exists $i_0 \in I$ such that $M$ is not entirely contained in $S_{i_0}$. Thanks to previous statement 1, we can then assert that $\dim (M \cap S_{i_0}) < \dim M$. 

Write $M \cap S_{i_0}$ as a finite disjoint union of connected Nash submanifolds $T_1,\ldots,T_l$ of~$\R^n$ (cf. Proposition 2.9.10 of \cite{BCR}): since the dimension of $M \cap S_{i_0}$ is smaller than the dimension of~$M$, so are the dimensions of $T_1,\ldots,T_l$. We can therefore apply the induction hypothesis: for any $j \in \{1,\ldots, l\}$, there exists a finite subset $I_j$ of $I$ such that $T_j  \cap \bigcap_{i \in I} S_i = T_j \cap  \bigcap_{i \in I_j} S_i$. Denoting $J := \bigcup_{j =1}^l I_j$, we then have $T_j  \cap \bigcap_{i \in I} S_i = T_j \cap  \bigcap_{i \in J} S_i$ for any~$j \in \{1,\ldots, l\}$ and, consequently, 
$$M \cap  \bigcap_{i \in I} S_i = M \cap S_{i_0} \cap \bigcap_{i \in I} S_i = \bigcup_{j=1}^l \left(T_j \cap \bigcap_{i \in I} S_i\right) =  \bigcup_{j=1}^l \left(T_j \cap \bigcap_{i \in J} S_i\right) = M \cap S_{i_0} \cap \bigcap_{i \in J} S_i.$$
\end{enumerate}
\end{proof}

\begin{rem} \label{remafterlemprooftopar} If $S$ and $T$ are semialgebraic subsets of $\R^n$ such that $S \subset T$ and $\dim S = \dim T$, then there exists a (non-empty) open subset of $T$ which is contained in $S$ (i.e. the interior of $S$ in $T$ is non-empty). Indeed, as in the proof of proposition \ref{propsemialggroupunionofccofzarclos}, consider a stratification of $T$ adapted to $S$ (cf. Corollary 3.8 of \cite{CosteSA}): $T$ is the finite disjoint union of semialgebraic subsets $C_1,\ldots, C_N$ such that $S$ is a union of some $C_k$'s, $k \in \{1,\ldots,N\}$, and, if $k \in \{1,\ldots,N\}$, the Euclidean closure $\overline{C_k}$ of $C_k$ in $T$ is the (disjoint) union of $C_k$ and of some $C_l$'s, $l \in \{1,\ldots, N\} \setminus \{k\}$, of smaller dimension. Let $k \in \{1,\ldots,N\}$ such that $C_k \subset S$ and $\dim C_k = \dim S = \dim T$, and denote $C:=C_k$. We have $C = T \setminus \left(\bigcup_{1\leq l \leq N,\, l \neq k} \overline{C_l}\right)$ (since for any $l \in \{1,\ldots, N\} \setminus \{k\}$, $\overline{C_l} \setminus C_l$ is a union of strata of dimension smaller than $\dim C = \dim T$) and then $C$ is an open subset of $T$ contained in $S$.
\end{rem}

We now use lemma \ref{lemprooftopar} to prove theorem \ref{theotopar}:

\begin{proof}[Proof of theorem \ref{theotopar}] If $\left(S_i\right)_{i \in I}$ is a family of semialgebraic arc-symmetric subsets of $\R^n$, we can apply lemma \ref{lemprooftopar} with $M = \R^n$ to assert that the intersection $\bigcap_{i \in I} S_i$ is a finite one and is therefore a semialgebraic arc-symmetric subset of $\R^n$. 

In order to show the Noetherianity of the topology whose closed sets are exactly the semialgebraic arc-symmetric subsets of $\R^n$, consider a sequence $\left(T_n\right)_{n \in \mathbb{N}}$ of semialgebraic arc-symmetric subsets of $\R^n$ such that for all $n \in \mathbb{N}$, $T_{n+1} \subset T_n$. By lemma \ref{lemprooftopar}, there exist $n_1,\ldots,n_k \in \mathbb{N}$ with $n_1 < \cdots < n_k$ such that 
$$\bigcap_{n \in \mathbb{N}} T_n = T_{n_1} \cap \cdots \cap T_{n_k} = T_{n_k}.$$
As a consequence, $T_n \supset T_{n_k}$ for any $n \in \N$ and then, if $m \in \mathbb{N}$ satisfies $m \geq n_k$, we have $T_m = T_{n_k}$ (we have $T_m \subset T_{n_k}$ since the sequence $\left(T_n\right)_{n \in \mathbb{N}}$ is decreasing).
\end{proof}

\begin{de} The topology on $\R^n$ whose closed sets are exactly the semialgebraic arc-symmetric subsets of $\R^n$ is called the \emph{$\mathcal{AR}$ topology} on $\R^n$. If $S$ is a subset of $\R^n$, the closure of~$S$ with respect to $\mathcal{AR}$ topology is denoted by $\overline{S}^{\mathcal{AR}}$.
\end{de}

\begin{rem} \label{remeqdimarclosuresasa} If $S$ is a semialgebraic subset of $\R^n$, we have $\dim \overline{S}^{\mathcal{AR}} = \dim S$. Indeed, since any real algebraic set is arc-symmetric (example \ref{exrealalgsetisarcsym}), we have $S \subset \overline{S}^{\mathcal{AR}} \subset \overline{S}^{\mathcal{Z}}$ while $\dim S = \dim \overline{S}^{\mathcal{Z}}$.

\end{rem}

Since the $\mathcal{AR}$ topology of $\R^n$ is Noetherian, any non-empty $\mathcal{AR}$-closed subset of $\R^n$ (by definition, a subset of $\R^n$ is $\mathcal{AR}$-closed if and only if it is a semialgebraic arc-symmetric subset of $\R^n$) has a unique finite decomposition into irreducible $\mathcal{AR}$-closed components that we will call its~\emph{$\mathcal{AR}$-irreducible components}. Furthermore:

\begin{prop} \label{proparirredisconnected} Let $S$ be an irreducible $\mathcal{AR}$-closed set of $\R^n$. Then $S$ is connected with respect to Euclidean topology. Furthermore, the Zariski closure $\overline{S}^{\mathcal{Z}}$ of $S$ in $\R^n$ is Zariski-irreducible.
\end{prop}

\begin{proof} Let $S_1,\ldots,S_k$ be the connected components of the $\mathcal{AR}$-irreducible semialgebraic arc-symmetric set $S$. The sets $S_1,\ldots,S_k$ are semialgebraic (see Theorem 2.4.5 of \cite{BCR}) and arc-symmetric (by remark \ref{remconncomparcsymarcsym}). Since $S$ writes as the disjoint union $\bigsqcup_{i=1}^k S_i$ of $\mathcal{AR}$-closed sets and is $\mathcal{AR}$-irreducible, we have $k=1$ and $S = S_1$ is connected. 

As for the second statement, suppose that $\overline{S}^{\mathcal{Z}} = X_1 \cup X_2$, where $X_1$ and $X_2$ are two real algebraic subsets of $\R^n$. Then $S$ writes as the union of semialgebraic arc-symmetric sets $(X_1 \cap S) \cup (X_2 \cap S)$: since $S$ is $\mathcal{AR}$-irreducible, we have, without loss of generality, $S = X_1 \cap S$ i.e. $S \subset X_1$ and then $\overline{S}^{\mathcal{Z}} \subset X_1$.
\end{proof}

\subsubsection{$\mathcal{AR}$-irreducible components of real algebraic sets}

Let us dedicate one part of this section to the description of the $\mathcal{AR}$-irreducible components of maximal dimension of the real algebraic subsets of $\R^n$. First, consider the nonsingular case:

\begin{lem} \label{lemconncompnonsingrasisarirred} Let $X \subset \R^n$ be a nonsingular real algebraic set. Any connected component of~$X$ is an irreducible $\mathcal{AR}$-closed subset of $\R^n$. As a consequence, the decomposition of $X$ into connected components is the decomposition of $X$ into its $\mathcal{AR}$-irreducible components.
\end{lem}

\begin{proof} Let $C$ be a connected component of $X$ ($C$ is then an $\mathcal{AR}$-closed set of $\R^n$ by example~\ref{exanalyticsubsetarcsym}~2) and suppose that $C$ is a union $S_1 \cup S_2$ of semialgebraic arc-symmetric subsets $S_1$ and $S_2$ of $\R^n$ such that $C \neq S_2$. Since $X$ is nonsingular, $C$ is a Nash submanifold of dimension $d := \dim X$ of $\R^n$ (see Propositions 3.3.11 and 2.8.14 of \cite{BCR}) and then, because $C$ is furthermore connected, we have $\dim S_2 < d$ by lemma \ref{lemprooftopar} 1. The dimension of $S_1$ is therefore necessarily equal to~$d$ and, by the same statement of lemma \ref{lemprooftopar}, $C = S_1$. As a result, $C$ is~$\mathcal{AR}$-irreducible.
\end{proof}

\begin{rem} The statement of lemma \ref{lemconncompnonsingrasisarirred} is not true for singular real algebraic sets: if~$X$ is for instance the Cartan umbrella of example \ref{excartanumbrella}, $X$ is connected and writes as the union of the canopy $S$ and the handle $T$ of the umbrella (cf. example \ref{excartanumbrella}), which are both semialgebraic arc-symmetric in $\R^3$ (recall also that $X$ is Zariski-irreducible).
\end{rem}  

The general case is a consequence of the following result. For its proof, we refer to the proof of Theorem 2.21 of \cite{KurPar}.

\begin{theo} \label{theokurpararirredconncompresofsing} Let $X$ be a real algebraic set of $\R^n$ and $S$ be an irreducible $\mathcal{AR}$-closed set contained in $X$ such that $\dim S = \dim X$. Now, let $\pi : \widetilde{X} \rightarrow X$ be a resolution of singularities of $X$, i.e. a birational regular map from a nonsingular real algebraic set $\widetilde{X}$ to $X$. There exists a unique connected component $\widetilde{S}$ of $\widetilde{X}$ such that 
	\begin{itemize}
		\item $\pi^{-1}(S)$ is the union of $\widetilde{S}$ and of irreducible $\mathcal{AR}$-closed sets of smaller dimension,
		\item $\overline{\pi\left(\widetilde{S}\right)}^{\mathcal{AR}} = S$.
	\end{itemize}
\end{theo}
 
\begin{cor} \label{cortheokurpararirredconncompresofsing} Keep the real algebraic set $X$ and the resolution of singularities $\pi : \widetilde{X} \rightarrow X$ of $X$ of theorem \ref{theokurpararirredconncompresofsing}. There is a bijective correspondence between the $\mathcal{AR}$-irreducible components of $X$ of dimension $\dim X$ and the connected components of $\widetilde{X}$.
\end{cor}

\begin{proof} Let $S_1,\ldots,S_r$ be the $\mathcal{AR}$-irreducible components of $X$ of dimension $d := \dim X$ and, for $i \in \{1,\ldots,r\}$, let $\widetilde{S}_i$ denote the unique connected component of $\widetilde{X}$ such that $\pi^{-1}(S_i)$ is the union of $\widetilde{S}$ and of irreducible $\mathcal{AR}$-closed sets of smaller dimension and $\overline{\pi\left(\widetilde{S}_i\right)}^{\mathcal{AR}} = S_i$ ($\widetilde{S}_i$ is provided by theorem~\ref{theokurpararirredconncompresofsing}). On the one hand, notice that, if $i,j$ are distinct indices of $\{1,\ldots,r\}$, we have $\widetilde{S}_i \neq \widetilde{S}_j$ since, otherwise, we would have $S_i = \overline{\pi\left(\widetilde{S}_i\right)}^{\mathcal{AR}} = \overline{\pi\left(\widetilde{S}_j\right)}^{\mathcal{AR}} = S_j$.

On the other hand, let $C_1,\ldots,C_s$ be the connected components of the nonsingular real algebraic set $\widetilde{X}$ (in particular, $C_1,\ldots,C_s$ are all of dimension $d$ and, by lemma \ref{lemconncompnonsingrasisarirred}, the equality $\widetilde{X} = \bigcup_{j=1}^s C_j$ is the decomposition of $\widetilde{X}$ into its $\mathcal{AR}$-irreducible components). Because $\widetilde{X} = \bigcup_{i=1}^r \pi^{-1}(S_i)$, the set $\widetilde{X}$ writes as the union of the sets $\widetilde{S}_i$, $i \in \{1,\ldots,r\}$, and of irreducible $\mathcal{AR}$-closed sets of dimension smaller than $d$. By a dimension argument and since the sets $\widetilde{S}_1,\ldots,\widetilde{S}_r$ are connected components of $\widetilde{X}$, we can then assert that $\bigcup_{j=1}^s C_j = \bigcup_{i=1}^r \widetilde{S}_i$ and, because the decomposition of $\widetilde{X}$ into its connected components (or irreducible $\mathcal{AR}$-closed components) is unique, we have $s = r$ and, for all $j \in \{1,\ldots,r\}$, there is a unique $i \in \{1,\ldots,r\}$ such that $C_j = \widetilde{S}_i$.
\end{proof}

\begin{ex} Consider again the Cartan umbrella $X$ of example \ref{excartanumbrella}. Blowing-up $\R^3$ at the intersection point $(0,0,0)$ of the canopy $S$ and the handle $T$ of $X$ provides, via the birational regular map $\sigma_1 : (x,y,z) \in \R^3 \mapsto (xz,yz,z) \in \R^3$, the strict transform 
$$X_1 := \left\{(x,y,z) \in \R^3~|~x^2+y^2 = x^3\right\}$$
of $X$. Now, blowing-up $\R^3$ along the singular line $\left\{(x,y,z) \in \R^3~|~x=y=0\right\}$ of $X_1$ provides, via the birational regular (polynomial) map $\sigma_2 : (x,y,z) \in \R^3 \mapsto (x,xy,z) \in \R^3$, the strict transform 
$$X_2 := \left\{(x,y,z) \in \R^3~|~1+y^2 = x\right\}.$$ 
The composition $\pi := \sigma_1 \circ \sigma_2 : X_2 \rightarrow X~;~(x,y,z) \mapsto (xz,xyz,z)$ is then a resolution of singularities of $X$ (the restriction $X_2 \setminus \{z = 0\} \rightarrow X_1 \setminus \{x=0\}$ of $\pi$ is a biregular map with inverse $(x,y,z) \in X_1 \setminus \{x=0\} \mapsto \left(\frac{x}{z}, \frac{y}{x},z\right) \in X_2 \setminus \{z = 0\}$) and the connected nonsingular real algebraic set $X_2$ corresponds to the $\mathcal{AR}$-irreducible canopy $S = \overline{\pi(X_2)}^{\mathcal{AR}}$ of $X$ (we have $\overline{\pi(X_2)}^{\mathcal{AR}} \subset S = \overline{\pi(X_2)} \subset \overline{\pi(X_2)}^{\mathcal{AR}}$), while $\pi(X_2) \cap T = \{(0,0,0)\}$. Since the real algebraic handle $T$ is nonsingular and connected, $T$ is $\mathcal{AR}$-irreducible and $X = S \cup T$ is the decomposition of $X$ into $\mathcal{AR}$-irreducible components.
\end{ex}

Previous statements also imply Lemma 2.25 of \cite{KurPar}, which is the $\mathcal{AR}$ analog of the fact that a proper algebraic subset of an irreducible algebraic set has smaller dimension (see for instance Lemma 2.2.9 of \cite{AK}):

\begin{prop} \label{proparirreddim} Let $S, T \subset \R^n$ be semialgebraic arc-symmetric sets of $\R^n$ such that $S \subset T$.
	\begin{enumerate}
		\item If $S$ is $\mathcal{AR}$-irreducible and $\dim S = \dim T$, then $S$ is an $\mathcal{AR}$-irreducible component of $T$.
		\item If $T$ is $\mathcal{AR}$-irreducible and $\dim S = \dim T$, then $S = T$.
	\end{enumerate}
\end{prop}
 
\begin{proof}[Proof of proposition \ref{proparirreddim}] Suppose that $S$ is $\mathcal{AR}$-irreducible and that $\dim S = \dim T$. First assume that $T$ is a real algebraic subset of $\R^n$ and consider a resolution of singularities $\pi : \widetilde{T} \rightarrow T$ of $T$: by theorem \ref{theokurpararirredconncompresofsing}, there exists a connected component $\widetilde{S}$ of $\widetilde{T}$ such that $\overline{\pi\left(\widetilde{S}\right)}^{\mathcal{AR}} = S$ and, by corollary \ref{cortheokurpararirredconncompresofsing} (and its proof), $\overline{\pi\left(\widetilde{S}\right)}^{\mathcal{AR}}$ is an $\mathcal{AR}$-irreducible component of $T$. If now~$T$ is any $\mathcal{AR}$-closed set of $\R^n$, let $T_0$ be an $\mathcal{AR}$-irreducible component of $T$ containing the $\mathcal{AR}$-irreducible subset $S$ of $T$: in particular, $\dim S = \dim T_0 = \dim T$. Let $X$ denote the Zariski closure of $T_0$. Since $S$ is an $\mathcal{AR}$-irreducible subset of the real algebraic set $X$ and $\dim S = \dim X$, we can assert, by the above case, that $S$ is an $\mathcal{AR}$-irreducible component of $X$ i.e. a maximal element, with respect to inclusion, among $\mathcal{AR}$-irreducible subsets of $X$. Because $S \subset T_0$ and $T_0$ is an $\mathcal{AR}$-irreducible subset of $X$, we therefore have $S = T_0$ and $S$ is then an~$\mathcal{AR}$-irreducible component of $T$.

As for the second statement of the proposition, suppose that $T$ is $\mathcal{AR}$-irreducible and satisfies $\dim S = \dim T$, and let $S_0$ be an $\mathcal{AR}$-irreducible component of $S$ of dimension $\dim S = \dim T$: by the first statement, $S_0$ is then an $\mathcal{AR}$-irreducible component of $T$. Because $T$ is~$\mathcal{AR}$-irreducible, we deduce that $S_0 = T$ and then $T \subset S$.
\end{proof} 
 
\begin{rem} If $S$ is an $\mathcal{AR}$-closed set of $\R^n$ of dimension $d$ and if $S_1$ and $S_2$ are two distinct~$\mathcal{AR}$-irreducible components of $S$ of dimension $d$, then $\dim(S_1 \cap S_2) < d$. Indeed, if $\dim(S_1 \cap S_2)$ was equal to $d$, we would have $S_1 = S_1 \cap S_2 = S_2$ by proposition \ref{proparirreddim} 2.
\end{rem} 
 
\subsubsection{Arc-symmetric sets and arc-analytic maps} 

In this part, we study the links between arc-symmetric sets and arc-analytic maps: let $m \in \Nstar$.
 
\begin{de} \label{defarcanalyticmap} Let $S$ be a subset of $\R^n$, $T$ be a subset of $\R^m$ and $f : S \rightarrow T$ be a map. We say that $f$ is \emph{arc-analytic} if for any analytic arc $\gamma : ]-1;1[ \rightarrow S$, the composition $f \circ \gamma : ]-1;1[ \rightarrow T$ is analytic.
\end{de}
 
\begin{exs}
~
\begin{enumerate}
	\item Any analytic map from an open set of $\R^n$ to a subset of $\R^m$ is arc-analytic. 
	\item The class of arc-analytic maps is larger than the class of analytic maps. For instance, the map 
$$f : \begin{array}{ccc}\R^2 & \rightarrow & \R\\(x,y) & \mapsto & \begin{cases}\frac{x^3}{x^2+y^2} & \mbox{if $(x,y) \neq (0,0)$,}\\0 & \mbox{if $(x,y) = (0,0)$,}\end{cases}\end{array}$$	
is arc-analytic and not differentiable at $(0,0)$. Indeed, notice that $f$ is analytic at any point $(x,y) \neq (0,0)$ and consider a non-identically zero analytic arc $\gamma : ]-1;1[ \rightarrow \R^2$ and $\xi \in ]-1;1[$ such that $\gamma(\xi) = (0,0)$: composing $\gamma$ with an analytic diffeomorphism as we did in the proof of lemma \ref{lemequivalentchararcsymrn}, we can suppose that $\xi = 0$. We can then write, in a neighborhood of $0$, $\gamma(t) = t^{\nu} \eta(t)$ with $\nu \in \N$ and $\eta = (\eta_1,\eta_2)$ analytic such that $\eta(0) \neq(0,0)$, and we have, if $t$ is a non-zero real number in the considered neighborhood of $0$,
	$$f \circ \gamma(t) = \frac{t^{3\nu} (\eta_1(t))^3}{t^{2\nu}\big((\eta_1(t))^2 + (\eta_2(t))^2\big)} = t^{\nu} \frac{(\eta_1(t))^3}{(\eta_1(t))^2 + (\eta_2(t))^2},$$
so that the map $f \circ \gamma$ coincide, on a neighborhood of $0$, with the analytic function defined by the latter expression. On the other hand, $f$ is not differentiable at $(0,0)$: if there existed two real numbers $\alpha$ and $\beta$ such that 
$$\frac{f(x,y) - \alpha x - \beta y}{\sqrt{x^2+y^2}} \underset{(x,y) \rightarrow (0,0)}{\longrightarrow} 0,$$
we would have $\frac{f(0,y) - \beta y}{\sqrt{y^2}} \underset{y \rightarrow 0}{\longrightarrow} 0$ i.e. $\beta = 0$ and $\frac{f(x,0) - \alpha x}{\sqrt{x^2}} \underset{x \rightarrow 0}{\longrightarrow} 0$ i.e. $\alpha = 1$, but
$$\frac{f(x,x) - x}{\sqrt{x^2+x^2}} = - \frac{1}{2 \sqrt{2}} \frac{x}{|x|}\underset{x \rightarrow 0\hspace{4pt}}{\not\hspace{-4pt}\longrightarrow} 0.$$
	\item The compositions of arc-analytic maps is arc-analytic.
\end{enumerate} 
\end{exs}
 
\begin{rem} By arguments similar to the ones we used in the proof of lemma \ref{lemequivalentchararcsymrn}, a map~$f$ from a subset $S$ of $\R^n$ to a subset $T$ of $\R^m$ is arc-analytic if and only if for any analytic arc $\gamma : ]\epsilon;\epsilon'[ \rightarrow S$, the composition $f \circ \gamma : ]\epsilon;\epsilon'[ \rightarrow T$ is analytic.
\end{rem} 
 
Since analytic arcs come up in the definition of both notions, we can expect that arc-symmetric sets and arc-analytic maps are somehow ``compatible''. The proposition below is a first result in this direction (see also Proposition 5.1 of \cite{KurAR} and Proposition 2.31 of \cite{KurPar}):

\begin{prop} \label{proparcanalyticmaps} Let $S$ be a subset of $\R^n$, $T$ be a subset of $\R^m$ and $f : S \rightarrow T$ be an arc-analytic map. 
\begin{enumerate}
	\item If $S$ is an arc-symmetric set of $\R^n$ and $T$ is an arc-symmetric set of $\R^m$, the graph $\Gamma_f := \left\{(x,y) \in S \times T~|~y = f(x)\right\}$ of $f$ is an arc-symmetric set of $\R^{n+m}$. 
	\item If $S$ is an arc-symmetric set of $\R^n$ and $A$ is any arc-symmetric set of $\R^m$ included in $T$, then $f^{-1}(A)$ is an arc-symmetric set of $\R^n$.
\end{enumerate}
\end{prop} 
 
\begin{proof}
\begin{enumerate}
	\item Suppose that $S$ is an arc-symmetric set of $\R^n$ and $T$ is an arc-symmetric set of $\R^m$. Let $\gamma = (\gamma_1,\gamma_2) : ]-1;1[ \rightarrow \R^{n+m} = \R^n \times \R^m$ be an analytic arc such that $\gamma(]-1;0[) \subset \Gamma_f$ i.e. for all $t \in ]-1;0[$, $\gamma_1(t) \in S$, $\gamma_2(t) \in T$ and $\gamma_2(t) = f \circ \gamma_1(t)$. First, we have $\gamma_1(]-1;1[) \subset S$, resp. $\gamma_2(]-1;1[) \subset T$, since $S$, resp. $T$, is arc-symmetric. Moreover, since $f \circ \gamma_1 : ]-1;1[ \rightarrow \R^m$ is an analytic function (because $f$ is arc-analytic), we have $\gamma_2 = f \circ \gamma_1$ on the entire interval $]-1;1[$ (by the identify theorem for analytic functions) and then $\gamma(]-1;1[) \subset \Gamma_f$.
	\item Suppose that $S$ is an arc-symmetric set of $\R^n$, let $A$ be an arc-symmetric set of~$\R^m$ included in $T$ and let $\gamma : ]-1;1[ \rightarrow \R^n$ be an analytic arc such that $\gamma(]-1;0[) \subset f^{-1}(A)$ i.e. $\gamma(]-1;0[) \subset S$ and $f \circ \gamma(]-1;0[) \subset A$. Since $S$ is arc-symmetric, $\gamma(]-1;1[) \subset S$ and, since $f \circ \gamma : ]-1;1[ \rightarrow \R^m$ is an analytic arc (because $f$ is arc-analytic) and $A$ is arc-symmetric, we have $f \circ \gamma(]-1;1[) \subset A$ i.e. $\gamma(]-1;1[) \subset f^{-1}(A)$.
\end{enumerate}
\end{proof}

As in \cite{KurAR} and \cite{KurPar}, let us now focus on semialgebraic arc-symmetric sets and semialgebraic arc-analytic maps:

\begin{prop} \label{proparcanalyticmapsandarclosedsets} Let $S$ be an $\mathcal{AR}$-closed set of $\R^n$, $T$ be a semialgebraic subset of $\R^m$ and $f : S \rightarrow T$ be a semialgebraic arc-analytic map.
\begin{enumerate}
	\item If $T$ is furthermore arc-symmetric in $\R^m$, then the graph $\Gamma_f := \left\{(x,y) \in S \times T~|~y = f(x)\right\}$ of $f$ is $\mathcal{AR}$-closed in $\R^{n+m}$ and $f$ is furthermore continuous with respect to (the topologies induced on $S$ and $T$ by) Euclidean topology.
	\item For any $\mathcal{AR}$-closed set $A$ of $\R^m$ included in $T$, the set $f^{-1}(A)$ is $\mathcal{AR}$-closed in $\R^n$, i.e.~$f$ is continuous with respect to (the topologies induced on $S$ and $T$ by) $\mathcal{AR}$ topology.
\end{enumerate}
\end{prop}

\begin{proof} The semialgebraic arc-symmetric counterparts of the statements of proposition \ref{proparcanalyticmaps} come from the latter together with the very definition of a semialgebraic map and the fact that the inverse image of a semialgebraic set by a semialgebraic map is a semialgebraic (see for instance Proposition 2.2.7 of \cite{BCR}). Assuming that $T$ is arc-symmetric in $\R^n$, it then remains to show that $f$ is continuous with respect to Euclidean topology.
 
Since the graph $\Gamma_f$ of $f$ is $\mathcal{AR}$-closed in $\R^{n+m}$, it is in particular closed with respect to Euclidean topology by proposition \ref{propsemialgarcsymisclosed}. Let us now show, using the arguments of the proof of Proposition 3.5 of \cite{DFP}, that $f$ is locally bounded on $S$, which will imply that $f$ is continuous by remark \ref{remproofclosedgraphlocallyboundimpliescontinuous} below.

Suppose by absurd that there exists $x \in S$ such that for any open (Euclidean) neighborhood~$U$ of $x$ in $\R^n$, $f$ is not bounded on $S \cap U$. In particular, for any $N \in \Nstar$, there exists $x_N$ in the intersection of $S$ with the open ball $B\left(x,\frac{1}{N}\right) \subset \R^n$ of radius $\frac{1}{N}$ centered at $x$ such that $\| f(x_N) \| \geq N$. Now, consider the semialgebraic subset $A := \left\{(z,t) \in S \times ]0;+\infty[~|~\| f(z) \| \geq \frac{1}{t}\right\}$ of $\R^{n+1}$ ($f$ is semialgebraic): for all $N \in \Nstar$, we have $\left(x_N, \frac{1}{N}\right) \in A$ and 
$$\left\| \left(x_N,\frac{1}{N}\right) - (x,0)\right\|^2 = \left\| x_N - x \right\|^2 + \frac{1}{N^2} \leq \frac{2}{N^2} \underset{N \rightarrow +\infty}{\longrightarrow} 0,$$ 
so that $(x,0)$ belongs to the Euclidean closure $\overline{A}$ of $A$ in $\R^{n+1}$. By Proposition 8.1.13 of \cite{BCR}, there then exists a Nash arc $\gamma : ]-1;1[ \rightarrow \R^{n+1}$ such that $\gamma(0) = (x,0)$ and $\gamma(]0;1[) \subset A$. Denote by $\gamma_1 : ]-1;1[ \rightarrow \R^n$ and $\gamma_2 : ]-1;1[ \rightarrow \R$ the respective compositions of $\gamma$ with the (Nash) projections $\R^{n+1} \rightarrow \R^n~;~(z,t) \mapsto z$ and $\R^{n+1} \rightarrow \R~;~(z,t) \mapsto t$. Notice that since $\gamma_1(]0;1[) \subset S$ and $S$ is arc-symmetric, we have $\gamma_1(]-1;1[) \subset S$ and we can consider the composition $f \circ \gamma_1 : ]-1;1[ \rightarrow T$, which is analytic since $f$ is arc-analytic. In particular, $f \circ \gamma_1$ is continuous at $0$. But, by continuity of $\gamma$ at $0$, for any $M \in \Nstar$, there exists $t_M \in \left]0;\frac{1}{M}\right[$ such that $\left\| \gamma_1(t_M) - x\right\|^2 + \left(\gamma_2(t_M)\right)^2 \leq \frac{1}{M^2}$, in particular $\left| \gamma_2(t_M)\right| \leq \frac{1}{M}$ and then, because $\gamma(t_M) = \left(\gamma_1(t_M), \gamma_2(t_M) \right) \in A$,
$$\left\| f \circ \gamma_1(t_M) \right\| \geq \frac{1}{\gamma_2(t_M)} \geq M,$$
which contradicts the fact that $f \circ \gamma_1(t_M) \underset{M \rightarrow +\infty}{\longrightarrow} f \circ \gamma_1(0) = f(x)$ ($t_M \underset{M \rightarrow +\infty}{\longrightarrow} 0$ and $f \circ \gamma_1$ is continuous at $0$).
\end{proof}
 
\begin{rem} \label{remproofclosedgraphlocallyboundimpliescontinuous} If $S$ is any subset of $\R^n$ and $f : S \rightarrow \R^m$ is a map such that $\Gamma_f$ is closed and~$f$ is locally bounded, then $f$ is continuous (with respect to Euclidean topology). Indeed, let $x \in S$ and let $(x_k)_{k \in \N}$ be a sequence of elements of $S$ converging to $x$. Let also $\delta$ and~$K$ be positive real numbers such that for any $z$ in $S \cap \overline{B(x,\delta)}$, $\| f(z)\| \leq K$. Since $(x_k)_{k \in \N}$ converges to $x$, there exist $k_0 \in \N$ such that for any $k$ greater than $k_0$, $x_{k_0}$ is in $S \cap B(x,\delta)$ and then $f(x_k)$ is in the compact closed ball $B$ of $\R^m$ centered at the origin and of radius $K$. Consider the restriction $h : S \cap \overline{B(x,\delta)} \rightarrow B$ of $f$ and recall that a map $\varphi$ from a topological space $X$ to a \emph{compact} topological space $Y$ is continuous as soon as its graph $\Gamma_{\varphi}$ is closed in $X \times Y$ (if~$A$ is any closed subset of $Y$, $\varphi^{-1}(A)$ is the image of the closed subset $\Gamma_{\varphi} \cap (X \times A)$ of $X \times Y$ by the projection map $X \times Y \rightarrow X$ which is closed because $Y$ is compact): since the graph $\Gamma_h = \Gamma_f \cap (\overline{B(x,\delta)} \times B)$ of $h$ is closed and $B$ is compact, the map $h$ is continuous. As a consequence, $f(x_k) = h(x_k) \underset{k \rightarrow +\infty}{\longrightarrow} h(x) = f(x)$.
\end{rem} 
 
In particular, if $f : \R^n \rightarrow \R$ is a semialgebraic arc-analytic map, then $f^{-1}(0)$ is a semialgebraic arc-symmetric set of $\R^n$. In \cite{AdSey} is  proved the converse implication (Theorem 1) so that:
 
\begin{theo} \label{theoadsey} Let $S$ be a subset of $\R^n$. Then $S$ is $\mathcal{AR}$-closed in $\R^n$ if and only if there exists a semialgebraic arc-analytic function $f : \R^n \rightarrow \R$ such that $S = f^{-1}(0)$.
\end{theo}   

In \cite{AdSey} are stated several consequences of this result. Let us emphasize the following one. As in \cite{KurAR} and \cite{AdSey}, we denote by $\mathcal{A}_a(\R^n)$ the ring of semialgebraic arc-analytic functions from~$\R^n$ to~$\R$ and, if $S$ is any subset of $\R^n$, by $\mathcal{I}(S)$ the set of functions of $\mathcal{A}_a(\R^n)$ that vanish on~$S$:~$\mathcal{I}(S)$ is an ideal of $\mathcal{A}_a(\R^n)$.

\begin{prop} \label{propvisarclosure}Let $S$ be any subset of $\R^n$. We have $V\left(\mathcal{I}(S)\right) = \overline{S}^{\mathcal{AR}}$, in other words,~$\overline{S}^{\mathcal{AR}}$ is the set of points of $\R^n$ at which vanish the semialgebraic arc-analytic functions of $\mathcal{A}_a(\R^n)$ that vanish on $S$.
\end{prop}

\begin{proof} By definition, we have $S \subset V\left(\mathcal{I}(S)\right)$. Moreover, $V\left(\mathcal{I}(S)\right) = \bigcap_{f \in \mathcal{I}(S)} f^{-1}(0)$ and, for $f \in \mathcal{I}(S)$, $f^{-1}(0)$ is an $\mathcal{AR}$-closed set of $\R^n$ since the function $f : \R^n \rightarrow \R$ is semialgebraic arc-analytic (proposition \ref{proparcanalyticmapsandarclosedsets} 2). Since any intersection of $\mathcal{AR}$-closed sets of $\R^n$ is an $\mathcal{AR}$-closed set by theorem \ref{theotopar}, $V\left(\mathcal{I}(S)\right)$ is an $\mathcal{AR}$-closed set of $\R^n$ containing $S$ and therefore $\overline{S}^{\mathcal{AR}} \subset V\left(\mathcal{I}(S)\right)$.

Now, notice that $V\left(\mathcal{I}(S)\right) = S$ if and only if $S$ is $\mathcal{AR}$-closed in $\R^n$. Indeed, $V\left(\mathcal{I}(S)\right)$ is~$\mathcal{AR}$-closed and if we suppose that $S$ is~$\mathcal{AR}$-closed, i.e. $S = f^{-1}(0) = V(f)$ with $f \in \mathcal{A}_a(\R^n)$ by theorem \ref{theoadsey}, then $S = \overline{S}^{\mathcal{AR}} \subset V\left(\mathcal{I}(S)\right) \subset V(f) = S.$

Finally, remark that $\mathcal{I}\left(\overline{S}^{\mathcal{AR}}\right) \subset \mathcal{I}(S)$, so that $V\left(\mathcal{I}(S)\right) \subset V\left(\mathcal{I}\left(\overline{S}^{\mathcal{AR}}\right)\right) = \overline{S}^{\mathcal{AR}}$ (the set~$\overline{S}^{\mathcal{AR}}$ is $\mathcal{AR}$-closed).
\end{proof}
 
\begin{cor} \label{corarclosureprodisproductarclos} Let $S$ be any subset of $\R^n$ and $T$ be any subset of $\R^m$. Then the $\mathcal{AR}$-closure of $S \times T$ in $\R^{n+m}$ is the Cartesian product $\overline{S}^{\mathcal{AR}} \times \overline{T}^{\mathcal{AR}}$.
\end{cor} 

\begin{proof} The proof is similar to the proof of the equality of the Zariski closure of a Cartesian product of subsets of affine spaces with the Cartesian product of their Zariski closures.

First, $\overline{S}^{\mathcal{AR}} \times \overline{T}^{\mathcal{AR}}$ is an $\mathcal{AR}$-closed set of $\R^{n+m}$ (a Cartesian product of arc-symmetric, resp. semialgebraic, sets is arc-symmetric, resp. semialgebraic) containing $S \times T$, so that $\overline{S \times T}^{\mathcal{AR}} \subset \overline{S}^{\mathcal{AR}} \times \overline{T}^{\mathcal{AR}}$. In order to show the converse inclusion, let $f \in \mathcal{I}(S \times T)$ and let $x \in S$: the function $y \in \R^m \mapsto f(x,y) \in \R$ is semialgebraic arc-analytic (as the composition of $f$ with the semialgebraic arc-analytic map $y \in \R^m \mapsto (x,y) \in \R^n \times \R^m$) and vanishes on~$T$, hence on $\overline{T}^{\mathcal{AR}}$ (recall that, by previous proposition \ref{propvisarclosure}, $\overline{T}^{\mathcal{AR}}$ is the set of points of~$\R^m$ on which vanish the semialgebraic arc-analytic functions of $\mathcal{A}_a(\R^m)$ vanishing on $T$). As a consequence, $f$ vanishes on $S \times \overline{T}^{\mathcal{AR}}$. Now fix $y \in \overline{T}^{\mathcal{AR}}$ and consider the semialgebraic arc-analytic function $x \in \R^n \mapsto f(x,y) \in \R$, which vanishes on $S$ hence on $\overline{S}^{\mathcal{AR}}$. Therefore $f$ vanishes on $\overline{S}^{\mathcal{AR}} \times \overline{T}^{\mathcal{AR}}$, so that $\mathcal{I}(S \times T) \subset \mathcal{I}\left(\overline{S}^{\mathcal{AR}} \times \overline{T}^{\mathcal{AR}}\right)$ and then
$$\overline{S}^{\mathcal{AR}} \times \overline{T}^{\mathcal{AR}} = V\left(\mathcal{I}\left(\overline{S}^{\mathcal{AR}} \times \overline{T}^{\mathcal{AR}}\right)\right) \subset V\left( \mathcal{I}(S \times T) \right) = \overline{S}^{\mathcal{AR}} \times \overline{T}^{\mathcal{AR}}.$$
\end{proof}

\subsubsection{Arc-symmetric sets of real analytic manifolds}

The notion of arc-symmetric set can be generalized to any real analytic manifold: let $M$ be a real analytic manifold of dimension $n \in \mathbb{N}$ and recall that a continuous arc $f : I \rightarrow M$, where~$I$ is an open interval of $\R$, is analytic if and only if for any $t \in I$, there exists an analytic chart $\phi : U \subset M  \rightarrow \R^n$ of $M$ around $f(t)$ such that the arc $\phi \circ f : f^{-1}(U) \rightarrow \R^n$ is analytic.

\begin{de} \label{defarcsymsetsrealanman} Let $S$ be a subset of $M$. We say that $S$ is an \emph{arc-symmetric set of $M$} if for any analytic arc $\gamma : ]-1;1[ \rightarrow M$, if $\gamma(]-1;0[) \subset S$, then $\gamma(]-1;1[) \subset S$.
\end{de}

Remark that the characterizations of lemma \ref{lemequivalentchararcsymrn} remain valid if we replace $\R^n$ by $M$ (the proof of lemma \ref{lemequivalentchararcsymrn} amounts to deal with analytic diffeomorphisms on open intervals). Let us also highlight a first class of examples of arc-symmetric sets of $M$, which generalizes example~\ref{exanalyticsubsetarcsym}~1:

\begin{lem} \label{lemanalyticsubsetmanifoldarcsym} Let $N$ be an \emph{analytic subset of $M$} i.e. for any $x \in M$, there exists an open neighborhood $U$ of $x$ in $M$ and an analytic function $f : U \rightarrow \R$ such that $N \cap U = f^{-1}(0)$. If $N$ is closed in $M$, then $N$ is arc-symmetric in $M$. In particular, any closed analytic submanifold of $M$ is arc-symmetric in $M$.
\end{lem}

\begin{proof} Suppose that $N$ is closed in $M$ and let $\gamma : ]-1;1[ \rightarrow M$ be an analytic arc such that $\gamma(]-1;0[) \subset N$: since $N$ is closed and $\gamma$ is continuous, we have $\gamma(]-1;0]) \subset N$. Considering an open neighborhood $U$ of $x$ in $M$ and an analytic function $f : U \rightarrow \R$ such that $N \cap U = f^{-1}(0)$, we can then follow the arguments of example \ref{exanalyticsubsetarcsym} 1 to show that $N$ is arc-symmetric in $M$.

In order to prove the second statement of the lemma, let us show that any analytic submanifold of $M$ is an analytic subset of $M$. Let $S$ be an analytic submanifold of dimension $d \in \mathbb{N}$ of~$M$: for any point $x \in S$, there exists an analytic diffeomorphism $\phi$ from an open neighborhood $U$ of $x$ in $M$ onto an open neighborhood $V$ of the origin in $\R^n$ such that $\phi(x) = \bf{0}$ and 
$$\phi(S \cap U) = \{(x_1,\ldots,x_n) \in V~|~x_{d+1} = \ldots = x_n = 0\}.$$
If we denote by $\pi_i : V \rightarrow \R$ the projection on the $i$\textsuperscript{th} coordinate, $i \in \{1,\ldots,n\}$, we have $S \cap U = f^{-1}(0)$ where $f$ is the analytic function $U \rightarrow \R~;~y \mapsto \left(\pi^2_{d+1} + \cdots + \pi^2_{n}\right)(\phi(y))$ (the argument is identical to example \ref{exanalyticsubsetarcsym} 3).
\end{proof}

\begin{rem} \label{remremconncomparcsymnashmanarcsym} The topological arguments of remark \ref{remconncomparcsymarcsym} are also valid in this context: the connected components of an arc-symmetric set of $M$ are arc-symmetric in $M$.
\end{rem}

Considering arc-symmetric sets of real analytic manifolds, it is natural to study the compatibility of the notion with analytic diffeomorphisms. But actually, we can already consider the following generalizations of definition \ref{defarcanalyticmap} and proposition \ref{proparcanalyticmaps}: 

\begin{de} Let $S$ be a subset of $M$, $T$ be a subset of a real analytic manifold $N$ and $f : S \rightarrow T$ be a map. We say that $f$ is \emph{arc-analytic} if for any analytic arc $\gamma : ]-1;1[ \rightarrow S$, the composition $f \circ \gamma : ]-1;1[ \rightarrow T$ is analytic.
\end{de}

\begin{prop} \label{propcomparcanalyticmapwitharcsymsetrealanalmanifold} Let $S$ be a subset of $M$, $T$ be a subset of a real analytic manifold $N$ of dimension~$m$ and $f : S \rightarrow T$ be an arc-analytic map. 
\begin{enumerate}
	\item If $S$ is an arc-symmetric set of $M$ and $T$ is an arc-symmetric set of $N$, then the graph $\Gamma_f := \left\{(x,y) \in S \times T~|~y = f(x)\right\}$ of~$f$ is an arc-symmetric set of the product analytic manifold $M \times N$. 
	\item If $S$ is an arc-symmetric set of $M$ and $A$ is any arc-symmetric set of $N$ included in $T$, then $f^{-1}(A)$ is an arc-symmetric set of $M$.
\end{enumerate}
In particular, if $\Phi : N \rightarrow M$ is an analytic diffeomorphism, $S \subset M$ is arc-symmetric in $M$ if and only if $\Phi^{-1}(S)$ is arc-symmetric in $N$ (an analytic map between real analytic manifolds is arc-analytic).
\end{prop} 
 
\begin{proof}
\begin{enumerate}
	\item Suppose that $S$ is an arc-symmetric set of $M$. Let $\gamma = (\gamma_1,\gamma_2) : ]-1;1[ \rightarrow M \times N$ be an analytic arc such that $\gamma(]-1;0[) \subset \Gamma_f$ i.e. for all $t \in ]-1;0[$, $\gamma_1(t) \in S$, $\gamma_2(t) \in T$ and $\gamma_2(t) = f \circ \gamma_1(t)$. First, we have $\gamma_1(]-1;1[) \subset S$, resp. $\gamma_2(]-1;1[) \subset T$, since $S$ is arc-symmetric in $M$, resp. $T$ is arc-symmetric in $N$. Furthermore, since the maps $\gamma_2 : ]-1;1[ \rightarrow N$ and $f \circ \gamma_1 : ]-1;1[ \rightarrow N$ are both analytic maps ($f$ is arc-analytic), they are in particular continuous so that $\gamma_2(0) = f \circ \gamma_1(0)$. Moreover, if $\psi : W \rightarrow \R^m$ is an analytic chart on $N$ such that $\gamma_2(0) = f \circ \gamma_1(0) \in W$ and if $\epsilon$ is a positive real number such that $]-\epsilon; \epsilon[$ is included is the non-empty open subset $\gamma_2^{-1}(W) \cap \left(f \circ \gamma_1\right)^{-1}(W)$ of~$]-1;1[$, we have $\psi \circ \gamma_2 = \psi \circ f \circ \gamma_1$ on $]-\epsilon;0]$, hence on the entire interval $]-\epsilon;\epsilon[$ (by the identify theorem for analytic functions) and then $\gamma(]-1;\epsilon[) \subset \Gamma_f$.	
	\item The proof is identical to the proof of proposition \ref{proparcanalyticmaps} 2.
\end{enumerate}
\end{proof}

As a particular case of the last statement of proposition \ref{propcomparcanalyticmapwitharcsymsetrealanalmanifold}, consider the inclusion map $i : N \hookrightarrow M$ of an analytic submanifold $N$ of $M$: if a subset $T$ of $N$ is arc-symmetric in $M$, then $T=i^{-1}(T)$ is arc-symmetric in $N$. The converse is true if $N$ is a \emph{closed} analytic submanifold of $M$:

\begin{lem} \label{corarcsymsubsetclosedsubman} Let $N$ be a closed analytic submanifold of $M$ and $T$ be a subset of $N$. Then~$T$ is arc-symmetric in $N$ if and only if $T$ is arc-symmetric in $M$.
\end{lem}

\begin{proof} Suppose that $T$ is arc-symmetric in $N$ and let $\gamma : ]-1;1[ \rightarrow M$ be an analytic arc such that $\gamma(]-1;0[) \subset T \subset N$. Since, by lemma \ref{lemanalyticsubsetmanifoldarcsym}, $N$ is arc-symmetric in $M$, we have $\gamma(]-1;1[) \subset N$ and we can restrict $\gamma$ into an arc $\widetilde{\gamma} : ]-1;1[ \rightarrow N$, which is analytic with respect to the real analytic manifold $N$ since $N$ is an analytic submanifold of $M$. As a consequence, because $T$ is arc-symmetric in $N$ and $\widetilde{\gamma}(]-1;0[) = \gamma(]-1;0[) \subset T$, we have $\widetilde{\gamma}(]-1;1[) \subset T$ i.e. $\gamma(]-1;1[) \subset T$.
\end{proof}

\begin{rem} \label{remexarcsymr2notarcsymp2} The above-proved implication is no longer true if $N$ is not closed in $M$. Consider for instance the analytic submanifold $\R^2$ of $\mathbb{P}^2(\R)$, identified with the open subset 
$$U_0 := \left\{[x_0 : x_1 : x_2] \in \mathbb{P}^2(\R)~|~x_0 \neq 0\right\}$$ 
of $\mathbb{P}^2(\R)$ ($\R^2$ is not closed in $\mathbb{P}^2(\R)$), and let $C := \left\{(x,y) \in \R^2~|~xy = 1, x> 0\right\}$ be the positive half-branch of the hyperbola of equation $xy=1$ in $\R^2$. As a connected component of the real algebraic hyperbola, $C$ is arc-symmetric in $\R^2$ (see example \ref{exanalyticsubsetarcsym}). But $C$ is not arc-symmetric in $\mathbb{P}^2(\R)$. Indeed, consider the map $\gamma : t \in ]-1;1[ \mapsto [-t:1:t^2] \in \mathbb{P}^2(\R)$: the arc $\gamma$ is analytic since, if we denote by $\phi$ the chart
	$$[x_0 : x_1 : x_2] \in U_1 := \left\{[x_0 : x_1 : x_2] \in \mathbb{P}^2(\R)~|~x_1 \neq 0\right\} \mapsto (x_0,x_2) \in \R^2,$$
the map $\phi \circ \gamma : t \in ]-1;1[ \mapsto (-t,t^2) \in \R^2$ is analytic. We have furthermore $\gamma(]-1;0[) \subset C$ because, if $t \in ]-1;0[$, $\gamma(t) = [-t:1:t^2] = \left[1 : -\frac{1}{t} : -t\right] \in C$. However, the point $\gamma(0) =  [0:1:0]$ does not belong to $U_0$, in particular does not belong to $C$, so that $C$ is not arc-symmetric in $\mathbb{P}^2(\R)$ (remark that the analytic arc $\gamma$ also allows to show that $\R^2$ is not arc-symmetric in~$\mathbb{P}^2(\R)$).
\end{rem}

Let us end this paragraph with the general real analytic counterpart of proposition \ref{propclosedarcsymtop}: the proof is identical.

\begin{prop} \label{propintersecfiniteunionclarcsymisclarcsym} Any intersection or finite union of closed arc-symmetric subsets of $M$ is a closed arc-symmetric set of $M$.
\end{prop}

\subsubsection{Semialgebraic sets of Nash manifolds} \label{subsubsecsaofnashmfds}

As we did in real affine spaces, we would like to consider semialgebraic arc-symmetric sets in this more general context of real analytic manifolds and we are then led to consider Nash manifolds. We therefore suppose from now that $M$ is a \emph{Nash manifold} of dimension $n$ i.e. a real analytic manifold with respect to a finite atlas of charts $\phi : U \rightarrow \R^n$ with semialgebraic open image $\phi(U)$ in $\R^n$ such that the analytic transition maps are Nash diffeomorphisms, i.e. semialgebraic diffeomorphisms, between semialgebraic open subsets of $\R^n$ (such an atlas will be called a \emph{Nash atlas}). For instance, the real projective space $\mathbb{P}^n(\R)$, with its usual analytic atlas, is a Nash manifold of dimension $n$. We refer to \cite{Shiota} for a thorough study of Nash manifolds and we begin this part by recalling below the definition and some properties of semialgebraic sets of Nash manifold. A \emph{Nash chart} on $M$ is a continuous map $\varphi : V \rightarrow \R^n$ from an open subset $V$ of $M$ which restricts to an homeomorphism $\varphi : V \rightarrow \varphi(V)$ onto a semialgebraic open subset of $\R^n$ such that for any chart $\phi : U \rightarrow \R^n$ of the Nash atlas of $M$, the map $\varphi \circ \phi^{-1} : \phi(U \cap V) \rightarrow \varphi(U \cap V)$ is a Nash diffeomorphism between semialgebraic open subsets of $\R^n$.

\begin{de} Let $S$ be a subset of $M$. We say that $S$ is \emph{semialgebraic} if for any chart $\phi : U \rightarrow \R^n$ of the Nash atlas of $M$, the set $\phi(S \cap U)$ is a semialgebraic subset of $\R^n$.
\end{de}

Notice that, in particular, $M$ is a semialgebraic subset of $M$ and that, if $\varphi : V \rightarrow \R^n$ is a Nash chart on $M$, $V$ is a semialgebraic subset of $M$. Remark also that $S \subset M$ is semialgebraic in $M$ if and only if for any Nash chart $\varphi : V \rightarrow \R^n$ on $M$, $\varphi(S \cap V)$ is a semialgebraic subset of $\R^n$.

Let us now state and show the following stability properties that extend the usual stability properties of semialgebraic sets of real affine spaces:

\begin{lem} \label{lembooleancomsaisasclosureissa} Any finite Boolean combination of semialgebraic subsets of $M$ is semialgebraic in $M$. Moreover, if $S$ is a semialgebraic set of $M$, the closure and interior of $S$ in $M$ are semialgebraic in $M$ and, if $T$ is a semialgebraic subset of a Nash manifold $N$, the Cartesian product $S \times T$ is a semialgebraic subset of the Nash manifold $M \times N$. 
\end{lem}

\begin{proof} Let $A$ and $B$ be semialgebraic subsets of $M$ and let $\varphi : V \rightarrow \R^n$ be a Nash chart on $M$. The subsets $\varphi((A \cap B) \cap V) = \varphi(A \cap V) \cap \varphi(B \cap V)$, $\varphi((A \cup B) \cap V) = \varphi(A \cap V) \cup \varphi(B \cap V)$ and $\varphi((M \setminus A) \cap V) = \varphi(V) \setminus \varphi(A \cap V)$ of $\R^n$ are all semialgebraic, so that the subsets $A \cap B$, $A \cup B$ and $M \subset A$ are all semialgebraic in $M$.

As for the second part of the statement, keep the Nash chart $\varphi : V \rightarrow \R^n$ on $M$, denote by $\overline{S}$ the closure of $S$ in $M$ and remark that $\overline{S} \cap V$ is the closure of $S \cap V$ in $V$ ($\overline{S} \cap V$ is a closed subset of $V$ containing $S \cap V$ and if $x \in \overline{S} \cap V$ and $W$ is any open subset of~$V$ containing $x$, we have $W \cap S \cap V = W \cap S \neq \emptyset$, since $W$ is also an open neighborhood of~$x$ in $M$). As a consequence, the set $\varphi\left(\overline{S} \cap V\right)$ is the closure of the semialgebraic set $\varphi(S \cap V)$ in the semialgebraic open subset $\varphi(V)$ of $\R^n$, which is the intersection of $\varphi(V)$ with the closure of $\varphi(S \cap V)$ in $\R^n$ and is therefore semialgebraic (see Proposition 2.2.2 of \cite{BCR}): $\overline{S}$ is then a semialgebraic subset of $M$. The interior of $S$ in $M$, being the set $M \setminus \left(\overline{M \setminus S}\right)$, is consequently also semialgebraic in $M$.

In order to show the final statement, let $\left(\phi_i : U_i \rightarrow \R^n\right)_{i \in \{1,\ldots,r\}}$ be the Nash atlas of $M$, $\left(\psi_j : W_j \rightarrow \R^m\right)_{j \in \{1,\ldots,s\}}$ be the Nash atlas of $N$ and consider the product Nash manifold $M \times N$ with Nash atlas 
$$\phi_i \times \psi_j : \begin{array}{ccc}U_i \times W_j & \rightarrow & \R^n \times \R^m\\(x,y) & \mapsto & \left( \phi_i(x), \psi_j(y)\right)\end{array},~i \in \{1,\ldots,r\},~j \in \{1,\ldots,s\}.$$
If $i \in \{1,\ldots,r\}$ and $j \in \{1,\ldots,s\}$, we have 
$$(\phi_i \times \psi_j)\big( (S \times T) \cap (U_i \times W_j)\big) = \phi_i\left(S \cap U_i\right) \times \psi_j\left(T \cap W_j\right)$$ 
and the latter set is then a semialgebraic subset of $\R^n \times \R^m$, so that $S \times T$ is semialgebraic in~$M \times N$.
\end{proof}

Some classes of maps are natural to be considered in this context: if $N$ is a Nash manifold of dimension $m$ and $\Phi : M \rightarrow N$ is a map, we say that $\Phi$ is a \emph{semialgebraic map}, resp. a \emph{Nash map}, if for any Nash chart $\phi : U \rightarrow \R^n$ on $M$ and any Nash chart $\psi : W \rightarrow \R^m$ on $N$, the set $\Phi^{-1}(W)$ is semialgebraic in $M$, resp. open semialgebraic in $M$, and the map $\psi \circ \Phi \circ \phi^{-1} : \phi\left(U \cap \Phi^{-1}(W)\right) \rightarrow \R^m$ is a semialgebraic map, resp. a Nash (i.e. semialgebraic and analytic) map. A \emph{Nash diffeomorphism} is a bijective Nash map between Nash manifolds whose inverse is also Nash. 

Notice that, in the above definition, we can consider only the (Nash) charts of the respective Nash atlases of $M$ and $N$.
 
\begin{lem} \label{lemdirectiminverseimsaissa} Let $N$ be a Nash manifold and $\Phi : M \rightarrow N$ be a map.
\begin{enumerate}
	\item The map $\Phi$ is semialgebraic if and only if the graph $\Gamma_{\Phi}$ of $\Phi$ is semialgebraic in $M \times N$.
	\item The map $\Phi$ is a Nash map if and only if $\Phi$ is a semialgebraic analytic map.
	\item If $S$ is a semialgebraic subset of $M$, $T$ is a semialgebraic subset of $N$ and $\Phi$ is a semialgebraic map, then $\Phi^{-1}(T)$ is semialgebraic in $M$ and $\Phi(S)$ is semialgebraic in $N$.
\end{enumerate}
\end{lem}

\begin{proof} Let $\left(\phi_i : U_i \rightarrow \R^n\right)_{i \in \{1,\ldots,r\}}$ be the Nash atlas of $M$, $\left(\psi_j : W_j \rightarrow \R^m\right)_{j \in \{1,\ldots,s\}}$ be the Nash atlas of $N$ and let $i \in \{1,\ldots,r\}$ and $j \in \{1,\ldots,s\}$
\begin{enumerate}
	\item We have 
\begin{eqnarray*}
\phi_i \times \psi_j\big(\Gamma_{\Phi} \cap (U_i \times W_j) \big) & = & \left\{\left(\phi_i(x), \psi_j \circ \Phi(x)\right)~|~x \in U_i \cap \Phi^{-1}(W_j)\right\} \\
& = & \left\{\left(z, \psi_j \circ \Phi \circ \phi_i^{-1}(z)\right)~|~z \in \phi_i\left(U_i \cap \Phi^{-1}(W_j)\right)\right\} = \Gamma_{\psi_j \circ \Phi \circ \phi_i^{-1}},
\end{eqnarray*}	
hence the graph $\Gamma_{\Phi}$ is semialgebraic in $M \times N$ if and only if $\Phi$ is a semialgebraic map.
	\item By definition, a Nash map is an analytic map together with a semialgebraic map, since a Nash map from an open semialgebraic set of an affine space to an affine space is analytic and semialgebraic. Conversely, suppose that $\Phi$ is semialgebraic and analytic, let $i,j \in \{1,\ldots,n\}$ and consider the map $\psi_j \circ \Phi \circ \phi_i^{-1} : \phi_i\left(U \cap \Phi^{-1}(W_j)\right) \rightarrow \R^m$: the set $\Phi^{-1}(W_j)$ is semialgebraic in~$M$, since $\Phi$ is semialgebraic, and open in $M$, since $\Phi$ is analytic, and the latter map is semialgebraic and analytic, hence Nash.
	\item We use the argument of the proof of Proposition 2.2.7 of \cite{BCR}: if $i \in \{1,\ldots,r\}$, we have 
$$\phi_i\left(\Phi^{-1}(T) \cap U_i\right) = \bigcup_{j = 1}^s \phi_i\left(\Phi^{-1}(T \cap W_j) \cap U_i\right)$$
and, for any $j \in \{1,\ldots,s\}$, the set $\phi_i\left(\Phi^{-1}(T\cap W_j) \cap U_i\right)$ is the image of the semialgebraic set $\big(\phi_i\left(U_i\right) \times \psi_j(T \cap W_j) \big) \cap \Gamma_{\psi_j \circ \Phi \circ \phi_i^{-1}}$ of $\R^n \times \R^m$ by the semialgebraic projection map $\phi_i\left(U_i\right) \times \psi_j(W_j) \rightarrow \phi_i\left(U_i\right)$. On the other hand, if $j \in \{1,\ldots,s\}$, we have
$$\psi_j\left(\Phi(S) \cap W_j\right) = \bigcup_{i=1}^r \psi_j\left(\Phi(S \cap U_i) \cap W_j\right)$$
and, for any $i \in \{1,\ldots,r\}$, the set $\psi_j\left(\Phi(S \cap U_i) \cap W_j\right)$ is the image of the semialgebraic set $\left(\phi_i\left( S \cap U_i\right) \times \psi_j(W_j) \right) \cap \Gamma_{\psi_j \circ \Phi \circ \phi_i^{-1}}$ of $\R^n \times \R^m$ by the semialgebraic projection map $\phi_i\left(U_i\right) \times \psi_j(W_j) \rightarrow \psi_j(W_j)$.
\end{enumerate}
\end{proof}

We can more generally consider a notion of semialgebraic map between semialgebraic subsets of Nash manifolds:

\begin{de} Let $S$ be a semialgebraic subset of $M$, $T$ be a semialgebraic subset of a Nash manifold $N$ of dimension $m$ and $f : S \rightarrow T$ be a map. We say that $f$ is a \emph{semialgebraic map} if for any Nash chart $\phi : U \rightarrow \R^n$ on $M$ and any Nash chart $\psi : W \rightarrow \R^m$ on $N$, the set $f^{-1}(T \cap W)$ is semialgebraic in $M$ and the map $\psi \circ f \circ \phi^{-1} : \phi\left(U \cap f^{-1}(T \cap W)\right) \rightarrow \R^m$ is a semialgebraic map.
\end{de}

\begin{lem} \label{lemsamapbetweensasetsofnashmanifolds} Let $S$ be a semialgebraic subset of $M$, $T$ be a semialgebraic subset of a Nash manifold $N$ of dimension $m$ and $f : S \rightarrow T$ be a map.
\begin{enumerate}
	\item The map $f$ is semialgebraic if and only if the graph $\Gamma_{f}$ of $f$ is semialgebraic in $M \times N$.
	\item If $B$ is a semialgebraic subset of $M$ included in $S$, $A$ is a semialgebraic subset of $N$ included in $T$ and $f$ is a semialgebraic map, then $f^{-1}(A)$ is semialgebraic in $M$ and~$f(B)$ is semialgebraic in $N$.
	\item If $f$ is a bijective semialgebraic map, then its inverse $f^{-1}$ is a semialgebraic map as well.
\end{enumerate}
\end{lem}

\begin{proof} The same arguments as in the proof of lemma \ref{lemdirectiminverseimsaissa} show statements 1 and 2. As for the third statement, remark that, if $f$ is bijective, the graph $\Gamma_{f^{-1}}$ of the inverse of $f$ is the image of $\Gamma_f$ by the Nash diffeomorphism $M \times N \rightarrow N \times M~;~(x,y) \mapsto (y,x)$. 
\end{proof}

\begin{rems} \label{remsrestrsamapcompasmapissa} Keep the notations of previous lemma \ref{lemsamapbetweensasetsofnashmanifolds}.
	\begin{enumerate}
		\item If $B$ is a semialgebraic subset of $M$ included in $S$, $A$ is a semialgebraic subset of $N$ included in $T$ such that $f(B) \subset A$ and $f$ is a semialgebraic map, then the map $\widetilde{f} : B \rightarrow A~;~x \mapsto f(x)$ induced by $f$ is a semialgebraic map as well, since $\Gamma_{\widetilde{f}} = \Gamma_f \cap (B \times A)$ is a semialgebraic set of $M \times N$ by lemma \ref{lemsamapbetweensasetsofnashmanifolds} 1 and lemma \ref{lembooleancomsaisasclosureissa}.
		\item If $X$ is a semialgebraic subset of a Nash manifold $L$ and $h : T \rightarrow X$ is a semialgebraic map, the composition $h \circ f : S \rightarrow X$ is a semialgebraic map as well, since the graph $\Gamma_{h \circ f}$ of $h \circ f$ is the image of the semialgebraic intersection $(\Gamma_f \times Z) \cap (X \times \Gamma_h)$ (lemma \ref{lemsamapbetweensasetsofnashmanifolds} 1 and lemma \ref{lembooleancomsaisasclosureissa}) by the Nash (and then semialgebraic: cf. lemma \ref{lemdirectiminverseimsaissa} 2) projection map $M \times N \times L \rightarrow M \times L~;~(x,y,z) \mapsto (x,z)$ (this argument is identical to the proof of Proposition 2.2.6 (i) of \cite{BCR}).
	\end{enumerate}
\end{rems}

Let us also mention the following generalization of Proposition 2.9.10 of \cite{BCR}. A semialgebraic subset $N$ of $M$ is a \emph{Nash submanifold of dimension $d$ of $M$} if for any $x \in N$, there exists a Nash chart $\phi : U \rightarrow \R^n$ on $M$ such that $x \in U$ and 
$$\phi(N \cap U) = \left\{(x_1,\ldots,x_n) \in \phi(U)~|~x_{d+1} = \ldots = x_n = 0\right\}$$ 
(this is equivalent to say that for any $x \in N$, there exists a Nash diffeomorphism $\phi$ from a semialgebraic open neighborhood $U$ of $x$ in $M$ onto a semialgebraic open neighborhood $V$ of the origin in $\R^n$ such that $\phi(x) = \bf{0}$ and $\phi(N \cap U) = \{(x_1,\ldots,x_n) \in V~|~x_{d+1} = \ldots = x_n = 0\}$: notice that any open semialgebraic subset of $M$ is itself a Nash manifold and that any Nash chart on $M$ is a Nash map). Remark that the inverse image of a Nash submanifold of dimension~$d$ by a Nash diffeomorphism is a Nash submanifold of dimension $d$. 

\begin{prop} \label{propsanashmanifisdisjunionofconashsub} Let $S$ be a semialgebraic subset of $M$. Then $S$ is the disjoint union of connected Nash submanifolds of $M$.
\end{prop}

\begin{proof} Consider the Nash atlas $\left(\phi_i : U_i \rightarrow \R^n\right)_{i \in \{1,\ldots,m\}}$ of $M$ and write $M$ as the disjoint union of semialgebraic subsets $U_I := \left(\bigcap_{i \in I} U_{i}\right) \setminus \left(\bigcup_{j \notin I} U_j\right)$, $I \subset \{1,\ldots,m\}$. 

Let $I \subset \{1,\ldots,m\}$ and $i \in I$: since $S$ is a semialgebraic subset of $M$, $\phi_i\left(S \cap U_I\right)$ is a semialgebraic subset of $\R^n$. By Proposition 2.9.10 of \cite{BCR}, $\phi_i\left(S \cap U_I\right)$ is a finite disjoint union of connected Nash submanifolds $T_1,\ldots,T_l$ of $\R^n$. Consider one set $T \subset \phi_i(U_i)$ among the latter and denote $\widetilde{T} := \phi_i^{-1}(T)$: since $\phi_i$ is a Nash diffeomorphism onto its image, $\widetilde{T}$ is a connected Nash submanifold of dimension $\dim T$ of $U_i$ and then of $M$ ($U_i$ is open semialgebraic in $M$).

As a result, all sets $S \cap U_I$, $I \subset \{1,\ldots,m\}$, of the disjoint union $S := \bigcup_{I \subset \{1,\ldots,m\}} S \cap U_I$ are disjoint unions of connected Nash submanifolds of $M$, hence the statement.
\end{proof}

\begin{cor} \label{remconnectedcompsemialgnashman} Let $S$ be a semialgebraic subset of $M$. All the connected components of $S$ are also semialgebraic in $M$.
\end{cor}

\begin{proof} Write $S$ as a disjoint union of connected Nash submanifolds $T_1,\ldots,T_s$ of $M$ and let $S_1,\ldots,S_r$ be the connected components of $S$. Then, for $i \in \{1,\ldots,r\}$, $S_i$ is the union of some~$T_j$'s and is therefore semialgebraic in $M$.
\end{proof}

\begin{rem} \label{remcompactnashsubisnashman} Let us mention that a Nash submanifold $N$ of $M$ is not necessarily a Nash manifold itself: the point is that, even if $N$ can be equipped with an atlas of Nash charts induced from Nash charts on $M$ making $N$ into a Nash submanifold ($N$ is therefore a real analytic manifold), this atlas is not necessarily finite, hence is not a Nash atlas. On the other hand, if $N$ is furthermore compact, a finite atlas can be extracted from the previous one and~$N$ is therefore a Nash manifold in this case.
\end{rem}

We can also consider the following natural notion of dimension for a semialgebraic subset of $M$: let $S$ be such a set.

\begin{de} The \emph{dimension} of $S$, denoted by $\dim S$, is the maximum of the dimensions of the semialgebraic sets $\phi(S \cap U)$ over the Nash charts $\phi : U \rightarrow \R^n$ of the Nash atlas of $M$. 
\end{de}

Notice that the dimension of $S$ is equal to the maximum of the dimensions of the semialgebraic sets $\varphi(S \cap V)$ over all the Nash charts $\varphi : V \rightarrow \R^n$ on $M$. Remark also that the dimension of a finite union of semialgebraic subsets of $M$ is the maximum of their dimensions, using and extending Proposition 2.8.5 i) of \cite{BCR}. In particular, if $T$ is a semialgebraic subset of~$M$ included in $S$, since $S = T \cup (S \setminus T)$, we have $\dim T \leq \dim S$.

To end this paragraph, let us show three further properties involving dimension:

\begin{lem} \label{lemdimclosminussa} We have $\dim \overline{S} \setminus S < \dim S$.
\end{lem} 

\begin{proof}
If $\varphi : V \rightarrow \R^n$ is any Nash chart on $M$, since $\varphi\left(\overline{S} \cap V\right)$ is the intersection of $\varphi(V)$ with the closure $\overline{\varphi(S \cap V)}$ of $\varphi(S \cap V)$ in $\R^n$ (see the proof of lemma \ref{lembooleancomsaisasclosureissa}), we have 
$$\dim \varphi\big(\left(\overline{S} \setminus S\right) \cap V\big) = \dim \varphi\left(\overline{S} \cap V\right) \setminus \varphi(S \cap V) \leq \dim \overline{\varphi(S \cap V)} \setminus \varphi(S \cap V) < \dim \varphi(S \cap V)$$
(Proposition 2.8.13 of \cite{BCR}). As a consequence, $\dim \varphi\big(\left(\overline{S} \setminus S\right) \cap V\big) < \dim S$ and then $\dim \overline{S} \setminus S < \dim S$.
\end{proof}

\begin{lem}
If $S$ is a Nash submanifold of dimension $d$ of $M$, we have $d = \dim S$.
\end{lem}

\begin{proof}
Let $\varphi : V \rightarrow \R^n$ be a Nash chart on $M$. The set $S \cap V$ is a then a Nash submanifold of dimension $d$ of $M$ and $\varphi(S \cap V)$ is a Nash submanifold of dimension $d$ of $\R^n$ ($V$ is open semialgebraic in $M$ and $\varphi$ is a Nash diffeomorphism onto its image which is open semialgebraic in $\R^n$). Consequently, by Proposition 2.8.14 of \cite{BCR}, $\dim \varphi(S \cap V) = d$, so that $\dim S =d$.
\end{proof}

\begin{lem} \label{lemdimimagesamap} Let $f : S \rightarrow T$ be a semialgebraic map from $S$ to a semialgebraic subset of a Nash manifold $N$ of dimension $m$. We have $\dim f(S) \leq \dim S$ and, if $f$ is bijective, $\dim f(S) = \dim S$.
\end{lem}

\begin{proof} Let $\psi : W \rightarrow \R^m$ be a Nash chart on $N$ and consider the semialgebraic set $\psi\left(f(S) \cap W\right) = \psi \circ f \left( f^{-1}(T \cap W)\right)$. If $\left(\phi_i : U_i \rightarrow \R^n\right)_{i \in \{1,\ldots,r\}}$ is the Nash atlas of $M$, the set $\psi \circ f \left( f^{-1}(T \cap W)\right)$ is the union of the sets $\psi \circ f \left( U_i \cap f^{-1}(T \cap W)\right)$, $i \in \{1,\ldots,r\}$, so that the dimension of the former is the maximum of the dimensions of the latter. But, if $i \in \{1,\ldots,r\}$, since $f$ is a semialgebraic map, the map $\psi \circ f \circ \phi^{-1} : \phi\left(U_i \cap f^{-1}(T \cap W)\right) \rightarrow \R^m$ is semialgebraic and then, according to Theorem 2.8.8 of \cite{BCR},
$$\dim \psi \circ f \left(U_i \cap f^{-1}(T \cap W)\right) \leq \dim U_i \cap f^{-1}(T \cap W) \leq \dim S.$$
As a consequence, $\dim \psi\left(f(S) \cap W\right) \leq \dim S$ and then $\dim f(S) \leq \dim S$.

Furthermore, if $f$ is bijective, then the inverse map $f^{-1}$ is a semialgebraic map as well (lemma \ref{lemsamapbetweensasetsofnashmanifolds} 3) and we have $\dim f(S) = \dim T \geq \dim f^{-1}(T) = \dim S$. 
\end{proof}

\subsubsection{The $\mathcal{AS}$ topology} \label{subsubsecastop}

We are now ready to consider semialgebraic arc-symmetric subsets of $M$ and generalize theorem~\ref{theotopar} to Nash manifolds. First:

\begin{prop} \label{propsemialgsarcsymnashmanclosed} Let $S$ be a semialgebraic arc-symmetric subset of $M$. Then $S$ is closed in~$M$.
\end{prop}

\begin{proof} Let $x \in \overline{S}$ and let $\phi : U \rightarrow \R^n$ be a Nash chart on $M$ such that $x \in U$: since $x \in \overline{S}$ and $U$ is an open neighborhood of $x$ in $M$, $S \cap U \neq \emptyset$. Now, consider the semialgebraic subset $\phi(S \cap U)$ of $\R^n$ ($S$ is semialgebraic in $M$). The point $\phi(x)$ is in the closure of $\phi(S \cap U)$ in $\R^n$ (see the proof of lemma \ref{lembooleancomsaisasclosureissa}) and we can then apply Nash curve selection lemma (Proposition~8.1.13 of \cite{BCR}) to $\phi(S \cap U)$ and $x$: there exists a Nash arc $\gamma : ]-1;1[ \rightarrow \R^n$ such that $\gamma(0) = \phi(x)$ and $\gamma(]0;1[) \subset \phi(S \cap U)$. In particular, $0$ belongs to the open subset $\gamma^{-1}(\phi(U))$ of $]-1;1[$ and there exists $\epsilon >0$ such that $]-\epsilon,\epsilon[ \subset \gamma^{-1}(\phi(U))$, so that $\gamma(]-\epsilon;1[) \subset \phi(U)$. We can therefore consider the analytic arc $\widetilde{\gamma} := \phi^{-1} \circ \gamma : ]-\epsilon ; 1[ \rightarrow M$ on $M$. Since $\gamma(]0;1[) \subset \phi(S \cap U)$, we have $\widetilde{\gamma}(]0;1[) \subset S \cap U \subset S$ and then, because $S$ is arc-symmetric in $M$, $\widetilde{\gamma}(]-\epsilon;1[) \subset S$. In particular, $x = \widetilde{\gamma}(0) \in S$.
\end{proof}

As in the case of $\R^n$, the semialgebraic arc-symmetric subsets of the Nash manifold $M$ are the closed sets of a Noetherian topology on $M$:

\begin{theo} \label{theotopas} Any intersection or finite union of semialgebraic arc-symmetric subsets of~$M$ is a semialgebraic arc-symmetric subset of $M$. Furthermore, the topology on $M$ whose closed sets are exactly the semialgebraic arc-symmetric subsets of $M$ is a Noetherian topology.
\end{theo}

Analogously to theorem \ref{theotopar}, because of proposition \ref{propsemialgsarcsymnashmanclosed}, proposition \ref{propintersecfiniteunionclarcsymisclarcsym} and lemma \ref{lembooleancomsaisasclosureissa}, the following generalization of lemma \ref{lemprooftopar} will prove theorem \ref{theotopas}:

\begin{lem} \label{lemtoprovethtopas} Let $N$ be a connected Nash submanifold of $M$.
	\begin{enumerate}
		\item If $S$ is a semialgebraic arc-symmetric subset of $M$, then $N \subset S$ or $\dim (N \cap S) < \dim N$.
		\item If $\left(S_i\right)_{i \in I}$ is a family of semialgebraic arc-symmetric subsets of $M$, then there exist $i_1,\ldots,i_k \in I$ such that
		$$N \cap \bigcap_{i \in I} S_i = N \cap S_{i_1} \cap \cdots \cap S_{i_k}.$$	
	\end{enumerate}
\end{lem}

\begin{proof} 
	\begin{enumerate}
		\item Denote $d := \dim N$ and let $S$ be a semialgebraic arc-symmetric subset of $M$ such that $\dim (N \cap S) = d$. Consider a Nash chart $\phi : U \rightarrow \R^n$ on $M$ such that $\phi(N \cap S \cap U)$ is a semialgebraic set of dimension $d$ de $\R^n$. Since $\phi(N \cap S \cap U) \subset \phi(N \cap U)$ and $\dim N = d$, we necessarily have $\dim \phi(N \cap U) = d$. Therefore, by remark \ref{remafterlemprooftopar}, there exists a non-empty open subset $\widetilde{U}$ of $\phi(N \cap U)$ which is contained in $\phi(N \cap S \cap U)$.  The inverse image~$\phi^{-1}\left(\widetilde{U}\right)$ is then a non-empty open subset of $N \cap U$, hence of $N$ (since $U$ is open in $M$), which is contained in $N \cap S \cap U$, hence in $N \cap S$. As a consequence, the interior~$T$ of $N \cap S$ in $N$ is non-empty. It can then be shown that $T$ is furthermore closed in $N$ by using the same arguments as in the proof of lemma \ref{lemprooftopar} and the fact that $S$ is arc-symmetric in~$M$. Since $T$ is an open and closed subset of the connected topological space $N$, we have $N = T$ and then $N = N \cap S$ i.e. $N \subset S$.
	\item The statement can be proved by the same reasoning as in the proof of lemma \ref{lemprooftopar}, since a connected Nash submanifold of dimension $0$ is a point and any semialgebraic subset of~$M$ is a finite disjoint union of connected Nash submanifolds of~$M$ by proposition \ref{propsanashmanifisdisjunionofconashsub}.
\end{enumerate}
\end{proof}

\begin{proof}[Proof of theorem \ref{theotopas}]  Let us prove that, if $\left(S_i\right)_{i \in I}$ is a family of semialgebraic arc-symmetric subsets of $M$, the intersection $\bigcap_{i \in I} S_i$ is a finite one. Write $M$ as the union of its connected components $M_1,\ldots,M_s$, which are all connected Nash submanifolds of dimension $\dim M$ of $M$ (as open subsets of $M$). If $j \in \{1,\ldots,s\}$, there exist, according to previous lemma \ref{lemtoprovethtopas} 2, indices $i_{j,1},\ldots,i_{j,k_j} \in I$ such that 
$$M_j \cap \bigcap_{i \in I} S_i = M_j \cap S_{i_{j,1}} \cap \cdots \cap S_{i_{j,k_j}}.$$ 
Therefore, if we denote $J := \bigcup_{j =1}^s \left\{i_{j,1},\ldots,i_{j,k_j}\right\}$, the previous intersection is equal to $M_j \cap \bigcap_{i \in J} S_i$ and
$$\bigcap_{i \in I} S_i  = M \cap \bigcap_{i \in I} S_i = \bigcup_{j=1}^s \left(M_j \cap \bigcap_{i \in I} S_i\right)  =  \bigcup_{j=1}^s \left(M_j \cap \bigcap_{i \in J} S_i\right) = M \cap \bigcap_{i \in J} S_i.$$

In particular, if $\left(T_n\right)_{n \in \mathbb{N}}$ is a decreasing sequence of semialgebraic arc-symmetric subsets of~$M$, there exist $n_1,\ldots,n_k \in \mathbb{N}$ with $n_1 < \cdots < n_k$ such that 
$$\bigcap_{n \in \mathbb{N}} T_n = T_{n_1} \cap \cdots \cap T_{n_k} = T_{n_k}$$
and then, if $m \in \mathbb{N}$ satisfies $m \geq n_k$, we have $T_m = T_{n_k}$.
\end{proof}

\begin{de} The topology on $M$ whose closed sets are exactly the semialgebraic arc-symmetric subsets of $M$ is called the \emph{$\mathcal{AS}$ topology} on $M$. If $S$ is a subset of $M$, the closure of~$S$ with respect to $\mathcal{AS}$ topology is denoted by $\overline{S}^{\mathcal{AS}}$.
\end{de}

\begin{rem} The $\mathcal{AS}$ topology on the Nash manifold $\R^n$ is the $\mathcal{AR}$ topology on $\R^n$.
\end{rem}

Any irreducible $\mathcal{AS}$-closed set of $M$ will be called an \emph{$\mathcal{AS}$-irreducible set} and, since the~$\mathcal{AS}$ topology is Noetherian, any $\mathcal{AS}$-closed set of $M$ has a unique finite decomposition into $\mathcal{AS}$-irreducible sets that we will call its \emph{$\mathcal{AS}$-irreducible components}.

Like the $\mathcal{AR}$ topology on real affine spaces (proposition \ref{proparirredisconnected}), we have:

\begin{lem} Any $\mathcal{AS}$-irreducible set of $M$ is connected with respect to the (strong) topology of $M$.
\end{lem}

\begin{proof} Let $S$ be a $\mathcal{AS}$-closed irreducible set of $M$ and $S_1,\ldots,S_k$ be its connected components, which are semialgebraic (by corollary \ref{remconnectedcompsemialgnashman}) and arc-symmetric (see remark \ref{remremconncomparcsymnashmanarcsym}), i.e. $\mathcal{AS}$-closed in $M$. Therefore, since $S = \sqcup_{i=1}^k S_i$ is $\mathcal{AS}$-irreducible, we have $k=1$ and $S = S_1$ is connected.
\end{proof}

\begin{rem} Since there is no well-defined notion of ``real algebraic subset'' for a Nash manifold, we cannot consider the ``Zariski closure'' of a subset of $M$, whereas this was crucial to prove that the dimension of the $\mathcal{AR}$-closure of a semialgebraic subset of $\R^n$ is equal to the dimension of the latter.
\end{rem}

Recall that a compact Nash submanifold of $M$ is itself a Nash manifold (remark \ref{remcompactnashsubisnashman}) and it is therefore natural to investigate the relation between the two associated $\mathcal{AS}$ topology.

\begin{prop} \label{propcompactnashsubastopcoincide} Let $N$ be a compact Nash submanifold of $M$. The $\mathcal{AS}$ topology on the induced Nash manifold $N$ coincides with the topology induced on $N$ by the $\mathcal{AS}$ topology on $M$.
\end{prop}

\begin{proof} If $S$ is a subset of $N$, then $S$ is semialgebraic in $N$ if and only if $S$ is semialgebraic in~$M$ by lemma \ref{lemdirectiminverseimsaissa} 3 (applied to the inclusion map $N \hookrightarrow M$ which is Nash) and, because~$N$ is in particular a closed analytic submanifold of $M$, $S$ is arc-symmetric in $N$ if and only if $S$ is arc-symmetric in $M$ by lemma \ref{corarcsymsubsetclosedsubman}. As a consequence, the $\mathcal{AS}$-closed sets of the Nash manifold $N$ are exactly the subsets of $N$ obtained by intersecting the $\mathcal{AS}$-closed sets of $M$ with~$N$ ($N$ is itself an $\mathcal{AS}$-closed set of $M$). 
\end{proof}

As in the case of real affine spaces (proposition \ref{proparcanalyticmapsandarclosedsets}), we can establish the ``compatibility'' of $\mathcal{AS}$-closed sets of Nash manifolds with semialgebraic arc-analytic maps:

\begin{prop} \label{propsaarcanalmapcontastop} Let $S$ be an $\mathcal{AS}$-closed set of $M$, $T$ be a semialgebraic subset of a Nash manifold $N$ and $f : S \rightarrow T$ be a semialgebraic arc-analytic map.
\begin{enumerate}
	\item If $T$ is $\mathcal{AS}$-closed in $N$, the graph $\Gamma_f := \left\{(x,y) \in S \times T~|~y = f(x)\right\}$ of $f$ is $\mathcal{AS}$-closed in $M \times N$.
	\item For any $\mathcal{AS}$-closed set $A$ of $N$ included in $T$, the set $f^{-1}(A)$ is $\mathcal{AS}$-closed in $M$, i.e. $f$ is continuous with respect to (the topologies induced on $S$ and $T$ by) $\mathcal{AS}$ topology.
\end{enumerate}
In particular, a Nash diffeomorphism between Nash manifolds is an homeomorphism with respect to $\mathcal{AS}$ topology.
\end{prop}

\begin{proof} 
\begin{enumerate}
	\item By lemma \ref{lemsamapbetweensasetsofnashmanifolds} 1, the graph $\Gamma_f$ is semialgebraic in $M \times N$ since $f$ is semialgebraic. If $T$ is arc-symmetric in $N$, $\Gamma_f$ is furthermore arc-symmetric in $M \times N$ by proposition~\ref{propcomparcanalyticmapwitharcsymsetrealanalmanifold}~1.
	\item If $A$ is an $\mathcal{AS}$-closed set of $N$ included in $T$, $f^{-1}(A)$ is semialgebraic in $M$ by lemma~\ref{lemsamapbetweensasetsofnashmanifolds}~2 and arc-symmetric in $M$ by proposition \ref{propcomparcanalyticmapwitharcsymsetrealanalmanifold} 2.
\end{enumerate}
\end{proof}

\begin{rems}
~
\begin{enumerate}
	\item Suppose that $M$ is an \emph{affine} Nash manifold, i.e. that there exists a Nash diffeomorphism $\Psi$ from $M$ to a Nash submanifold $L$ of a real affine space $\R^s$ (we will call such a map an \emph{affine Nash embedding} of $M$). Suppose furthermore that $M$ is compact ($L$ is therefore a compact Nash submanifold of $\R^s$). Then, according to proposition \ref{propsaarcanalmapcontastop} and proposition \ref{propcompactnashsubastopcoincide}, if $S$ is a subset of $M$, the set~$S$ is $\mathcal{AS}$-closed in $M$ if and only if $\Psi(S)$ is $\mathcal{AR}$-closed in $\R^s$.
	\item Any semialgebraic arc-analytic map between $\mathcal{AS}$-closed sets of compact affine Nash manifolds is continuous with respect to the strong topologies induced by the latter (compose the considered map with respective affine Nash embeddings of the considered compact Nash manifolds and use proposition \ref{proparcanalyticmapsandarclosedsets} 1).
\end{enumerate}
\end{rems}

\subsubsection{Quasi-arc-symmetric sets and $\mathcal{AS}$-sets} \label{subsubsecqarcsymassets}

In this paragraph, we will consider, together with the semialgebraic arc-symmetric sets of the Nash manifold $M$, two more general classes of objets, namely the \emph{quasi-arc-symmetric sets} of~$M$ and the \emph{$\mathcal{AS}$-sets} of $M$. We will show that when $M$ is compact and affine, the two notions coincide. From now on, we will focus on semialgebraic subsets of $M$.

\begin{de} \label{defquasiarcsym} Let $S$ be a semialgebraic subset of $M$. We say that $S$ is a \emph{quasi-arc-symmetric set of $M$} if for any analytic arc $\gamma : ]-1;1[ \rightarrow M$, if $\gamma(]-1;0[) \subset S$, then there exists $\epsilon \in ]0;1[$ such that $\gamma(]0;\epsilon[) \subset S$.
\end{de}

\begin{rem} By arguments similar to the proof of lemma \ref{lemequivalentchararcsymrn}, a semialgebraic subset $S$ of~$M$ is quasi-arc-symmetric if and only if for any analytic arc $\gamma : ]\epsilon;\epsilon'[ \rightarrow M$,
$$\begin{cases}\mbox{ if $\gamma(]a,b[) \subset S$ with $\epsilon \leq a < b < \epsilon'$, then there exists $\xi \in ]b;\epsilon'[$ such that $\gamma(]b;\xi[) \subset S$,}\\ \mbox{ if $\gamma(]a,b[) \subset S$ with $\epsilon < a < b \leq \epsilon'$, then there exists $\xi \in ]\epsilon;a[$ such that $\gamma(]\xi;a[) \subset S$.}\end{cases}$$
\end{rem}

Of course, an $\mathcal{AS}$-closed set of $M$ is a quasi-arc-symmetric set, which is furthermore closed (for the strong topology of $M$) by proposition \ref{propsemialgsarcsymnashmanclosed}. The converse is also true:

\begin{lem} \label{lemquasiarcsymclosedas} A quasi-arc-symmetric set $S$ of $M$ is arc-symmetric in $M$ if and only if it is closed for the (strong) topology of $M$.
\end{lem}

\begin{proof} Suppose that $S$ is a closed subset of $M$ and let $\gamma : ]-1;1[ \rightarrow M$ be an analytic arc such that $\gamma(]-1;0[) \subset S$: since $S$ is quasi-arc-symmetric in $M$, there exists $\epsilon \in ]0;1[$ such that $\gamma(]0;\epsilon[) \subset S$. But $S$ is closed in $M$ so that $\gamma(0) \in \gamma\left(\overline{]-1;0[}\right) \subset \overline{\gamma(]-1;0[)} \subset S$. As a consequence, $\gamma(]-1;\epsilon[) \subset S$ and $S$ is then arc-symmetric in $M$.
\end{proof}

Like $\mathcal{AS}$-closed sets, the class of quasi-arc-symmetric sets of $M$ is stable under finite intersection and finite union. But it is stable under set-theoritic difference as well:

\begin{prop} \label{propasisquasiarcsym} Any finite Boolean combination of quasi-arc-symmetric subsets of $M$ is quasi-arc-symmetric in $M$.
\end{prop}

\begin{proof} Recall lemma \ref{lembooleancomsaisasclosureissa} and let $A$ and $B$ be quasi-arc-symmetric subsets of $M$. Let $\gamma : ]-1;1[ \rightarrow M$ be an analytic arc.

First, suppose that $\gamma(]-1;0[) \subset A \cap B$. Since $A$ and $B$ are quasi-arc-symmetric, there exist $\epsilon, \epsilon' \in ]0;1[$ such that $\gamma(]0;\epsilon[) \subset A$ and $\gamma(]0;\epsilon'[) \subset B$ and then, if we denote by $\epsilon''$ the minimum of $\epsilon$ and $\epsilon'$, we have $\gamma(]0;\epsilon''[) \subset A \cap B$. As a consequence, $A \cap B$ is quasi-arc-symmetric in $M$.

Now, suppose that $\gamma(]-1;0[) \subset A \cup B$. By lemma \ref{lemforunionquasiarcsymisquasiarcsym} below, there exists $\xi \in ]-1;0[$ such that $\gamma(]\xi;0[) \subset A$ or there exists $\xi \in ]-1;0[$ such that $\gamma(]\xi;0[) \subset B$. As a consequence, since $A$ and $B$ are quasi-arc-symmetric in $M$, there exists $\epsilon \in ]0;1[$ such that $\gamma(]0;\epsilon[) \subset A \subset A \cup B$ or there exists $\epsilon \in ]0;1[$ such that $\gamma(]0;\epsilon[) \subset B \subset A \cup B$: the union $A \cup B$ is therefore quasi-arc-symmetric in $M$.

Finally, suppose that $\gamma(]-1;0[) \subset A \setminus B$. Since $A$ is quasi-arc-symmetric in $M$, there exists $\epsilon_0 \in ]0;1[$ such that $\gamma(]0;\epsilon_0[) \subset A$. Now, suppose by absurd that for any $\epsilon \in ]0;\epsilon_0[$, there exists $t_{\epsilon} \in ]0;\epsilon[$ such that $\gamma(t_{\epsilon}) \in B$: we can adapt the arguments of the proof of lemma \ref{lemforunionquasiarcsymisquasiarcsym} below to assert that there then exists $\epsilon' \in ]0;\epsilon_0[$ such that $]0;\epsilon'[ \subset \gamma^{-1}(B)$. Therefore, since $B$ is quasi-arc-symmetric in $M$, there is $\xi \in ]-1;0[$ such that $\gamma(]\xi;0[) \subset B$. But this contradicts the inclusion $\gamma(]-1;0[) \subset A \setminus B$. Consequently, there exists $\epsilon \in ]0;\epsilon_0[$ such that $\gamma(]0;\epsilon[) \subset A \setminus B$ and, as a result, $A \setminus B$ is quasi-arc-symmetric in $M$.
\end{proof}

\begin{lem} \label{lemforunionquasiarcsymisquasiarcsym} Let $\gamma : ]-1;1[ \rightarrow M$ be an analytic arc and let $A$ and $B$ be distinct semialgebraic subsets of $M$ such that $\gamma(]-1;0[) \subset A \cup B$. Then there exists $\xi \in ]-1;0[$ such that $\gamma(]\xi;0[) \subset A$ or there exists $\xi \in ]-1;0[$ such that $\gamma(]\xi;0[) \subset B$.
\end{lem}

\begin{proof} Suppose that the former is not satisfied: for any $\xi \in ]-1;0[$, there exists $t_{\xi} \in ]\xi;0[$ such that $\gamma(t_{\xi}) \in B \setminus A$ (we have $\gamma(]-1;0[) \subset A \cup B$). For $n \in \Nstar$, let $t_n \in ]-\frac{1}{n};0[$ such that $\gamma(t_n) \in B \setminus A$, and consider a Nash chart $\phi : U \rightarrow \R^n$ on $M$ such that $\gamma(0) \in U$. Since $0$ is contained in the open subset $\gamma^{-1}(U)$ of $]-1;1[$, there exists $\alpha \in ]0;1[$ such that $]-\alpha;\alpha[ \subset \gamma^{-1}(U)$ and, since $(t_n)_{n \in \Nstar}$ converges to~$0$, there exists $n_0 \in \Nstar$ such that for any $n \geq n_0$, $t_{n} \in ]-\alpha;\alpha[$ and then $\gamma(t_n) \in U$: we have $(t_n)_{n \geq n_0} \subset \gamma^{-1}(B \cap U)$. 

Denote by $T$ the semialgebraic subset $\phi(B \cap U)$ of $\R^n$: $T$ is a finite union of semialgebraic subsets of $\R^n$ of the form $\{p = 0, q_1 > 0, \ldots, q_p > 0\}$ where $p, q_1,\ldots,q_p$ are polynomial functions on $\R^n$, and there exists a subsequence $(t_{n_k})_{k \in \N}$ of $(t_{n})_{n \geq n_0}$ such that $(\gamma(t_{n_k}))_{k \in \N}$ is included in one such set. Keeping the above notations, we therefore have, for any $k \in \N$, $p \circ \phi \circ \gamma(t_{n_k}) = 0$ and $q_i \circ \phi \circ \gamma(t_{n_k}) >0$, $i \in \{1,\ldots,p\}$. Since the analytic map $p \circ \phi \circ \gamma : ]-\alpha;\alpha[ \rightarrow \R$ vanishes on the sequence $(t_{n_k})_{k \in \N}$ with accumulation point $0$, $p \circ \phi \circ \gamma$ vanishes on the entire interval $]-\alpha;\alpha[$. 

Furthermore, if $i \in \{1,\ldots,p\}$, we clame that there exists $k_i \in \N$ such that for all $t \in ]t_{n_{k_i}};0[$, $q_i \circ \phi \circ \gamma(t) >0$. Indeed, if for all $k \in \N$, there was $\tau_k \in ]t_{n_k};0[$ such that $q_i \circ \phi \circ \gamma(\tau_k) = 0$ (the map $q_i \circ \phi \circ \gamma$ is continuous and is positive at $t_{n_k}$), since the sequence $(\tau_{k})_{k \in \N}$ converges to the accumulation point $0$, the analytic map $q_i \circ \phi \circ \gamma$ would vanish on the entire interval $]-\alpha;\alpha[$, which is not the case. As a consequence, if we let $\xi_0$ denote the maximum of the numbers $t_{n_{k_i}}$, $i \in \{1,\ldots,p\}$, we have 
$$]\xi_0;0[ \subset (\phi \circ \gamma)^{-1}\big(\{p = 0, q_1 > 0, \ldots, q_p > 0\}\big) \subset (\phi \circ \gamma)^{-1}(T) = \gamma^{-1}(B \cap U).$$
In particular, $]\xi_0;0[ \subset \gamma^{-1}(B)$.
\end{proof}

\begin{rem} \label{remcartprodquasiarcsymsisquasiarcsym} If $S$ is a quasi-arc-symmetric set of $M$ and $T$ is a quasi-arc-symmetric set of a Nash manifold $N$, then the cartesian product $S \times T$ is a quasi-arc-symmetric set of the Nash manifold $M \times N$. Indeed, we have lemma \ref{lembooleancomsaisasclosureissa} and, if $\gamma = (\gamma_1,\gamma_2) : ]-1;1[ \rightarrow M \times N$ is an analytic arc satisfying $\gamma(]-1;0[) \subset S \times T$ i.e. $\gamma_1(]-1;0[) \subset S$ and $\gamma_2(]-1;0[) \subset T$, then, since $S$ and $T$ are quasi-arc-symmetric, there exist $\epsilon_1 >0$, $\epsilon_2 > 0$ such that $\gamma_1(]0;\epsilon_1[) \subset S$ and $\gamma_2(]0;\epsilon_2[) \subset T$. As a consequence, if we set $\epsilon$ to be the minimum of $\epsilon_1$ and $\epsilon_2$, we have $\gamma(]0;\epsilon[) \subset S \times T$.
\end{rem}

Quasi-arc-symmetric sets are also stable under inverse image by semialgebraic arc-analytic maps:

\begin{lem} \label{leminverseimagequasiarcsymbysaarcanalisquasiarcsym} Let $S$ be an $\mathcal{AS}$-closed set of $M$, $T$ be a semialgebraic set of a Nash manifold~$N$ and $f : S \rightarrow T$ be a semialgebraic arc-analytic map. If $A$ is a quasi-arc-symmetric set of $N$ included in $T$, then $f^{-1}(A)$ is quasi-arc-symmetric in $M$.
\end{lem}

\begin{proof} If $\gamma : ]-1;1[ \rightarrow M$ is an analytic arc such that $\gamma(]-1;0[) \subset f^{-1}(A) \subset S$, then $\gamma(]-1;1[) \subset S$ (since $S$ is arc-symmetric in $M$), the analytic arc $f \circ \gamma : ]-1;1[ \rightarrow N$ satisfies $f \circ \gamma(]-1;0[) \subset A$ and, since $A$ is quasi-arc-symmetric in $N$, there exists $\epsilon \in ]0;1[$ such that $f \circ \gamma(]0;\epsilon[) \subset A$ i.e. $\gamma(]0;\epsilon[) \subset f^{-1}(A)$.
\end{proof}

We also have the quasi-arc-symmetric analog of proposition \ref{propcompactnashsubastopcoincide}:

\begin{prop} \label{remquasiarcsymcompactnashsubsemialaamap} Let $N$ be a compact Nash submanifold of $M$ and $S$ be a subset of $N$. The set $S$ is quasi-arc-symmetric in the Nash manifold $N$ if and only if $S$ is quasi-arc-symmetric in $M$.
\end{prop}

\begin{proof} Suppose that $S$ is quasi-arc-symmetric in $N$ and let $\gamma : ]-1;1[ \rightarrow M$ be an analytic arc such that $\gamma(]-1;0[) \subset S \subset N$: since $N$ is $\mathcal{AS}$-closed in $M$, we have $\gamma(]-1;1[) \subset N$ and we can restrict $\gamma$ into an arc $\widetilde{\gamma} : ]-1;1[ \rightarrow N$, which is analytic with respect to the real analytic manifold~$N$ since $N$ is an analytic submanifold of $M$.  As a consequence, because $S$ is quasi-arc-symmetric in $N$ and $\widetilde{\gamma}(]-1;0[) = \gamma(]-1;0[) \subset S$, there exists $\epsilon \in ]0;1[$ such that $\widetilde{\gamma}(]0;\epsilon[) \subset S$ i.e. $\gamma(]0;\epsilon[) \subset S$.

The converse implication is given by previous lemma \ref{leminverseimagequasiarcsymbysaarcanalisquasiarcsym} applied to the inclusion map $N \hookrightarrow M$.
\end{proof}

We now state a result which will reveal crucial to show that the notion of quasi-arc-symmetric set and the notion of $\mathcal{AS}$-set defined below coincide when $M$ is a \emph{compact affine} Nash manifold.

\begin{lem} \label{lemcomplementassetinasclos} Let $S$ be a quasi-arc-symmetric subset of $M$. We have $\overline{S}^{\mathcal{AS}} = S \cup \overline{\overline{S} \setminus S}^{\mathcal{AS}}$, so that $\overline{S}^{\mathcal{AS}} \setminus S \subset \overline{\overline{S} \setminus S}^{\mathcal{AS}}$.
\end{lem}

\begin{proof} The proof is the same as in \cite{FMI} (Proposition 1.6). The quasi-arc-symmetric set 
$$S \cup \overline{\overline{S} \setminus S}^{\mathcal{AS}} = S \cup (\overline{S} \setminus S) \cup  \overline{\overline{S} \setminus S}^{\mathcal{AS}} = \overline{S} \cup \overline{\overline{S} \setminus S}^{\mathcal{AS}}$$
is closed for the strong topology of $M$ (proposition \ref{propsemialgsarcsymnashmanclosed}): by lemma \ref{lemquasiarcsymclosedas}, it is therefore an~$\mathcal{AS}$-closed set of $M$. Since it contains $S$, we have $\overline{S}^{\mathcal{AS}} \subset S \cup \overline{\overline{S} \setminus S}^{\mathcal{AS}}$.

Conversely, since $\overline{S} \setminus S \subset \overline{S}$ and $\overline{S} \subset \overline{S}^{\mathcal{AS}}$ (because $\overline{S}^{\mathcal{AS}}$ is closed in $M$ and contains $S$), we have $\overline{S} \setminus S \subset \overline{S}^{\mathcal{AS}}$ and then $\overline{\overline{S} \setminus S}^{\mathcal{AS}} \subset \overline{S}^{\mathcal{AS}}$, so that $S \cup \overline{\overline{S} \setminus S}^{\mathcal{AS}} \subset \overline{S}^{\mathcal{AS}}$.
\end{proof}

\begin{de} Let $S$ be a subset of $M$. We say that $S$ is an \emph{$\mathcal{AS}$-set of $M$} if $S$ is a finite Boolean combination of semialgebraic arc-symmetric subsets of $M$.
\end{de}

Remark that the class of $\mathcal{AS}$-sets of $M$ is stable under finite Boolean combination by definition and that, by proposition \ref{propasisquasiarcsym}, an $\mathcal{AS}$-set of $M$ is quasi-arc-symmetric in $M$. We are going to show that the converse is true when $M$ is a real affine space. Thereafter, we will show the equivalence of the two notions when $M$ is supposed compact and affine, using the following analogs of lemma \ref{leminverseimagequasiarcsymbysaarcanalisquasiarcsym} and proposition \ref{remquasiarcsymcompactnashsubsemialaamap}: 

\begin{lem} \label{remassetcompactsubnashsaaamap}
~
\begin{enumerate}
	\item If $f$ is a semialgebraic arc-analytic map from an $\mathcal{AS}$-closed set $S$ of $M$ to a semialgebraic set $T$ of a Nash manifold $N$ and $A$ is an $\mathcal{AS}$-set of $N$ included in $T$, then the set $f^{-1}(A)$ is an $\mathcal{AS}$-set of $M$.
	\item If $N$ is a compact Nash submanifold of $M$ and $S$ is a subset of $N$, then $S$ is an $\mathcal{AS}$-set of $N$ if and only if $S$ is an $\mathcal{AS}$-set of $M$.
\end{enumerate}
\end{lem}

\begin{proof} The two statements are direct consequences of proposition \ref{propcompactnashsubastopcoincide} and proposition \ref{propsaarcanalmapcontastop}.
\end{proof}

\begin{prop} \label{propquasiarcsymrniffasrn} A subset of $\mathbb{R}^n$ is quasi-arc-symmetric in $\mathbb{R}^n$ if and only if it is an $\mathcal{AS}$-set of $\mathbb{R}^n$.
\end{prop}

\begin{proof} We show that any quasi-arc-symmetric set of $\R^n$ is an $\mathcal{AS}$-set of $\R^n$ by induction on dimension.

A quasi-arc-symmetric set of $\mathbb{R}^n$ of dimension $0$ is a finite union of points, so let $S$ be a quasi-arc-symmetric subset of $\mathbb{R}^n$ of dimension larger than $0$ and write
$$S = \overline{S}^{\mathcal{AR}} \setminus \left(\overline{S}^{\mathcal{AR}} \setminus S\right).$$
The set-theoritic difference of quasi-arc-symmetric sets of $\R^n$ is quasi-arc-symmetric in $\R^n$ (proposition \ref{propasisquasiarcsym}) and, by lemma \ref{lemcomplementassetinasclos}, we have $\overline{S}^{\mathcal{AR}} \setminus S \subset \overline{\overline{S} \setminus S}^{\mathcal{AR}}$ (recall that the $\mathcal{AS}$ topology of $\R^n$ is its $\mathcal{AR}$ topology), so that
$$\dim \overline{S}^{\mathcal{AR}} \setminus S \leq \dim \overline{\overline{S} \setminus S}^{\mathcal{AR}} = \dim \overline{S} \setminus S < \dim S:$$
the equality is given by remark \ref{remeqdimarclosuresasa} (that we do not have in the general Nash manifold context) while we refer for instance to Proposition 2.8.13 of \cite{BCR} for the right-hand side inequality (which is also true in the general Nash manifold context: see lemma \ref{lemdimclosminussa}). We can therefore assert by induction that $\overline{S}^{\mathcal{AR}} \setminus S$ is a finite Boolean combination of $\mathcal{AR}$-closed sets of $\mathbb{R}^n$, and then so is $S$.
\end{proof}

\begin{cor} \label{corcompactnashmaniembquasiarcsymas} Suppose that $M$ is a compact affine Nash manifold. Then a subset of $M$ is quasi-arc-symmetric in $M$ if and only if it is an $\mathcal{AS}$-set of $M$.
\end{cor}

\begin{proof} Let $\Phi$ be an affine Nash embedding of $M$ onto a Nash submanifold $N$ of a real affine space $\R^s$ (notice that $N$ is compact as well) and let $S$ be a subset of $M$. Then $S$ is quasi-arc-symmetric in $M$ if and only if $\Phi(S)$ is quasi-arc-symmetric in $N$ (lemma \ref{leminverseimagequasiarcsymbysaarcanalisquasiarcsym}) if and only if~$\Phi(S)$ is quasi-arc-symmetric in $\R^s$ (proposition \ref{remquasiarcsymcompactnashsubsemialaamap}) if and only if $\Phi(S)$ is an $\mathcal{AS}$-set of $\R^s$ (proposition \ref{propquasiarcsymrniffasrn}) if and only if $\Phi(S)$ is an $\mathcal{AS}$-set of $N$ (lemma \ref{remassetcompactsubnashsaaamap} 2) if and only if $S$ is an~$\mathcal{AS}$-set $M$ (lemma \ref{remassetcompactsubnashsaaamap} 1).
\end{proof}

\begin{rems} \label{remsafterquasiarcsymcoincideascompactaffinenashmanifold}
~
	\begin{enumerate}
		\item According to lemmas \ref{leminverseimagequasiarcsymbysaarcanalisquasiarcsym} and \ref{remassetcompactsubnashsaaamap} 1 and the arguments that we used in the above proof, the statement of corollary \ref{corcompactnashmaniembquasiarcsymas} is more generally true if there exists a bijective semialgebraic arc-analytic map with arc-analytic inverse (recall lemma~\ref{lemsamapbetweensasetsofnashmanifolds}~3) from $M$ onto a compact Nash submanifold of a real affine space.
		\item Suppose that $M$ is a compact affine Nash manifold. Then, if $S$ is any semialgebraic subset of $M$, we have $\dim \overline{S}^{\mathcal{AS}} = \dim S$, so that we could have used the same arguments as in the proof of proposition \ref{propquasiarcsymrniffasrn} to prove the statement of corollary \ref{corcompactnashmaniembquasiarcsymas}. Indeed, if~$\Phi$ is an affine Nash embedding of $M$ onto a Nash submanifold $N$ of a real affine space,~$\Phi$ is in particular an homeomorphism with respect to $\mathcal{AS}$ topology and we have, using lemma~\ref{lemdimimagesamap},
$$\dim \overline{S}^{\mathcal{AS}} = \dim \Phi\left(\overline{S}^{\mathcal{AS}}\right) = \dim \overline{\Phi\left(S\right)}^{\mathcal{AR}} = \dim \Phi(S) = \dim S.$$
	\end{enumerate}
\end{rems}

\begin{ex} \label{exrealprojspacecompactaffnashmanifasquasiarcsymcoincide} The real projective space $\mathbb{P}^n(\R)$ is a compact affine Nash manifold (any real projective space, with its usual atlas, is actually a compact affine nonsingular real algebraic variety: see section 3.4 of \cite{BCR}). As a consequence, the quasi-arc-symmetric subsets of $\mathbb{P}^n(\R)$ are~$\mathcal{AS}$-sets and we recover the characterization of Proposition 3.2 of \cite{KurPar}, as well as Proposition~3.3 (by lemma \ref{lemcomplementassetinasclos}, remark \ref{remsafterquasiarcsymcoincideascompactaffinenashmanifold} 2 above and lemma \ref{lemdimclosminussa}).
\end{ex}

Let us also show the following characterization of $\mathcal{AS}$-sets of real affine spaces or compact affine Nash manifolds that we will use in next paragraph:

\begin{lem} \label{lemcharacassetnashinjectarcs} Suppose that $M$ is $\R^n$ or a compact affine Nash manifold and let $S$ be a semialgebraic subset of $M$. Then $S$ is an $\mathcal{AS}$-set of $M$ if and only if for any \emph{injective Nash} arc $\gamma : ]-1;1[ \rightarrow M$, if $\gamma(]-1;0[) \subset S$, then there exists $\epsilon \in ]0;1[$ such that $\gamma(]0;\epsilon[) \subset S$.
\end{lem}

\begin{proof} Remark that, if we replace analytic arcs with injective Nash arcs in the definition of quasi-arc-symmetric sets, proposition \ref{propasisquasiarcsym} remains true as well as lemma \ref{leminverseimagequasiarcsymbysaarcanalisquasiarcsym} for injective semialgebraic arc-analytic maps, and then proposition \ref{remquasiarcsymcompactnashsubsemialaamap} is also true.

Now, let $B$ be a closed semialgebraic subset of $M$ satisfying the condition of the lemma. If $M = \R^n$, let $\Psi$ denote the identity map on $\R^n$, and if $M$ is a compact affine Nash manifold, let~$\Psi$ be an affine Nash embedding of $M$ onto a compact Nash submanifold of an affine space~$\R^s$. In both cases, the closed semialgebraic subset $\Psi(B)$ of $\R^n$, resp. $\R^s$, satisfies, by the above considerations, the same condition as a subset of $\R^n$, resp. $\R^s$, and is then, according to the proof of lemma~\ref{lemquasiarcsymclosedas} and the characterization of $\mathcal{AR}$-closed sets of $\R^s$ of lemma \ref{lemcharsaarcsymrninjnasharcs}, an $\mathcal{AR}$-closed set of $\R^n$, resp.~$\R^s$. By propositions \ref{propcompactnashsubastopcoincide} and \ref{propsaarcanalmapcontastop}, $B$ is therefore an $\mathcal{AS}$-closed set of $M$. This proves the ``injective Nash arcs'' version of lemma \ref{lemquasiarcsymclosedas} which implies the corresponding analog of lemma~\ref{lemcomplementassetinasclos}. 

As a result, we can use the same arguments as in the proof of proposition \ref{propquasiarcsymrniffasrn} to assert the equivalence of the lemma: see remark \ref{remsafterquasiarcsymcoincideascompactaffinenashmanifold} 2.
\end{proof}

In the last subsection, we will focus on $\mathcal{AS}$-sets of real projective spaces, as in section 3 of \cite{KurPar}. Before, let us consider maps between $\mathcal{AS}$-sets of compact affine Nash manifolds which have $\mathcal{AS}$-graph.

\subsubsection{$\mathcal{AS}$-maps}

In the following, we will consider only \emph{compact affine} Nash manifolds, on which we proved that quasi-arc-symmetric sets and $\mathcal{AS}$-sets form one single notion: in accordance with \cite{KurPar}, we will just speak of $\mathcal{AS}$-sets in this context.

So suppose that $M$ is a compact affine Nash manifold and let $N$ be a compact affine Nash manifold. We are going to study maps between~$\mathcal{AS}$-sets of $M$ and $N$ which generalize semialgebraic arc-analytic maps between $\mathcal{AS}$-closed sets:

\begin{de} \label{defasmap} Let $S$ be an $\mathcal{AS}$-set of $M$, $T$ be an $\mathcal{AS}$-set of $N$ and $f : S \rightarrow T$. We say that $f$ is an \emph{$\mathcal{AS}$-map} if its graph $\Gamma_f$ is an $\mathcal{AS}$-set of $M \times N$ (notice that the Cartesian product~$M \times N$ is a compact affine Nash manifold as well).
\end{de}

\begin{exs} \label{exsasmaps}
~
	\begin{enumerate}
		\item By proposition \ref{propsaarcanalmapcontastop} 1, a semialgebraic arc-analytic map from an $\mathcal{AS}$-closed set of $M$ to an $\mathcal{AS}$-closed set of $N$ is an $\mathcal{AS}$-map.
		\item If $f : S \rightarrow T$ is an $\mathcal{AS}$-map between respective $\mathcal{AS}$-sets $S$ and $T$ of $M$ and $N$, if $B$ is an $\mathcal{AS}$-set of $M$ included in $S$ and $A$ is an $\mathcal{AS}$-set of $N$ included in $T$ such that $f(B) \subset A$, then the map $\widetilde{f} : B \rightarrow A~;~x \mapsto f(x)$ induced by $f$ is an $\mathcal{AS}$-map as well, since $\Gamma_{\widetilde{f}} = \Gamma_f \cap (B \times A)$ (see remark \ref{remcartprodquasiarcsymsisquasiarcsym}).
	\end{enumerate}
\end{exs}

\begin{rems} \label{remasmapembeddings} Keep the notations of definition \ref{defasmap}. 
	\begin{enumerate}
		\item If $f$ is a bijective $\mathcal{AS}$-map, then the inverse $f^{-1} : T \rightarrow S$ of $f$ is an $\mathcal{AS}$-map as well, since the graph $\Gamma_{f^{-1}}$ of $f^{-1}$ is the image of $\Gamma_f$ by the Nash diffeomorphism $M \times N \rightarrow N \times M~;~(x,y) \mapsto (y,x)$ (see lemma \ref{leminverseimagequasiarcsymbysaarcanalisquasiarcsym} or lemma \ref{remassetcompactsubnashsaaamap} 1).
		\item If $\Phi$ is a Nash embedding of $M$ into a real affine space $\R^n$ and $\Psi$ is a Nash embedding of $N$ into an affine space $\R^m$, then the map $f$ is an $\mathcal{AS}$-map if and only if the map $\widetilde{f} := \Psi \circ f \circ \Phi^{-1} : \Phi(S) \rightarrow \Psi(T)$ is an $\mathcal{AS}$-map if and only if the graph $\Gamma_{\widetilde{f}} = (\Phi \times \Psi)\left(\Gamma_f \right)$ of $\widetilde{f}$ is an $\mathcal{AS}$-set of $\R^{n+m}$ (see lemma \ref{leminverseimagequasiarcsymbysaarcanalisquasiarcsym} or lemma \ref{remassetcompactsubnashsaaamap} 1).
	\end{enumerate}
\end{rems}

We established in previous paragraph that the class of $\mathcal{AS}$-sets of compact affine Nash manifolds is stable under inverse image by semialgebraic arc-analytic maps from an $\mathcal{AS}$-closed set to a semialgebraic set. In fact, it is more generally stable under inverse image by $\mathcal{AS}$-maps: we will show this result after having proved that the image of an $\mathcal{AS}$-set by an \emph{injective} $\mathcal{AS}$-map is an $\mathcal{AS}$-set as well.

\begin{theo} \label{thimageofassetbyinjasmapisas} Let $S$ be an $\mathcal{AS}$-set of $M$, $T$ be an $\mathcal{AS}$-set of $N$ and $f : S \rightarrow T$ be an $\mathcal{AS}$-map. If $f$ is injective, then $f(S)$ is an $\mathcal{AS}$-set of $N$.
\end{theo}

\begin{proof} According to remark \ref{remasmapembeddings} 2 and lemma \ref{remassetcompactsubnashsaaamap}, we can suppose that $M$, resp. $N$, is a compact Nash submanifold of an affine space $\R^n$, resp. $\R^m$, and we have to prove that $f(S)$ is an $\mathcal{AS}$-set of $\R^m$.

By virtue of lemma \ref{lemcharacassetnashinjectarcs}, let $\gamma : ]-1;1[ \rightarrow \R^m$ be an injective Nash arc such that $\gamma(]-1;0[) \subset f(S)$. Denote by $\widetilde{f}$ the bijective semialgebraic map $S \rightarrow f(S)~;~x \mapsto f(x)$ induced by $f$, by $\widetilde{\gamma}$ the injective semialgebraic arc $\widetilde{\gamma} := \widetilde{f}^{-1} \circ \gamma : ]-1;0[ \rightarrow S$ and consider the injective semialgebraic arc 
$$\big(\widetilde{\gamma}, \gamma\big) : \begin{array}{ccc}]-1;0[ & \rightarrow & S \times T\\ t & \mapsto &\big(\widetilde{\gamma}(t), \gamma(t)\big).\end{array}$$
Since $S \times T$ is included in a compact subset of $\R^{n+m}$, there exist a sequence $\left(t_{k}\right)_{k \in \N}$ and a point $(a,b) \in \R^{n+m}$ such that $\big(\widetilde{\gamma}(t_k), \gamma(t_k)\big) \xrightarrow[k \rightarrow + \infty]{} (a,b)$. By the curve-selection lemma (see for instance Theorem 2.5.5 of \cite{BCR}), there then exists a continuous semialgebraic map $\lambda = (\lambda_1,\lambda_2) : [0;1[ \rightarrow \R^{n} \times \R^m$ such that $\lambda(0) = (a,b)$ and $\lambda(]0;1[) \subset (\widetilde{\gamma}, \gamma)(]-1;0[)$. Furthermore, by Lemma 2.6 of \cite{KurPar} and its proof, there exist a positive rational number $q$, a positive real number $\epsilon$ lower than~$1$ and an injective analytic curve $\eta = (\eta_1,\eta_2) : ]-\epsilon;\epsilon[ \rightarrow \R^{n} \times \R^m$ such that, for all $t \in [0;\epsilon[$, $\lambda(t^q) = \eta(t)$.

We then have $\eta(]0;\epsilon[) \subset \lambda(]0;1[) \subset (\widetilde{\gamma}, \gamma)(]-1;0[) \subset \Gamma_f$ and consequently, since $\Gamma_f$ is an~$\mathcal{AS}$-set, there exists a negative number $\delta > -\epsilon$ such that $\eta(]\delta;0[) \subset \Gamma_f$, i.e. for any $t \in ]\delta;0[$, $\eta_2(t) = f(\eta_1(t))$. 

On the other hand, $\eta_2([0;-\delta[) \subset \lambda_2([0;1[) \subset \gamma(]-1;0])$ (notice that $\gamma(0) = b = \lambda_2(0) = \eta_2(0)$ since $\gamma$ is Nash) so that the $1$-dimensional irreducible analytic germs $\big(\eta_2(]\delta;-\delta[),0\big)$ and $\big(\gamma(]-1;1[),0\big)$ (see Lemma 2.6 (b) of $\cite{KurPar}$)  have an intersection germ of positive dimension and are therefore equal (see also the proof of Théorème 4.2 in \cite{KurAR}, as well as the proof of Theorem 2.33 of \cite{KurPar}): there exists an open neighborhood $U$ of $b$ in $\R^m$ such that $\eta_2(]\delta;-\delta[) \cap U = \gamma(]-1;1[) \cap U$. In particular, there exists a positive real number $\alpha < 1$ such that $\gamma(]-\alpha;\alpha[) \subset \eta_2(]\delta;-\delta[)$.

As a result, if $t \in ]0;\alpha[$, there exists $s \in ]\delta;-\delta[ \setminus \{0\}$ such that $\gamma(t) = \eta_2(s) = f(\eta_1(s)) \in f(S)$ (if $t \neq 0$, $\gamma(t)$ cannot be equal to $\eta_2(0) = \lambda_2(0) = b = \gamma(0)$ by injectivity of $\gamma$).
\end{proof}

\begin{cor} \label{corinverseimageassetbyasmapisasset} Let $S$ be an $\mathcal{AS}$-set of $M$, $T$ be an $\mathcal{AS}$-set of $N$ and $f : S \rightarrow T$ be an~$\mathcal{AS}$-map. If $A$ is an $\mathcal{AS}$-set of $N$ included in $T$, then $f^{-1}(A)$ is an $\mathcal{AS}$-set of $M$.
\end{cor}

\begin{proof} We adapt the proof of Proposition 2.2.7 of \cite{BCR}: the set $f^{-1}(A)$ is the image of $(S \times A) \cap \Gamma_f$ by the map $\pi : (S \times A) \cap \Gamma_f \rightarrow T~;~(x,y) \mapsto x$. On the one hand, the former set is an $\mathcal{AS}$-set of~$M \times N$ by definition and by remark \ref{remcartprodquasiarcsymsisquasiarcsym}. On the other hand, the map $\pi$ is injective and is moreover an $\mathcal{AS}$-map, as the restriction of the Nash projection map $M \times N \rightarrow N~;~(x,y) \mapsto x$ (see examples \ref{exsasmaps}). Therefore, $f^{-1}(A)$ is an $\mathcal{AS}$-set of $M$ by theorem \ref{thimageofassetbyinjasmapisas}.
\end{proof}

Another corollary of theorem \ref{thimageofassetbyinjasmapisas} is that the composition of $\mathcal{AS}$-maps is an $\mathcal{AS}$-map:

\begin{cor} \label{corcompasmapsisasmap} Keeping the notations of corollary \ref{corinverseimageassetbyasmapisasset}, let $L$ be a compact affine Nash manifold, $X$ be an $\mathcal{AS}$-set of $L$ and $h : T \rightarrow X$ be an $\mathcal{AS}$-map. The composition $h \circ f : S \rightarrow X$ is an $\mathcal{AS}$-map.
\end{cor}

\begin{proof} We adapt the argument of remark \ref{remsrestrsamapcompasmapissa} 2: the graph $\Gamma_{h \circ f}$ of $h \circ f$ is the image of the~$\mathcal{AS}$ intersection $(\Gamma_f \times Z) \cap (X \times \Gamma_h)$ by the injective $\mathcal{AS}$-map 
$$(\Gamma_f \times Z) \cap (X \times \Gamma_h) \rightarrow S \times X~;~(x,y,z) \mapsto (x,z)$$
(induced by the Nash projection map $M \times N \times L \rightarrow M \times L~;~(x,y,z) \mapsto (x,z)$), so that we can apply theorem \ref{thimageofassetbyinjasmapisas} to conclude.
\end{proof}

\begin{rems} \label{remafterthcorasmaps}
~
\begin{enumerate}
	\item Keeping again the notations of corollary \ref{corinverseimageassetbyasmapisasset}, if $L$ is a compact affine Nash manifold, $X$ is an $\mathcal{AS}$-set of $L$ and $h : S \rightarrow X$ is an $\mathcal{AS}$-map, then the map $S \rightarrow T \times X~;~x \mapsto \big(f(x),h(x)\big)$ is an $\mathcal{AS}$-map. Indeed, the graph of the latter map is the intersection of $\Gamma_f \times X$ with the image of $\Gamma_h \times T$ by the Nash diffeomorphism $M \times L \times N \rightarrow M \times N \times L~;~(x,z,y) \mapsto (x,y,z)$.
	\item Theorem \ref{thimageofassetbyinjasmapisas} is not true for $\mathcal{AS}$-sets of real affine spaces (recall that the classes of $\mathcal{AS}$-sets and quasi-arc-symmetric sets of $\R^n$ coincide by proposition \ref{propquasiarcsymrniffasrn}). If we consider for instance (see Remark 3.6 of \cite{KurPar}) the arc-symmetric half-hyperbola $C := \{(x,y) \in \R^2~|~xy = 1, x> 0\}$ of $\R^2$ and the injective map with arc-symmetric graph $\pi : C \rightarrow \R~;~(x,y) \mapsto x$, the image $\pi(C) = ]0;+\infty[$ of~$C$ by $\pi$ is not an $\mathcal{AS}$-set of $\R$ (the analytic arc $\gamma : t \in ]-1;1[ \mapsto -t \in \R$ satisfies $\gamma(]-1;0[) \subset ]0;+\infty[$ but for any $\epsilon \in ]0;1[$, $\gamma(\epsilon) \notin ]0;+\infty[$).
\end{enumerate}
\end{rems}

\subsubsection{$\mathcal{AS}$-sets and $\mathcal{AS}$-maps of real projective spaces} \label{subsectionarcsymasprojspaces}

In this final part, we will consider the compact affine Nash manifold $\mathbb{P}^n(\R)$ (see example \ref{exrealprojspacecompactaffnashmanifasquasiarcsymcoincide}): this is the framework of section 3 of \cite{KurPar}. 

First, notice the following examples of $\mathcal{AS}$-sets of $\mathbb{P}^n(\R)$: recall that the Zariski-closed sets of the real algebraic variety $\mathbb{P}^n(\R)$ are the projective algebraic sets of $\mathbb{P}^n(\R)$ (see for instance section 3.4.1 of \cite{BCR}) and denote by $\Phi : \mathbb{P}^n(\R) \rightarrow P$ the biregular isomorphism given in section~3.4.2 of \cite{BCR} sending $\mathbb{P}^n(\R)$ onto a (connected irreducible) compact nonsingular real algebraic set (and then compact Nash submanifold) of $\R^{(n+1)^2}$.

\begin{lem} \label{lemprojalgsetisarcsymrasisasset} Any projective algebraic subset of $\mathbb{P}^n(\R)$ is $\mathcal{AS}$-closed in $\mathbb{P}^n(\R)$. Any real algebraic set of $\R^n$ is an $\mathcal{AS}$-set of $\mathbb{P}^n(\R)$.
\end{lem}

\begin{proof} Let $Z$ be a projective algebraic subset of $\mathbb{P}^n(\R)$, then $\Phi(Z)$ is a real algebraic set of~$\R^{(n+1)^2}$ contained in the compact Nash submanifold $P$: by example \ref{exrealalgsetisarcsym} and propositions~\ref{propcompactnashsubastopcoincide} and \ref{propsaarcanalmapcontastop}, $Z$ is therefore a semialgebraic arc-symmetric set of $\mathbb{P}^n(\R)$.

Now, let $X$ be a real algebraic subset of $\R^n$: if $(x_1,\ldots,x_n)$ denotes the coordinates of $\R^n$ and $[x_0:x_1:\ldots:x_n]$ denotes the homogeneous coordinates of $\mathbb{P}^n(\R)$, we identify $\R^n$ with the Zariski-open subset $\{x_0 \neq 0\}$ of $\mathbb{P}^n(\R)$. The set $X$ is then a Zariski-open set of $\mathbb{P}^n(\R)$ as the set-theoretic difference of the projective algebraic set obtained by homogenizing generators of the ideal $I(X)$ and the hyperplane at infinity of equation $x_0 = 0$. Therefore, by the above showed statement, $X$ is an $\mathcal{AS}$-set of $\mathbb{P}^n(\R)$ as the set-theoretic difference of two $\mathcal{AS}$-closed sets of $\mathbb{P}^n(\R)$.
\end{proof}

\begin{rem} More generally, any connected component of a projective algebraic subset of~$\mathbb{P}^n(\R)$ is a semialgebraic arc-symmetric set of $\mathbb{P}^n(\R)$, using the same argument as above together with example \ref{exanalyticsubsetarcsym} 2 (recall that the connected components of a semialgebraic subset of a real affine space are semialgebraic: cf. Theorem 2.4.5 of \cite{BCR}).
\end{rem}

Even if all real algebraic sets, as well as all Zariski-open sets, of $\R^n$ are $\mathcal{AS}$-sets of $\mathbb{P}^n(\R)$, this is not the case for any $\mathcal{AS}$-set of $\R^n$:

\begin{ex}
Recall the example of remark \ref{remexarcsymr2notarcsymp2}: the positive half-branch $C$ of the hyperbola of equation $x y = 1$ in $\R^2$ is $\mathcal{AR}$-closed in $\R^2$ but if we consider the analytic arc $\gamma : t \in ]-1;1[ \mapsto [-t:1:t^2] \in \mathbb{P}^2(\R)$, we have $\gamma(]-1;0[) \subset C$ and, for all $t \in ]0;1[$,  $\gamma(t) = [-t:1:t^2] = \left[1 : -\frac{1}{t} : -t\right] \in \left\{(x,y) \in \R^2~|~xy = 1, x< 0\right\}$, so that $C$ is not an $\mathcal{AS}$-set of $\mathbb{P}^2(\R)$.

In addition, remark that the $\mathcal{AS}$-closure of $C$ in $\mathbb{P}(\R^2)$ is the projective algebraic set 
$$H := \left\{[x_0 : x_1 : x_2] \in \mathbb{P}^2(\R)~|~x_1 x_2 = x_0^2\right\} = C \cup \{[0:1:0],[0:0:1]\}.$$ 
Indeed, if $S$ is any arc-symmetric subset of $\mathbb{P}(\R^2)$ containing $C$, then $S$ must contain the image of the above arc $\gamma$ (since $\gamma(]-1;0[) \subset C \subset S$ and $S$ is arc-symmetric), but $\gamma(]-1;1[) = H \setminus \{[0:0:1]\}$ (we have $[0:1:0] = \gamma(0)$ and, if $[1:x:y] \in H$, $[1:x:y] = \gamma(-y)$). Analogously, $S$ must contain the image of the analytic arc $\eta : t \in ]-1;1[ \mapsto [-t:t^2:1]$, which is $H \setminus \{[0:1:0\}$. As a consequence, $S$ contains $H$ and therefore $\overline{C}^{\mathcal{AS}}$ contains $H$: since $H$ is a projective algebraic set of $\mathbb{P}(\R^2)$ containing $C$, we finally have $\overline{C}^{\mathcal{AS}} = H$ (see also Remark 3.6 of \cite{KurPar}).
\end{ex}

On the other hand, a \emph{compact} subset of $\R^n$ is $\mathcal{AS}$-closed in $\R^n$ if and only if it is $\mathcal{AS}$-closed in $\mathbb{P}^n(\R)$ (see also Remark 3.5 of \cite{KurPar}):

\begin{prop} \label{propcompactsubrnequivsaarcsympnsaarcsymrn} Let $S$ be compact subset of $\R^n$. Then $S$ is semialgebraic arc-symmetric in $\mathbb{P}^n(\R)$ if and only if $S$ is semialgebraic arc-symmetric in $\R^n$.
\end{prop}

\begin{proof} Notice first that $S$ is semialgebraic in $\mathbb{P}^n(\R)$ if and only if $S$ is semialgebraic in $\R^n$: consider the Nash inclusion map $i : \R^n \hookrightarrow \mathbb{P}^n(\R)$ ($\R^n$ is a Nash submanifold of the Nash manifold $\mathbb{P}^n(\R)$) and apply for instance lemma \ref{lemdirectiminverseimsaissa} 3 (this does not use the compactness hypothesis on $S$).

Now, suppose that $S$ is arc-symmetric in $\mathbb{P}^n(\R)$, then $S = i^{-1}(S)$ is arc-symmetric in $\R^n$ by proposition \ref{propcomparcanalyticmapwitharcsymsetrealanalmanifold} 2 (this is also independent of the compactness hypothesis on $S$).

Conversely, suppose that $S$ is arc-symmetric in $\R^n$ and let $\gamma : ]-1;1[ \rightarrow \mathbb{P}^n(\R)$ be an analytic arc such that $\gamma(]-1;0[) \subset S \subset \R^n$. On the one hand, since $\R^n$ is an $\mathcal{AS}$-set of~$\mathbb{P}^n(\R)$ (lemma \ref{lemprojalgsetisarcsymrasisasset}), there exists $\epsilon \in ]0;1[$ such that $\gamma(]0;\epsilon[) \subset \R^n$. On the other hand, since~$S$ is compact,~$S$ is closed in $\mathbb{P}^n(\R)$ and then $\gamma(0) \in S \subset \R^n$. As a consequence, $\gamma$ induces the analytic arc $\widetilde{\gamma} : t \in ]-1;\epsilon[ \mapsto \gamma(t) \in \R^n$ and, since $\widetilde{\gamma}(]-1;0[) = \gamma(]-1;0[) \subset S$ and $S$ is arc-symmetric in $\R^n$, we have $\gamma(]-1;\epsilon[) = \widetilde{\gamma}(]-1;\epsilon[) \subset S$, so that $S$ is arc-symmetric in~$\mathbb{P}^n(\R)$.
\end{proof}

\begin{ex} Any compact connected component of an $\mathcal{AR}$-closed set of $\R^n$ (e.g. of a real algebraic set of $\R^n$) is $\mathcal{AS}$-closed in $\mathbb{P}^n(\R)$: cf. remark \ref{remconncomparcsymarcsym} (see also example \ref{exanalyticsubsetarcsym} 2) as well as Theorem 2.4.5 of \cite{BCR}.
\end{ex}

We deduce a sufficient condition for an $\mathcal{AS}$-set of $\R^n$ to be an $\mathcal{AS}$-set of $\mathbb{P}^n(\R)$:

\begin{cor} \label{corascompactarcloasinproj} Let $S$ be a subset of $\R^n$ such that the $\mathcal{AR}$-closure $\overline{S}^{\mathcal{AR}}$ of $S$ in $\R^n$ is compact. Then $S$ is an $\mathcal{AS}$-set of $\mathbb{P}^n(\R)$ if and only if $S$ is an $\mathcal{AS}$-set of $\R^n$.
\end{cor}

\begin{proof} Remark that, as in the proof of proposition \ref{propcompactsubrnequivsaarcsympnsaarcsymrn}, $S$ is semialgebraic in $\mathbb{P}^n(\R)$ if and only if $S$ is semialgebraic in $\R^n$ and that, if $S$ is an $\mathcal{AS}$-set of $\mathbb{P}^n(\R)$, $S$ is an $\mathcal{AS}$-set of $\R^n$ by lemma \ref{leminverseimagequasiarcsymbysaarcanalisquasiarcsym}.

Now, suppose that $S$ is an $\mathcal{AS}$-set in $\R^n$ and let $\gamma : ]-1;1[ \rightarrow \mathbb{P}^n(\R)$ be an analytic arc such that $\gamma(]-1;0[) \subset S \subset \overline{S}^{\mathcal{AR}}$. Since $\overline{S}^{\mathcal{AR}}$ is a compact $\mathcal{AR}$-closed subset of $\R^n$, by previous proposition \ref{propcompactsubrnequivsaarcsympnsaarcsymrn}, $\overline{S}^{\mathcal{AR}}$ is $\mathcal{AS}$-closed in $\mathbb{P}^n(\R)$ and then $\gamma(]-1;1[) \subset \overline{S}^{\mathcal{AR}} \subset \R^n$. Consequently,~$\gamma$ induces the analytic arc $\widetilde{\gamma} : t \in ]-1;1[ \mapsto \gamma(t) \in \R^n$ and, since $\widetilde{\gamma}(]-1;0[) = \gamma(]-1;0[) \subset S$ and $S$ is an $\mathcal{AS}$-set of $\R^n$, there exists $\epsilon \in ]0;1[$ such that $\gamma(]0;\epsilon[) = \widetilde{\gamma}(]0;\epsilon[) \subset S$, so that $S$ is an $\mathcal{AS}$-set of $\mathbb{P}^n(\R)$.
\end{proof}

\begin{rem} Up to Nash diffeomorphism, any $\mathcal{AS}$-set of a compact affine Nash manifold can be considered as an $\mathcal{AS}$-set of a real projective space: if $S$ is a subset of a compact Nash manifold~$M$ with affine Nash embedding $\Psi$ onto a Nash submanifold $L$ of $\R^n$, then $S$ is an $\mathcal{AS}$-set of~$M$ if and only if $\Psi(S)$ is an $\mathcal{AS}$-set of $\R^n$ (cf. for instance lemma \ref{remassetcompactsubnashsaaamap}) if and only if $\Psi(S)$ is an~$\mathcal{AS}$-set of~$\mathbb{P}^n(\R)$ (the Euclidean-closed $\mathcal{AR}$-closure of $\Psi(S)$ in $\R^n$ is contained in the compact $\mathcal{AR}$-closed set $L$ of $\R^n$, so that $\overline{\Psi(S)}^{\mathcal{AR}}$ is compact in $\R^n$ as well).
\end{rem}

If $m \in \Nstar$, the real affine space $\R^{n+m} = \R^n \times \R^m$ can also be considered as a Nash submanifold of the product Nash manifold $\mathbb{P}^n(\R) \times \mathbb{P}^m(\R)$. We then prove the following analog of previous corollary \ref{corascompactarcloasinproj} that gives, in particular, a sufficient condition for a map between~$\mathcal{AS}$-sets of real projectives spaces included in affine spaces to be an $\mathcal{AS}$-map.

\begin{prop} \label{propasinprodprojective} Let $S$ be a subset of $\R^{n+m}$ such that the $\mathcal{AR}$-closure $\overline{S}^{\mathcal{AR}}$ of $S$ in $\R^{n+m}$ is compact. Then $S$ is an $\mathcal{AS}$-set of $\mathbb{P}^n(\R) \times \mathbb{P}^m(\R)$ if and only if $S$ is an $\mathcal{AS}$-set of $\R^{n+m}$.
\end{prop}

\begin{proof} Considering the Nash inclusion map $\R^{n+m} \hookrightarrow \mathbb{P}^n(\R) \times \mathbb{P}^m(\R)$, the set $S$ is semialgebraic in $\mathbb{P}^n(\R) \times \mathbb{P}^m(\R)$ if and only if $S$ is semialgebraic in $\R^{n+m}$ (cf. proposition \ref{propcomparcanalyticmapwitharcsymsetrealanalmanifold} 2). Moreover, if $(y_1,\ldots,y_m)$ denotes the coordinates of $\R^m$ and $[y_0:y_1:\ldots:y_m]$ denotes the homogeneous coordinates of $\mathbb{P}^m(\R)$, we have
$$\R^{n+m} = \big(\mathbb{P}^n(\R) \times \mathbb{P}^m(\R)\big) \setminus \big( \{x_0 = 0\} \times \mathbb{P}^n(\R) \cup \mathbb{P}^m(\R) \times \{y_0 = 0\}\big)$$
and $\R^{n+m}$ is therefore an $\mathcal{AS}$-set of $\mathbb{P}^n(\R) \times \mathbb{P}^m(\R)$ (use lemma \ref{lemprojalgsetisarcsymrasisasset} and, for instance, remark~\ref{remcartprodquasiarcsymsisquasiarcsym}). As a consequence, we can use the same arguments as in the proofs of proposition~\ref{propcompactsubrnequivsaarcsympnsaarcsymrn} and corollary \ref{corascompactarcloasinproj} to assert the desired equivalence.
\end{proof}

\begin{cor} \label{corsufficientconditionasmapbetweenrealassets} Let $S$ be an $\mathcal{AS}$-set of $\mathbb{P}^n(\R)$ included in $\R^n$, $T$ be an $\mathcal{AS}$-set of $\mathbb{P}^m(\R)$ included in $\R^m$ and $f : S \rightarrow T$ be a map. If the $\mathcal{AR}$-closure of $\Gamma_f$ in $\R^{n+m}$ is compact, then $f$ is an $\mathcal{AS}$-map if and only if $\Gamma_f$ is an $\mathcal{AS}$-set of $\R^{n+m}$.
\end{cor}

\begin{proof} It is a direct application of proposition \ref{propasinprodprojective}.
\end{proof}

\begin{ex} \label{exaftercorsufficientconditionasmapbetweenrealassets} Keep the notations of previous corollary \ref{corsufficientconditionasmapbetweenrealassets}. If the respective $\mathcal{AR}$-closures of~$S$ and $T$ in $\R^n$ and $\R^m$ are compact, the graph $\Gamma_f$ of $f$ satisfies the condition of the statement since $\overline{\Gamma_f}^{\mathcal{AR}}$ is a (Euclidean) closed set of $\R^{n+m}$ (recall for instance proposition \ref{propsemialgarcsymisclosed}) included in $\overline{S}^{\mathcal{AR}} \times \overline{T}^{\mathcal{AR}}$.
\end{ex}

\begin{rem} \label{rempolynomialmaphasasgraphinproj} Notice that any real algebraic set $X$ of $\R^{n+m}$ is an $\mathcal{AS}$-set of $\mathbb{P}^n(\R) \times \mathbb{P}^m(\R)$. Indeed, write $X = V(P)$ with $P \in \R[x_1,\ldots,x_n,y_1,\ldots,y_m]$ and 
$$P = \sum_{\alpha \in \N^n, \beta \in \N^m} a_{\alpha,\beta} x^{\alpha} y^{\beta}.$$ 
Denote by $d$, resp. $e$, the maximum of the quantities $|\alpha| := \alpha_1 + \cdots \alpha_n$, resp. $|\beta|$, for $\alpha \in \N^n, \beta \in \N^m$ such that $a_{\alpha,\beta} \neq 0$, and consider the homogeneous polynomial
$$\widetilde{P} := \sum_{\alpha \in \N^n, \beta \in \N^m} a_{\alpha,\beta} x^{\alpha} x_0^{d - |\alpha|} y^{\beta} y_0^{e - |\beta|}$$
of $\R[x_0,x_1,\ldots,x_n,y_0,y_1,\ldots,y_m]$. If we denote $\widetilde{X} := \left\{\widetilde{P} = 0\right\} \subset \mathbb{P}^n(\R) \times \mathbb{P}^m(\R)$, the set $\widetilde{X}$ is then a semialgebraic arc-symmetric set of $\mathbb{P}^n(\R) \times \mathbb{P}^m(\R)$ (lemma \ref{lemanalyticsubsetmanifoldarcsym}: notice that $\widetilde{X}$ is closed with respect to any of the real affine spaces covering $\mathbb{P}^n(\R) \times \mathbb{P}^m(\R)$ hence is closed in $\mathbb{P}^n(\R) \times \mathbb{P}^m(\R)$) and we have 
$$X = \widetilde{X} \setminus \big(\{x_0 = 0\} \times \mathbb{P}^n(\R) \cup \mathbb{P}^m(\R) \times \{y_0 = 0\}\big),$$
so that $X$ is an $\mathcal{AS}$-set of $\mathbb{P}^n(\R) \times \mathbb{P}^m(\R)$.

In particular, if $B$ and $A$ are respective real algebraic subsets of $\R^n$ and $\R^m$ and $f : B \rightarrow A$ is a map with algebraic graph (i.e. the graph $G_f$ of $f$ is a real algebraic set of $\R^{n+m}$), then $f$ is an~$\mathcal{AS}$-map between the $\mathcal{AS}$-sets $B \subset \mathbb{P}^n(\R)$ and $A \subset \mathbb{P}^m(\R)$.
\end{rem}

\end{document}